\documentclass[11pt]{article}

\usepackage{amssymb,epsfig}
\usepackage{amsmath}
\usepackage{comment}
\usepackage[dvips]{color}

\usepackage{mathrsfs}
\renewcommand\mathcal{\mathscr}

\usepackage[T1]{fontenc}
\usepackage[dvips]{color}
\usepackage{amssymb}

\newtheorem{defi}{Definition}[section]
\newtheorem{prop}[defi]{Proposition}
\newtheorem{theo}[defi]{Theorem}
\newtheorem{conj}[defi]{Conjecture}
\newtheorem{lemm}[defi]{Lemma}
\newtheorem{coro}[defi]{Corollary}
\newtheorem{rema}[defi]{Remark}
\newtheorem{exem}[defi]{Example}
\newtheorem{exems}[defi]{Examples}

\newcommand{\bdefi}{\begin{defi}}
\newcommand{\edefi}{\end{defi}}
\newcommand{\bprop}{\begin{prop}}
\newcommand{\eprop}{\end{prop}}
\newcommand{\btheo}{\begin{theo}}
\newcommand{\etheo}{\end{theo}}
\newcommand{\blemm}{\begin{lemm}}
\newcommand{\brema}{\begin{rema}}
\newcommand{\erema}{\end{rema}}
\newcommand{\bexer}{\begin{exem}}
\newcommand{\eexer}{\end{exem}}
\newcommand{\bexems}{\begin{exems}}
\newcommand{\eexems}{\end{exems}}
\newcommand{\bconj}{\begin{conj}}
\newcommand{\econj}{\end{conj}}
\newcommand{\elemm}{\end{lemm}}
\newcommand{\bcoro}{\begin{coro}}
\newcommand{\ecoro}{\end{coro}}
\newcommand{\dem}{\noindent{\bf Proof. }}
\newcommand{\rem}{\noindent{\bf Remark. }}

\newcommand{\A}{{\cal A}}

\newcommand{\M}{{\cal M}}
\newcommand{\N}{{\cal N}}

\newcommand{\E}{{\cal E}}

\renewcommand{\H}{{\cal H}}
\renewcommand{\O}{{\cal O}}
\newcommand{\C}{{\cal C}}
\newcommand{\I}{{\cal I}}

\newcommand{\maths}[1]{{\mathbb #1}}  

\newcommand{\RR}{\maths{R}}
\newcommand{\NN}{\maths{N}}
\newcommand{\CC}{\maths{C}}
\newcommand{\QQ}{\maths{Q}}

\renewcommand{\SS}{\maths{S}}

\newcommand{\HH}{\maths{H}}

\newcommand{\ZZ}{\maths{Z}}
\newcommand{\PP}{\maths{P}}

\newcommand{\ra}{\rightarrow}

\newcommand{\ov}[1]{{\overline #1}} 
\newcommand{\wt}[1]{{\widetilde{#1}}}

\newcommand{\ga}{\gamma}
\newcommand{\Ga}{\Gamma}

\newcommand{\cqfd}{\hfill$\Box$}

\newcommand{\PSLZ}{\mbox{${\rm{ PSL}}_{2}(\ZZ)$}}

\newcommand{\hdr}{{\HH\,}^2_\RR}

\newcommand{\hnr}{{\HH\,}^n_\RR}
\newcommand{\hnc}{{\HH\,}^n_\CC}

\newcounter{fig}

\title{Prescribing the behaviour of geodesics 
in negative curvature} 
\author{Jouni Parkkonen \and Fr\'ed\'eric Paulin}
\date{\today}

\newcommand{\Isom}{{\operatorname{Isom}}}
\newcommand{\conv}{{\operatorname{Conv}}}
\newcommand{\diam}{{\operatorname{diam}}}
\newcommand{\arcosh}{{\operatorname{arcosh}}}
\newcommand{\arsinh}{{\operatorname{arsinh}}}
\renewcommand{\Im}{{\operatorname{Im}}}
\newcommand{\height}{\operatorname{ht}}
\newcommand{\maxht}{\operatorname{maxht}}
\newcommand{\maxsp}{\operatorname{MaxSp}}
\newcommand{\limsupsp}{\operatorname{LimsupSp}}
\newcommand{\inj}{\operatorname{inj}}

\newcommand{\gf}{\operatorname{\mathcal{GF}}}
\def\norm#1{{\lVert#1\rVert_\infty}}

\renewcommand{\P}{{\cal P}}
\renewcommand{\tilde}{\widetilde}
\renewcommand{\ov}{\overline}
\newcommand{\Q}{{\cal Q}}

\newcommand{\bp}{{\operatorname{\mathfrak b\mathfrak p}}}
\newcommand{\ph}{{\operatorname{\mathfrak p\mathfrak h}}}
\newcommand{\ipp}{{\operatorname{\mathfrak i\mathfrak p\mathfrak p}}}
\newcommand{\cp}{{\operatorname{\mathfrak c\mathfrak r\mathfrak p}}}
\newcommand{\ftp}{{\operatorname{\mathfrak f\mathfrak t\mathfrak p}}}
\newcommand{\codim}{\operatorname{codim}}

\begin{document}
\maketitle

\begin{abstract}
Given a family of (almost) disjoint strictly convex subsets of a
complete negatively curved Riemannian manifold $M$, such as balls,
horoballs, tubular neighborhoods of totally geodesic submanifolds,
etc, the aim of this paper is to construct geodesic rays or lines in
$M$ which have exactly once an exactly prescribed (big enough)
penetration in one of them, and otherwise avoid (or do not enter too
much in) them. Several applications are given, including a definite
improvement of the  unclouding problem of \cite{PP}, the
prescription of heights of geodesic lines in a finite volume
such $M$, or of  spiraling times around a closed geodesic in a
closed such $M$. We also prove that the Hall ray phenomenon described
by Hall in special arithmetic situations and by Schmidt-Sheingorn for
hyperbolic surfaces is in fact only a negative curvature property.
\footnote{ {\bf AMS codes: } 53 C 22, 11 J 06, 52 A 55, 53 D 25. {\bf
Keywords:} geodesics, negative curvature, horoballs, Lagrange
spectrum, Hall ray.  }
\end{abstract}

\section{Introduction}
\label{sec:intro}

%
%
%

The problem of constructing obstacle-avoiding geodesic rays or lines
in negatively curved Riemannian manifolds has been studied in various
different contexts.  For example, Dani \cite{Dan} and others
\cite{Str,AL,KW} have constructed (many) geodesic rays that are
bounded (i.e.~avoid a neighborhood of infinity) in noncompact
Riemannian manifolds. This work has deep connections with Diophantine
approximation problems, see for instance the papers by Sullivan
\cite{Sul}, Kleinbock-Margulis \cite{KM} and Hersonsky-Paulin
\cite{HP2}.  Hill and Velani \cite{HV} and others (see for instance
\cite{HP3}) have studied the shrinking target problem for the geodesic
flow.  Schroeder \cite{Schr} and others \cite{BSW} have worked on the
construction of geodesic lines avoiding given subsets, see also the
previous work \cite{PP} of the authors on the construction of geodesic
rays and lines avoiding a uniformly shrunk family of horoballs.

In this paper, we are interested in constructing geodesic rays or
lines in negatively curved Riemannian manifolds which, given some
family of obstacles, have exactly once an exactly prescribed (big
enough) penetration in one of them, and otherwise avoid (or do not
enter too much in) them. We also study an asymptotic version of this
problem.  This introduction contains a sample of our results (see also
\cite{PP2}).

\medskip %
Given a horoball $H$ of center $\xi$ or a ball of center $x$ and
radius $r$ in a CAT$(-1)$ metric space (such as a complete simply
connected Riemannian manifold of sectional curvature at most 
$-1$
), for every $t\geq 0$, let $H[t]$ be the
concentric horoball or ball contained in $H$, whose boundary is at
distance $t$ from the boundary of $H$ (with $H[t]$ empty if $H$ is a
ball of radius $r$ and $t>r$).  The following result (see Section
\ref{subsec:uncloudingshort}) greatly improves the main results,
Theorem 1.1 and Theorem 4.5, of \cite{PP}. The fact that the constant
$\mu_0$ is universal (and not very big, though not optimal) is indeed
remarkable.

\btheo \label{theo:introbetterGAFA} %
Let $X$ be a proper geodesic CAT$(-1)$ metric space with arcwise
connected boundary $\partial_\infty X$ and extendible geodesics, let
$(H_\alpha)_{\alpha\in\A}$ be any family of balls or horoballs with
pairwise disjoint interiors, and let $\mu_0=1.534$. For every $x$ in
$X-\bigcup_{\alpha\in\A} H_\alpha$, there exists a geodesic ray
starting from $x$ and avoiding $H_\alpha[\mu_0]$ for every $\alpha$.
\etheo

\medskip %
From now on, we denote by $M$ a complete connected Riemannian manifold
with sectional curvature at most $-1$ and dimension at least $3$.

If $M$ has finite volume and is non compact, let $e$ be an end of $M$.
Let $V_e$ be the maximal Margulis neighborhood of $e$ (see for
instance \cite{BK,Bow,HP2} and Section \ref{sec:ballsandhoroballs}).
If $\rho_e$ is a minimizing geodesic ray in $M$ starting from a point
in the boundary of $V_e$ and converging to $e$, let $\height_e:M\ra
\RR$ be the height map defined by $\height_e(x)= \lim_{t\ra\infty}
\;(t-d(\rho_e(t),x))$.  The {\it maximum height spectrum}
$\maxsp(M,e)$ of the pair $(M,e)$ is the subset of $]-\infty,+\infty]$
consisting of elements of the form $\sup_{t\in\RR}
\height_e(\gamma(t))$ where $\ga$ is a locally geodesic line in $M$.

As a consequence of Theorem \ref{theo:introbetterGAFA} (see Corollary
\ref{coro:lowclosedgeod}), we prove that in any (noncompact) finite
volume complete Riemannian manifold with dimension at least $2$ and
sectional curvature at most $-1$, there exists universally low closed
geodesics.  Furthermore, we have the following result on the upper
part of the spectrum:

\btheo\label{theo:introhallhoro} %
If $M$ has finite volume and $e$ is an end of $M$, then $\maxsp(M,e)$
contains the interval $[4.2,+\infty]$.  
\etheo

For more precise analogous statements when $M$ is geometrically
finite, and for finite subsets of cusps of $M$, see Section
\ref{sec:ballsandhoroballs}.  Schmidt and Sheingorn \cite{SS} proved
the two-dimensional analog of Theorem \ref{theo:introhallhoro} in
constant curvature $-1$. They showed that the maximum height spectrum
of a finite area hyperbolic surface with respect to any cusp contains
the interval $[4.61,+\infty]$.

\medskip
The previous result is obtained by studying the penetration properties
of geodesic lines in a family of horoballs. Our next theorem  concerns
families of balls (see Section \ref{sec:ballsandhoroballs} for
generalizations).

\medskip %
\btheo \label{theo:balinjrad} %
Let $x$ be point in $M$ with $r=\inj_M x\geq 56$. Then, for every
$d\in [2,r-54]$, there exists a locally geodesic line $\ga$ passing at
distance exactly $d$ from $x$ at time $0$ and remaining at distance
greater than $d$ from $x$ at any nonzero time.  
\etheo

Given a closed geodesic $L$ in $M$, the behavior of a locally geodesic
ray $\ga$ in $M$ with respect to $L$ is typically that $\ga$ spirals
around $L$ for some time, then wanders away from $L$, then spirals
again for some time around $L$, then wanders away, etc. Our next aim
is to construct such a $\ga$ which has exactly one (big enough)
exactly prescribed spiraling length, and all of whose other spiraling
lengths are bounded above by some uniform constant. Let us make this
precise (see the figure in Section \ref{subsec:totgeod}).

Let $L$ be an embedded compact totally geodesic submanifold in $M$
with $1\leq{\rm dim}\; L\leq {\rm dim}\; M-1$, and $\epsilon>0$ small
enough so that the (closed) $\epsilon$-neighborhood $\N_\epsilon L$ of
$L$ is a tubular neighborhood. For every locally geodesic line $\ga$
in $M$, the set of $t\in\RR$ such that $\ga(t)$ belongs to
$\N_\epsilon L$ is the disjoint union of maximal closed intervals
$[s_n,t_n]$, with $s_n\leq t_n<s_{n+1}$. For every such interval, let
$x_n$ (resp.~$y_n$) be the origin of the locally geodesic ray starting
from a point of $L$ perpendicularly to $L$, which is properly
homotopic, by a homotopy with origins in $L$, to the piecewise
geodesic path starting with the geodesic segment from the closest
point of $\ga(s_n)$ (resp.~$\ga(t_n)$) on $L$ to $\ga(s_n)$
(resp.~$\ga(t_n)$), and then following $\ga$ at times less than $s_n$
(resp.~more than $t_n$). The length of the locally geodesic path in
$L$ between $x_n, y_n$ in the obvious homotopy class will be called a
{\it fellow-traveling time} of $\ga$ along $L$ (see Section
\ref{subsec:totgeod}).

\btheo 
\label{theo:introftp} %
Let $L$ be as above.  There exist constants $c,c'>0$, depending only
on $\epsilon$, such that for every $h\geq c$, there exists a locally
geodesic line in $M$, having one fellow-traveling time exactly $h$,
all others being at most $c'$.  
\etheo

See Section \ref{subsec:totgeod} for an extension of Theorem
\ref{theo:introftp} when $L$ is not necessarily embedded, and to
finitely many disjoint such neighborhoods $\N_\epsilon L$. If $M$ has
finite volume, we also construct bounded locally geodesic lines with
the above property (with a control of the heights uniform in
$\epsilon$). In constant curvature, we can also prescribe one of the
penetration lengths $|t_n-s_n|$ at least $c$, while keeping all the
other ones at most $c'$.  Schmidt and Sheingorn \cite{SS} sketch the
proof of a result for hyperbolic surfaces which is analogous to
Theorem \ref{theo:introftp} with a different way of measuring the
affinity of locally geodesic lines. Other results about the spiraling
properties of geodesic lines around closed geodesics will be given in
\cite{HP6,HPP}.

For our next result, we specialize to the case where $M$ is a
hyperbolic $3$-manifold. See Section \ref{subsec:recur} for a more
general statement, and for instance \cite{MT} for references on
$3$-manifolds and Kleinian groups.

\btheo \label{theo:intromaxrec} Let $N$ be a compact, connected,
orientable, irreducible, acylindrical, atoroi\-dal, boundary
incompressible $3$-manifold with boundary, with $\partial N$ having
exactly one torus component $e$. For every compact subset $K$ in the
space $\gf(N,e)$ of (isotopy classes of) complete geometrically finite
hyperbolic metrics in the interior of $N$ with one cusp, there exists
a constant $c\geq 0$ such that for every $h\ge c$ and every $\sigma\in
K$, there exists a locally geodesic line $\ga$ contained in the
convex core of $\sigma$ such that the maximum height of $\ga$ is
exactly $h$.  
\etheo

If $M$ has finite volume and $e$ is an end of $M$, define the {\it
asymptotic height spectrum} $\limsupsp(M,e)$ of the pair $(M,e)$ to be
the subset of $]-\infty,+\infty]$ consisting of elements of the form
$\limsup_{t\in\RR} \height_e(\gamma(t))$ where $\ga$ is a locally
geodesic line in $M$.

\btheo[The ubiquity of Hall rays]\label{theo:introhallquoi} %
If $M$ has finite volume and $e$ is an end of $M$, then
$\limsupsp(M,e)$ contains $[6.8, \infty]$.  
\etheo

The interval given by Theorem \ref{theo:introhallquoi} is called a
{\it Hall ray}. Note that the value $6.8$ is uniform on all couples
$(M,e)$, but we do not know the optimal value.  If $M$ is the
one-ended hyperbolic $2$-orbifold ${\rm PSL}_2(\ZZ)\backslash
\HH^2_\RR$ where $\HH^2_\RR$ is the real hyperbolic plane with
sectional curvature $-1$, then the existence of a Hall ray follows
from the work of Hall \cite{Hall1,Hall2} on continued fractions.
Freiman \cite{Freiman} (see also \cite{Slo}) has determined the
maximal Hall ray of ${\rm PSL}_2(\ZZ)\backslash \HH^2_\RR$, which is
approximately $[3.02,+\infty]$. The generality of Theorem
\ref{theo:introhallquoi} proves in particular that the Hall ray
phenomenon is neither an arithmetic nor a constant curvature property.
See Section \ref{subsec:limsupspec} for a more precise version of
Theorem \ref{theo:introhallquoi}, which is valid also in the
geometrically finite case.

\medskip %
The results of Hall and Freiman cited above were originally formulated
in terms of Diophantine approximation of real numbers by rationals.
The projective action of the modular group ${\rm PSL}_2(\ZZ)$ on the
upper halfplane provides a way to obtain the geometric interpretation.
We conclude this sample of our results by giving applications of our
methods to Diophantine approximation problems (see Section
\ref{sec:diophapp} for generalizations in the framework of
Diophantine approximation on negatively curved manifolds, developped
in \cite{HP3,HP4,HP2}).  These results were announced in \cite{PP2}.

\btheo \label{theo:applibianchi} 
Let $m$ be a squarefree positive integer, and $\I$ a non-zero ideal in
an order $\O$ in the ring of integers $\O_{-m}$ of the imaginary
quadratic number field $\QQ(i\sqrt{m})$. For every $x\in\CC-
\QQ(i\sqrt{d})$, let
$$
c(x)=\liminf_{(p,q)\in \O\times\I,\,
\langle p,q\rangle=\O, \,|q|\ra\infty}\;\;|q|^2\Big|x-\frac{p}{q}\Big|
$$ 
be the approximation constant of the complex number $x$ by elements
of $\O\I^{-1}$, and ${\rm Sp}_{Lag}$ the Lagrange spectrum consisting
of the real numbers of the form $c(x)$ for some
$x\in\CC-\QQ(i\sqrt{m})$. Then ${\rm Sp}_{Lag}$ contains the interval
$[0,0.033]$.  
\etheo

Theorem \ref{theo:applibianchi} follows from Hall's result and from
the work of Poitou \cite{Poi} in the particular case $\I=\O=\O_{-m}$.
Other arithmetic applications of our geometric methods can be obtained
by varying the (non uniform) arithmetic lattice in the isometry group
of a negatively curved symmetric space. We only state the following
result in this introduction (with the notation of Section
\ref{subsec:comphypgeom}), see Section \ref{subsec:ubiqhallways} and
\cite{PP2} for other ones.

\btheo \label{theo:applicygan} Let $\Q(\RR)$ be the real quadric
$\{(z,w)\in\CC^2\;:\;2\,{\rm Re}\; z - |w|^2=0\}$ endowed with the Lie
group law $(z,w)\cdot(z',w')=(z+z'+w'\overline{w},w+w')$ and
$\Q(\QQ)=\Q(\RR)\cap\QQ(i)^2$ be its rational points. If
$r=(p/q,p'/q)\in\Q(\QQ)$ with $p,p',q\in\ZZ[i]$ relatively prime, let
$h(r)=|q|$. Let $d'_{\rm Cyg}$ be the left-invariant distance on
$\Q(\RR)$ such that $d'_{\rm Cyg}((z,w),(0,0))=
\sqrt{2|z|+|w|^2}$. For every $x\in\Q(\RR)-\Q(\QQ)$, let
$$
c(x)=
\liminf_{r\in\Q(\QQ)\,,\; h(r)\ra\infty}h(r)\;d'_{\rm Cyg}(x,r)
$$ be the approximation constant of $x$ by rational points, and ${\rm
Sp}_{Lag}$ the Lagrange spectrum consisting of the real numbers of the
form $c(x)$ for some $x\in\Q(\RR)-\Q(\QQ)$.  Then ${\rm Sp}_{Lag}$
contains the interval $[0,0.047]$.  
\etheo

The paper is organized as follows. In Section \ref{sec:strictconvex},
we define a class of uniformely strictly convex subsets of metric
spaces, that we call {\it $\epsilon$-convex subsets}. We study the
interaction of geodesic rays and lines with $\epsilon$-convex sets in
CAT($-1$)-spaces.  In particular, we give various estimates on the
distance between the entering and exiting points in an
$\epsilon$-convex set of two geodesic rays starting from a fixed point
in the space and of two geodesic lines starting from a fixed point in
the boundary at infinity. Section \ref{sec:proppene} is devoted to
defining and studying several penetration maps which are used to
measure the penetration of geodesic rays and lines in an
$\epsilon$-convex set. We emphasize the case of penetration maps
in horoballs, balls and tubular neighborhoods of totally geodesic
submanifolds. We show that in a number of geometrically interesting
cases, it is possible to adjust the penetration of a geodesic line or
ray in one $\epsilon$-convex set while keeping the penetration in
another set fixed. Section \ref{sec:main} contains the inductive
construction that gives geodesic rays and lines with prescribed
maximal penetration with respect to a given collection of
$\epsilon$-convex sets. As a warm-up for the construction, we prove
Theorem \ref{theo:introbetterGAFA} in Subsection
\ref{subsec:uncloudingshort}.  
The other theorems in the introduction
besides the last two and a number of others are proved in Section
\ref{sec:results} where the results of Section \ref{sec:main} are
applied in the cases studied in Section \ref{sec:proppene}. Finally,
we give our arithmetic applications in Section \ref{sec:diophapp}.

\medskip \noindent {\small {\it Acknowledgments. } Each author
  acknowledges the support of the other author's institution, where
  part of this work was done. This research was supported by the
  Center of Excellence "Geometric analysis and mathematical physics"
  of the Academy of Finland. We thank P.~Pansu, Y.~Bugeaud,
  A.~Schmidt, P.~Gilles, A.~Guilloux, D.~Harari for various
  discussions and comments on this paper. }

\section{On strict convexity in CAT$(-1)$ spaces}
\label{sec:strictconvex}

%
%
%

\subsection{Notations and background}
\label{subsec:notations}

In this section, we introduce some of the objects which are central in
this paper.  We refer to \cite{BH,GH} for the definitions and basic
properties of CAT($-1$) spaces.  Our reference for hyperbolic geometry
is \cite{Bea}.

\medskip %
Let $(X,d)$ be a proper geodesic CAT$(-1)$ metric space, and
$X\cup\partial_\infty X$ be its compactification by the asymptotic
classes of geodesic rays.  By a {\it geodesic line} (resp.~{\it ray}
or {\em segment}) in $X$, we mean an isometric map $\ga:\RR\ra X$
(resp.~$\ga:[\iota_0,+\infty[\;\ra X$ with $\iota_0\in\RR$ or
$\ga:[a,b]\;\ra X$, with $a\leq b$). We sometimes also denote by $\ga$
the image of this map. For $x,y$ in $X$, we denote by $[x,y]$ the
(unique) closed geodesic segment between $x,y$, with the obvious
extension to open and half-open geodesic segments, rays and lines
(with one or two endpoints in $\partial_\infty X$). We say that $X$
{\it has extendible geodesics} if every geodesic segment can be
extended to a geodesic line.

We denote by $T^1X$ the space of geodesic lines in $X$, endowed with
the compact-open topology. When $X$ is a Riemannian manifold, the
space $T^1X$ coincides with the usual definition of the unit tangent
bundle, upon identifying a geodesic line $\ga$ and its (unit) tangent
vector $\dot{\ga}(0)$ at time $t=0$.  For every geodesic ray or line
$\ga$, we denote by $\ga(+ \infty)$ the point of $\partial_\infty X$
to which $\ga ( t)$ converges as $t\ra+\infty$, and we define $\ga(-
\infty)$ similarly when $\ga$ is a geodesic line. We say that a
geodesic line (resp.~ray) $\ga$ {\em starts} from a point $\xi\in
\partial_\infty X$ (resp.~$\xi\in X$) if $\xi=\ga(-\infty)$
(resp.~$\ga(\iota_0)=\xi$). For every $\xi$ in $X\cup\partial_\infty
X$, we denote by $T^1_\xi X$ the space of geodesic lines (if
$\xi\in\partial_\infty X$) or rays (if $\xi\in X$) starting from
$\xi$, endowed with the compact-open topology.

If $Y$ is a subset of $X$ and $\xi$ a point in $X\cup\partial_\infty
X$, the {\it shadow of $Y$ seen from $\xi$} is the set 
$\O_\xi Y$ 
of points $\gamma(+\infty)$ where $\ga$ is a geodesic ray
or line starting from $\xi$ and meeting $Y$.

The {\it Busemann function} $\beta_{\xi}:X\times X\ra \RR$
at a point $\xi$ in $\partial_\infty X$ is defined by 
$$
\beta_{\xi}(x,y) =
\lim_{t\ra+\infty} \big(d(x,\rho(t))-d(y,\rho(t))\big),
$$ 
where $\rho$ is any geodesic ray ending at $\xi$.  The function
$x\mapsto\beta_\xi(x,y)$ can be thought of as a normalized signed
distance to $\xi\in\partial_\infty X$, or as the {\it height} of the
point $x$ with respect to $\xi$ (relative to $y$).  Accordingly, if
$\beta_\xi(x,y)=\beta_\xi(x',y)$, then the points $x$ and $x'$ are
said to be {\it equidistant} to $\xi$.  If $\xi\in X$, we define
$\beta_{\xi}(x,y)= d(x,\xi)-d(y,\xi)$. This is convenient in Section
\ref{sec:induction} and in the proof of Corollary \ref{coro:ballline}.
For every $x,y,z$ in $X$ and $\xi\in X\cup \partial_\infty X$, we have
$$
\beta_{\xi}(x,y)+\beta_{\xi}(y,z)=\beta_{\xi}(x,z),
$$
$\beta_{\xi}(x,x)=0$, and $|\beta_\xi(x,y)|\leq d(x,y)$.

A {\it horoball} in $X$ centered at $\xi\in\partial_\infty X$ is the
preimage of $[s,+\infty[$ for some $s$ in $\RR$ by the map
$y\mapsto\beta_\xi(x,y)$ for some $x$ in $X$.  If
$$
H=\{y\in X\;:\; \beta_\xi(x,y)\geq s\}
$$ 
is a horoball, we define its {\it boundary horosphere} by 
$$
\partial H=\{y\in X\;:\;\beta_\xi(x,y)= s\},
$$ 
and for every $t\geq 0$, its {\it $t$-shrunk horoball} by 
$$
H[t]=\{y\in X\;:\; \beta_\xi(x,y)\geq s+t\}.
$$
(In \cite{PP}, we denoted $H[t]$ by $H(t)$.)  Similarly, if $B$ is a
ball of center $x$ and radius $r$, for every $t\leq r$, we denote by
$B[t]$ the ball of center $x$ and radius $r-t$.  By convention, if
$t>r$, define $B[t]=\emptyset$. Note that for every ball or horoball
$H$, we have $H[t']\subset H[t]$ if $t'\geq t$. The point at infinity
of an horoball $H$ is denoted by $H[\infty]$.  Note that, in this
paper, all balls and horoballs in $X$ are assumed to be closed.

Recall that a subset $C$ in a CAT$(-1)$ metric space is {\it convex}
if $C$ contains the geodesic segment between any two points in $C$.
Let $C$ be a convex subset in $X$. We denote by $\partial_\infty C$
its set of points at infinity, and by $\partial C$ its boundary in
$X$.  If $C$ is nonempty and closed, for every $\xi$ in
$\partial_\infty X$, we define {\it the closest point to $\xi$ on the
  convex set $C$} to be the following point $p$ in
$C\cup\partial_\infty C$: if $\xi\notin\partial_\infty C$, then $p$
belongs to $C$ and maximizes the map $y\mapsto \beta_{\xi}(x_0,y)$ for
some (hence any) given point $x_0$ in $X$; if $\xi\in\partial_\infty
C$, then we define $p=\xi$. This $p$ exists, is unique, and depends
continuously on $\xi$, by the properties of CAT$(-1)$-spaces.

If $x,y,z\in X\cup\partial_\infty X$, we denote by $(x,y,z)$ the
triangle formed by the three geodesic segments, rays or lines with
endpoints in $\{x,y,z\}$.  Recall that if $\alpha:t\mapsto \alpha_t$
and $\beta:t\mapsto\beta_t$ are two (germs of) geodesic segments
starting from a point $x_0$ in $X$ at time $t=0$, if $(\overline{x}_0,
\overline{\alpha}_t, \overline{\beta}_t)$ for $t>0$ small enough is a
comparison triangle to $(x_0,{\alpha_t}, {\beta_t})$ in the real
hyperbolic plane $\hdr$, then the {\it comparison angle} between
$\alpha$ and $\beta$ at $x_0$ is the limit, which exists, of the angle
$\angle_{\ov {x_0}}(\ov{\alpha_t},\ov{\beta_t})$ as $t$ tends to $0$.

\bigskip %
We end this section with the following (well known) exercises in
hyperbolic geometry.

\blemm\label{lem:expopinch}%
For all points $x,y$ in $X$ and $z$ in $X\cup\partial X$, and every
$t$ in $[0,d(x,z)]$, if $x_t$ is the point on $[x,z]$ at distance $t$
from $x$, then $$d(x_t,[y,z])\leq e^{-t}\sinh d(x,y)\leq \frac{1}{2}\;
e^{-t+d(x,y)}\;.$$
\elemm

\dem %
By comparison, we may assume that $X=\HH^2_\RR$. As it does not
decrease $d(x_t,[y,z])$ to replace $z$ by the point at infinity of the
geodesic ray starting from $x$ and passing through $z$, we may assume
that $z$ is the point at infinity in the upper halfspace model of
$\HH^2_\RR$.  Let $p$ be the orthogonal projection of $x_t$ on the
geodesic line $\ga$ through $y$ and $z$. Assume first that $p$ belongs
to $[y,z[$.

\noindent
\begin{minipage}{10.5cm}
  ~~~If we replace $y$ by the orthogonal projection of $x$ on $\ga$,
  then we decrease $d(x,y)$, and do not change $t$ and $d(x_t,[y,z])$.
  Hence we may assume that $y=i$ and $x$ is on the (Euclidean) circle
  of center $0$ and radius $1$. If $\alpha$ is the (Euclidean) angle
  at $0$ between the ho\-ri\-zontal axis and the (Euclidean) line from
  $0$ passing through $x$, then an easy computation in hyperbolic
  geometry (see also \cite{Bea}, page 145) gives $\sinh d(x,y) =\cos
  \alpha/\sin\alpha$.  Similarly, $\sinh d(x_t,[y,z]) =\cos
  \alpha/(e^t\sin\alpha)$. So that
$$d(x_t,[y,z])\leq \sinh d(x_t,[y,z])=e^{-t}\sinh
d(x,y)\;.$$
\end{minipage}
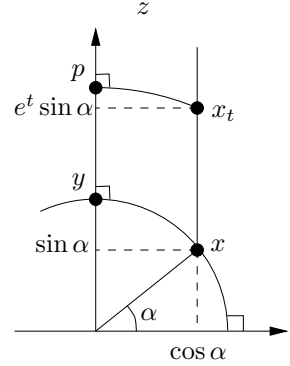
\begin{minipage}{4.3cm}
\begin{center}
\input{fig_expopinch.pstex_t}
\end{center}
\end{minipage}
\medskip

Assume now that $p$ does not belong to $[y,z[$. In particular, $y\neq
x_t$.  Let $\overline{x_t}$ be the point at same distance from $y$ as
$x_t$ (and on the same side) such that $y$ is the orthogonal
projection of $\overline{x_t}$ on $\ga$, so that
$$
d(\overline{x_t},[y,z])=d(\overline{x_t},y)=d(x_t,y)=d(x_t,[y,z]).
$$

\noindent
\begin{minipage}{9.7cm}
  ~~~Let $\overline{x}$ be the intersection of the geodesic line from
  $z$ through $\overline{x_t}$ with the (hyperbolic) circle of center
  $y$ and radius $d(x,y)$, so that $d(\overline{x},y)=d(x,y)$. Then,
  with $\overline{t}= d(\overline{x_t},\overline{x})$, we have
  $\overline{t}\geq t$, as the angle at $\overline{x_t}$ of
  $[\overline{x_t},\overline{x}]$ with the outgoing unit vector of the
  geodesic ray from $y$ through $\overline{x_t}$ is bigger than the
  corresponding one for $x_t$ and $x$. Hence we may assume that
  $x_t=\overline{x_t}$ and $x=\overline{x}$. As then the orthogonal
  projection of $x_t$ on the geodesic line through $y$ and $z$ is $y$,
  this reduces the situation to the first case treated above. \cqfd
\end{minipage}
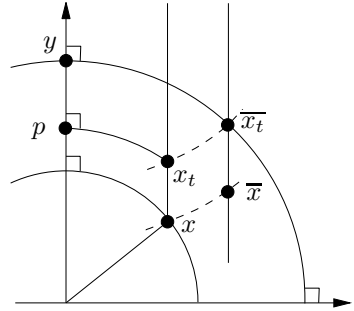
\begin{minipage}{4.7cm}
\begin{center}
\input{fig_expopinch2.pstex_t}
\end{center}
\end{minipage}

\blemm\label{lem:midpointclose}%
For every $\epsilon>0$, if
\begin{equation}\label{eq:c0}
c_0(\epsilon)=
2\log\Big(\frac{2(1+e^{\epsilon/2})\sinh\epsilon}\epsilon\Big)\;, 
\end{equation}
then for all points $a,b,a',b'$ in $X$ 
such that
$$
d(a,a')\leq \epsilon\;,\;\;\;d(b,b')\leq \epsilon\;, 
\;\;\;d(a,b)\geq c_0(\epsilon)\;,
$$ 
if $m$ is the midpoint of the geodesic segment $[a,b]$, then
$d(m,[a',b'])\leq\frac{\epsilon}{2}$.  
\elemm

\dem %
Let $p$ be the point in $[a,b']$ the closest to $m$, and $q$ the point
of $[a',b']$ the closest to $p$. Let $t=d(a,m)=d(b,m)=d(a,b)/2$. By
Lemma \ref{lem:expopinch}, we have
$$
d(m,p)\leq e^{-d(b,m)}\sinh d(b,b')\leq e^{-t}\sinh \epsilon
$$ 
and, as $d(m,p)\leq\epsilon/2$ by convexity, $$d(p,q)\leq
e^{-d(a,p)}\sinh d(a,a')\leq e^{-d(a,m)+d(m,p)}\sinh d(a,a')\leq
e^{-t+\epsilon/2}\sinh \epsilon\;.$$ Hence $d(m,q)\leq d(m,p)+
d(p,q)\leq e^{-t}(1+e^{\epsilon/2})\sinh \epsilon $, and the result
follows by the assumption on $d(a,b)$.  
\cqfd

\medskip 
\rem 
If we want a simpler expression, we can also take
$c_0(\epsilon)=3\epsilon+4\log 2$.

%
%
%

\subsection{Entering and exiting  $\epsilon$-convex subsets}
\label{subsec:enteringexiting}

For every subset $A$ in $X$ and $\epsilon>0$, we denote by
$\N_\epsilon A$ the closed $\epsilon$-neighborhood of $A$ in $X$.  For
every $\epsilon>0$, a subset $C$ of $X$ will be called {\it
  $\epsilon$-convex} if there exists a convex subset $C'$ in $X$ such
that $C=\N_\epsilon C'$. As the metric space $X$ is CAT$(-1)$, it is
easy to see that an $\epsilon$-convex subset $C$ is closed, convex,
equal to the closure of its interior, and {\it strictly convex} in the
sense that for every geodesic line $\ga$ meeting $C$ in at least two
points, the segment $\ga\cap C$ is the closure of $\ga\;\cap
\stackrel{\circ}{C}$.  If $X$ is a smooth Riemannian manifold, then an
$\epsilon$-convex subset has a ${\rm C}^{1,1}$-smooth boundary, see
\cite{Wal}.

\medskip\noindent{\bf Examples. } %
(1) For every $\epsilon>0$, any ball of radius at least $\epsilon$ is
$\epsilon$-convex, and any horoball is $\epsilon$-con\-vex.
Conversely, as proved below, if a subset $C\subset X$ is
$\epsilon$-convex for every $\epsilon>0$, then $C$ is $X$, $\emptyset$
or a horoball.  Accordingly, we will sometimes refer to horoballs as
$\infty$-convex subsets.

To prove the above statement, assume that $C\neq X,\emptyset$ and that
for all $\epsilon>0$, there exists a convex subset $C_{-\epsilon}$ in
$X$ such that $C=\N_\epsilon C_{-\epsilon}$.  For every $x$ in
$\partial C$ (note that $\partial C$ is non empty as $C\neq
X,\emptyset$) and every $t\geq 0$, let $x_t$ be the point of the
closed convex subset $\overline{C_{-t}}$ which is the closest to $x$.
Then $t\mapsto x_t$ is a geodesic ray, which converges to a point
called $x_\infty$. We claim that $x_\infty=y_\infty$ for every $x,y$
in $\partial C$. Otherwise, the geodesic segment between $x_t$ and
$y_t$, contained in $\overline{C_{-t}}$ by convexity,  converges
to the geodesic line between $x_\infty=y_\infty$. Hence, the point
$x_t$ would not be the closest one to $x$, for $t$ big enough.
Therefore $\partial C$ is an horosphere whose point at infinity is
$x_\infty$, and by convexity, $C$ is an horoball.

(2) For manifolds, the property for a closed convex subset with
nonempty interior to be $ \epsilon$-convex is related with extrinsic 
curvature properties of its boundary.  Let us explain this
relationship.

Let $(X,\langle\cdot,\cdot\rangle)$ be a complete simply connected
smooth Riemannian manifold of dimension $m\geq 3$ with pinched
negative sectional curvature $-b^2\leq K\leq -a^2<0$. Let $C$ be a
compact convex subset of $X$ with nonempty interior, and with ${\rm
  C}^\infty$-smooth boundary $S=\partial C$. (Compactness and ${\rm
  C}^\infty$ instead of ${\rm C}^{1,1}$ are not really necessary, but
as statements and proofs are then simpler, we will work only under
these hypotheses). Let ${\rm II}_S:TS\oplus TS\ra\RR$ be the second
fundamental form of $S$ associated to the inward normal unit vector
field $\vec{n}$ along $S$, that is
$$
{\rm  II}_S(V,W)=\langle\nabla_VW,\vec{n}\rangle=
-\langle\nabla_V\vec{n},W\rangle\;,
$$
where $V,W$ are tangent vectors to $S$ at the same point, extended to
vector fields tangent to $S$ at every point of $S$ (the definition of
${\rm II}_S$ depends on the choice between $\vec{n}$ and $-\vec{n}$,
and the various references differ on that point, see for instance
\cite[page 217]{GHL}, \cite[page 36]{Pet} vs \cite[page 37]{Gra}; we
have chosen the inward pointing vector field in order for the
symmetric bilinear form ${\rm II}_S$ to be nonnegative, by convexity
of $C$).  Let $\overline{\rm II}_S$ be the upper bound of ${\rm
  II}_S(V,V)$ for every unit tangent vector $V$ to $S$.

\bprop
\begin{itemize}
\item[$\bullet$] %
If $\overline{\rm II}_S\leq a\coth(a\epsilon)$, then $C$ is
$\epsilon$-convex.
\item[$\bullet$] %
If $C$ is $\epsilon$-convex, then $\overline{\rm II}_S\leq
    b\coth(b\epsilon)$.
\end{itemize}
\eprop

\dem %
Assume first that $C$ is $\epsilon$-convex. Let $x$ be a point in $S$
and $y=\exp_x(\epsilon \,\vec{n}(x))$. Note that the sphere
$S_X(y,\epsilon)$ of center $y$ and radius $\epsilon$ in $X$ is
contained in $C$, as $C$ is $\epsilon$-convex. Locally over the
tangent space $T_xS=T_x(S_X(y,\epsilon))$, the graph of
$S_X(y,\epsilon)$ is above the graph of $S$ (when $\vec{n}$ points
upwards).  Hence, for every $V$ in $T_xS$, we have ${\rm
  II}_S(V,V)\leq {\rm II}_{S_X(y,\epsilon)}(V,V)$.  As the sectional
curvature of $X$ is at least $-b^2$, we have by comparison
$\overline{\rm II}_{S_X(y,\epsilon)}\leq b\coth(b\epsilon)$ (see for
instance \cite[page 145]{Pet}). The second result follows.  Note that
$b\coth(b\epsilon)$ is the (value on unit tangent vectors of) the
second fundamental form (with respect to the inward normal unit vector
field) of a sphere of radius $\epsilon$ in the real hyperbolic
$3$-space with constant curvature $-b^2$).

Assume now that $\overline{\rm II}_S\leq a\coth(a\epsilon)$.  For
every $t\geq 0$ and $x$ in $S$, let $x_t=\exp_x(t \,\vec{n}(x))$.
Identify by parallel transport $|\!|_{x_t}^{x}:T_{x_t}X\ra T_{x}X$
along $t\mapsto x_t$ the tangent spaces $T_{x_t}X$ with $T_{x}X$.  For
every $x$ in $S$, let $A(t)$ and $R(t)$ be the symmetric endomorphisms
of $T_{x}S$ defined by $v\mapsto -|\!|_{x_t}^{x}
\big(\nabla_{|\!|^{x_t}_{x}(v)}\dot{x}_t\big)$ and $v\mapsto
R(|\!|_{x}^{x_t}(v), \dot{x}_t) \dot{x}_t$, respectively, where $R$ is
the curvature tensor of $X$. For every $v\in T^1S$, note that $\langle
R(t)v,v\rangle$, being the curvature of the plane generated by the
orthonormal tangent vectors $||_{x}^{x_t}v$ and $\dot{x}_t$ at $x_t$,
is at most $-a^2$. It is well known that $t\mapsto A(t)$ satisfies the
following matrix Riccati equation
$$
\dot{A}(t) + A(t)^2 + R(t)=0
$$
(see for instance Theorem 3.6 on page 37 of \cite{Pet}).  The
following standard result, controlling the time of explosion of the
solution to the matrix Riccati equation, follows for instance from (a
minor modification of) Theorem 2.3 on page 144 of \cite{Pet}.

\blemm %
Let $t\mapsto A(t)$ be a smooth map from a neighborhood of $0$ in
$\RR$ to the space of symmetric matrices on an Euclidean space, such
that $\dot{A}(t) + A(t)^2 -a^2{\rm Id}$ is nonnegative, and the
biggest eigenvalue of $A(0)$ is at most $a\coth a\epsilon$. Then
$t\mapsto A(t)$ is defined and smooth at least on
$]-\epsilon,+\epsilon\,[$, and $A(t)$ is nonnegative for every $t$ in
$]-\epsilon,+\epsilon\,[$.  \cqfd \elemm

As $\langle A(t)v,v\rangle = \overline{\rm II}_S(v,v)$ for every $v$
in $T^1S$ and $\overline{\rm II}_S\leq a\coth(a\epsilon)$, it follows
that for every $t \in[0,\epsilon[$, the map $x\mapsto x_t$ is a smooth
immersion of the compact $(m-1)$-manifold $S$ into $X$, such that, if
$S_t$ is the image, then ${\rm II}_{S_t}$ is everywhere nonnegative.
As $m\geq 3$, it follows from the Hadamard-Alexander theorem (see for
instance \cite{Ale}) that $S_t$ bounds a convex subset $C_t$, with
$C_t\subset C_{t'}$ if $t'\leq t$. If $C'=\cap_{t\in[0,\epsilon[}
C_t$, then $C'$ is compact, nonempty, convex, and $C=\N_\epsilon C'$,
so that $C$ is $\epsilon$-convex.  \cqfd

\medskip %
In particular, the above result implies that if $M$ has constant
curvature $-a^2$, then $C$ is $\epsilon$-convex in our sense if and
only if $\overline{\rm II}_S\leq a\coth(a\epsilon)$.  Our notion of
$\epsilon$-convexity is hence related to, but different from, the
opposite of the notion of $\lambda$-convexity studied for instance in
\cite{GR,BGR,BM}.  Alexander and Bishop \cite{AB} (see also
\cite{Lytchak}), have introduced a natural notion of an ``extrinsic
curvature bounded from above'' for subspaces of ${\rm
  CAT}(\kappa)$-spaces, extending the notion of having a bounded
(absolute value of the) second fundamental form for submanifolds of
Riemannian manifolds.  Thus, the opposite of this concept of \cite{AB}
is related to our notion of $\epsilon$-convex subsets (see in
particular Proposition 6.1 in \cite{AB}).

\bigskip %
The rest of this section is devoted to several lemmas concerning
the relative distances between entering points and exiting points, in
and out of an $\epsilon$-convex subset of $X$, of two geodesic rays or
lines starting from the same point.  The asymptotic behaviour of the
various constants appearing in this section is described in Remark
\ref{rem:asymptotconst}.

\blemm\label{lemma:c1} %
Let $C$ be a convex subset in $X$, let $\epsilon>0$ and let $\xi_0\in
(X\cup\partial_\infty X) -(\N_\epsilon C\cup\partial_\infty C)$. If
two geodesic segments, rays or lines $\gamma,\gamma'$ which start from
$\xi_0$ intersect $\N_\epsilon C$, then the first intersection points
$x,x'$ of $\gamma,\gamma'$ respectively with $\N_\epsilon C$ are at a
distance at most
$$
c'_1(\epsilon)= 2\,\arsinh(\coth\epsilon).
$$
\elemm

\dem Let $y$ and $y'$ be the closest points in $C$ to $x$ and $x'$
respectively.  As $x,x'\in\partial\N_\epsilon[y,y']$, it is sufficient
to prove the result when $C=[y,y']$.

\medskip\noindent
\begin{minipage}{7.5cm}
  ~~~~  We may assume that $x\neq x'$,
  and, by a continuity argument, that $y\neq y'$. Let us construct a
  pentagon in $\hdr$ with vertices $\overline{\xi_0},\overline
  x,\overline y,\overline{ x'},\overline{y'}$ by gluing together the
  comparison triangles of $(\xi_0,x,x')$, $(x,x',y')$ and $(x,y',y)$.
  By comparison (see for instance \cite[Prop.~1.7.(4)]{BH}), the
  comparison angles at $\overline x,\overline y,\overline
  {x'},\overline {y'}$ are at least $\pi/2$. Hence, the segments or
  rays $]\overline {\xi_0},\overline {x}[$ and $]\overline
  {\xi_0},\overline {x'}[$ do not meet $\N_\epsilon[\overline
  y,\overline{y'}]$, and the point $\overline y$ is the closest point
  on $[\overline y,\overline{y'}]$ to $\overline x$.
\end{minipage}
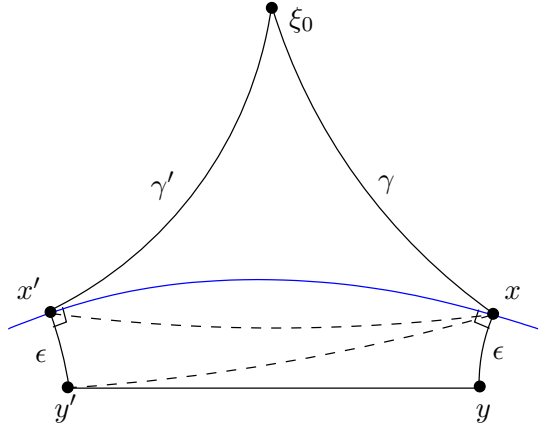
\begin{minipage}{7.3cm}
\begin{center} 
\input{fig_pentanew.pstex_t}
\end{center}
\end{minipage}
\smallskip

Furthermore, $\overline{y'}$ is the closest point on $[\overline
y,\overline{y'}]$ to $\overline{x'}$. Indeed, the angle at
$\overline{y'}$ of the pentagon is at most $3\pi/2$ since
$\angle_{\overline{y'}}(\overline y,\overline x)\leq\pi/2$ and
$\angle_{\overline{y'}}(\overline x, \overline{x'})\leq\pi$.
Therefore, if by absurd $\overline z\in[\overline y,\overline{y'}[$ is
closest to $\overline{x'}$, the geodesic segment $[\ov{x'},\ov{z}]$
intersects $[\ov{y'},\ov x]$ at a point $\ov u$.  If $z\in[y',y]$ and
$u\in[y',x]$ are such that $d(y',z)=d(\ov{y'},\ov z)$ and
$d(y',u)=d(\ov{y'},\ov u)$, then by comparison
$$
d(x',z)  \le d(x',u)+d(u,z)
 \le d(\ov{x'},\ov u)+d(\ov u,\ov
z) =d(\ov{x'},\ov{z})<d(\ov{x'},\ov{y'})=d(x',y')\;,
$$
a contradiction.

As $d(x,x')=d(\overline x,\overline {x'})$, we only have to prove that
$d(\overline x,\overline {x'})\leq c'_1(\epsilon)$, i.e.~we may assume
that $X=\hdr$.  Up to replacing $\xi_0$ by the point at infinity of
the geodesic ray starting at $x$ and passing through $\xi_0$, we may
assume that $\xi_0$ is at infinity. By homogeneity, we may assume that
$\xi_0$ is the point at infinity $\infty$ in the upper halfplane model
of $\hdr$.  As a geodesic line starting from $\infty$ and meeting the
$\epsilon$-neighborhood of a vertical geodesic segment enters it in
the sphere of radius $\epsilon$ centered at its highest point, we may
assume that $C$ is a segment of the geodesic line $\ell$ between the
points $-1$ and $1$ of the real line.  As there are geodesic lines
starting at $\infty$ whose first intersection points with
$\N_\epsilon\ell$ are at distance at least $d(x,x')$, we may assume
that $C=\ell$.
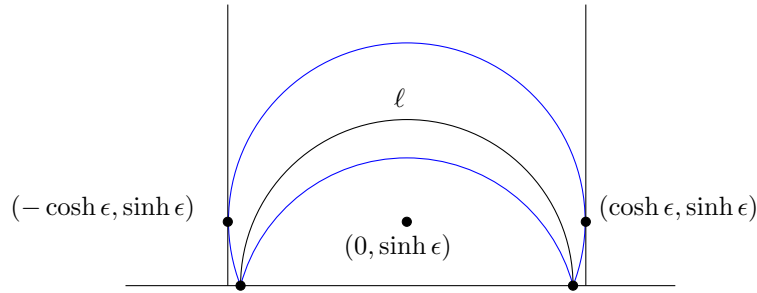
\begin{center}
\input{fig_intersectbis.pstex_t}
\end{center}

The distance $d(x,x')$ is then maximized when the geodesic lines are
tangent to $\N_\epsilon\ell$ on both sides (see the above figure). The
upper component of the boundary of $\N_\epsilon\ell$ is the
intersection with $\hdr$ of the Euclidean circle through $(\pm 1,0)$
and $(0,e^\epsilon)$, hence centered at $(0,\sinh\epsilon)$. Thus, we
may assume that the points $x,x'$ are $(\pm\cosh\epsilon,
\sinh\epsilon)$. A computation then yields the result.  
\cqfd

\bigskip
The following technical result will be used in Lemma \ref{lemma:c2}.
Define, for every $\epsilon>0$,
$$
c''(\epsilon)=
\frac{2}{\epsilon}\,{\rm arcosh}(2\cosh (\epsilon/2)) \;.
$$

\blemm\label{lem:subh0} %
For every $\epsilon>0$, for every convex subset $C$ in $X$, for 
every $a,b$ in $\N_\epsilon C$ and for every $a_0$ in $[a,b]$, if
$d(a,b)\geq c_0(\epsilon)$ and
 $$
\eta=\frac{1}{c''(\epsilon)}
\min\{d(a_0,a),d(a_0,b)\}\leq \frac{\epsilon}{2}\;,
$$  
then $d(a_0,C)\leq \epsilon-\eta$.  
\elemm

\dem %
Let $\epsilon>0$.  Let $C,a,b,a_0,\eta$ be as in the statement, and
let us prove that $d(a_0,C)\leq \epsilon-\eta$.  By an easy
computation, we have $c''(\epsilon)\epsilon\leq c_0(\epsilon)$. By
symmetry, we may assume that $d(a,a_0)\leq d(b,a_0)$, so that our
assumptions give the following inequalities:
\begin{equation}\label{eq:majomoit}
d(a,a_0)=c''(\epsilon)\eta\leq  c''(\epsilon)\epsilon/2\leq
c_0(\epsilon)/2\leq d(a,b)/2\;.
\end{equation}

Let $a',b'$ be the points in $C$ the closest to $a,b$ respectively.
As $[a',b']$ is contained in $C$, we may assume that $C=[a',b']$. Let
$m$ be the midpoint of $[a,b]$, and $m'$ its closest point on
$[a',b']$. By Lemma \ref{lem:midpointclose}, we have
$d(m,m')\leq\frac{\epsilon}{2}$.

As $\eta\leq \epsilon/2$, if $d(a,a')\leq \epsilon-\eta$, then every
point in $[a,m]$ is at distance at most $\epsilon-\eta$ from $C$. In
particular, this is true for $a_0$, since $d(a,a_0)\leq d(a,m)$.
Hence, we may assume that $d(a,a')> \epsilon-\eta$.

\medskip %
Consider the quadruple $(a,a',m,m')$ of points of $X$, which satisfies 
\begin{itemize}
\item $\epsilon-\eta<d(a,a')\leq \epsilon$, 
\item $d(m,m')\leq \frac{\epsilon}{2}$, 
\item $a'$ is the point in $[a',m']$ the closest to $a$, and 
\item $m'$ is the point in
$[a',m']$ the closest to $m$.  
\end{itemize}
Define $t=t(a,a',m,m')$ as the distance
between $a$ and the point $z=z(a,a',m,m')$ in $[a,m]$ at distance
$\epsilon-\eta$ from $[a',m']$ (which exists and is unique by
convexity), see the figure below.

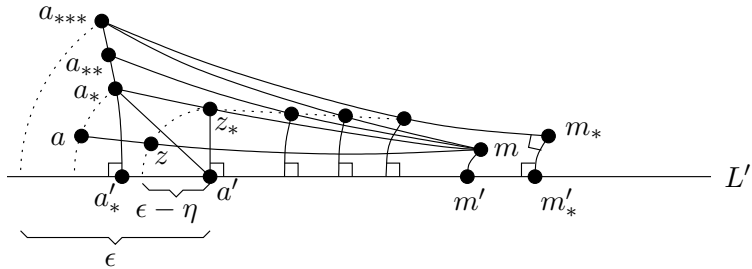
\begin{center}
\input{fig_reducnew.pstex_t}
\end{center}

\noindent We claim that $t\leq c''(\epsilon)\eta=d(a,a_0)$. Before
proving this claim, we note that it implies by Equation
(\ref{eq:majomoit}) that $t\leq d(a,a_0)\leq d(a,b)/2$, hence, by
convexity, $d(a_0,[a',m'])\leq \epsilon-\eta$, and Lemma
\ref{lem:subh0} will follow.

We will make several reductions, in order to reach
a situation where easy computations will be possible.

First we may assume, by comparison, that $X=\hdr$.  By an
approximation argument, we may assume that $a'\neq m'\neq m$. The
assumptions on the quadruple $(a,a',m,m')$ imply that the angles
$\angle_{a'}(a,m')$ and $\angle_{m'}(m,a')$ are at least
$\frac{\pi}{2}$.

If the segment $[a,m]$ cuts the segment $[a',b']$ in a point $u$, then
replacing $m$ and $m'$ by the intersection point $u$ gives a new
quadruple with the same $t$. Hence, we may assume that $a$ and $m$ are
on the same side of the geodesic line $L'$ through $a'$ and $m'$.

If $[a,m]$ does not enter $\N_{\epsilon-\eta}C$ in the sphere
$\partial B(a',\epsilon-\eta)$, then define $a_*=a$. Otherwise,
replace $a$ by the point $a_*$ at distance equal to $d(a,a')$ from
$a'$, such that the geodesic segment between $a_*$ and $m$ goes
through the point $z_*\in\partial B(a',\epsilon-\eta)
\cap\partial\N_{\epsilon-\eta}L'$ (on the same side of $L'$ as $a$).
This gives a new quadruple $(a_*,a',m,m')$ satisfying the same
properties, whose $t$ has not decreased, by convexity.  

Replace $a'$ by $a'_*$ and $a_*$ by $a_{**}$ such that
$\angle_{a'_*}(a_{**},m') =\frac\pi 2$, $d(a_{**},a'_*)=d(a,a')$, and
$a_*\in[a_{**},a'_*]$. Clearly, this does not decrease $t$. Now
replace $a_{**}$ by the point $a_{***}$ such that
$d(a_{***},a'_*)=\epsilon$ and $[a_{**},a'_*]\subset[a_{***},a'_*]$.
Let $m_*$ be the point on $\N_{\epsilon/2}C$ such that there is a
geodesic line through $a_{***}$ and $m_*$ which is tangent to
$\N_{\epsilon/2}C$ at $m_*$. Let $m'_*$ be its closest point in
$L'$. Again, the value of $t$ for the quadruple
$(a_{***},a'_{*},m_{*},m'_{*})$ has not decreased.

Hence, after these reductions, we may assume that $X=\hdr$, that the
quadrilateral $(a,a',m,m')$ has right angles at $a',m',m$, and that
$d(a,a')=2d(m,m')=\epsilon$.

\medskip %
Now, let $\ell=d(m,a)-t$ be the distance between $m$ and the point at
distance $\epsilon-\eta$ from $[a',m']$. An easy computation (see
\cite[page 157]{Bea}) shows that
$$
\cosh (t+\ell) =\frac{\sinh (\epsilon)}{\sinh(\epsilon/2)}
\;\;\;{\rm and}\;\;\; \cosh \ell =\frac{\sinh
  (\epsilon-\eta)}{\sinh(\epsilon/2)}\;.
$$ 
Consider the map $f_\epsilon:s\mapsto {\rm arcosh}\;
\frac{\sinh(\epsilon +s)} {\sinh (\epsilon/2)}$. This function is
increasing and concave on $[-\epsilon/2,0]$, with
$f_\epsilon(-\epsilon/2)=0$.  By concavity, the graph of $f_\epsilon$
on $[-\epsilon/2,0]$ is above the line passing through its endpoints
$(-\epsilon/2,0)$ and $(0,f_\epsilon(0))$.  Hence, for every $s$ in
$[0,\epsilon/2]$, we have $f_\epsilon(0)-f_\epsilon(-s)\leq
c''(\epsilon) s$.  Therefore $t = f_\epsilon(0)-f_\epsilon(-\eta)\leq
c''(\epsilon)\eta$ as $\eta\leq\epsilon/2$. This proves our claim, and
ends the proof of Lemma \ref{lem:subh0}.
\cqfd

\medskip %
Here is a finer version of Lemma \ref{lemma:c1} which shows that the
entry points of a geodesic which enters an $\epsilon$-convex set for a
long enough time and that of any nearby geodesic are close.  For every
$\epsilon>0$, we define
\begin{equation}
c'_2(\epsilon)=\max\Big\{\;c''(\epsilon)+1\;,\;
\frac{2c'_1(\epsilon)}{\epsilon}\;,\;
\sqrt{\frac{\cosh\epsilon}{\cosh\epsilon-1}}\;\frac{\sinh
  c'_1(\epsilon)}{c'_1(\epsilon)}\;\Big\}\;.  \label{eq:c2}
\end{equation}

\blemm\label{lemma:c2} %
For every $\epsilon>0$, every $\xi_0$ in $X\cup\partial_\infty X$, 
every convex subset $C$ in $X$, and all geodesic rays or lines
$\ga,\ga'$ in $X$ which start at $\xi_0$ and enter $\N_\epsilon C$ at
the points $x,x'$ in $X$ respectively, if the length of
$\ga'\cap\N_\epsilon C$ is at least $c_0(\epsilon)$, then we have
$$d(x,x')\leq c'_2(\epsilon) d(x,\ga')\;.$$
\elemm

\noindent{\bf Remarks.} (1) Without assuming that the geodesic ray or
line $\ga'$ has a sufficiently big penetration distance inside
$\N_\epsilon C$, the result is false.

(2) The curvature assumption is necessary, as can be seen by
considering geodesics which enter a half-plane in $\RR^2$ almost
parallel to the boundary.

\medskip \dem 
Let $\epsilon>0$ and assume that $\xi_0,C,\ga,\ga',x,x'$ are as in the
statement. We may assume that $x\neq x'$. In particular $\xi_0\notin
C$.  Let $p'$ be the point of $\ga'$ the closest to $x$.  Let
$[x',y']$ be the intersection of $\ga'$ with $\N_\epsilon C$ (or
$[x',y'[$ with $y'\in\partial_\infty X$ if $\ga'\cap\N_\epsilon C$ is
unbounded).

\medskip \noindent %
{\bf Case 1: } Assume first that $p'$ does not belong to $[\xi_0,x']$.
If $d(x',p')\leq \frac{\epsilon}{2}c''(\epsilon)$, then let $a_0$ be
the point $p'$.  Otherwise let $a_0$ be the point in $[x',y'[$ at
distance $\frac{\epsilon}{2}c''(\epsilon)$ from $x'$.  This point
exists and is at distance at least $\frac{\epsilon}{2}
c''(\epsilon)\geq d(a_0,x')$ from $y'$, as $d(x',y') \geq
c_0(\epsilon)\geq \epsilon\,c''(\epsilon)$. By Lemma \ref{lem:subh0},
we have $d(a_0,C)\leq \epsilon-\frac{1}{c''(\epsilon)}d(a_0,x')$.

Hence, if $a_0=p'$, then, as $d(x,C)=\epsilon$, we have $d(x,p')\geq
\frac{1}{c''(\epsilon)}d(p',x')$.  So that
$$
d(x,x')\leq d(x,p')+d(p',x')\leq (1+c''(\epsilon))\;d(x,p')\;,
$$
which proves the result, by the definition of $c'_2(\epsilon)$.

If $a_0\neq p'$, then $p'\notin [a_0,\xi_0[$. Let us prove that
$d(x,p')\geq \frac{\epsilon}{2}$. This implies, by Lemma
\ref{lemma:c1}, that
$$
d(x,x')\leq c'_1(\epsilon)\leq
\frac{2c'_1(\epsilon)} {\epsilon}\;d(p',x)\;,
$$ 
which proves the result, by the definition of $c'_2(\epsilon)$.
Let $b_0$ be the point in $[x',y']$ at distance $\frac{\epsilon}{2}
c''(\epsilon)$ from $y'$ (or $b_0=y'$ if $y'$ is at infinity).  By
Lemma \ref{lem:subh0}, we have
$$
\max\{d(a_0,C),d(b_0,C)\} \leq\epsilon-\frac{1}{c''(\epsilon)} 
\min\{d(a_0,x'),d(b_0,y')\}= \frac{\epsilon}{2} \;.
$$ 
If $p'\in [a_0,b_0]$, then by convexity $d(p',C)\leq \frac
{\epsilon}{2}$.  As $d(x,C)=\epsilon$, this implies that $d(x,p')
\geq\frac{\epsilon}{2}$, as wanted. If otherwise $p'\notin [a_0,b_0]$,
then assume by absurd that $d(x,p')< \frac{\epsilon}{2}$.  Let $z$ be
the point in $[x,\xi_0[$ whose closest point to $[p',\xi_0[$ is $b_0$.
By convexity, $d(z,b_0)< \frac{\epsilon}{2}$. Hence $d(z,C)\leq
d(z,b_0)+d(b_0,C) <\epsilon$, which contradicts the fact that $\ga$
enters $\N_\epsilon C$ at $x$.

\medskip\noindent %
{\bf Case 2: } Assume now that $p'$ belongs to $[\xi_0,x']$. Let $a'$
and $b'$ be the points of $C$ the closest to $x'$ and $y'$
respectively. They are at distance $\epsilon>0$ from $x'$ and $y'$
respectively (except that $b'=y'$ if $y'$ is at infinity). Let $\phi$
be the comparison angle at $x'$ between the geodesic segments
$[x',a']$ and $[x',y'[$. We claim that $\sin{\phi}\le
\frac{1}{\sqrt{\cosh \epsilon}}$.

To prove this claim, if $y'\in X$, we construct a comparison
quadrilateral with vertices $\overline{x'},\overline{a'},
\overline{b'}, \overline{y'} \in\hdr$ by gluing together the
comparison triangles $(\overline{x'},\overline{a'},\overline{y'})$ of
$({x'},{a'},{y'})$ and $(\overline{a'},\overline{b'},\overline{y'})$
of $({a'},{b'}, {y'})$ along their isometric edges $[\overline{a'},
\overline{y'}]$.  If $y'\notin X$, then $b'=y'$, and the above
quadrilateral is replaced by the comparison triangle with vertices
$\overline{x'},\overline{a'}, \overline{y'} \in \hdr \cup \{\infty\}$.
By comparison, all angles in the quadrilateral in $\hdr$ with vertices
$\overline{x'}, \overline{a'}, \overline{b'},\overline{y'}$ are
greater than or equal to those in the quadrilateral in $X$ with
vertices $x',a',b',y'$. In particular, if the angle at $\overline{x'}$
is $\overline\phi$, we have $\phi\le\overline{\phi}$.  If the
quadrilateral with vertices $\overline{x'}, \overline{a'},
\overline{b'}, \overline{y'}$ is replaced by the one with vertices
$\overline{x'},\overline{a'}_*, \overline{b'}_*, \overline{y'}$ with
$d(\overline{x'}, \overline{a'}_*) = \epsilon =
d(\overline{y'},\overline{b'}_*)$ and right angles at
$\overline{a'}_*$ and $\overline{b'}_*$, the angle $\overline{\phi}_*$
at $\overline{x'}$ of this quadrilateral is at least $\overline\phi$.
Furthermore, this quadrilateral is symmetric: the angle at
$\overline{y'}$ is also $\overline{\phi}_*$. Thus, we get an upper
bound for $\phi$ by estimating $\overline{\phi}_*$.

Let $[m,m']$ be the common perpendicular segment between
$[\overline{x'}, \overline{y'}]$ and $[\overline{a'}_*,
\overline{b'}_*]$, with $m\in[\overline{x'}, \overline{y'}]$. We have
(see for instance \cite[page 157]{Bea}),
$$
\sin\ov{\phi}_* =\frac{\cosh d(m,m')}{\cosh \epsilon}
\;\;\;{\rm and}\;\;\; \cosh d(\overline{x'},{m})= \frac{\sinh
  \epsilon}{\sinh d(m,m')}\;.
$$
Hence, as $d(x',y')\geq c_0(\epsilon)$,
{\small
$$\sin\ov{\phi}_*= \frac{\sqrt{1+(\sinh^2 \epsilon)/
(\cosh^2 d(\overline{x'},{m}))}}{\cosh\epsilon}
 \leq \frac{\sqrt{1+(\sinh^2 \epsilon)/(\cosh^2
    (c_0(\epsilon)/2))}} {\cosh\epsilon}
\leq \frac{1}{\sqrt{\cosh \epsilon}}\;,
$$}
as $c_0(\epsilon)\ge \epsilon c''(\epsilon)\ge 2\;\arcosh(\sqrt
2\cosh(\epsilon/2))$.  This proves the claim.

\medskip %
By convexity, the comparison angle at $x'$ between the geodesic
segments $[x',x]$ and $[x',a']$ is at most $\frac{\pi}{2}$.  Hence the
comparison angle $\theta$ at $x'$ between $[x',x]$ and $[x',\xi_0[$ is
at least $\pi-\frac{\pi}{2} -\phi= \frac{\pi}{2} -\phi$. In
particular,
$$
\frac{1}{\sin \theta}\leq
\frac{1}{\sin(\frac{\pi}{2}-\phi)}
=\frac{1}{\sqrt{1-\sin^2\phi}}\leq
\sqrt{\frac{\cosh\epsilon}{\cosh\epsilon-1}}\;.
$$ 
With $\overline{x},\overline{x'}$ as above, consider
$\overline{\xi_0}$ in $\hdr\cup\partial_\infty\hdr$ such that
$\angle_{\overline{x'}}(\overline{x},\overline{\xi_0})=\theta$ 
and $d(\overline{x'},\overline{\xi_0})=d({x'},{\xi_0})$.  By
comparison, the point $\overline{p'}$ on $[\,\overline{x'},
\overline{\xi_0}\,]$ the closest to $\overline{x}$ is at distance from
$\overline{x}$ at most equal to $d(x,p')$. If $\overline{p'}=
\overline{x'}$, then $d(x,p')\geq d(\overline{x},\overline{p'})=
d(x,x')$, which implies the result, as $c'_2(\epsilon)\geq 1$.
Otherwise, the angle $\angle_{\overline{p'}} (\overline{x},
\overline{x'})$ is at least $\frac{\pi}{2}$ (equality holds if
$\overline{p'}\neq \overline{\xi_0}$). By the formulae in right-angled
hyperbolic triangles, we have $\sinh d(\overline{x},\overline{p'})\geq
\sinh d(\overline{x},\overline{x'}) \sin\theta$.

As closest point maps do not increase distances, we have $d(x,p')\leq
d(x,x')\leq c'_1(\epsilon)$. In particular $$\sinh d(x,p') \leq
\frac{\sinh c'_1(\epsilon)}{c'_1(\epsilon)} \;d(x,p')$$ by convexity
of the map $t\mapsto \sinh t$ on $[0,+\infty[\,$.  Hence
$$
d(x,x')=d(\overline{x},\overline{x'})\leq \sinh
d(\overline{x},\overline{x'}) \leq \frac{\sinh
  d(\overline{x},\overline{p'})}{\sin \theta} \leq \frac{\sinh
  d({x},{p'})}{\sin \theta} \leq c'_2(\epsilon)\; d(x,p') \;,
$$ 
by the definition of $c'_2(\epsilon)$.  
\cqfd

\medskip %
In general, there is no estimate analogous to Lemma \ref{lemma:c2} for
the distance between the points $y,y'$ where two geodesic rays or
lines $\ga,\ga'$ starting from a point $\xi_0$ exit an
$\epsilon$-convex subset $\N_\epsilon C$.  For instance, the geodesic
line $\ga$ could be tangent to $\N_\epsilon C$, and $\ga'$ could enter
for a long time in $\N_\epsilon C$, so that $y$ and $y'$ would not be
close.  But the result is not true even if we assume that both $\ga$
and $\ga'$ meet $\N_\epsilon C$ in a long segment. Here is a
counterexample when $X$ is a tree (but this phenomenon is not specific
to trees).

\medskip\noindent
\begin{minipage}{7.8cm}
  ~~~~ Let $\ga,\ga'$ be two geodesic lines in a tree $X$, coinciding
  on their negative subrays, starting at $\xi_0\in\partial_\infty X$,
  and with disjoint positive subrays. Let $\epsilon=\eta=1$, and
  $C=\ga'([-\ell,+\ell])$. Then the entering points of $\ga,\ga'$ in
  $\N_\epsilon C$ are $x=x'=\ga'(-\ell-1)$.  Besides, $y=\ga(1)$,
  $y'=\ga'(\ell+1)$ and $d(y,\ga')\leq 1$.  But we have $d(y,y')=\ell
  +2$, which goes to $+\infty$ as $\ell\ra+\infty$.
\end{minipage}
\begin{minipage}{7cm}
\begin{center} 
\input{fig_contrextree.pstex_t}
\end{center}
\end{minipage}

\medskip %
This explains the dichotomy in the following result on the exiting
points from an $\epsilon$-convex sets of two geodesic lines which
start from the same point at infinity.  For every $\epsilon,\eta>0$,
we define
\begin{equation}
h'(\epsilon,\eta)= \max\big\{\;
2\eta+\max\{0,-2\log\frac{\epsilon}{2}\}\;,\;
\eta+c'_1(\epsilon)+c_0(\epsilon) 
\;\big\}  \label{eq:h}
\end{equation}
and
\begin{equation}
c'_3(\epsilon)= 3+\frac{2 c_1'(\epsilon)}{\epsilon}\;. \label{eq:c3}
\end{equation}

\blemm\label{lem:h0} %
Let $\epsilon,\eta>0$. Let $C$ be a convex subset in $X$, %
$\xi_0\in X\cup\partial_\infty X$, and $\ga,\ga'$ geodesic rays or
lines starting from $\xi_0$. If $\ga$ enters $\N_\epsilon C$ at a
point $x\in X$ and exits $\N_\epsilon C$ at a point $y\in X$ such that
$d(x,y)\geq h'(\epsilon,\eta)$ and $d(y,\ga')\leq \eta$, then $\ga'$
meets $\N_\epsilon C$, entering it at a point $x'\in X$, exiting it at
a point $y'\in X\cup\partial_\infty X$ such that
$$
d(y,y')\leq c'_3(\epsilon)d(y,\ga')\;\;{\rm or} \;\; 
d(x',y')> d(x,y)\;.
$$ 
\elemm

\dem %
Let $p'$ be the closest point on $\ga'$ to $y$. Let $q$ be the closest
point on $\ga$ to $p'$. The point $q$ belongs to  $[y,\xi_0]$ and
satisfies $d(y,q)\leq d(y,p')\leq \eta$, as closest point maps do not
increase the distances.  By the properties of 
geodesic triangles in CAT$(-1)$ spaces, we have
$$
d(p',q)\leq {\rm arsinh~} 1=\log(1+\sqrt{2}).
$$

Let us first prove that $\ga'$ meets $\N_\epsilon C$. Let $m$ be the
midpoint of $[x,y]$. As 
$$
d(y,m)=d(x,y)/2\geq h'(\epsilon,\eta)/2\geq
\eta\geq d(y,q)\;,
$$ 
the point $q$ belongs to $[m,y]$.  Furthermore,
$$
d(q,m)=d(y,m)-d(y,q)\geq h'(\epsilon,\eta)/2- \eta\geq
-\log\frac{\epsilon}{2}\;,
$$ 
by the definition of $h'(\epsilon,\eta)$.
By Lemma \ref{lem:expopinch}, we have 
$$
d(m,\ga')\leq e^{-d(q,m)} \sinh
d(q,p')\leq \frac{\epsilon}{2}\;.
$$  
By Lemma \ref{lem:midpointclose}, as $d(x,y)\geq h'(\epsilon,\eta)
\geq c_0(\epsilon)$ by the definition of $h'(\epsilon,\eta)$, we have
$d(m,C)\leq\frac{\epsilon}{2}$.  Hence the point $m'$ of $\ga'$ the
closest to $m$ belongs to $\N_\epsilon C$, which is what we wanted.

Let $x'$ and $y'$ be the entering point in $\N_\epsilon C$ and exiting
point out of $\N_\epsilon C$ of $\ga'$ respectively, where $y'$ could
for the moment be at infinity, in which case the second possibility
below would hold.

\medskip
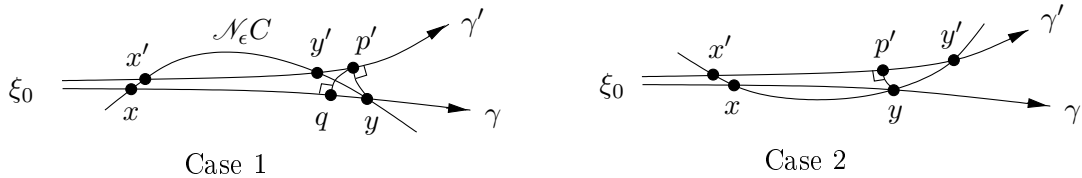
\begin{center}
\input{fig_etudcaslemh0.pstex_t}
\end{center}

\medskip\noindent {\bf Case 1 :} Assume that $p'\notin[y',\xi_0]$.
Let $\eta_\epsilon=\epsilon \,c''(\epsilon)/2$.  There are two
subcases. First assume that $d(y,p')\geq \epsilon/2$. Let 
$$
t_\epsilon=\max\big\{\eta_\epsilon, -\log \frac{\epsilon}{2}\big\}.
$$  
Note that
$h'(\epsilon,\eta)\geq \eta_\epsilon+t_\epsilon+\eta$ by the
definition of $h'(\epsilon,\eta)$, as $c_0(\epsilon)\geq \epsilon
\,c''(\epsilon)=2\eta_\epsilon$. Hence we have
$$
d(y,x)-d(y,q)-t_\epsilon\geq h'(\epsilon, \eta)
-\eta-t_\epsilon\geq\eta_\epsilon\geq 0\;.
$$ 
Therefore, the point $y_0$ in $[x,q]$ at distance $t_\epsilon$ of $q$
exists and satisfies $d(y_0,x)\geq \eta_ \epsilon$. Furthermore,
$d(y_0,q)=t_\epsilon\geq \eta_\epsilon$ and
$$
d(x,q)=d(x,y)-d(y,q)\geq h'(\epsilon, \eta) -\eta \geq
c_0(\epsilon)\geq 2\eta_\epsilon\;,
$$ 
by the definition of $h'(\epsilon,\eta)$. Let $a_0$ and $b_0$ be the
points in $[x,q]$ at distance $\eta_\epsilon$ from $x$ and $q$
respectively, which are at distance at least $\eta_\epsilon$ from $q$
and $x$ respectively. By Lemma \ref{lem:subh0}, we have $d(a_0,C)\leq
\epsilon- \eta_\epsilon/ c''(\epsilon) = \epsilon/2$, and similarly
$d(b_0,C)\leq \epsilon/2$. Note that $y_0$ belongs to $[a_0,b_0]$.
Hence by convexity, we have $d(y_0,C)\leq \epsilon/2$. By Lemma
\ref{lem:expopinch}, we have
$$
d(y_0,\ga')\leq e^{-t_\epsilon}\sinh 
d(q,p')\leq \frac{\epsilon}{2}\;.
$$
Therefore the point $q'$ on $\ga'$ the closest to $y_0$ belongs to
$\N_\epsilon C$.  As $y'$ is the exiting point of $\ga'$ from
$\N_\epsilon C$, it belongs to $[q',p']$. As closest point maps do not
increase the distances, we have $d(p',q')\leq d(y,y_0)$.  Hence
\begin{align*}
d(y,y') & \leq d(y,p')+d(p',y')\leq d(y,p')+d(p',q')\leq 
d(y,p')+d(y,y_0) \\
& \leq d(y,p')+d(y,q)+d(q,y_0) \leq 2d(y,p')+t_\epsilon \\
& \leq (2+2t_\epsilon/\epsilon)d(y,p')\le c'_3(\epsilon)\,d(y,p')\;,
\end{align*}
as it can be checked that $2c'_1(\epsilon)/\epsilon +1 \geq
2t_\epsilon/\epsilon$.

Assume now that $d(y,p')\leq \epsilon/2$. Since
$$
d(x,y)\geq h'(\epsilon,\eta)\geq c_0(\epsilon)\geq
2\eta_\epsilon\geq 2c''(\epsilon)d(y,p')\;,
$$ 
the point $y_0$ in $[x,y]$ at distance $c''(\epsilon)d(y,p')$ from $y$
exists and $d(y_0,x)\geq c''(\epsilon)d(y,p')$. Hence by Lemma
\ref{lem:subh0}, we have $d(y_0,C)\leq \epsilon-d(y,p')$. Let $q'$ be
the point on $\ga'$ the closest to $y_0$.  By convexity, $q'$ is at
distance at most $d(y,p')$ from $y_0$, hence belongs to $\N_\epsilon
C$. As $y'$ is the exiting point of $\ga'$ from $\N_\epsilon C$, it
belongs to $[q',p']$. As closest point maps do not increase distances,
we have $d(q',p')\leq d(y_0,y)$. Hence, as above,
$$
d(y,y')\leq  
d(y,p')+d(y_0,y)\leq \big(1+c''(\epsilon)\big)\:d(y,p')\;,
$$ 
which proves the result, by the definition of $c'_3(\epsilon)$,
as $2c'_1(\epsilon)/\epsilon +1 \geq c''(\epsilon)$.

\bigskip
\noindent 
{\bf Case 2 :} Assume that $p'\in\;]y',\xi_0]$. Lemma \ref{lemma:c1}
implies that $d(x,x')\le c_1'(\epsilon)$. Note that $p'\notin
[x',\xi_0]$. Otherwise, with $q$ and $s$ the closest points to $p'$
and $x'$ on $\ga$ respectively, we would have $s\notin\;]q,\xi_0]$ by
convexity. As $q\in [y,\xi_0]$, we would then have 
$$
d(x,y)\leq d(x,s)+d(q,y)\leq d(x,x')+d(p',y)\leq 
c'_1(\epsilon)+\eta< h'(\epsilon,\eta)\;,
$$ 
by the definition of $h'(\epsilon,\eta)$, a contradiction.

\medskip
Assume first that $d(y,p')< \epsilon/2$. We start by observing that
$d(p',y')\le c''(\epsilon)d(y,p')$. Indeed, suppose by absurd that
$d(p',y')>c''(\epsilon)d(y,p')$. By continuity of the closest point
maps, let $y_0$ be a point on $\ga$ that does not belong to
$\N_\epsilon C$, but is close enough to $y$, so that the closest point
$q'$ to $y_0$ on $\ga'$ belongs to $[p',y']$ and satisfies
$d(y_0,q')\leq\epsilon/2$ and $d(q',y')\geq c''(\epsilon)d(q',y_0)$.
Hence, using the definition of $h'(\epsilon,\eta)$, we have
\begin{equation}
\label{formula}
\begin{aligned}
d(y',x') & \geq d(q',x')\geq d(p',x')\geq  d(x,y)-d(p',y)-d(x,x')\\
& \geq h'(\epsilon,\eta)-\eta-c'_1(\epsilon)
\geq c_0(\epsilon)\geq \epsilon\,c''(\epsilon) \geq
2c''(\epsilon)d(q',y_0)\;.
\end{aligned}
\end{equation}
Let $a_0$ and $b_0$ be the points in $[x',y']$ at distance
$c''(\epsilon)d(q',y_0)\leq\epsilon\,c''(\epsilon)/2 $ from $x'$ and
$y'$ respectively. The estimate \eqref{formula} implies that $a_0$ and
$b_0$ are at distance at least $c''(\epsilon) d(q',y_0)$ from $y'$ and
$x'$ respectively. By Lemma \ref{lem:subh0}, we have $d(a_0,C)\leq
\epsilon -d(q',y_0)$ and $d(b_0,C)\leq \epsilon -d(q',y_0)$. Hence,
the point $q'$, which belongs to $[a_0,b_0]$ by Formula
\eqref{formula} and the construction of $q'$, is by convexity at
distance at most $\epsilon-d(q',y_0)$ from $C$.  Therefore by the
triangular inequality, $d(y_0,C)\leq\epsilon$, which is a
contradiction.  Hence $d(p',y')\leq c''(\epsilon) d(y,p')$, and
$$
d(y,y')\leq d(y,p')+d(p',y')\leq (1+c''(\epsilon))d(y,p')\;,
$$
which proves the result, as in Case 1.

\medskip %
Assume now that $d(y,p')\geq \epsilon/2$.  Suppose first that
$d(p',y')> d(y,p')+c_1'(\epsilon)$.  Then, as $p'\in[x',y']$,
$$
d(x',y')=d(x',p')+d(p',y')\geq d(p',y')+d(x,y) -d(y,p')-d(x,x')> 
d(x,y)\;,
$$ 
which is one of the two possible conclusions. 
Otherwise, 
$$
d(y,y')\leq d(y,p')+d(p',y')\leq 2d(y,p') +c_1'(\epsilon)\leq
c'_3(\epsilon)d(y,p')\;, 
$$ 
by the definition of $c'_3(\epsilon)$. This is the other
possible conclusion.
\cqfd

\bigskip 
\brema \label{rem:asymptotconst}
{\rm The asymptotic behaviour of the constants when $\epsilon$ is very
big or very small is as follows.
\begin{itemize}
\item 
$c_0(\epsilon)\sim 3\epsilon$ as $\epsilon\ra +\infty$ and
$\lim_{\epsilon\ra 0} c_0(\epsilon)= 4\log 2\approx 2.77$.
\item
$\lim_{\epsilon\to+\infty}c'_1(\epsilon)=c'_1(\infty)=
2\log(1+\sqrt2)\approx 1.76, $ and $c'_1(\epsilon)\sim-2\log\epsilon$
as $\epsilon\ra 0$. Note that $\epsilon\mapsto c'_1(\epsilon)$ is
decreasing.
\item
$\lim_{\epsilon\ra +\infty} c''(\epsilon)=1$, and $c''(\epsilon)
\sim \frac{2}{\epsilon}\log(2+\sqrt3)$ as $\epsilon\ra 0$.  
\item 
For $\epsilon$ big, $c'_2(\epsilon)=c''(\epsilon)+1$, hence
$\lim_{\epsilon\to+\infty}c'_2(\epsilon)=2$.  For $\epsilon>0$ small,
$$
c'_2(\epsilon)=\sqrt{\frac{\cosh\epsilon}{\cosh\epsilon-1}}\;
\frac{\sinh c'_1(\epsilon)}{c'_1(\epsilon)}\sim 
\frac{\sqrt2} {4\epsilon^3\log(1/\epsilon)}\;.
$$
\item 
$\lim _{\epsilon\ra +\infty} c'_3(\epsilon)= 3$, and $c'_3(\epsilon)
\sim -\frac{4}{\epsilon}\log\epsilon$ as $\epsilon\ra 0$.
\item 
$h'(\epsilon,\eta) \sim 3\epsilon$ as $\epsilon\ra+\infty$, and
$h'(\epsilon,\eta)\sim -2\log\epsilon$ as $\epsilon\ra 0$, uniformly
on compact subsets of $\eta$'s.
\end{itemize}
When $\epsilon$ goes to $+\infty$, 
$c'_1(\epsilon)$  and $c'_3(\epsilon)$ have finite limits, and the
limiting values apply for the horoball case, see Lemmas
\ref{lem:c1pinfty} and  \ref{lem:c3pinfty} below.
On the other hand,  the constants
$c_0(\epsilon)$ and  $h'(\epsilon,\eta)$ behave badly as
$\epsilon\to\infty$, and we will improve them in Section
\ref{subsec:hittinghoroball}.  }
\erema

When $X$ is a tree, the constants $c'_3(\epsilon)$ and
$h'(\epsilon,\eta)$ can be simplified, we can take $c'_3(\epsilon)=2$
and any $h'(\epsilon,\eta) >2\eta$, as the following more precise
result shows, improving Lemma \ref{lem:h0} for trees. Note that the
versions of Lemmas \ref{lemma:c1} and \ref{lemma:c2} for trees simply
say that we can take $c_0(\epsilon)=\epsilon$, and $c'_1(\epsilon)
=c'_2(\epsilon)=0$, since for every point or end $\xi_0$ of a (real)
tree, for every convex subset $C$, for all geodesic rays or lines
$\ga,\ga''$ starting from $\xi_0$ and entering $C$ in $x,x'$
respectively, we have $x=x'$.

\brema %
Let $X$ be an $\RR$-tree and $\epsilon>0$. Let $C$ be a convex subset 
in $X$, $\xi_0\in X\cup \partial_\infty X$,
and $\ga,\ga'$ geodesic rays or lines starting from $\xi_0$. If $\ga$
enters $\N_\epsilon C$ at a point $x\in X$ and exits $\N_\epsilon C$
at a point $y\in X$ such that $d(x,y)>2d(y,\ga')$, then $\ga'$ meets
$\N_\epsilon C$, entering it at $x'=x$, exiting it at a point $y'$
(possibly at infinity) such that
$$
d(y,y')\leq 2\;d(y,\ga')\;\;{\rm or} \;\; d(x',y')> d(x,y).
$$
\erema

\dem %
Let $p'$ be the closest point to $y$ on $\ga'$. Note that $p'$ belongs
to $]\xi_0,y]$, as $X$ is a tree and $\ga'$ also starts from $\xi_0$.
If $p'\in\;]\xi_0,x[$, then $d(y,\ga')> d(x,y)$, a contradiction.
Hence $p'\in[x,y]\subset \N_\epsilon C$, and $\ga'$ enters
$\N_\epsilon C$ at $x'=x$.

Suppose first that $d(x,y)<2\epsilon$. Then the closest point $z$ to
$y$ in $C$ does not belong to $[x,y]$. Let $q$ be the midpoint of
$[x,y]$, which is also the closest point to $z$ on $[x,y]$. As
$d(x,y)>2d(y,\ga')$, the point $p'$ belongs to $]q,y]$, hence
$d(y,y')=2d(y,\ga')$, which is fine.

Assume now that $d(x,y)\geq 2\epsilon$. If $x_\epsilon$ and
$y_\epsilon$ are the points in $[x,y]$ at distance $\epsilon$ from $x$
and $y$ respectively, then $[x,y]\cap C=[x_\epsilon, y_\epsilon]$.  If
$p'$ belongs to $]y_\epsilon, y]$, then $y_\epsilon$ is also the
closest point to $y'$ in $C$, and $d(p',y)=d(p',y')$, so that
$d(y,y')= 2d(y,\ga')$, which is fine. Otherwise, we have $d(y,\ga')
\geq \epsilon$. If $d(x',y')\leq d(x,y)$, then $d(p',y')\leq
d(p',y)$. Hence
$$
d(y,y')= d(y,p')+ d(p',y')\leq 2d(y,p')= 2d(y,\ga')
\;.\;\;\;\mbox{\cqfd}
$$

%
%
%

\subsection{Hitting horoballs}
\label{subsec:hittinghoroball}

As shown in Remark \ref{rem:asymptotconst}, the constants
$c_0(\epsilon)$ and $h'(\epsilon,\eta)$, used to describe the
penetration of geodesic lines inside $\epsilon$-convex subsets, do not
have a finite limit as $\epsilon$ goes to $+\infty$. Horoballs are
$\epsilon$-convex subsets for every $\epsilon$, and we could use for
instance $\epsilon=1$ in these constants to get numerical values.  But
in order to get better values, we will prove analogs for horoballs of
the lemmas \ref{lemma:c1},  \ref{lem:subh0}, \ref{lemma:c2} and
\ref{lem:h0}. The proofs of the lemmas below follow the same lines
as the ones for the general case of
$\epsilon$-convex subsets given in Section
\ref{subsec:enteringexiting}, with many simplifications.

\medskip %
As $c'_1(\epsilon)$ tends to 
$$
c'_1(\infty)= 2\log(1+\sqrt{2}),
$$ 
the next lemma follows by passing to the limit in Lemma
\ref{lemma:c1}. It is not hard to see (for instance by considering  the
real hyperbolic plane) that the constant $c'_1(\infty)$ is
optimal.

\blemm \label{lem:c1pinfty} 
For every horoball $H$ in $X$, for every $\xi_0$ in $(X\cup
\partial_\infty X) -(H\cup H[\infty])$, for all geodesic rays or lines
$\ga$ and $\ga'$ starting from $\xi_0$ and entering $H$ in $x$ and
$x'$ respectively, we have
$$
d(x,x')\leq c'_1(\infty)=2\log(1+\sqrt{2})\;.\;\;\;\;\mbox{\cqfd}
$$
\elemm

The following result, Lemma \ref{lem:cppinfty}, improves Lemma
\ref{lem:subh0} for horoballs, and says that when the
$\epsilon$-convex subset under consideration is a horoball, we can
replace $c_0(\epsilon)$ by
\begin{equation}
c_0(\infty)=4.056\;, \label{eq:c0i}
\end{equation}
and $c''(\epsilon)$ by $c''(\infty) =\frac{3}{2}$. Lemma
\ref{lem:c2pinfty} below is the analog of Lemma \ref{lemma:c2} for
horoballs, and says that when the $\epsilon$-convex subset under
consideration is a horoball, we can replace $c_0(\epsilon)$ by
$c_0(\infty)=4.056$ and $c'_2(\epsilon)$ by
\begin{equation}
c'_2(\infty)=\frac{5}{2}\;. \label{eq:c2i}
\end{equation}
Note that $c''(\infty), c_0(\infty)$ and $c'_2(\infty)$ are not limits
as $\epsilon$ goes to $\infty$ of $c_0(\epsilon)$ and
$c'_2(\epsilon)$, but this notation will be useful in Section
\ref{sec:main}.

\blemm \label{lem:cppinfty} 
For every horoball $H$, for every $a$ and $b$ in $\partial H$ with 
$d(a,b)\geq c_0(\infty)$, for every $a_0$ in $[a,b]$, we have
$$a_0\in H[\,\mbox{$\frac{2}{3}$}\min \{d(a_0,a),d(a_0,b)\}\,]\;.$$
\elemm

\dem 
Let $\xi=H[\infty]$ be the point at infinity of $H$.  By symmetry, we
may assume that $\ell=d(a_0,a)=\min \{d(a_0,a),d(a_0,b)\}$.

Let $(\ov a, \ov b,\ov{\xi}=\infty)$ be a comparison triangle of
$(a,b,\xi)$ in $\hdr$. By comparison, the difference $\ell'$ of the
heights of $a_0$ and $a$ with respect to $\xi$ is bigger than the
corresponding quantity $\ov{\ell '}$ for the comparison points
$\ov{a_0}$ and $\ov a$.  Thus, in order to show that $\ell'\ge
\frac{2}{3}\ell$, it suffices to show that $\ov{\ell'}\ge
\frac{2}{3}\ell$, and the question reduces to the case $X=\hdr$.

\medskip
\noindent
\begin{minipage}{8.8cm}
  ~~~ We assume that $[b,a]$ lies on the unit circle, with $a$ (hence
  $a_0$, as $a$ and $b$ have the same (Euclidean) vertical coordinate)
  in the closed positive quadrant. Let $s$ be the (Euclidean) vertical
  coordinate of $a_0$ and $t$ the one of $a$, with $0<t\leq s\leq 1$.
  An easy computation in hyperbolic geometry (see also the proof of
  Lemma \ref{lem:expopinch}) gives $ \ell'=\log \frac{s}{t}$ and
$$\ell={\rm ~arsinh~}\frac{\sqrt{1-t^2}}{t}-
{\rm ~arsinh~}\frac{\sqrt{1-s^2}}{s}$$
$$= \log \frac{s}{t}+ \log
\frac{1+\sqrt{1-t^2}}{1+\sqrt{1-s^2}}\;.$$
\end{minipage}
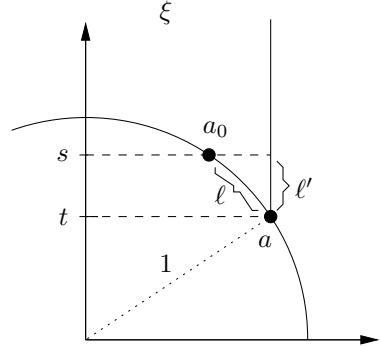
\begin{minipage}{6cm}
\begin{center}
\input{fig_subhzerhoronew.pstex_t}
\end{center}
\end{minipage}
\medskip

Hence, to prove that $\ell\leq \frac{3}{2}\ell'$, we only have to show
that $\log \frac{1+\sqrt{1-t^2}}{1+\sqrt{1-s^2}}\leq \frac{1}{2}\log
\frac{s}{t}$, which is equivalent to $\sqrt{t}(1+\sqrt{1-t^2})\leq
\sqrt{s}(1+\sqrt{1-s^2})$. The map $f:x\mapsto\sqrt{x}(1+
\sqrt{1-x^2})$ on $[0,1]$ is increasing from $f(0)=0$ to 
$f(\frac{\sqrt{5}}{3})$, and then decreasing to $f(1)=1$.  Let
$t'=0.25873$. As $f(t')<1$ and $s\geq t$, to prove that $f(t)\leq
f(s)$, it is sufficient to show that $t\leq t'$.  Let $a'$ and $b'$ be
the two points of the unit circle at (Euclidean) height $t'$. As $a$
and $b$ are at the same (Euclidean) height $t$ on the unit circle, to
prove that $t\leq t'$, we only have to show that $d(a',b')\leq
d(a,b)$.  By the definition of $c_0(\infty)$, we have
$$
d(a',b')= 2\,{\rm arsinh}\; \frac{\sqrt{1-{t'}^2}}{t'}\leq
c_0(\infty)\leq d(a,b)\;.   
$$
Hence  the result follows. 
\cqfd

\blemm \label{lem:c2pinfty} 
For every horoball $H$ in $X$, for every $\xi_0$ in  
$X\cup \partial_\infty X$, for all geodesic rays or lines $\ga$ and
$\ga'$ starting from $\xi_0$ and entering $H$ in $x\in X$ and $x'\in
X$ respectively, if the length of $\ga'\cap H$ is at least
$c_0(\infty)$, then
$$
d(x,x')\leq \frac{5}{2}\,d(x,\ga')\;.
$$
\elemm

\dem 
Let $p'$ be the point of $\ga'$ the closest to $x$. Let $\xi$ be the
point at infinity of $H$. Define $y'$ by $[x',y']=\ga'\cap H$ if this
intersection is bounded, and $y'=\xi$ otherwise. We may assume that
$x\neq x'$. In particular, $\xi_0\notin H\cup\{\xi\}$.

Assume first that $p'$ does not belong to $[\xi_0,x'[$. As closest
point projections do not increase distances and by Lemma
\ref{lem:c1pinfty}, we have $d(x',p')\leq c'_1(\infty)$, and since
$$
d(x',y')\geq c_0(\infty)\geq 2c'_1(\infty),
$$
the point $p'$ belongs to $H$, and $d(p',y')\geq d(p',x')$. Let $z$
be the point of intersection of $]\xi,x']$ with the horosphere
centered at $\xi$ passing through $p'$, so that in particular
$d(x,p')\geq d(x',z)$. By Lemma \ref{lem:cppinfty}, we have
$d(x',z)\geq \frac{2}{3}\,d(x',p')$. Hence
$$
d(x,x')\leq d(x,p')+d(p',x')\leq d(x,p')+\frac{3}{2}\,d(x',z)\leq 
\frac{5}{2}\,d(x,p')\;.
$$

\medskip
Assume now that $p'$ belongs to $[\xi_0,x'[$. Let $\beta$ be the
comparison angle at $x'$ between the (nontrivial) geodesic segments
or rays $[x',y'[$ and $[x',\xi[$. By comparison, $\beta$ is at most
the angle $\overline \beta$ between $[\,\overline{x'},\overline{y'}[$
and $[\,\overline{x'},\overline{\xi}[$, where $\overline{x'}$ and
$\overline{y'}$ are two points, at distance $d(x',y')$, on an
horosphere in $\HH^2_\RR$ centered at $\overline{\xi}$. An easy
computation in the upper half space model shows that
$$
\tan \overline\beta = \left(\sinh \frac{1}{2} 
d(\overline{x'},\overline{y'}) \right)^{-1}\leq
\left(\sinh\big(2\log(1+\sqrt{2})\big)\right)^{-1}\leq 
\frac{1}{\sqrt{3}}\;.
$$
As $0\leq \overline\beta\leq \frac{\pi}{2}$, this implies that
$\beta\leq\overline\beta\leq \frac{\pi}{6}$.

Let $\alpha$ be the comparison angle at $x'$ between the (nontrivial)
geodesic segments $[x',p']$ and $[x',x]$, which is at most
$\frac{\pi}{2}$, as $p'$ is the closest point to $x$ on
$]x',\xi_0]$. As the geodesic segment $[x,x']$ lies in $H$, we have
$\alpha\geq \pi-\frac{\pi}{2}- \beta\geq\frac{\pi}{3}$. By Lemma
\ref{lem:c1pinfty}, we have $d(x,p')\leq d(x,x')\leq
c'_1(\infty)$. Using the formulae for right-angled hyperbolic
triangles (see \cite{Bea}) and the comparison triangle in $\HH^2_\RR$
to the triangle $(x,x',p')$ in $X$, we have, by convexity of $t\mapsto
\sinh t$,
$$
d(x,x')\leq \sinh d(x,x')\leq \frac{1}{\sin\alpha}\sinh d(x,p')\leq 
\frac{2}{\sqrt{3}}\frac{\sinh c'_1(\infty)}{c'_1(\infty)}d(x,p')\leq 
\frac{5}{2}\,d(x,p')\;.
$$
This proves the result.
\cqfd

\medskip 
The following Lemma is the analog of Lemma \ref{lem:h0} for
horoballs. It says that when the $\epsilon$-convex subset under
consideration is a horoball, we can replace $c'_3(\epsilon)$ and
$h'(\epsilon,\eta)$ by
\begin{equation}
c'_3(\infty)=\frac{5}{2}
\quad{\rm and}\quad
h'(\infty,\eta)=3\eta+c_0(\infty)+c'_1(\infty)\approx 3\eta +
5.8188\;, \label{eq:c3i}
\end{equation}
and that the first of the two possible conclusions of Lemma
\ref{lem:h0} always holds. Note that $c'_3(\infty)$ is not the limit as
$\epsilon$ goes to $+\infty$ of $c'_3(\epsilon)$, and that
$h'(\epsilon,\eta)$ diverges as $\epsilon\to\infty$. However, in both
cases, this notation will be useful in Section \ref{sec:main}.

\blemm \label{lem:c3pinfty} 
For every horoball $H$ in $X$, for every $\xi_0$ in 
$X\cup\partial_\infty X$, for all geodesic rays or lines $\ga,\ga'$
starting from $\xi_0$, if $\ga$ enters $H$ at a point $x\in X$ and
exits $H$ at a point $y\in X$, and if $d(x,y)\geq
h'(\infty,d(y,\ga'))$, then $\ga'$ meets $H$, exiting it at a point
$y'\in X$ such that
$$
d(y,y')\leq \frac{5}{2}\;d(y,\ga')\;.
$$
\elemm

\smallskip\noindent
\begin{minipage}{10cm}
  \dem Let $\xi$ be the point at infinity of $H$, let $p$ be the
  closest point on $[x,y]$ from $\xi$, and let $p_x$ and $p_y$ be the
  intersection of the horosphere $\partial H_p$ centered at $\xi$
  passing through $p$ with the geodesic rays $[x,\xi[$ and $[y,\xi[$
  respectively.  By comparison, we have $d(p_x,p_y)\le 2\log(1+\sqrt
  2)=c'_1(\infty)$.  Thus, the triangle inequality, along with the
  fact that $p_y$ is the closest point to $y$ on $\partial H_p$ and
  the assumption on $d(x,y)$, gives
$$ 
2\min\{d(y,p),d(x,p)\}\ge d(y,p_y)+d(x,p_x)\ge 
$$ 
$$ 
d(x,y)-2\log(1+\sqrt 2)\geq c_0(\infty)\geq 3\;. 
$$
\end{minipage}
\begin{minipage}{4.8cm}
\begin{center}
\input{fig_horopenepropnewer.pstex_t}
\end{center}
\end{minipage}

\medskip
In particular, as $d(x,y)\geq c_0(\infty)$, Lemma \ref{lem:cppinfty}
implies that $p$ belongs to $H[1]$.  By Lemma \ref{lem:expopinch} and
the assumption on $d(x,y)$, we have
$$ 
d(p,\ga')\leq \frac12\; e^{-d(y,p) +d(y,\ga')}\leq \frac12\; 
e^{-\frac12 d(x,y)+\log(1+\sqrt{2})+d(y,\ga')}\leq \frac12\;.
$$ 
This implies that $\ga'$ meets $H$, because $\N_{\frac 12}(H[1])$ is 
contained in $H$.

Let $x'$ and $y'$ be the entering point in $H$ and the exiting point
out of $H$ of $\ga'$, respectively.  Let $p'$ be the point on $\ga'$
the closest to $y$.

\medskip\noindent 
{\bf Case 1 :} %
Assume that $p'\notin[y',\xi_0]$. Note that 
$$
d(x,y)-\frac32 \,d(y,p')\geq \frac32\,d(y,p')\geq 0\;,
$$ 
as $d(x,y)\ge 3 \,d(y,\gamma')$.  Hence, there is a point $y_0$ in
$[x,y]$ at distance $\frac32 \,d(y,p')$ of $y$ which satisfies
$d(x,y_0)\geq \frac32\,d(y,p')$. By Lemma \ref{lem:cppinfty}, we have
$y_0\in H[d(y,p')]$. Let $q'$ be the point of $\ga'$ the closest to
$y_0$. By convexity, we have $d(y_0,q')\leq d(y,p')$. Hence $q'$
belongs to $H$.  By the intermediate value theorem, the point $y'$
belongs to $[q',p']$.  As closest point maps do not increase the
distances, we have $d(p',q')\leq d(y,y_0)=\frac32 \,d(y,p')$.
Therefore,
$$ 
d(y,y')\leq d(y,p')+d(p',y')\leq d(y,p')+d(p',q')
\leq \frac52 \,d(y,p')\;, 
$$
which  proves the result.

\medskip\noindent 
{\bf Case 2 :} Assume that $p'\in\;]y',\xi_0]$.  By the same argument
as in Case 2 of the proof of Lemma \ref{lem:h0}, we have
$p'\notin[x',\xi_0]$. If $d(p',y')\leq \frac32 \,d(y,p')$, then
$d(y,y')\leq \frac52 \,d(y,p')$, and the result is proved. Therefore,
assume by absurd that $d(p',y')>\frac32 \,d(y,p')$. By the continuity
of the closest point maps, there exists a point $y_0$ in $\ga$ that
does not belong to $H$, whose closest point $q'$ on $\ga'$, which lies
in $\ga'-\;]p',\xi_0]$, satisfies $d(q',y')\geq \frac32 \,d(y_0,q')$
and $d(y_0,q')\leq d(y,p')+ \frac12 \,c_0(\infty)$. Lemma
\ref{lem:c1pinfty} implies that $d(x,x')\le c_1'(\infty)$. Thus, by
the assumption on $d(x,y)$,
$$
\begin{aligned}
d(x',y') & \geq d(x',q')\geq d(p',x')\geq d(x,y)-d(x,x')-d(y,p')\\
& \geq 3\,d(y,p')+c_0(\infty)+c'_1(\infty)-c_1'(\infty)-d(y,p') \\
& \ge 2d(y,p')+ c_0(\infty) \geq \max\{2d(y_0,q'),c_0(\infty)\}\;.
\end{aligned}
$$
In particular, $d(y_0,q')\leq \frac23 d(q',x')$ and we already had
$d(y_0,q')\leq \frac23 d(q',y')$.  Hence, by Lemma \ref{lem:cppinfty},
we have $q'\in H[d(y_0,q')]$. This implies that $y_0$ belongs to $H$,
a contradiction.  \cqfd

\section{Properties of penetration in $\epsilon$-convex 
sets}
\label{sec:proppene}

%
%
%

\subsection{Penetration maps}
\label{subsec:penetrationmaps}

Let $X$ be a proper geodesic CAT$(-1)$ space, and $\xi_0\in
X\cup\partial_\infty X$.  We are interested in controlling the
penetration of geodesic rays or lines starting from $\xi_0$ in
$\epsilon$-convex subsets of $X$.  One way to measure this
penetration is the intersection length.
If $C$ is a closed convex subset in $X$ such
that $\xi_0\notin C\cup\partial_\infty C$, we define a map
$\ell_C:T^1_{\xi_0} X\ra [0,+\infty]$, called the {\it penetration
  length map}, which associates to every $\ga$ in $T^1_{\xi_0} X$ the
length of the intersection $\ga\cap C$ (which is connected by
convexity).

When we study specific geometric situations, such as collections of
horoballs and $\epsilon$-neighbourhoods of geodesics, there are
further natural ways of measuring the penetration. These will be used
in many applications in Section \ref{sec:results} and in \cite{HPP}.
If $C$ is an $\epsilon$-convex subset of $X$ such that $\xi_0\notin
C\cup\partial_\infty C$, we will require our penetration maps
$f:T^1_{\xi_0} X\ra [0,+\infty]$ in $C$ to have one or two of the
following properties, the first one depending on a constant
$\kappa\geq 0$.  The $\sup$--norm of a real valued function $f$ on
$T^1_{\xi_0} X$ is denoted by $\norm{f}$.
\begin{enumerate}
\item[$(i)$] {\bf (Penetration property) } %
$\norm{f-\ell_C}\le\kappa\;.$
\item[$(ii)$] {\bf (Lipschitz property) } %
  For every $\ga,\ga'$ in $T^1_{\xi_0} X$ which intersect $C$, if
  $\ga\cap C=[a,b]$ and  
  $\ga'\cap C=[a',b']$ with $a,b,a',b'$ in $X$, then
  $$|f(\ga)-f(\ga')|\leq 2\max\big\{d(a,a'),d(b,b')\big\}\;.$$
\end{enumerate}
If $C$ is an $\epsilon$-convex subset of $X$ such that $\xi_0\notin
C\cup\partial_\infty C$, and $f:T^1_{\xi_0} X\to[0,+\infty[$ is a map
which satisfies $(i)$ for some $\kappa\ge 0$, we say that $f$ is a
$\kappa$-{\em penetration map in (the $\epsilon$-convex set)} $C$.  We
also say that $(C,f)$ is an $(\epsilon,\kappa)$-{\em penetration
  pair}.  In the condition $(ii)$, we could have replaced $2$ by some
$\lambda\geq 2$, but if $f$ also satisfies the property $(i)$, then
only $\lambda =2$ is really relevant in the large scale.

Note that if $(C,f)$ is an $(\epsilon',\kappa')$-penetration pair, if
$\epsilon'\geq \epsilon$ and $\kappa'\leq \kappa$, then $(C,f)$ is an
$(\epsilon,\kappa)$-penetration pair.  If $C$ is $\infty$-convex and
$(C,f)$ is an $(\epsilon,\kappa)$-penetration pair in every
$\epsilon>0$ then $f$ is a $\kappa${\em -penetration map in (the
  $\infty$-convex set)} $C$.

\medskip\noindent{\bf Penetration maps in general $\epsilon$-convex
  subsets. } %
If $C$ is a closed convex subset of $X$, the map $\ell_C$ is in
general not continuous on $T^1_{\xi_0} X$, as can be seen by taking
$C$ to be a geodesic segment of positive length.  The following result
shows that the situation is nicer for $\epsilon$-convex subsets.  Note
that the statement of Lemma \ref{lem:proppenegene} is not true in
$\RR^n$ (which is not a CAT$(-1)$ space).

\blemm \label{lem:proppenegene} %
Let $\epsilon>0$ and let $C$ be an $\epsilon$-convex subset of $X$
such that $\xi_0\notin C\cup\partial_\infty C$. The map
$\ell_C:T^1_{\xi_0} X\ra [0,+\infty]$ is a continuous $0$-penetration
map in $C$ satisfying the Lipschitz property $(ii)$.  \elemm

\dem %
The Lipschitz property $(ii)$ of the penetration length map $\ell_C$
follows from the triangular inequality. It remains to show the
continuity of the map.

Choose a convex subset $C'$ such that $C=\N_\epsilon (C')$, and note
that by the definition of the topology of $X\cup \partial_\infty X$,
the subsets $C$ and $C'$ have the same points at infinity. Let
$\ga_0\in T^1_{\xi_0} X$, and let us prove that $\ell_C$ is continuous
at $\ga_0$.

Assume first that $\ga_0(+\infty)$ is a point at infinity of $C$. Then
there exists a geodesic ray contained in $C'$ ending at this point at
infinity. As geodesic rays converging to the same point at infinity
become exponentially close, this implies that $\ell_C(\ga_0)=\infty$.
Let $A>0$.  As $\ga_0\cap C$ is the closure of $\ga_0\;\cap
\stackrel{\circ}{C}$, let $[x,y]$ be a geodesic segment of length
$A+2$ contained in $\ga_0\;\cap \stackrel{\circ}{C}$.  Let $\eta \in
\;]0,1]$ be such that the balls $B$ and $B'$ of radius $\eta$ and of
center $x$ and $y$ respectively are contained in
$\stackrel{\circ}{C}$. If $\ga\in T^1_{\xi_0} X$ is close enough to
$\ga_0$, then $\ga$ meets $B$ and $B'$, and by convexity,
$\ell_C(\ga)\geq A$, which proves the result.

Assume now that $\ga_0(+\infty)$ is not a point at infinity of $C$,
but that $\ga_0$ does meet $C$. Then $\ga_0\cap C$ is a nonempty
compact segment $[a,b]$. For every $\eta>0$, let $a_+,b_+$ be points
in $\ga_0-[a,b]$, at distance at most $\eta/4$ from $a,b$
respectively, and, if $d(a,b)>0$, let $a_-,b_-$ be points in $]a,b[$
at distance at most $\eta/4$ from $a,b$ respectively. As $C$ is closed
and $\ga_0\cap C$ is the closure of $\ga_0\;\cap \stackrel{\circ}{C}$
if $a\neq b$, there exists $\eta'\in\;]0,\eta/4]$ such that the balls
$B(a_+), B(b_+)$ of radius $\eta'$ and centers $a_+,b_+$ respectively
are contained in $X-C$ and, if $d(a,b)>0$, the balls $B(a_-), B(b_-)$
of radius $\eta'$ and centers $a_-,b_-$ respectively are contained in
the interior of $C$.  If $\ga\in T^1_{\xi_0} X$ is close enough to
$\ga_0$, then $\ga$ meets $B(a_+), B(b_+)$ (and hence $B(a_-), B(b_-)$
by convexity, if $d(a,b)>0$). It is easy to see then that
$|\ell_C(\ga)-\ell_C(\ga_0)|\leq \eta$.

Assume now that $\ga_0$ does not meet $C$. Let $U,V$ be neighborhoods
of the endpoints of $\ga_0$ in $X\cup\partial_\infty X$ that are
disjoint from $C\cup\partial_\infty C$. Let $\eta>0$ be such that the
$\eta$-neighborhood of $\ga_0$ is disjoint from $C$, which exists, as
$\inf_{x\in\ga_0}d(x,C) >0$. If $\ga\in T^1_{\xi_0} X$ is close enough
to $\ga_0$, then (the image of) $\ga$ lies in $U\cup V\cup
\N_\eta\ga_0$, hence does not meet $C$. So that
$\ell_C(\ga)=\ell_C(\ga_0)=0$.
\cqfd

\medskip %
In particular, if $H$ is a horoball such that $\xi_0\notin
H\cup\partial_\infty H$, then $\ell_H$ is a continuous $0$-penetration
map for $H$ satisfying the Lipschitz property $(ii)$.

\medskip %
Let $C$ be a convex subset of $X$ such that $\xi_0\notin
C\cup\partial_\infty C$.  For every $\ga$ in $T^1_{\xi_0} X$, let
$\ga_-=\xi_0$ and $\ga_+=\ga(+\infty)$, and let $q_{\ga_\pm}$ be the
closest point on $C$ to $\ga_\pm$. Define the {\it boundary-projection
  penetration map} $\bp_C: T^1_{\xi_0} X \ra [0,+\infty]$ by
$$
\bp_C(\ga)=d(q_{\ga_-},q_{\ga_+})\;,
$$
with the obvious convention that $\bp_C(\ga)=+\infty$ if $q_{\ga_+}$
is at infinity.

\blemm\label{lem:boundaryprojpenemap} %
Let $C$ be an $\epsilon$-convex subset of $X$ such that $\xi_0\notin
C\cup\partial_\infty C$. The map $\bp_C$ is a continuous
$2c'_1(\epsilon)$-penetration map in $C$.
\elemm

\dem %
The continuity of $\bp_C$ follows from the continuity of the
projection maps and the endpoint maps. Let $\ga\in T^1_{\xi_0} X$, and
let us prove that $|\bp_C(\ga)-\ell_C(\ga)|\leq 2c'_1(\epsilon)$. If
$\ga_+$ is a point at infinity of $C$, then
$\bp_C(\ga)=\ell_C(\ga)=+\infty$, and the result is true. Otherwise,
if $\ga$ meets $C$, then $\ga$ enters $C$ at $x$ and exits $C$ at $y$,
with $x,y$ in $X$. By Lemma
\ref{lemma:c1}, we hence have
$$
|d(x,y)-d(q_{\ga_-},q_{\ga_+})|\leq d(x,q_{\ga_-})+d(q_{\ga_+},y)\leq
 2c'_1(\epsilon)\;,
$$
and the result follows.

If $\ga$ does not meet $C$, let $[p,q]$ be the shortest connecting
segment between a point $p$ in $\ga$ and a point $q$ in $C$. By angle
comparison, the geodesic segment or ray between $q$ and $\ga_\pm$
meets $C$ exactly in $q$. Hence, by Lemma \ref{lemma:c1},
$$
d(q_{\ga_-},q_{\ga_+})\leq d(q_{\ga_-},q)+d(q,q_{\ga_+})
\leq 2c'_1(\epsilon)\;.
$$
As $\ell_C(\ga)=0$, the result follows.
\cqfd

\medskip\noindent{\bf Penetration maps in horoballs. } %
If $H$ is a horoball in $X$, with $\xi$ its point at infinity, such
that $\xi_0\notin H\cup\{\xi\}$, and if $x_0$ is any point in the
boundary of $H$ in $X$, define a $1$-Lipschitz map $\beta_H:X\ra
[0,+\infty[$, called the {\it height map} of $H$ by
$$
\beta_H:x\mapsto \max\{\beta_\xi(x_0,x),0\},
$$ 
whose values are positive in the interior of $H$, and $0$ outside $H$.
By convention, define $\beta_H(\xi)=+\infty$. For every $\ga$ in
$T^1_{\xi_0} X$, let $p_\ga$ be the closest point to $\ga(+\infty)$ on
the geodesic line between $\xi_0$ and $\xi$, with $p_\ga=\xi$ if
$\ga(+\infty)=\xi$.

We will study two penetration maps associated with the height map.
The map $\ph_H:T^1_{\xi_0} X\ra [0,+\infty]$ defined by
$$
\ph_H(\ga)=2\;\sup_{t\in\RR} \beta_H(\ga(t))
$$ 
will be called the {\it penetration height map} inside $H$.  The map
$\ipp_H:T^1_{\xi_0} X\ra [0,+\infty]$ defined by
$$
\ipp_H(\ga)= 2\;\beta_H(p_\ga)
$$ 
will be called the {\it inner-projection penetration map} inside $H$.
Note that for every $t\geq 0$ and $\ga\in T^1_{\xi_0} X$, we have
$\ph_{H[t]}(\ga) = \max\{0,\ph_{H}(\ga)-2t\}$ and
$\ipp_{H[t]}(\ga) = \max\{0,\ipp_{H}(\ga)-2t\}$.

\blemm\label{lem:proppenehoro} %
Let $H$ be an horoball in $X$, such that $\xi_0\notin
H\cup\partial_\infty H$. The maps $\ph_H,\ipp_H: T^1_{\xi_0} X\ra
[0,+\infty]$ are continuous $2\log(1+\sqrt2)$-penetration maps for
$H$, and $\ph_H$ has the Lipschitz property $(ii)$.  Furthermore,
$$
\norm{\ph_H-\ipp_H}\leq 2\log(1+\sqrt{2}).
$$  
\elemm

\noindent{\bf Remark.} %
In $\hnr$, the equality $\ipp_H(\ga)=\ph_H(\ga)+\log 2$ holds for any
horoball $H$ with $\xi_0\notin H\cup H[\infty]$ and $\ga\in T^1_{\xi_0}
X$ meeting $H$. Thus, in $\hnr$ the map $\ipp_H$ satisfies the
Lipschitz property $(ii)$. We do not know whether $(H,\ipp_H)$
satisfies the Lipschitz property $(ii)$ in general.

\medskip \dem %
Let us prove that $(H,\ph_H)$ satisfies the Lipschitz property $(ii)$.
Let $\ga,\ga'$ be elements in $T^1_{\xi_0} X$ such that $\ga\cap
H=[a,b]$ and $\ga'\cap H=[a',b']$. Then, for every $x$ in $[a,b]$, if
$x'$ is the point on $[a',b']$ the closest to $x$, we have, with $\xi$
the point at infinity of $H$,
$$ 
|\beta_\xi(x,a)-\beta_\xi(x',a)|=|\beta_\xi(x,x')|\leq d(x,x')
\leq \max\{d(a,a'),d(b,b')\}
$$ 
by convexity. Taking $x$ the highest point in $[a,b]$, and using a
symmetry argument, the result follows.

\medskip %
Let us prove that $(H,\ph_H)$ satisfies the Penetration property $(i)$
with $\kappa = c'_1(\infty)= 2\log(1+\sqrt2)$. Let $\ga\in T^1_{\xi_0}
X$. Note that $\ga$ enters the interior of $H$ if and only if
$\ell_H(\ga)>0$, and if and only if $\ph_H(\ga)>0$.  Hence we may
assume that $\ga$ meets $H$ in a segment $[x,y]$. By the first
paragraph of the proof of Lemma \ref{lem:c3pinfty}, we have
$$
\ph_H(\ga)\leq\ell_H(\ga)\leq \ph_H(\ga) + 2\log(1+\sqrt2)\;.
$$

\medskip %
Let us prove that $(H,\ipp_H)$ satisfies the Penetration property
$(i)$ with $\kappa =2\log(1+\sqrt{2})$. Let $\ga\in T^1_{\xi_0} X$. If
$p_\ga=\xi$, then $\ipp_H(\ga)=\ell_H(\ga)= +\infty$, and the result
holds, hence we may assume that $p_\ga$ belongs to $X$. If $p_\ga$
does not belong to $H$, as the closest point projection of $\gamma$ on
the geodesic line $]\gamma_-,\xi[$ is $]\gamma_-,p_\gamma[$, then
$\ga$ does not enter $H$, and hence $\ipp_H(\ga)= \ell_H(\ga)=0$, and
the result is proven.

\medskip\noindent
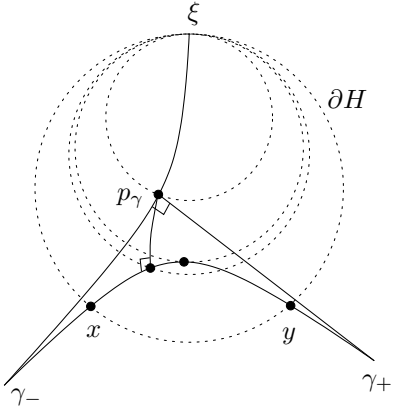
\begin{minipage}{5.6cm}
\begin{center}
\input{fig_horopenprojprop.pstex_t}
\end{center}
\end{minipage}
\begin{minipage}{9.2cm}
  Assume that $p_\ga$ belongs to $H$, and note that by comparison and
  an easy hyperbolic estimate, we have $d(p_\ga,\ga)\leq
  \log(1+\sqrt{2})$.  In particular, if $\ga$ does not enter $H$, then
  $0\leq \beta_H(p_\ga)\leq d(p_\ga,\ga)\leq \log(1+\sqrt{2})$, and
  $|\ipp_H(\ga)- \ell_H(\ga)|\leq 2\log(1+\sqrt{2})$, hence the result
  holds. Therefore we may assume that $\ga$ enters $H$ at the point
  $x$ and exits $H$ at the point $y$. We then have $\ph_H(\ga) \leq
  \ipp_H(\ga)\leq\ph_H(\ga)+2d(p_\ga,\ga)$. Hence,
  $$
  \ell_H(\ga)-2\log(1+\sqrt2)\leq\ipp_H(\ga)\leq
  \ell_H(\ga)+2\log(1+\sqrt{2})\;,
  $$
  and the result is proven.
\end{minipage}

\medskip %
The continuity of $\ipp_H$ follows from the continuity of the endpoint
maps, of the closest point projection maps and of $\beta_H:X\cup\{\xi\} \ra
[0,+\infty]$. To prove the continuity of $\ph_H$ at a point $\ga_0$ of
$T^1_{\xi_0} X$, note that if $\ga_0(+\infty)=\xi$, then $\ph_H(\ga_0)
= +\infty$, and the continuity follows from the Penetration property
$(i)$ of $(H,\ph_H)$ and the continuity of $ \ell_H$. Otherwise,
$\ga_0\cap H$ is a compact segment. If it is nonempty, then if $\ga$
is close enough to $\ga'$, the argument in the proof of Lemma
\ref{lem:proppenegene} shows that the Hausdorff distance between
$\ga\cap H$ and $\ga_0\cap H$ is as small as wanted.  The result
follows then since $\beta_H$ is $1$-Lipschitz. If $\ga_0$ does not
meet $ H$, then if $\ga$ is close enough to $\ga'$, the argument in
the proof of Lemma \ref{lem:proppenegene} shows that $\ga$ also avoids
$H$, hence $\ph_H(\ga)=\ell_H(\ga)=0$.  \cqfd

\bigskip\noindent{\bf Penetration maps in balls. } %
If $B$ is a ball of center $x_0$ and radius $r_0$ in $X$ with
$\xi_0\notin B$, define a $1$-Lipschitz map $\beta_B:X\ra
[0,+\infty[$, called the {\it height map} by $\beta_B:x\mapsto
\max\{r_0-d(x_0,x),0\}$, whose values are positive in the interior of
$B$, and $0$ outside $B$. For every geodesic line $\ga$ in
$T^1_{\xi_0} X$, let $p_\ga$ be the closest point to $\ga(+\infty)$ on
the geodesic segment or ray between $\xi_0$ and $x_0$.

The map $\ph_B:T^1_{\xi_0} X\ra [0,2r_0]$ defined by 
$$
\ph_B(\ga)= 2\;\sup_{t\in\RR} \beta_B(\ga(t))
$$ 
will be called the {\it penetration height map} inside $B$. The map
$\ipp_B:T^1_{\xi_0} X\ra [0,2\,r_0]$ defined by
$$
\ipp_B(\ga)=2\;\beta_B(p_\ga)
$$ 
will be called the {\it inner-projection penetration map} inside $B$.

As above, the maps $\ph_B,\ipp_B$ are continuous
$2\log(1+\sqrt2)$-penetration maps, and $\ph_B$ has the Lipschitz
property {\em (ii)}.  Furthermore,
$$
\norm{\ph_B-\ipp_B}\leq
2\log(1+\sqrt{2}).
$$
In the proof of the Penetration
property of $\ipp_B$, note that if $p_\gamma=x_0$, then, by
comparison, $d(\gamma,x_0)\le\log(1+\sqrt 2)$, and the claim follows
in this case as $\ipp_B(\gamma)=2r_0$.

If a sequence of balls $(B_i)_{i\in\NN}$ converges to an horoball $H$
(for the Hausdorff distance on compact subsets of $X$), then the maps
$\ph_{B_i}, \ipp_{B_i}$ converge, uniformly on compact subsets of
$T^1_{\xi_0} X$, to $\ph_H,\ipp_H$ respectively.

\bigskip
\noindent{\bf Penetration maps in tubular neighborhoods 
of totally geodesic subspaces. } 

We define two functions on $T^1_{\xi_0}X$ which describe the closeness
of a geodesic line to a totally geodesic subspace $L$. If $\xi_0$ is
in the boundary at infinity of $X$, then these functions are defined
without reference to an $\epsilon$-neighbourhood of $L$. However, we
show that they are penetration maps in the $\epsilon$-neighbourhood
of $L$, with explicit constants which depend only on $\epsilon$.

Let $\epsilon>0$, and let $L$ be a complete totally geodesic subspace
of $X$, with set of points at infinity $\partial_\infty L$, such that
$\xi_0\notin \N_\epsilon L\cup \partial_\infty L$. For every $\ga$ in
$T^1_{\xi_0} X$, let $\ga_-=\xi_0$ and $\ga_+=\ga(+\infty)$, and let
$p_{\ga_\pm}$ be the point on $L$ the closest to $\ga_\pm$.

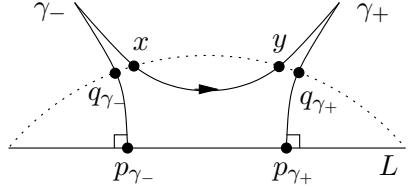
\begin{center} 
\input{fig_verifiitub.pstex_t}
\end{center}

\noindent %
We define the {\it fellow-traveller penetration map} 
${\ftp}_L:T^1_{\xi_0} X\ra[0,+\infty]$ by
$$
{\ftp}_L(\ga)=d(p_{\ga_-},p_{\ga_+})\;,
$$ 
with the convention that this distance is $+\infty$ if $p_{\ga_+}$ is
in $\partial_\infty X$.

\blemm \label{fellowtravelpenemap} %
Let $\epsilon>0$, and let $L$ be a complete totally geodesic subspace
of $X$ such that $\xi_0\notin \N_\epsilon L\cup \partial_\infty L$.
The map ${\ftp}_L$ is a continuous
$(2c'_1(\epsilon)+2\epsilon)$-penetration map in $ \N_\epsilon L$ and
$\norm{{\ftp}_L-\bp_{\N_\epsilon L}} \leq 2\epsilon$.
\elemm

\dem %
The continuity of ${\ftp}_L$ follows from the continuity of the
projection maps and of the endpoint maps.  Note that, for every $\ga$
in $T^1_{\xi_0} X$, the geodesic ray from $p_{\ga_\pm}$ to $\ga_\pm$
exits $\N_\epsilon L$ at the closest point $q_{\ga_\pm}$ on
$\N_\epsilon L$ to $\ga_\pm$. Hence, by the triangular inequality, and
as closest point maps do not increase distances, we have
$$
0\leq \bp_{\N_\epsilon L}(\ga)-{\ftp}_L(\ga)\leq 2\epsilon\;.
$$
Therefore the fact that ${\ftp}_L(\ga)$ satisfies the Penetration
property $(i)$ with $\kappa=2c'_1(\epsilon)+2\epsilon$ follows from
Lemma \ref{lem:boundaryprojpenemap}. 
\cqfd

\medskip

If $L$ is one-dimensional and $\xi_0\in \partial_\infty
X-\partial_\infty L$, a natural penetration map is defined using the
crossratios of the endpoints of $L$ and $\gamma$. Let $\partial_4X$ be
the set of quadruples $(a,b,c,d)$ in $(\partial_\infty X)^4$ such that
$a\neq b$ and $c\neq d$.  The {\it crossratio}
$[a,b,c,d]\in[-\infty,+\infty]$ of a quadruple $(a,b,c,d)$ in
$\partial_4X$ is defined as follows (see for instance
\cite{Ota,Bou,Pau}). If $a_t,b_t,c_t,d_t$ are any geodesic rays
converging to $a,b,c,d$ respectively, then
$$
[a,b,c,d] = \frac{1}{2}\lim _{t\rightarrow+\infty}
               d(a_t,c_t)-d(c_t,b_t)+d(b_t,d_t)-d(d_t,a_t).
$$
Note that the order conventions differ in the references, we are using
the ones of \cite{Bou,HP}), and that our crossratio is the logarithm
of the crossratio used in \cite{Bou}.   

Let us give other formulae for the crossratio. The {\em visual
distance} of two points $a$ and $b$ in $\partial_\infty X$ with
respect to $x_0$ is
$$
d_{x_0}(a,b)=
\lim_{t\to\infty}e^{-\frac12\big(d(x_0,a_t)+d(x_0,b_t)-d(a_t,b_t)\big)}.  
$$

If $\xi\in\partial_\infty X$, if $H$ is a horosphere centered at
$\xi$, and $a,b$ are points in $\partial_\infty X-\{\xi\}$, and
$t\mapsto x_t$ is a geodesic ray with $x_0\in H$ which converges to
$\xi$, the {\em Hamenst\"adt distance} (defined in \cite{Ham},
\cite[Appendix]{HP}) of $a$ and $b$ in $\partial_\infty X- \{\xi\}$ 
normalized with respect to $H$ is
$$
d_{H}(a,b)=\lim_{t\to\infty}e^td_{x_t}(a,b).
$$ 
Note that if $H'$ is another horosphere centered at $\xi$, then there
exists a constant $c>0$ such that $d_{H'}= c \,d_H$. In particular,
for every $\xi'\in \partial_\infty X- \{\xi\}$ and $r>0$, the sphere
of center $\xi'$ and radius $r$ for $d_H$ coincides with the sphere of
center $\xi'$ and radius $r$ for $d_{H'}$.  

It is easy to see that for any $x_0\in X$ and any horoball $H$, we
have for every $(a,b,c,d)\in\partial_4 X$
$$
[a,b,c,d]   =   \log\;\frac{d_{x_0}(a,c)}{d_{x_0}(c,b)}\;
               \frac{d_{x_0}(b,d)}{d_{x_0}(d,a)}
	    =   \log\;\frac{d_{H}(a,c)}{d_{H}(c,b)}\;
               \frac{d_{H}(b,d)}{d_{H}(d,a)}\;,
$$
if, in the second equation, $a,b,c,d$ are in $\partial_\infty
X-\{\xi\}$. Note that each expression in the above two equalities is
$-\infty$ if $a=c$ or $b=d$, and $+\infty$ if $c=b$ or $a=d$. If the
points $\xi$ and $a$ coincide, the expression of the crossratio
simplifies to
$$
[\xi,b,c,d]  =   \log\;\frac{d_{H}(b,d)}{d_{H}(c,b)}\;.
$$
The crossratio is continuous on $\partial_4X$, it is invariant under
the diagonal action of the isometry group of $\Ga$, and it has the
following symmetries
$$
[c,d,a,b]= [a,b,c,d]\;\;{\rm and }\;\; [a,b,d,c]= [b,a,c,d]= -[a,b,c,d].
$$  

If $X=\hnr$ and $\xi$ is the point at infinity $\infty$ in the upper
halfspace model of $\hnr$ , then the Hamenst\"adt distance coincides
with a constant multiple of the Euclidean distance of
$\partial_\infty\hnr-\{\infty\}=\mathbb R^{n-1}$ (see for instance
\cite{HP3}). In particular, if $n=2$, then our crossratio is the
logarithm of the modulus of the classical crossratio of four points in
$\mathbb C\cup\{\infty\}$.

 \medskip
If $\xi_0\in\partial_\infty X-\partial_\infty L$, we define the {\it
crossratio penetration map} ${\cp}_{L}:T^1_{\xi_0} X\ra
[0,+\infty]$ as follows. Let $\ga$ be a geodesic line starting at
$\ga_-=\xi_0$, and ending at $\ga_+\in\partial X$.  Let $L_1,L_2$ be
the endpoints of $L$. Set
$$
{\cp}_L(\ga)= \max\big\{0,[\ga_-,L_1,\ga_+,L_2],
[\ga_-,L_2,\ga_+,L_1]\big\}\;,
$$ 
if $\ga_+\neq L_1,L_2$, and $ {\cp}_L(\ga)=+\infty$ otherwise.  The
map ${\cp}_L$ is clearly continuous, in particular since $[a,b,c,d]$
tends to $0$ when $a,b,d$ are pairwise distinct and $c$ tends to $d$,
and is independant of the ordering $L_1,L_2$ of the endpoints of $L$.

\bigskip
If $H$ is a horosphere centered at $\xi_0$, then
$$
[\xi_0,L_1,\ga_+,L_2]=\log\frac{d_H(L_1,L_2)}{d_H(\ga_+,L_1)},
$$
and the level sets for ${\cp}_L$ have a simple form:
$[\xi_0,L_1,\ga_+,L_2]=c$ if and only if $\ga_+$ is on the sphere of
radius $e^{-c}d_H(L_1,L_2)$ centered at $L_1$ with respect to the
Hamenst\"adt metric. Thus, in particular, the  
boundary of the zero set of $\cp_L$ is the boundary of the union of
the two balls of radius $d_H(L_1,L_2)$ centered at $L_1$ and $L_2$. 
Furthermore, if $c>\log 2$, then the level set  ${{\cp}_L}^{-1}(c)$ is
the union of two spheres for the Hamenstädt distance $d_H$ of
centers $L_1$ and  $L_2$ and radius $e^{-c}d_H(L_1,L_2)$. These two
spheres are disjoint by the triangle inequality. Each of them
separates $\xi_0$ from exactly one of the endpoints of $L$. We will use
this in the proof of Lemma \ref{lem:topopeneprescri}.

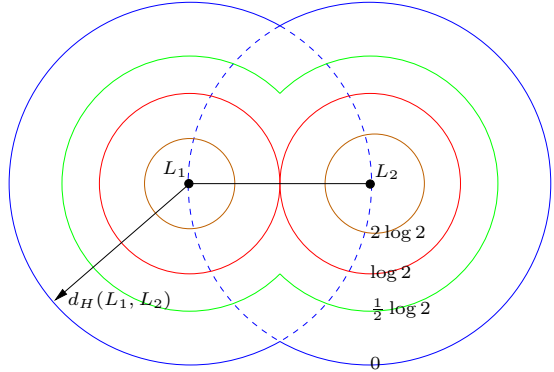
\begin{center}
\input{fig_crosspenfn.pstex_t}
\end{center}

Note that if $X$ is a negatively curved symmetric space, then the spheres and
balls of the Hamenst\"adt distance are topological spheres and balls
in the topological sphere $\partial _\infty X$ (see \cite{HP3} if
$X=\hnr$ and \cite{HP4} if $X=\hnc$). We don't know (and in
fact we doubt it) whether this always holds in the general variable
curvature case.

\blemm\label{lem:crossratio}
Let $(a,b,c,d)\in\partial_4 X$. If $b=d$, we define by convention
$p=q=b$ and $d(p,q)=0$. Otherwise, let $p$ and $q$ be the closest
points on $[b,d]$ of $a$ and $c$ respectively. 
\begin{enumerate} 
\item If $b,q,p,d$ are in this order on $[b,d]$ and $d(p,q)\ge
  c'_1(\infty)$, then $\big|[a,b,c,d]-d(p,q)\big|\le 2 c'_1(\infty)$.
\item If $b,p,q,d$ are in this order on $]b,d[$ and
$d(p,q)\ge c'_1(\infty)$, then $[a,b,c,d]\le c'_1(\infty)$.  
\item If $d(p,q)\le c'_1(\infty)$, then $[a,b,c,d]\le 2\;
c'_1(\infty)$.     
\end{enumerate}
\elemm 

\dem %
If $a=d$ or $c=b$, then $p=d$ or $q=b$, hence we are in case (1) and
$[a,b,c,d]=d(p,q)=\infty$, which proves the result. If $a=c$ or $b=d$,
then $p=q$, we are in case (3) and $[a,b,c,d]=-\infty$, which proves
the result. Hence we may assume that $a,b,c,d$ are pairwise
disjoint. 

Let $a_t,b_t,c_t,d_t$ be geodesic rays converging to respectively
$a,b,c,d$ as $t\to\infty$, and let $p_t$ and $q_t$ be the closest
points to $a_t$ and $c_t$ respectively on $[b_t,d_t]$.  Let
$p'\in[a_t,d_t]$ and $p''\in[a_t,c_t]$ be the closest points to $p_t$ on
$[a_t,d_t]$ and $[a_t,c_t]$, and let $q'\in[b_t,c_t]$ and
$q''\in[a_t,c_t]$ be the closest points to $q_t$ on $[b_t,c_t]$ and
$[a_t,c_t]$.

Recall that by an easy comparison argument, for pairwise distinct
points $u,v,w$ in $X\cup\partial_\infty X$, if $r$ is the closest
point to $w$ on $]u,v[$, then $r$ is at distance less than
$\delta=\log(1+\sqrt{2})$ from a point on $]u,w[\,$. We will apply
this remarkto $r=p_t$ and $r=q_t$.  Recall also that
$c'_1(\infty)=2\,\delta$.

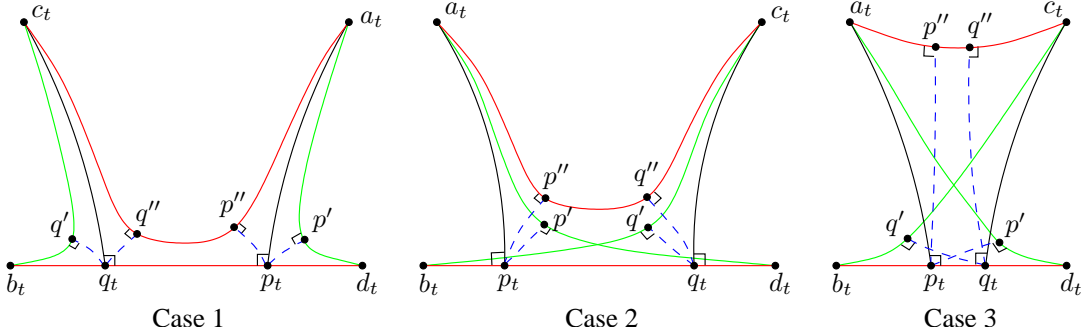
\begin{center}
\input{fig_crossratio3.pstex_t}
\end{center}

\noindent{\bf Case (1).} 
If $t$ is big enough, then the points
$b_t,q_t,p_t,d_t$ are in this order on $[b_t,d_t]$. 
Using the triangle inequality on $d(a_t, c_t)$ and
$d(b_t,d_t)$, and inserting the points $p'$ and $q'$, we have
\begin{align*}
 &d(a_t, c_t)-d(c_t,b_t)+d(b_t,d_t)-d(d_t,a_t)\\
 & \leq
d(a_t,p_t)+d(p_t,q_t)+d(q_t,c_t)-d(c_t,q_t)-d(q_t,b_t)+2\; d(q',q_t)\\
 & \ \ \ +
d(b_t,q_t)+d(q_t,p_t)+d(p_t,d_t)-
d(d_t,p_t)-d(p_t,a_t)+2\; d(p',p_t)\\ 
& \le  2\; d(p_t,q_t)+4\,\delta\;.
\end{align*}
By comparison and a standard argument on hyperbolic quadrilaterals with
three right angles (see \cite[page 157]{Bea}), for every $\epsilon
>0$, if $t$ is big enough,  we have that $d(q_t,q'')\le
2\,\delta +\epsilon/4$,  as $d(p_t,q_t)\ra d(p,q)\ge c'_1(\infty)$.
If we insert the points $p''$ and $q''$, we get, as above,
\begin{align*}
 &d(a_t, c_t)-d(c_t,b_t)+d(b_t,d_t)-d(d_t,a_t)\\
 & \geq
d(a_t,p_t)+d(p_t,q_t)+d(q_t,c_t)-2\;d(p_t,p'')-2\;d(q_t,q'')
-d(c_t,q_t)-d(q_t,b_t)\\
 & \ \ \ +
d(b_t,q_t)+d(q_t,p_t)+d(p_t,d_t)-
d(d_t,p_t)-d(p_t,a_t)\\ 
& \ge  2\; d(p_t,q_t)-8\,\delta +\epsilon\;.
\end{align*}

\noindent{\bf Case (2).} 
The proof is almost identical to the one of the upper bound in the
first inequality in Case 1. The different order of the points $p_t$
and $q_t$ now causes cancellations:
\begin{align*}
 &d(a_t, c_t)-d(c_t,b_t)+d(b_t,d_t)-d(d_t,a_t)\\
 & \leq
d(a_t,p_t)+d(p_t,q_t)+d(q_t,c_t)-d(c_t,q_t)-d(q_t,p_t)-d(p_t,b_t)+2\; d(q',q_t)\\
 & \ \ \ +
d(b_t,p_t)+d(p_t,q_t)+d(q_t,d_t)-
d(d_t,q_t)-d(q_t,p_t)-d(p_t,a_t)+2\; d(p',p_t)\le 4\,\delta\;.
\end{align*}

\noindent{\bf Case (3).} Let $\epsilon>0$. By taking $t$ big enough,
we can assume that $d(p_t,q_t)\le\;2\,\delta+\epsilon$.  Inserting the
points $p''$ and $q''$ and using the fact that closest point
maps do not increase distances, we have $d(p'',q'')\le
2\,\delta+\epsilon$, $d(c_t,q_t)\ge d(c_t,q'')$ and $d(a_t,p_t)\ge
d(a_t,p'')$. Thus, as in the cases above,
\begin{align*}
 &d(a_t, c_t)-d(c_t,b_t)+d(b_t,d_t)-d(d_t,a_t)\\
 & \leq
d(a_t,p'')+d(p'',q'')+d(q'',c_t)-d(c_t,q_t)-d(q_t,b_t)+2\; d(q',q_t)\\
 & \ \ \ +
d(b_t,q_t)+d(q_t,p_t)+d(p_t,d_t)-
d(d_t,p_t)-d(p_t,a_t)+2\; d(p',p_t)\\ 
& \le  2\; 8\,\delta+2\epsilon\;.
\end{align*}
As this holds for any $\epsilon>0$, the result follows.
\cqfd

\medskip \blemm \label{lem:crossratiopene} 
Let $\epsilon>0$, let $L$ be a geodesic line in $X$, and assume that
$\xi_0\in\partial_\infty X-\partial_\infty L$.  The map $\cp_L$ is a
continuous $(2c'_1(\epsilon)+2c'_1(\infty)+2\epsilon)$-penetration map in the
$\epsilon$-convex set $\N_\epsilon L$ and $\norm{\cp_L-\ftp_L}\le
2c'_1(\infty)$.
\elemm 

\dem 
Let $\ga\in T^1_{\xi_0}X$, let $\ga_-=\xi_0$ and $\ga_+$ be
the endpoints of $\ga$, and $L_1$ and $L_2$ be the endpoints of $L$.
Let $p$ and $q$ be the closest points to $\ga_-$ and $\ga_+$ on $L$
respectively.

If $d(p,q)\le c'_1(\infty)$, then Lemma \ref{lem:crossratio}
implies that $0\le\cp_L(\ga)\le2\;c'_1(\infty)$, and thus
$|\cp_L(\ga)-\ftp_L(\ga)|\le 2\;c'_1(\infty)$.

If $d(p,q)>c'_1(\infty)$, then up to renaming the endpoints of $L$, we
have $[\ga_-,L_2,\ga_+,L_1]\le c'_1(\infty)$ and
$-2\,c'_1(\infty)+\ftp_L(\ga)\le [\ga_-,L_1,\ga_+,L_2]\le
2\,c'_1(\infty)+\ftp_L(\ga)$, which implies the result.
\cqfd

\bigskip

\noindent{\bf Remark.} %
The penetration maps can be defined for any fixed starting point which
is outside the $\epsilon$-convex set $C$, except for $\cp_L$, and its
boundary at infinity.  Thus, the penetration maps
$\ell_C,\bp_C,\ph_H,\ipp_H,\ph_B,\ipp_B,\ftp_L$ considered in this
section are all restrictions to $T^1_{\xi_0} X$ of maps defined, and
continuous (as an inspection of the above proof shows) on
$\bigcup_{\xi\notin C'\cup\partial_\infty C'} T^1_\xi X\subset T^1X$
with $C'$ respectively $C,C,H,H,B,B$, $\N_\epsilon L$. The penetration
map $\cp_L$ is defined and continuous on $\bigcup_{\xi\in
\partial_\infty X-\partial_\infty C'} T^1_\xi X$. This point of view
is used in cases (3) and (4) of Proposition
\ref{prop:prescribe} below, and will be useful to apply Corollary
\ref{coro:main}.


%
%
%
%
\subsection{Prescribing the penetration}

In Section \ref{sec:main}, we will use the following operation
repeatedly: a geodesic ray or line $\gamma$ starting from a given
point $\xi_0$ is given that penetrates two $\epsilon$-convex sets $C$
and $C'$ with penetration maps $f$ and $f'$, first entering $C$
with $f(\gamma)=h$, and then $C'$ with $f'(\gamma)\geq h'$. We will
need to pick a new geodesic ray or line $\gamma'$ starting from
$\xi_0$ which intersects $C$ before $C'$, for which we still have
$f(\gamma')=h$, and for which we now have the equality
$f'(\gamma')=h'$.  In the following result, we show that this
operation is possible in a number of geometric cases.  These cases
will be used in Section \ref{sec:results} for various applications.

\bprop\label{prop:prescribe} %
Let $X$ be a complete, simply connected Riemannian manifold with 
sec\-tion\-al curvature at most $-1$ and dimension at least $3$.  Let
$\epsilon>0$ and $\delta,h,h'\ge 0$. Let $C$ and $C'$ be
$\epsilon$-convex subsets of $X$, and $\xi_0\in(X\cup \partial_\infty
X)- (C\cup \partial_\infty C)$. Let $f$ and $f'$ be maps $T^1_{\xi_0}
X\ra [0,+\infty]$, with $f'$ continuous and $\kappa'=
\norm{f'-\ell_{C'}} <+\infty$. Consider the following cases:
\begin{enumerate}
\item \label{1} %
  $C$ is a horoball with $\diam(C\cap C')\le\delta$; $f$ is either the
  penetration height map $\ph_C$ or the inner-projection penetration
  map $\ipp_C$; 
  $$
  h\ge h^{\rm min}=2 c'_1(\epsilon)+2\delta+
  \norm{f-\ph_C} 
  $$ 
  and $h'\geq h^{\rm min}_0=\kappa'+2\delta$; if $C'$ is also a horoball we
  may take $\epsilon=+\infty$ in the definition of $h^{\rm min}$.
\item \label{3} %
  $C$ is a ball of radius $R$ $(\ge\epsilon)$ with $\diam(C\cap
  C')\le\delta$; $f$ is either the penetration height map $\ph_C$ or
  the inner-projection penetration map $\ipp_C$;
  $$ 
  h^{\rm min}=2c'_1(\epsilon)+2\delta+ \norm{f-\ph_C}
  \;\le \;h\;\le\; 2R-2c'_1(\epsilon)- \norm{f-\ph_C} = h^{\rm max}
  $$ 
  and $h'\geq h^{\rm min}_0=\kappa'+2\delta$;
\item \label{2a} %
  $C$ is the $\epsilon$-neighbourhood of a complete totally geodesic
  subspace $L$ of dimension at least $2$, with $\diam(C\cap
  C')\le\delta$; either $f=\ell_C$ and $X$ has constant curvature, or
  $f$ is the fellow-traveller penetration map $\ftp_L$;
  $$
  h\ge h^{\rm min}=4c'_1(\epsilon)+2\epsilon+\delta+\norm{f-\ftp_L}
  $$ 
  and $h'> h^{\rm min}_0=\kappa'+\delta$;
\item \label{2b} %
  \begin{itemize}
  \item $C$ is the $\epsilon$-neighbourhood of a geodesic line $L$; 
  \item   $h\ge h^{\rm  min}=
  4c'_1(\epsilon)+2\epsilon+\delta+\norm{f-\ftp_L}\;$;  
\item either
  $f=\ell_C$ and $X$ has constant curvature, 
or $f=\cp_L$,  $\xi_0\in\partial _\infty X$ and the metric spheres of
the Hamenst\"adt distance on $\partial_\infty X-\{\xi_0\}$ are topological
spheres, 
  or $f$
  is the fellow-traveller penetration map $\ftp_L$; 
    \item either $C'$ is any
  $\epsilon$-convex subset that does not meet $C$ (in which case
  $\delta=0$) and $h'> h^{\rm min}_0=\kappa'$, or $C'$ is the
  $\epsilon$-neighbourhood of a totally geodesic subspace with
  codimension at least two such that $\diam(C\cap C') \le\delta$ and
  $$
  h'\geq h^{\rm min}_0=
  3c'_1(\epsilon)+3\epsilon+\delta+\norm{f'-\ftp_{L'}}\:.
  $$
  \end{itemize}
\end{enumerate}
Assume that one of the above cases holds. If there exists a geodesic
ray or line $\gamma$ starting from $\xi_0$ which meets first $C$ and
then $C'$ with $f(\gamma)=h$ and $f'(\gamma)\ge h'$, then there exists
a geodesic ray or line $\ov\ga$ starting from $\xi_0$ which meets first $C$ and
then $C'$ with $f(\ov\gamma)=h$ and $f'(\overline\gamma)=h'$.  
\eprop

\dem %
Let $\ga$ be as in the statement, and $x$ (resp.~$y$) be the point
where $\ga$ enters (resp.~exits) $C$ (with $y$ in $X$ as $f(\ga)= h
<+\infty$). Let $x'$ (resp.~$y'$) be the point where $\ga$ enters
(resp.~exits) $C'$ (with $x'\in X$ but possibly with $y'$ at
infinity). By convexity, $\xi_0\notin C'\cup \partial_\infty C'$. For
every $h\geq 0$, we define $A$ as the set of points $\alpha(+\infty)$
where $\alpha\in T^1_{\xi_0} X$  satisfies
$f(\alpha) =h $. Let $A_0$ be the arcwise connected component of $A$
containing $\ga(+\infty)$. By considering the various cases, we will
prove below the following two claims :
\begin{enumerate}
\item[a)] every geodesic ray or line, starting from $\xi_0$ and
  meeting $C'$, first meets $C$ and then $C'$;
\item[b)] there exists a geodesic ray or line $\overline\gamma_0$
  starting from $\xi_0$ with $\overline\gamma_0(+\infty)$ belonging to
  $A_0$, and $f'(\overline\gamma_0)\leq h^{\rm min}_0$.
\end{enumerate}
As $f'$ is continuous and $A_0$ is arcwise connected, the intermediate
value theorem implies the existence of a geodesic $\overline\gamma$
with the desired properties, and Proposition \ref{prop:prescribe} is
proven.

\medskip
\noindent {\bf Case \eqref{1}. } %
Let $\kappa= \norm{f-\ph_C}$. Let $\xi$ be the point at infinity of
$C$, which is different from $\xi_0$, and let $p_\xi$ be the closest
point to $\xi$ on $\gamma$. As $f(\ga)=h>0$, the point $p_\xi$ belongs
to the interior of the horoball $C$. Let $\ga_\xi$ be the geodesic
ray or line starting from $\xi_0$ with  $\ga_\xi(+\infty)=\xi$.

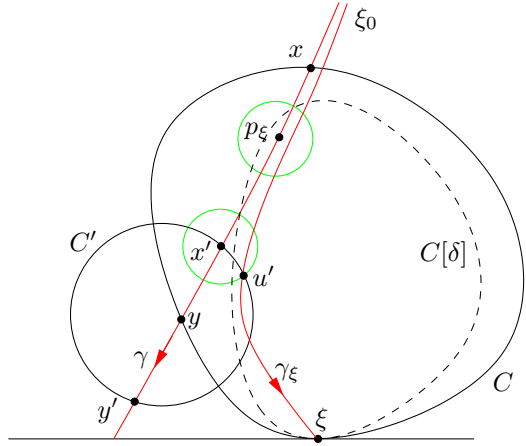
\begin{center}
\input{fig_case1new.pstex_t}
\end{center}

We start by proving the (stronger) first claim that every geodesic ray
or line starting from $\xi_0$ and meeting $C'$ meets $C[\delta]$ first
(hence meets $C$ before $C'$).  Note that
$$ 
d(y,p_\xi)\geq d(p_\xi,\partial C)= \frac{\ph_C(\ga)}{2}\geq
\frac{f(\ga)-\kappa}{2}= \frac{h-\kappa}{2}\geq \frac{h^{\rm
    min}-\kappa}{2}=c'_1(\epsilon)+\delta> \delta\;.
$$ 
As $f'(\ga)\geq h'\geq h_0^{\rm min}= \kappa'+2\delta$, we have
$\ell_{C'}(\ga)\geq f'(\ga)-\kappa'>\delta$, unless
$\ell_{C'}(\ga)=\delta=0$. Note that $\ga\cap C'$ is not contained in
the geodesic segment $[x,y]$. Otherwise, this would contradict the
assumption that $\diam(C\cap C')\leq \delta$ when
$\ell_{C'}(\ga)>\delta$. When $\ell_{C'}(\ga)=\delta=0$, as $\ga$
meets $C'$, the segment $\ga\cap C'$ would be reduced to a point by
the convexity of $C'$, which would be $\{x\}$ or $\{y\}$ (as $C'$ is
not a singleton). But then the tangent vector of $\ga$ at $x$ or its
opposite at $y$ would both enter strictly $C$ and be tangent to $C'$,
which contradicts the fact that $\delta=0$.

As $\ga$ meets $C$ before $C'$, this implies in particular that the
geodesic ray $[y,\ga(+\infty)[$ meets $C'$, and that the point $p_\xi$
belongs to $]x',\xi_0]$\;: otherwise $C\cap C'$ would contain a
segment of length at least $d(p_\xi,y)> \delta$, which is impossible.
Hence by convexity, any geodesic ray or line, starting from $\xi_0$
and meeting $B(x',c'_1(\epsilon))$, first meets
$B(p_\xi,c'_1(\epsilon))$.  By Lemma \ref{lemma:c1}, every geodesic
ray or line, starting from $\xi_0$ and meeting $C'$, meets the ball
$B(x',c'_1(\epsilon))$ before entering $C'$. This proves the first
claim, as the ball $B(p_\xi,c'_1(\epsilon))$ is contained in
$C[\delta]$, since $ d(p_\xi,\partial C)\geq c'_1(\epsilon)+\delta$,
as seen above.

\medskip %
Let us prove now the (stronger) second claim that there exists a geodesic
ray or 
line $\overline\gamma_0$ starting from $\xi_0$ with $\overline\gamma_0
(+\infty)$ belonging to $A_0$, and avoiding the interior of $C'$,
which implies the result, as then $f'(\overline\gamma_0)\leq
\ell_{C'}(\overline\ga_0) +\kappa'\le h^{\rm min}_0$.

The subspace $A$ of $\partial_\infty X$ is a codimension $1$
topological submanifold of the topological sphere $\partial_\infty X$,
which is homeomorphic to the sphere $\SS^{n-2}$, hence is arcwise
connected. Indeed, if $f=\ph_C$, then $A$ is the subset of endpoints
of the geodesic rays or lines starting from $\xi_0$ that are tangent
to $\partial\big(C[h/2]\big)$. If $f=\ipp_C$, the subset $A$ is the
preimage of a point in $]\xi_0,\xi[$ by the closest point map from
$\partial_\infty X$ to $[\xi_0,\xi]$, which is, over $]\xi_0,\xi[$, a
trivial topological bundle with fibers homeomorphic to $\SS^{n-2}$.

Note that $f$ is continuous, $f(\ga_\xi)=\infty>h$, and $f(\alpha)=0$
if $\alpha$ is a geodesic ray or line starting from $\xi_0$ with
$\alpha(+\infty)$ close enough to $\ga(-\infty)$. Therefore $A$
separates $\ga(-\infty)$ and $\ga_\xi(+\infty)$, as the connected
components of $\partial_\infty X-A$ are arcwise connected.

If the (stronger) second claim is not true, then the topological sphere $A_0=A$
of dimension $n-2$ is contained in the interior of the shadow
$\O_{\xi_0}C'$.  As $\xi_0\notin C'\cup \partial_\infty C'$, this
shadow is homeomorphic to a ball of dimension $n-1$. Thus, by Jordan's
theorem, one of the two connected components of $\partial_\infty X-A$
is contained in the interior of $\O_{\xi_0}C'$. As $\ga(-\infty)$ does
not belong to $\O_{\xi_0}C'$ and $A$ separates $\ga(-\infty)$ and
$\ga_\xi(+\infty)$, this implies that $\ga_\xi(+\infty)$ belongs to
the interior of $\O_{\xi_0}C'$. Hence $\ga_\xi$ meets the interior of
$C'$.

Therefore, by the first claim, the geodesic ray or line $\ga_\xi$ meets
$C[\delta]$ before meeting $C'$. Let $u'$ be the entering point of
$\ga_\xi$ in $C'$.  As $\xi$ is the point at infinity of $C[\delta]$,
the points $\xi_0,u',\xi$ are in this order on $\ga_\xi$. Hence by
convexity, this implies that $u'$ belongs to $C[\delta]$. As
$f'(\ga)\geq h'\geq h^{\rm min}_0=\kappa'+2\delta$, we have
$$
d(x',y')=\ell_{C'}(\ga)\geq f'(\ga)-\kappa'\geq 2\delta\;.
$$ 
Hence by
the triangular inequality, one of the two distances $d(u',x'),
d(u',y')$ is at least $\delta$, and by the strict convexity of the
distance, is strictly bigger than $\delta$ (as $u'$ does not belong to
$\ga$ (as $\ga\neq \ga_\xi$). Hence, if $u''$ is a point close enough to
$u$ in $]u,\xi[$, then $u''$ belongs to the interior of $C'$ and to the
interior of $C[\delta]$, and is at distance strictly greater than $\delta$
from either $x'$ or $y'$. Therefore, some   
geodesic segment of length
strictly bigger than $\delta$ is contained in the intersection $C\cap
C'$. This contradicts the assumption that $\diam (C\cap C')\le\delta$.

\bigskip
\noindent{\bf Case \eqref{3}. } The proof is completely similar to
Case \eqref{1}.  Let now $C=B(z,R)$, let $p_z$ be the point of $\ga$
the closest to $z$, and let $\ga_z$ be the geodesic ray or line
starting from $\xi_0$ and passing through $z$. Note that $R>h^{\rm
  max}/2 \geq h^{\rm min}/2\geq \delta$. We only have to replace $\xi$
by $z$, $p_\xi$ by $p_z$ and $\ga_\xi$ by $\ga_z$, and to replace two
arguments in the above proof, the one in order to show that $A$
separates $\ga(-\infty)$ from $\ga_z(+\infty)$, and the one in order
to show that $\xi_0,u',z$ are in this order on $\ga_z$, where $u'$ is
the entering point of $\ga_z$ in $C'$.

To prove that $A$ separates $\ga(-\infty)$ from $\ga_z(+\infty)$, we
simply use now that $f(\ga_z)=2R>h^{\rm max}\geq h$ instead of
$f(\ga_\xi)=\infty>h$. Let us prove that $\xi_0,u',z$ are in this
order on $\ga_z$. We have, with $\kappa=\norm{f-\ph_C}$,
$$
d(z,p_z)=R-\frac{\ph_C(\ga)}{2}\geq R-\frac{f(\ga)+\kappa}{2}\geq
R-\frac{h^{\rm max}+\kappa}{2} =c'_1(\epsilon)\;.
$$
By Lemma 1.3, we have $d(x',u')\leq c'_1(\epsilon)$. As $\ga_z$
meets the interior of $C'$, by the same argument as in Case 1, we even have
$d(u',\ga)< c'_1(\epsilon)$. 
Hence by strict convexity, we do have $u'\in \;]z,\xi_0[$. The rest of
the argument in the proof of Case \eqref{1} is unchanged.

\bigskip %
Before studying the last two cases, we start by proving two lemmas.
The first one implies the first of the two claims we need to prove in
Cases \eqref{2a}, \eqref{2b}, and the second one gives the topological
information on $A$ that we will need in these last two cases.

\blemm \label{lem:meetbeforelastcas}%
Let $L$ be a complete totally geodesic subspace with dimension at least $1$,
$\epsilon>0$, $C=\N_\epsilon L$, $\xi_0\in (X\cup\partial_\infty X) -
(C\cup\partial_\infty L$), and $C'$ be an $\epsilon$-convex subset of
$X$ such that $\diam(C\cap C')\leq \delta$. Let $f,f': T^1_{\xi_0}X
\ra [0,+\infty]$ be maps such that $\kappa= \norm{f-\ftp_L}<+\infty$,
$\kappa'=\norm{f'-\ell_{C'}}< +\infty$. Let $\ga$ be a geodesic ray or
line starting from $\xi_0$, entering $C$ before entering $C'$, such
that $ 4c'_1(\epsilon) + 2\epsilon + \delta+\kappa\le f(\ga)\le +\infty$ and
$f'(\ga)> \delta+\kappa'$. If $\tilde\gamma$ is a geodesic ray or line
starting from $\xi_0$ which meets $C'$, then $\tilde\gamma$ meets the
interior of $C$ before meeting $C'$.  
\elemm

\dem %
Note that $\xi_0\notin C'\cup\partial_\infty C'$, by convexity and the
assumptions on $\ga$, as $\xi\notin C\cup\partial_\infty C$. Let $L_0$ be the geodesic line passing through
the closest points $p_{\xi_0}, p_{\ga(+\infty)}$ on $L$ of
$\xi_0,\ga(+\infty)$, respectively. Note that
$$
d(p_{\xi_0},p_{\ga(+\infty)})=\ftp_L(\ga)\geq f(\ga)-\kappa 
\geq 4c'_1(\epsilon)+2\epsilon+\delta>0\;.
$$
Hence, by Lemma \ref{fellowtravelpenemap}, and as $\ftp_{L_0}(\gamma)
=\ftp_{L}(\gamma)$, we have 
$$
\ell_{\N_\epsilon L_0}(\ga) \geq
\ftp_{L_0}(\gamma)-2c'_1(\epsilon)- 2\epsilon>0.
$$
In particular, $\ga$
enters $\N_\epsilon L_0$ at a point $x_0$ and exits it at a point
$y_0$ in $X$ (as $\ga(+\infty)\notin \partial_\infty L$). Let
$u\mapsto p_u$ be the closest point map from $X\cup\partial_\infty X$
onto $L_0\cup\partial_\infty L_0$. Recall that this map does not increase
the distances (and even decreases them, unless the two points under
consideration are on $L_0$), and that it {\it preserves 
  betweenness}, that is, if $u''\in[u,u']$, then $p_{u''}\in
[p_u,p_{u'}]$. Let $x'$ (resp.~$\tilde x'$) be the point where $\ga$
(resp.~$\tilde\gamma$) enters $C'$, and $q_{\xi_0}$ and
$q_{\ga(+\infty)}$ be the closest point to $\xi_0$ and $\ga(+\infty)$
respectively on $\N_\epsilon L_0$.

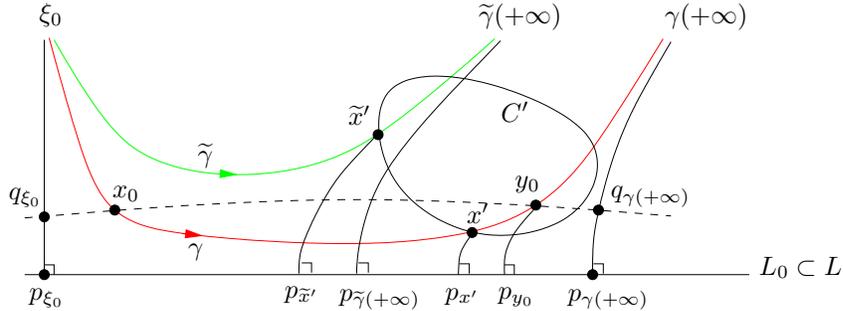
\begin{center}
\input{fig_case2a2bfm.pstex_t}
\end{center}

Recall that by Lemma \ref{lemma:c1}, the distances $d(\tilde x', x'),
d(x_0,q_{\xi_0}), d(y_0,q_{\ga(+\infty)})$ are at most $c'_1(\epsilon)$.
Note that $\tilde x'\in[\xi_0,\tilde\ga(+\infty)]$.  Hence, as
betweennes is preserved, 
\begin{align*}
\ftp_{L_0}(\tilde\gamma)& 
=d(p_{\xi_0},p_{\tilde \ga(+\infty)})\ge d(p_{\xi_0},p_{\tilde x'})
 \ge d(p_{\xi_0},p_{x'})-d(p_{x'},p_{\tilde x'})\\
& \ge d(p_{\xi_0},p_{x'})-d(x',\tilde x')\ge 
d(p_{\xi_0},p_{x'})-c'_1(\epsilon)\;.  
\end{align*}
Note that $d(p_{\xi_0},p_{x'})\ge d(p_{\xi_0},p_{y_0})$ when
$\xi_0,y_0,x'$ are in this order on $\ga$. When $\xi_0,y_0,x'$ are not
in this order on $\ga$, as $\ga$ enters in $C$ before $C'$, as
$\ell_{C'}(\ga)\geq f'(\ga)-\kappa'> \delta$ and as $\diam(C\cap
C')\leq\delta$, we have $d(x',y_0)\leq \delta$; hence
$$
d(p_{\xi_0},p_{x'})\geq d(p_{\xi_0},p_{y_0})-d(p_{y_0},p_{x'})\ge
d(p_{\xi_0},p_{y_0})-d(y_0,x')\geq d(p_{\xi_0},p_{y_0})-\delta\;.
$$
Therefore, in both cases, as $y_0\in[\xi_0,\ga(+\infty)]$ and
$u\mapsto p_u$ preserves the betweenness, and since
$p_{\ga(+\infty)}=p_{q_{\ga(+\infty)}}$, we have
\begin{align*}
\ftp_{L_0}(\tilde\gamma) & \geq d(p_{\xi_0},p_{\tilde x'}) 
\geq  d(p_{\xi_0},p_{y_0})-\delta-c'_1(\epsilon) \\ &\ge
d(p_{\xi_0},p_{\ga(+\infty)})- d(p_{q_{\ga(+\infty)}},p_{y_0})-
c'_1(\epsilon)-\delta\\ &> \ftp_L(\ga)- d(q_{\ga(+\infty)},y_0)-
c'_1(\epsilon)-\delta\ge \ftp_L(\ga)-2c'_1(\epsilon)-\delta
\geq 2c'_1(\epsilon)+2\epsilon \;.
\end{align*}
By Lemma \ref{fellowtravelpenemap}, we hence have 
$$
\ell_{\N_\epsilon
  L}(\tilde\ga) \geq\ell_{\N_\epsilon L_0}(\tilde\ga) \geq \ftp_{L_0}
(\tilde\gamma)-2c'_1(\epsilon)- 2\epsilon>0.
$$ 
In particular,
$\tilde\ga$ does enter the interior of $C$, at a point $\wt x$. Note
that the geodesic from $\xi_0$ through $p_{\xi_0}$ enters $C$ at
$q_{\xi_0}$. Now by absurd, if $\tilde\ga$ enters the interior of $C$
after it enters $C'$, then $\tilde x'\in[\xi_0,\tilde x]$, so that
$$
c'_1(\epsilon)\geq d(q_{\xi_0},\tilde x)\geq 
d(p_{\xi_0},p_{\tilde x})\geq d(p_{\xi_0},p_{\tilde x'})>
2c'_1(\epsilon)+2\epsilon\;,
$$
as seen above, a contradiction.
\cqfd

\blemm\label{lem:topopeneprescri} %
Let $X$ be a complete, simply connected Riemannian manifold with
sectional curvature at most $-1$ and dimension at least $3$.  Let
$\epsilon,h>0$. Let $L$ be a complete totally geodesic submanifold with
dimension at least $1$ and $\xi_0\in(X\cup \partial_\infty X)-(\N_\epsilon
L\cup \partial_\infty L)$.  Assume either that 
\begin{enumerate}
\item $f=\ell_{\N_\epsilon L}$, $X$
has constant curvature and $h\in [4c'_1(\epsilon) + 2\epsilon,
+\infty[$, or 
\item 
$f=\cp_{ L}$, 
$\xi_0\in\partial_\infty X$, $\dim L = 1$, $h\in\left]\log 2,
+\infty\right[$, 
and the metric spheres of
the Hamenst\"adt distance on $\partial_\infty X-\{\xi_0\}$ are
topological spheres, 
or 
\item $f=\ftp_L$. 
\end{enumerate}
Then 
$$
A=\{\alpha(+\infty) \;:\;
\alpha\in T^1_{\xi_0}X,\; f(\alpha)=h\}
$$ 
is a codimension $1$ topological submanifold of the topological sphere
$\partial_\infty X$, which is homeomorphic to the torus $\SS^{\dim
  L-1}\times \SS^{\codim L-1}$. Furthermore,
\begin{itemize}
\item[\rm(a)] if $\dim L= 1$, then $A$ has two arcwise connected
  components, homeomorphic to a sphere of dimension $n-2$.
  If $f=\cp_L$ or if $h>c'_1(\epsilon)$, then  each of
  them separates $\ga(-\infty)$ and exactly one of the two points at
  infinity of $L$, for every geodesic ray or line $\ga$ starting from
  $\xi_0$ if $f=\cp_L$, and for those meeting $\N_\epsilon L$ if $f\ne\cp_L$.
\item[\rm(b)] if $\codim L= 1$, then $A$ has two arcwise connected
  components, homeomorphic to a sphere of dimension $n-2$, separated
  by $\partial_\infty L$. 
\item[\rm(c)] if $\dim L\geq 2$ and $\codim L\ge 2$, then $A$ is arcwise
  connected. 
\end{itemize}
In cases (b) and (c), for every component $A'$ of $A$, for every
  geodesic ray $\rho$ in $L$ with $\rho(0)$ the closest point to
  $\xi_0$ on $L$, there exists $\eta\in A'$ such that $\rho(h)$ is at
  distance at most $\norm{f-\ftp_L}$ from the closest point to $\eta$
  on $\rho$.
\elemm

\dem %
Let $\pi_L: X\cup \partial_\infty X\ra L\cup \partial_\infty L$ be the
closest point map, $p_0=\pi_L(\xi_0)$ and $S_0=
\pi_L^{-1}(p_0)\cap \partial_\infty X$.

Assume first that $f=\ftp_L$. As $h>0$, the subspace $A$ of
$\partial_\infty X$ is the preimage of the sphere (of dimension $\dim
L-1$) of center $p_0$ and radius $h$ in $L$, by $\pi_L$. As
$\pi_L:\partial_\infty X\setminus\partial_\infty L\to L$ is a trivial
topological bundle whose fibers are spheres of dimension $\codim L-1$,
the topological structure of $A$ is immediate.  The final statement on
(b) and (c) is trivial as, by definition, $\rho(h)$ is the closest
point to some point in $A'$.

If $\dim L=1$ and $h>c'_1(\epsilon)$, if $\ga\in T^1_{\xi_0} X$ meets
$\N_\epsilon L$, then by Lemma \ref{lemma:c1} and by convexity,
$d(\pi_L(\ga(-\infty)), p_0) \leq c'_1(\epsilon)<h$. Hence, the
separation statement in (a) follows.
\medskip

Assume now that $f=\cp_{L}$, and that the hypotheses of (2) are satisfied.
The result in this case follows from the discussion  
before Lemma \ref{lem:crossratio}.

\medskip
Assume now that $f=\ell_{\N_\epsilon L}$, $X$ has constant curvature,
and $h\in[4c'_1(\epsilon) + 2\epsilon,+\infty[\,$. Using normal
coordinates along $L$, the topological sphere $\partial_\infty X$ is
homeomorphic to the topological join of the spheres
$\partial_\infty L$ of dimension $\dim L-1$ and $S_0$ of dimension
$\codim L-1$
$$
S_0\vee \partial_\infty L=
\big(S_0\times [0,+\infty]\times \partial_\infty L\big)/\sim\;,
$$
where $\sim$ is the equivalence relation generated by $(a,0,b)\sim
(a,0,b')$ and $(a,+\infty,b)\sim(a',+\infty,b)$, for every $a,a'$ in
$S_0$ and $b,b'$ in $\partial_\infty L$. We denote by $[a,t,b]$ the
equivalence class of $(a,t,b)$. We choose the parametrization of
$\partial_\infty X$ by $S_0\vee \partial_\infty L$ such that
$[a,0,b]=a$, $[a,+\infty,b]=b$, $d(\pi_L([a,t,b]),p_0)=t$, and the
geodesic rays $\big[\pi_L([a,t,b]),[a,t,b]\big[$ are parallel transports of
$[p_0,a[$ along the geodesic ray $[p_0,b[\,$, for $0<t<+\infty$.

For every $t$ in $]0,+\infty[$ and every $(a,b)$ in $S_0 \times
\partial_\infty L$, let $\ga_{[a,t,b]}$ be the geodesic ray or line
starting from $\xi_0$ and ending at $[a,t,b]$. By the proof of Lemma
\ref{fellowtravelpenemap} and by Lemma \ref{lem:boundaryprojpenemap}, we
have  
$$ 
-2c'_1(\epsilon)\leq
\ell_{\N_\epsilon L} (\alpha) -\ftp_L(\alpha) \leq 2c'_1(\epsilon) +
2\epsilon,
$$ 
for every $\alpha\in T^1_{\xi_0}X$.  In particular, if
$t=\ftp_L(\ga_{[a,t,b]}) > 2c'_1(\epsilon)$, then $\ell_{\N_\epsilon
  L} (\ga_{[a,t,b]})>0$, that is $\ga_{[a,t,b]}$ meets the interior of
$\N_\epsilon L$. By reducing to the case $X=\HH_\RR^3$ and $L$ a
geodesic line, it is easy to see 
that the map from $[2c'_1(\epsilon), +\infty[$ to $[0,+\infty[$
defined by $t\mapsto \ell_{\N_\epsilon L}(\ga_{[a,t,b]})$ is
continuous and strictly increasing, for every fixed $(a,b)$ in $S_0
\times
\partial_\infty L$. 
Hence, as 
$$
\ell_{\N_\epsilon L} (\ga_{[a,2c'_1(\epsilon),b]})\leq
4c'_1(\epsilon) +2\epsilon\leq h < +\infty,
$$ 
there exists a unique
$t_{a,b}\in[2c'_1(\epsilon),+\infty[$, depending continuously on
$(a,b)$, such that $\ell_{\N_\epsilon L}(\ga_{[a,t_{a,b},b]})=h$. In
particular, the subset of points of $\partial_\infty X$ of the form
$[a,t_{a,b},b]$ for some $(a,b)$ in $S_0 \times \partial_\infty L$ is
indeed a codimension $1$ topological submanifold of $\partial_\infty
X$, which is homeomorphic to the torus $\SS^{\dim
  L-1}\times\SS^{\codim L-1}$. The statements (b) and (c) follow.

If $L$ has dimension $1$, and if $\ga\in T^1_{\xi_0} X$ meets
$\N_\epsilon L$, then by Lemma \ref{lemma:c1} and by convexity,
$d(\pi_L(\ga(-\infty)), p_0) \leq c'_1(\epsilon)$. For every $\xi$ in
a component $A_0$ of $A$, if as above $\xi=[a,t_{a,b},b]$, then we
have $d(\pi_L(\xi), p_0)=t_{a,b}\geq 2c'_1(\epsilon)$, hence $A_0$
separates $\ga(-\infty)$ and $b$. This proves (a).

Let us prove the last assertion of the lemma.  Let $\kappa=
\norm{f-\ftp_L}$, and let $A'$ be a connected component of $A$.  For
every $u$ in $L$ such that $d(u,p_0) =h$, let $\eta_0=[a,h,b]$, on the
same side of $\partial_\infty L$ as $A'$ if $\codim L= 1$, be such
that $\pi_L (\eta_0) =u$.  Let $\eta_t=[a,h+t,b]$, which is on the
same side of $\partial_\infty L$ as $A'$ if $\codim L= 1$.  Note that
$$
f(\ga_{[a,h+\kappa,b]})\geq
\ftp_L(\ga_{[a,h+\kappa,b]})-\kappa=h,
$$ and similarly,
$f(\ga_{[a,h-\kappa,b]})\leq h$.  By the intermediate value theorem,
there exists $t\in[-\kappa, +\kappa]$ such that $\eta_t\in A'$. Hence
$d(u,\pi_L(\eta_t))=|t|\leq \kappa$.  
\cqfd

\bigskip

Now we proceed with the proof of the the remaining parts of
Proposition \ref{prop:prescribe}.

\medskip
\noindent{\bf Case \eqref{2a}. }  %
By Lemma \ref{lem:meetbeforelastcas}, we only have to prove the second
claim, that there exists a geodesic ray or line $\ov{\gamma_0}$ starting
from $\xi_0$ with $\ov{\ga_0}(+\infty)$ belonging to $A_0$, such that
$f'(\ov{\ga_0})\leq h_0^{\rm min}$.

Let $\kappa=\norm{f-\ftp_L}$. Let $p_0$ (resp.~$p_\ga$) be the point
of $L$ the closest to $\xi_0$ (resp.~$\ga(+\infty)$), so that, in
particular, 
$$
d(p_0,p_\ga)=\ftp_L(\ga)\geq f(\ga)-\kappa=h-\kappa>0.
$$
Let $p_{\ga}'$ be the point on the geodesic line $L_0$ (contained in
$L$) passing through $p_0$ and $p_\ga$ on the opposite side of
$p_{\ga}$ with respect to $p_0$, and at distance $h$ from $p_0$. By
Lemma \ref{lem:topopeneprescri}, there exists a geodesic line
$\ov{\gamma_0}$ starting from $\xi_0$ and ending at a point in
$A_0$ whose closest point $p_{\ov{\ga_0}}$ on $L$ is at distance at most
$\kappa$ from $p_{\ga}'$.

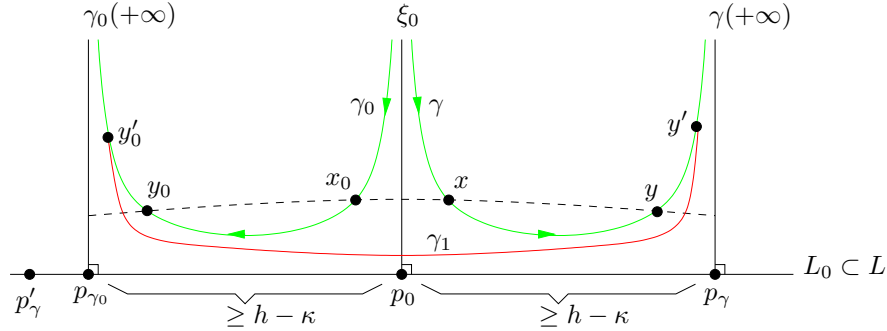
\begin{center}
\input{fig_case2a.pstex_t}
\end{center}

Assume by absurd that $f'(\ov{\ga_0})> h_0^{\rm min}$. We have
$$
\ell_{C'}(\ga)\geq f'(\ga)-\kappa'\geq h'-\kappa'\geq
h_0^{min}-\kappa'> \delta.
$$ 
Similarly $\ell_{C'}(\ov{\ga_0})> \delta$, and,
in particular, $\ov{\ga_0}$ enters $C'$. Let $y'$ (resp.~$y'_0$) be the
point, possibly at infinity, where $\ga$ (resp.~$\ov{\ga_0}$) exits $C'$.
By Lemma \ref{lem:meetbeforelastcas}, $\ov{\ga_0}$ meets $C$ before $C'$.
Let $x_0$ (resp.~$y_0$) be the point where $\ov{\ga_0}$ enters in
(resp.~exits) $C$. As $\diam(C\cap C')\leq \delta$, we have
$y'_0\in \; ]y_0,\ga_0(+\infty)[$ and $y'\in \;]y,\ga(+\infty)[$, so
that in particular $d(y'_0,L)> \epsilon$ and $d(y',L)> \epsilon$.

Let $\ga_1$ be the geodesic line through $y'$ and $y'_0$. The points
at infinity of $\gamma_1$ do not belong to $\partial_\infty L_0$, so that
$\ftp_{L_0}(\ga_1)$ and  $\ell_C(\ga_1)$ are finite. 
Note that by strict
convexity and by Lemma \ref{lemma:c1}, we have 
$$
d(y',[p_{\ga},
\ga(+\infty)[)< d(y,[p_{\ga}, \ga(+\infty)[) \leq c'_1(\epsilon),
$$ 
and
similarly $d(y'_0,[p_{\ov{\ga_0}}, \ov{\ga_0}(+\infty)[)$ $< c'_1(\epsilon)$.
Hence, with $\pi_{L_0}$ the closest point map to $L_0$, which
preserves the betweenness and does not increase the distances,
\begin{align*}
\ftp_{L_0}(\ga_1)& 
\geq d(\pi_{L_0}(y'_0),\pi_{L_0}(y'))> 
d(p_{\ov{\ga_0}},p_{\ga})-2c'_1(\epsilon)=
d(p_{\ov{\ga_0}},p_0)+d(p_0,p_{\ga})-2c'_1(\epsilon)\\ & \ge 
h-\kappa +h-\kappa-2c'_1(\epsilon)=2h -2\kappa-2c'_1(\epsilon)\;.  
\end{align*}
In particular, by Lemma \ref{fellowtravelpenemap},
$$
\ell_C(\ga_1)\geq \ell_{\N_{\epsilon} L_0}(\ga_1)\geq 
\ftp_{L_0}(\ga_1)- 2c'_1(\epsilon)-2\epsilon> 
2h -2\kappa-4c'_1(\epsilon)-2\epsilon\geq\delta\;,
$$ 
by the definition of $h^{\rm min}$. Hence $\ga_1$ meets $C$ in a
segment $I$ of length $>\delta$. But as $y'_0$ and $y'$ are at
a distance strictly bigger than $\epsilon$ of $L$, the segment $I$ is
contained in $[y',y'_0]$, which is contained in $C'$, by convexity.
This contradicts the assumption that $\diam (C\cap C')\leq\delta$.

\medskip
\noindent{\bf Case \eqref{2b}. }  %
Let $\kappa''=\norm{f'-\ftp_{L'}}$. Note that $f'(\ga)\ge h'\ge h_0^{\rm
  min}> \delta+ \kappa'$. This is true under both assumptions on the value
of $h^{\rm  min}_0$, as when $h^{\rm min}_0= 3c'_1(\epsilon) + 3\epsilon +
\delta + \kappa''$, we have, by Lemma \ref{fellowtravelpenemap},
$$
\delta+\kappa'\leq \delta+\kappa''+2c'_1(\epsilon)+2\epsilon < 
h^{\rm min}_0\;.
$$
By Lemma \ref{lem:meetbeforelastcas}, we only have to prove the second
claim that there exists a geodesic line $\overline{\gamma}_0$ starting
from $\xi_0$ with $\overline{\ga}_0(+\infty)$ belonging to $A_0$, such
that $f'(\overline{\ga}_0)\leq h_0^{\rm min}$.

\medskip %
We first consider the case $C'=\N_\epsilon L'$ where $L'$ is a totally
geodesic subspace of codimension at least $2$, with $\diam(C\cap
C')\leq \delta$, and $h'\geq h^{\rm min}_0= 3c'_1(\epsilon) +
3\epsilon + \delta + \kappa''$. Assume by absurd that every geodesic
ray or line $\alpha$ starting from $\xi_0$ with $\alpha(\infty)\in
A_0$ meets $C'$ with $f'(\alpha)> h_0^{\rm min}$.  Let
$$
B'=\big\{\beta(\infty):\beta\in T^1_{\xi_0}X,\; 
\ftp_{L'}(\beta)> h_0^{\rm min}-\kappa'' \big\}\;.
$$
By the absurdity hypothesis and the definition of $\kappa''$, we have
$A_0\subset B'$. Let $p'_0$ be the closest point to $\xi_0$ on
$L'$. Note that $B'$ is a (topological) open tubular neighbourhood of
$\partial_\infty L'$, whose fiber over a point $\xi$ in
$\partial_\infty L'$ is the preimage of $\rho_\xi (]h_0^{\rm min}-
\kappa'',+\infty])$ by the closest point map from $\partial_\infty X$
to $L'\cup\partial_\infty L'$, where $\rho_\xi$ is the geodesic ray
with $\rho(0)=p'_0$ and $\rho(+\infty)=\xi$.

By Lemma \ref{lem:topopeneprescri}(a), let $\xi_1$ be the point at
infinity of $L$ separated from $\ga(-\infty)$ by $A_0$. Note that as
$\ga$ enters $C'$ at $x'$, and $\xi_0\notin C'$, if
$p'_{\ga(-\infty)}$ is the closest point to $\ga(-\infty)$ on $L'$,
then by Lemma \ref{lemma:c1}, we have $d(p'_{\ga(-\infty)},p'_0)\leq
c'_1(\infty)<h_0^{\rm min}-\kappa''$ by the definition of $h_0^{\rm
  min}$. Hence the complement of $B'$ in 
$\partial_\infty X$, which is connected as $\codim L'\ge 2$, contains
 $\ga(-\infty)$. As $A_0$ separates $\ga(-\infty)$
from $\xi_1$ and is contained in $B'$, it follows that $B'$ contains $\xi_1$.

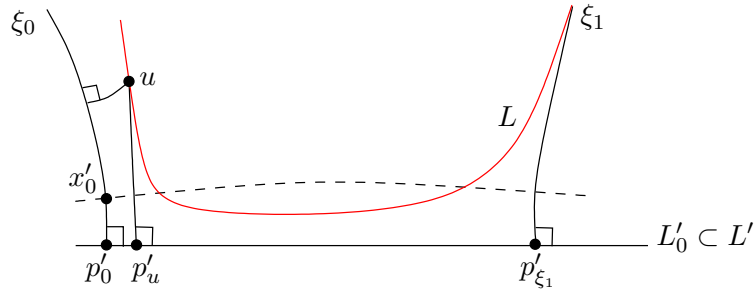
\begin{center}
\input{fig_totgeod.pstex_t}
\end{center}

Let $x'_0$ be the intersection point of $]\xi_0,p'_0]$ with $\partial
C'$. Lemma \ref{lemma:c1} implies that $d(x',x'_0)\le c'_1(\epsilon)$.
Hence, by convexity and as $\ga$ first meets $C$ and then $C'$, we
have $d\big(x,]\xi_0,p'_0]\big)\le c'_1(\epsilon)$, which implies that
there is a point $u$ in $L$ at distance at most $c'_1(\epsilon)
+\epsilon$ from $]\xi_0,p'_0]$. Let $p'_{\xi_1}$ be the closest point
to $\xi_1$ on $L'$, and $L'_0$ the geodesic line (contained in $L'$)
through $p'_0$ and $p'_{\xi_1}$. As the closest point map does
not increase distances, the closest point $p'_u$ to $u$ on $L'_0$
satisfies $d(p'_0,p'_u)\le c'_1(\epsilon)+\epsilon$.  Then, as the
closest point map to $L'_0$ preserves the betweenness and as
$\xi_1$ belongs to $B'$,
$$
\ftp_{L'_0}(L) \geq d(p'_u,p'_{\xi_1})\geq 
d(p'_{\xi_1},p'_0)-d(p'_0,p'_u)> h_0^{\rm min}-
\kappa''-c'_1(\epsilon)-\epsilon\;.
$$
Therefore, using Lemma \ref{fellowtravelpenemap},
\begin{align*}
\diam(C\cap C')&\geq \diam(C\cap \N_\epsilon L'_0)\geq
\ell_{\N_\epsilon L'_0}(L)\geq \ftp_{L'_0}(L) 
-2c'_1(\epsilon)- 2\epsilon\\ & 
> h_0^{\rm min}- \kappa''-3c'_1(\epsilon)- 3\epsilon=\delta\;,
\end{align*}
a contradiction.

\medskip %
Assume now that $C'$ is any $\epsilon$-convex subset such that $C\cap
C'=\emptyset$, and that $h'>h_0^{\rm min}=\kappa''$.  Let us
prove that there exists a geodesic ray or line $\overline\gamma_0$
starting from $\xi_0$ with $\overline\gamma_0(+\infty)$ in $A_0$, and
avoiding $C'$.  This implies the result as in Case \eqref{1}.

By absurd, suppose that for every $\xi$ in $A_0$, the geodesic ray or
line $\ga_\xi$ starting from $\xi_0$ and ending at $\xi$ meets $C'$.
By the first claim (see Lemma \ref{lem:meetbeforelastcas}), $\ga_\xi$
meets the interior of $C$ before meeting $C'$.  Let $x'_\xi$ be the
entering point of $\ga_\xi$ in $C'$ and $y_\xi$ be its exiting point
out of $C$.  As $C$ and $C'$ are disjoint, note that $\xi_0,
y_\xi,x'_\xi, \xi$ are in this order along $\ga_\xi$. The maps
$\xi\mapsto y_\xi$ and $\xi\mapsto x'_\xi$ are injective and
continuous on $a_0$ (by the strict convexity of $C$, as $\ga_\xi$ meets the 
interior of $C$). We know that $A_0$ is a topological sphere, by Lemma
\ref{lem:topopeneprescri}(a), separating the endpoints of $L$.  Hence the
subsets $A_0$ and $S'=\{x'_\xi\;:\xi\in A_0\}$ are spheres, that are
homotopic (by the homotopy along $\ga_\xi$) in the complement of $L$
in $X\cup\partial_\infty X$. By an homology argument, every disc with
boundary $S'$ in $X\cup\partial_\infty X$ has to meet $L$. But by
convexity of $C'$, there exists a disc contained in $C'$ with boundary
$S'$ (fix a point of $S'$ and take the union of the geodesic arcs from
this point to the other points of $S'$). This contradicts the fact
that $C\cap C'=\emptyset$.  
\cqfd

\medskip 
\noindent{\bf Remarks. } %
(1) In Case \eqref{3}, we have $h^{\rm max}\ge h^{\rm min}$ if $R$ is
big enough, as $c'_1(\epsilon)$ has a finite limit as
$\epsilon\to\infty$.

\medskip %
(2) In Case \eqref{2a}, if the codimension of $L$ is $1$, then we may
assume that $\overline{\ga}$ meets $L$ if $\ga$ meets $L$.  Indeed, as
we have seen in Lemma \ref{lem:topopeneprescri}(b), $L\cup\partial_\infty
L $ separates $X\cup \partial_\infty X$ into two connected components,
and $A$ (defined in the beginning of the proof) has exactly two
components separated by $L\cup\partial_\infty L$. If $A_0^+$ is the
component of $A$ on the same side of $\xi_0$ from
$L\cup\partial_\infty L$, and $A_0^-$ the component of $A$ on the
other side, then a geodesic ray or line starting from $\xi_0$ and
ending in $A_0^+$ does not meet $L$ (as $L$ is totally geodesic), and
any geodesic line starting from $\xi_0$ and ending in $A_0^-$ meets
$L$, by separation.  This observation on the crossing property will be
used in the proof of Corollary \ref{coro:max+recurrent} to make sure
that the locally geodesic ray or line constructed in the course of the
proof stays in the convex core.

\medskip %
(3) Case \eqref{2b} is not true if $C'$ is assumed to be any
$\epsilon$-convex subset, as shown by taking $X$ the real hyperbolic $3$-space,
and $C'$ the $\epsilon$-neighborhood of the (totally geodesic)
hyperbolic plane perpendicular to $L$ at a point at distance $h$ from
the closest point to $\xi_0$ on $L$: any geodesic ray or
line $\alpha$ starting from $\xi_0$, with $\ftp_L(\alpha) =h$ and
meeting $C'$ satisfies $f'(\alpha) =+\infty$ for every $f'$ which is a
$\kappa'$-penetration map in $C'$.

\section{The main construction}
\label{sec:main}

%
%
%

\subsection{Unclouding the sky}
\label{subsec:uncloudingshort}

The aim of this section is to prove the following result, improving on
our result in \cite{PP}. The first claim of Theorem
\ref{theo:betterthanGAFA} was stated as Theorem 
\ref{theo:introbetterGAFA}.

\btheo \label{theo:betterthanGAFA} %
Let $X$ be a proper geodesic CAT$(-1)$ metric space (having at least
two points), with arcwise connected boundary $\partial_\infty X$ and
extendible geodesics. Let $(H_\alpha)_{\alpha\in\A}$ be any family of
balls or horoballs with pairwise disjoint interiors. Let $\mu_0=1.534$. %
\begin{enumerate}
\item %
  For every $x$ in $X-\bigcup_{\alpha\in\A} H_\alpha$, there exists a
  geodesic ray starting from $x$ avoiding $H_\alpha[\mu_0]$ for every
  $\alpha$.
\item %
  For every $\alpha_0$ in $\A$, there exists a geodesic line starting
  from the point at infinity of $H_{\alpha_0}$ and avoiding
  $H_\alpha[\mu_0]$ for every $\alpha\neq \alpha_0$.
\end{enumerate}
\etheo

\noindent{\bf Remarks. } %
(1) Note that by its generality, Theorem \ref{theo:betterthanGAFA}
greatly improves the main results, Theorem 1.1 and Theorem 4.5, of
\cite{PP}, where (except for trees) $X$ was always assumed to be a
manifold, strict assumptions were made on the boundary of $X$, and no
definite value of $\mu_0$ was given except in special cases. But
besides this, an important point is that its proof is a much
simplified version of the upcoming main construction of Section
\ref{sec:main}, and hence could be welcome as a guide for reading
Section \ref{sec:induction}.

(2) Note that the constant $\mu_0$ is not optimal, but not by much.
For simplicial trees all of whose vertices have degree at least $3$,
the result is true, with any $\mu_0>1$ and this is optimal (though
they do not satisfy the hypotheses of the above result, the proof is
easy for them, see for instance \cite[Theo.~7.2 (3)]{PP}). We proved
in \cite{PP} that the optimal value for the second assertion of the
theorem, when $X=\HH^n_\RR$, is $\mu_0=-\log(4\sqrt{2}-5)\approx
0.42$. Hence Theorem \ref{theo:betterthanGAFA} (2) is not far from
optimal, despite its generality. Furthermore, when $X=\HH^n_\RR$, a
possible value of $\mu_0$ for the first assertion of the theorem that
was given in \cite[Theo.~7.1]{PP} was $\log(2+\sqrt
5)-\log(4\sqrt{2}-5)\approx 1.864$. Hence Theorem
\ref{theo:betterthanGAFA} (1) is even better than the corresponding
result in \cite{PP} when $X=\HH^n_\RR$, despite its generality.

\medskip \dem %
We start with the following geometric lemma. For every $\mu\geq 0$,
define
\begin{equation}
\nu(\mu)=\frac{2\;e^{-\mu}}{1+\sqrt{1-e^{-2\mu}}}\;, \label{eq:nu}
\end{equation}
which is positive and decreasing from $2$ to $0$ as $\mu$ goes from
$0$ to $+\infty$.

\blemm \label{lem:petitcalcul} %
Let $X$ be a proper geodesic CAT$(-1)$ space. Let $H$ be a ball or
horoball in $X$ and $\xi_0\in (X\cup\partial_\infty
X)-(H\cup H[\infty])$.  Let $\mu\geq\log 2$ be at most the radius of
$H$, and let $\ga$ and $\ga'$ be geodesic rays or lines starting at
$\xi_0$, meeting $H[\mu]$, parametrized such that $\ga'(s),\ga(s)$ are
equidistant to $\xi_0$ for some (hence every) $s$, and that $\ga$
enters $H$ at time $0$.
\begin{enumerate}
\item If $x=\gamma(0)$ and $x'$ are the points of entry in $H$ of $\ga$
and $\ga'$ respectively, then $d(x,x')\leq \nu(\mu)$.
\item 
For every $s\geq 0$, we have
$$
d(\ga(-s),\ga'(-s))\leq \nu(\mu)\;e^{-s}\;.
$$
\end{enumerate}

\elemm

\dem %
Let $\xi$ be the center or point at infinity of $H$, and let $t,t'$ be
the entrance times of $\ga,\ga'$ respectively in $H[\mu]$. Note that
$t\geq 0$ as $\mu\geq 0$. Let us prove first that $t'\geq 0$ too. We
refer to Section for the definition and properties of the map
$\beta_{\xi_0}$, especially when $\xi_0\in X$. Let $u$ be the point on
the geodesic $]\xi_0,\xi[$ such that $\beta_{\xi_0}(x,u)= 0$.  By the
convexity of the balls and horoballs, we have
$\beta_{\xi}(x,\ga'(0))\leq \beta_{\xi}(x,u)$. Let us prove that
$\beta_{\xi}(x,u)\leq \mu$, which will hence imply that $\ga'$ enters
$H[\mu]$ at a non-negative time (which is $t'$).

\medskip
\noindent
\begin{minipage}{10.7cm}
  ~~~Glue the two comparison triangles $(\,\overline{\xi_0},
  \overline{x}, \overline{\xi} \,)$ and $(\,\overline{\xi},
  \overline{x}, \overline{\ga(t)}\,)$ in $\HH^2_\RR$ for the geodesic
  triangles $({\xi_0},{x}, {\xi})$ and $({\xi},{x},{\ga(t)})$ along
  their sides $[\,\overline{x}, \overline{\xi}\,]$. Let $\overline{H}$
  be the ball or horoball centered at $\overline{\xi}$ such that
  $\overline{x}\in\partial \overline{H}$. By comparison, we have
  $\angle_{\overline{x}}\,(\,\overline{\xi_0}, \overline{\xi}\,) \leq
  \pi\leq \angle_{\overline{x}}\,(\,\overline{\xi_0},
  \overline{\xi}\,) + \angle_{\overline{x}}\,(\,\overline{\xi},
  \overline{\ga(t)}\,)$.  Hence the geodesic ray or line $\overline{\ga}$
  starting from $\overline{\xi_0}$ and passing through $\overline{x}$
  meets $[\,\overline{\ga(t)}, \overline{\xi}\,]$, therefore it enters
  $\overline{H}[\mu]$. Let $\overline u$ be the point on the geodesic
   $[\, \overline{\xi_0}, \overline{\xi}\,]$ such that
  $\beta_{\overline{\xi_0}} (\overline x,\overline u)= 0$. As
  $\beta_{\xi}(x,u)= \beta_{\overline \xi}(\overline x,\overline
  u)$, we only have to prove the result if $X$ is the upper halfspace
  model of the hyperbolic plane $\HH^2_\RR$. We may then assume that
  $\xi_0$ is the point at infinity $\infty$, and that $H$ is the
  horoball with point at infinity $0$ and Euclidean diameter $1$ (see
  the figure below). But then, the vertical coordinate of $\ga(0)$ is
  at least $\frac{1}{2}$, and as $e^{-\mu}\leq\frac{1}{2}$, the
  result follows: any geodesic line starting from $\xi_0$ meets the
  horizontal horosphere containing $\ga(0)$ before possibly meeting
  $H[\mu]$.
\end{minipage}
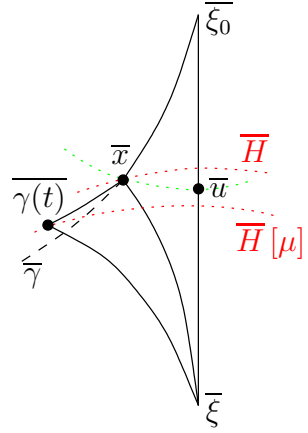
\begin{minipage}{4.1cm}
\begin{center}
\input{fig_positivitete.pstex_t}
\end{center}
\end{minipage}
\medskip

Now, in order to prove both assertions of Lemma \ref{lem:petitcalcul},
let us show that we may assume that $X=\HH^2_\RR$.

For the first one, glue the two comparison triangles
$(\,\overline{\xi_0}, \overline{x}, \overline{\xi} \,)$ and
$(\,\overline{\xi_0}, \overline{x'}, \overline{\xi} \,)$ for the
geodesic triangles $({\xi_0},{x}, {\xi})$ and $({\xi_0},{x'}, {\xi})$
along their sides $[\,\overline{\xi_0}, \overline{\xi}\,]$. As seen
above, the geodesic lines $\overline{\ga}$ (resp.~$\overline{\ga'}$)
starting from $\overline{\xi_0}$ and passing through $\overline{x}$
(resp.~$\overline{x'}$) enter $\overline{H}[\mu]$. And by comparison,
we have $d(x,x')\leq d(\overline{x},\overline{x'})$.

For the second assertion, we glue the two comparison triangles
$(\,\overline{\xi_0}, \overline{\ga(t)}, \overline{\ga'(t')}\,)$ and
$(\,\overline{\xi}, \overline{\ga(t)},\overline{\ga'(t')}\,)$ for the
geodesic triangles $({\xi_0},{\ga(t)},{\ga'(t')})$ and
$({\xi},{\ga(t)}, {\ga'(t')})$ along their isometric segments
$[\,\overline{\ga(t)}, \overline{\ga'(t')}\,]$.  As in the proof of
Lemma \ref{lemma:c1}, the geodesic segment or ray
$]\,\overline{\xi_0}, \overline{\ga(t)}\,[$ does not meet the ball or
horoball $\overline {H[\mu]}$ centered at $\overline{\xi}$ whose
boundary goes through $\overline{\ga(t)}$ and $\overline{\ga'(t')}$.
By comparison, if $H'$ is the ball or horoball centered at
$\overline{\xi}$ whose boundary passes through the point
$\overline{\ga(0)}$ on $]\,\overline{\xi_0}, \overline{\ga(t)}\,[$ at
distance $t$ from $\overline{\ga(t)}$, then $H'[\mu]$ contains
$\overline {H[\mu]}$, so that $]\,\overline{\xi_0},
\overline{\ga(t)}\,]$ and $]\,\overline{\xi_0},
\overline{\ga'(t')}\,]$ meet $H'[\mu]$. For every $s\geq 0$, as $t,t'\geq 0$, if
$\overline{\ga(-s)},\overline{\ga'(-s)}$ are the corresponding points
to ${\ga(-s)},{\ga'(-s)}$ on $]\,\overline{\xi_0},
\overline{\ga(t)}\,[$, $]\,\overline{\xi_0}, \overline{\ga'(t')}\,[$
respectively, then by comparison $d({\ga(-s)},{\ga'(-s)})\leq
d(\,\overline{\ga(-s)},\overline{\ga'(-s)}\,)$.

\medskip
\noindent
\begin{minipage}{9.5cm}
  ~~~Hence we may assume that $X$ is the upper half\-space model of
  the real hyperbolic plane $\HH^2_\RR$.  By homogeneity and
  monotonicity, it is sufficient to prove the result for $\xi_0$ the
  point at infinity $\infty$, for $H$ the horoball with point at
  infinity $0$ and Euclidean diameter $1$, and with $\ga$ and $\ga'$
  different and both tangent to $H[\mu]$. Then, by an easy
  computation, the Euclidean height of the point $\ga(0)$ is
  $\nu'(\mu)=\frac{1}{2}(1+\sqrt{1-e^{-2\mu}\,})$, so that the
  Euclidean height of the point $\ga(-s)$ is $\nu'(\mu)\;e^{s}$. The
  hyperbolic distance between $\ga(-s)$ and $\ga'(-s)$ is hence at
  most $\frac{e^{-\mu}}{\nu'(\mu)e^{s}} = \nu(\mu)\; e^{-s}$. With the
  case $s=0$, this proves both assertions. \cqfd
\end{minipage}
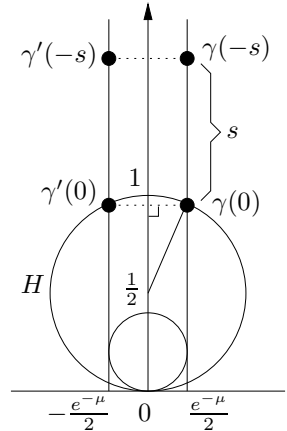
\begin{minipage}{5.3cm}
\begin{center}
\input{fig_horoexpopinch.pstex_t}
\end{center}
\end{minipage}
\bigskip

\noindent{\bf Proof of Theorem \ref{theo:betterthanGAFA}.} %
Let $X$ and $(H_\alpha)_{\alpha\in \A}$ be as in the statement. Let
$\xi_0$ be either a point in $X-\bigcup_{\alpha\in\A} H_\alpha$ or the
point at infinity of $H_{\alpha_0}$ for some $\alpha_0$ in $\A$.  For
every $\mu_1\geq \log 2$, define the following constants, with $\nu$ the
map introduced before Lemma \ref{lem:petitcalcul},
$$
\mu_2= \nu(\mu_1)>0\;,\;\;\;\;\mu_3=\mu_1+\mu_2>0\;,\;\;\;\;
\mu_4=2\mu_1- 2\mu_2\;.
$$
As $\mu_1\geq \log 2$, $\nu$ is decreasing and $\nu(\log 2)<\log 2$,
we have $\mu_4>0$. We define by induction an initial segment $\N$ in
$\NN$ and the following finite or infinite sequences
\begin{enumerate}
\item[$\bullet$] 
  $(\ga_k)_{k\in\N}$ of geodesic rays or lines starting from $\xi_0$,
\item[$\bullet$]
  $(\alpha_k)_{k\in \N-\{0\}}$ of elements in $\A$,
\item[$\bullet$]  
  $(t_k)_{k\in\N}$ of non-negative real numbers,
\item[$\bullet$] 
 $(u_k)_{k\in\N}$ of maps $u_k:[0,+\infty[\;\ra \;]0,+\infty[\,$,
\end{enumerate}
such that for every $k$ in $\N$, the following assertions hold:
\begin{enumerate}
\item [(1)] 
  If $\xi_0\in X$, then $\ga_k(0)=\xi_0$. Otherwise, $\ga_k$
  meets $\partial H_{\alpha_0}$ at time $0$.
\item [(2)]
  If $k\geq 1$, then $\ga_k$ enters $H_{\alpha_k}$ at the point
  $\ga_k(t_k)$ and meets $H_{\alpha_k}[\mu_1]$ in one and only one
  point.
\item [(3)]
  If $k\geq 1$, then $u_k(t)=u_{k-1}(t) + \mu_2\,e^{t-t_k}$ if
  $t\leq t_{k-1}$, and $u_k(t)=\mu_3$ if $t> t_{k-1}$.
\item [(4)]
  If $k\geq 1$, then $t_k\geq \mu_4+t_{k-1}$. 
\item [(5)]
  If  $t\in[0,t_k[\,$, then the point $\ga_k(t)$ does not
  belong to $\bigcup_{\alpha\in \A} H_\alpha[u_k(t)]$.
\end{enumerate}

If $\xi_0\in X$, let $\ga_0$ be a geodesic ray starting from $\xi_0$
at time $0$. Otherwise, let $\ga_0$ be a geodesic line starting from
$\xi_0$ and exiting $H_{\alpha_0}$ at time $0$. Such a $\ga_0$ exists
by the assumptions on $X$. Define $u_0$ as the constant map $t\mapsto
\mu_3$. Let $t_0=0$. The assertions (1)--(5) are satisfied for
$k=0$. Assume that $\ga_k, t_k,\alpha_k,u_k$ are constructed for
$0\leq k\leq n$ verifying the assertions (1)--(5).

If $\ga_n(]t_n,+\infty[)$ does not enter in the interior of any
element of the family $(H_\alpha[\mu_1])_{\alpha\in \A}$, then define
$\N=[0,n]\cap \NN$, and the construction terminates. Otherwise, let
$H_{\alpha_{n+1}}[\mu_1]$ be the first element of the family
$(H_\alpha[\mu_1])_{\alpha\in \A}$ such that the geodesic ray
$\ga_n(]t_n,+\infty[)$ enters in its interior.  Such an element exists
as the $H_\alpha$'s have disjoint interiors. Note that $\alpha_{n+1}
\neq \alpha_{n}$, as $\ga_n$ does not meet the interior of
$H_{\alpha_{n}} [\mu_1]$ by (2).

If $\xi_0\in X$, let $\ga_{n+1}$ be a geodesic ray starting from
$\xi_0$ at time $0$ and meeting $H_{\alpha_{n+1}}[\mu_1]$ in one and
only one point. This is possible as there exists a geodesic ray
starting from $\xi_0$ and avoiding $H_{\alpha_{n+1}}$ by the
properties of $X$ (consider for instance the extension to
$]-\infty,0]$ of $\ga_n$) and since $\partial_\infty X$ is arcwise
connected.  If $\xi_0\notin X$, let $\ga_{n+1}$ be a geodesic line
starting from $\xi_0$, and meeting $H_{\alpha_{n+1}}[\mu_1]$ in one
and only one point. Again, this is possible as $\partial_\infty X$ is
arcwise connected.  Parametrize $\ga_{n+1}$ such that $\ga_{n+1}$
exits $H_{\alpha_0}$ at time $0$. In particular, in both cases, the
assertion (1) for $k=n+1$ is satisfied.

\medskip %
Define $t_{n+1}\geq 0$ such that $\ga_{n+1}$ enters $H_{\alpha_{n+1}}$
at the point $\ga_{n+1}(t_{n+1})$, so that the assertion (2) for
$k=n+1$ is satisfied. As $\ga_{n}$ and $\ga_{n+1}$ both meet
$H_{\alpha_{n+1}}[\mu_1]$ and as $\mu_1\geq\log 2$, it follows from
Lemma \ref{lem:petitcalcul} (2) that, for every $t\leq t_{n+1}$,
\begin{equation}\label{decay}
d(\ga_{n+1}(t),\ga_{n}(t))\leq  \mu_2\;e^{t-t_{n+1}}\;.
\end{equation}
Define $\tau_n\geq t_n$ as the entrance time of $\ga_n$ in
$H_{\alpha_{n+1}}$. By Lemma \ref{lem:petitcalcul} (1), as both
$\ga_n$ and $\ga_{n+1}$ meet $H_{\alpha_{n+1}}[\mu_1]$ and
$\mu_1\geq\log 2$, we have
$$
d(\ga_{n+1}(t_{n+1}),\ga_{n}(\tau_n))\leq \mu_2\;.
$$ 
As $H_{\alpha_{n+1}}$ and $H_{\alpha_{n}}$ have disjoint
interiors, and since $H_{\alpha_{n}}$ and $H_{\alpha_{n}}[\mu_1]$ are at
distance $\mu_1$, we have $d(\ga_{n}(t_{n}),\ga_{n}(\tau_{n}))\geq
2\mu_1$.  Hence
$$
\begin{array}[b]{rl}\!\!\!d(\ga_{n}(t_{n}),\ga_{n}(t_{n+1}))\!\!\! & \geq
d(\ga_{n}(t_{n}),\ga_{n}(\tau_{n}))- 
d(\ga_{n}(\tau_{n}),\ga_{n+1}(t_{n+1}))-
 d(\ga_{n+1}(t_{n+1}),\ga_{n}(t_{n+1}))\smallskip\\ 
&\geq 2\mu_1- 2\mu_2=\mu_4>0\;.
\end{array}
$$
Hence $t_{n+1}- t_{n}$ is positive and at least $\mu_4$, which proves
the assertion (4) for $k=n+1$.

Define $t\mapsto u_{n+1}(t)$ by the induction formula in assertion
(3).  The only remaining assertion to verify is (5).  By absurd,
assume that there exist some $t$ in $[0,t_{n+1}[$ and some $\alpha\in
\A$ such that
$\ga_{n+1}(t)$ belongs to $H_\alpha[u_{n+1}(t)]$. As $u_{n+1}(t)>0$,
the element $\alpha$ is different from $\alpha_0$ if
$\xi_0\in\partial_\infty X$, and it is also different from
$\alpha_{n+1}$ by construction. By Equation \eqref{decay}, the point
$\ga_{n}(t)$ belongs to $H_\alpha[u_{n+1}(t)-\mu_2\,e^{t-t_{n+1}}]$.

Assume first that $t> t_n$, so that $u_{n+1}(t)=\mu_3$.  As $\mu_3-
\mu_2\,e^{t-t_{n+1}}>\mu_1$ (we cannot have $t=t_{n+1}$ as $H_\alpha$
and $H_{\alpha_{n+1}}$ have disjoint interiors), this implies that
$\ga_{n}(t)$ belongs to the interior of $H_\alpha[\mu_1]$. This
contradicts the fact that $H_{\alpha_{n+1}}[\mu_1]$ is the first
element of the family $(H_\alpha[\mu_1])_{\alpha\in \A}$ encountered
by $\ga_n(]t_n,+\infty[)$ in its interior.

Assume that $t\leq t_{n}$. Then $\ga_{n}(t)$ belongs to $H_\alpha
[u_{n}(t)]$. This contradicts the assertion (5) at step $n$. Thus, the
assertions (1)--(5) hold for all $k\in\N$.

\medskip %
Let us prove that the maps $u_n$ are uniformly bounded from above by
$$\mu_5=\mu_3+\frac{\mu_2}{e^{\mu_4}-1}\;.$$ As $\mu_4>0$, the
sequence $(t_k)_{k\in \NN}$ increases to $+\infty$. Fix $t\geq 0$. Let
$k=k(t)$ be the unique non-negative integer such that $t$ belongs to
$]t_{k-1},t_k]$ (by convention, $t_{-1}= -\infty$). Let us prove, by
induction on $n$, that
$$u_n(t)\leq \mu_3+\mu_2\;\sum_{j=1}^{n-k} \,e^{-\mu_4j}\;.$$ (Recall
that an empty sum is $0$). This implies that $u_n(t)\leq \mu_5$.

This is true if $n=0$, as $u_0(t)=\mu_3$. Assume that the result is
true for $n$. If $t>t_{n}$, then $u_{n+1}(t)=\mu_3$, and the result is
true.  Otherwise, by the property (3), we have $u_{n+1}(t)= u_n(t)
+\mu_2e^{t-t_{n+1}}$. Note that $t_k-t_{n+1}\leq -\mu_4(n+1-k)$ by the
property (4), and that $t\leq t_k$. Hence, by induction,
$$
u_{n+1}(t)\leq \mu_3+\mu_2\;\sum_{j=1}^{n-k} \,e^{-\mu_4j} +
\mu_2e^{-\mu_4(n+1-k)} =\mu_3+\mu_2\;\sum_{j=1}^{n+1-k} \,e^{-\mu_4j}\;.
$$
This proves the induction.

\medskip
Summarizing the above construction, there exist a sequence of geodesic
rays or lines $(\ga_n)_{n\in \NN}$ starting from $\xi_0$, and a
sequence of times $(t_n)_{n\in\NN}$ converging to $+\infty$, such that
for every $t$ in $[0,t_n]$, the point $\ga_n(t)$ does not belong to
$\bigcup_{\alpha\in\A} H_\alpha[\mu_5]$.  (Take an eventually constant
sequence $(\ga_n)_{n\in \NN}$ if the construction stops at a finite
stage, which is possible as $\mu_5>\mu_1$.)  As $(t_{n})_{n\in\NN}$
grows at least linearly, the formula \eqref{decay} implies that
$(\ga_{n}(t))_{n\in\NN}$ is a Cauchy sequence, uniformly on every
compact subset of non-negative $t$'s.  Hence, the geodesic rays or
lines $\ga_n$ converge to a geodesic ray or line avoiding
$\bigcup_{\alpha\in\A} H_\alpha[\mu_5-\epsilon]$, for every
$\epsilon>0$.  Taking $\mu_1=1.042\geq \log 2$, we can check that
$\mu_5< 1.5332$, hence the result follows.
\cqfd

\bcoro \label{coro:betterthanGAFA} %
Let $X$ and $(H_\alpha)_{\alpha\in\A}$ be as in Theorem
\ref{theo:betterthanGAFA}. 
For every $x\in X$, there exist $t>0$ and a geodesic ray $\gamma$
starting at $x$ such that $\gamma([t,\infty[)$ is contained in the
complement of $\bigcup_{\alpha\in\A}H_\alpha[\mu_0]$.
\ecoro

\dem %
We may assume that $x\in H_{\alpha_0}$ for some $\alpha_0\in\A$,
otherwise, Theorem \ref{theo:betterthanGAFA} (1) applies (with $t=0$).
Let $H'_\alpha=H_\alpha$ if $\alpha\neq\alpha_0$, and $H'_{\alpha_0}
=H_{\alpha_0}[d(x,\partial H_{\alpha_0})+1]$. Then $x\notin X -
\bigcup_{\alpha\in\A} H'_\alpha$. By Theorem \ref{theo:betterthanGAFA}
(1), let $\ga$ be a geodesic ray starting from $x$ and avoiding the
$H'_\alpha[\mu_0]$'s. Let $t=d(x,\partial H_{\alpha_0})+2+ 2\mu_0+
c'_1(\infty)$. As $\norm{\ell_{H_{\alpha_0}}-\ph_{H_{\alpha_0}}}\leq
c'_1(\infty)$ by Subsection \ref{subsec:penetrationmaps}, the geodesic
ray $\gamma([t,\infty[)$ does not meet $H_{\alpha_0}$. The result follows.
\cqfd

\medskip %
Let $e$ be an end of a finite volume complete negatively curved
Riemannian manifold $V$. Let ${\height}_e$ be the Busemann function of
$e$ normalized to be zero on the boundary of the maximal Margulis
neighbourhood of $e$ (see for instance \cite{BK,HP3,PP}, as well as
the paragraph above Corollary \ref{coro:maxsp}).  Our next result
improves Theorem 7.4 (hence Corollary 1.2) in \cite{PP}, with the same
proof as in [loc.~cit.], by removing the technical assumptions on the
manifold, and giving a universal upper bound on $h_e(V)$.

\bcoro \label{coro:lowclosedgeod}
Let $V$ be a finite volume complete Riemannian manifold with dimension
at least $2$ and sectional curvature $K\le -1$. Then there exists a
closed geodesic in $V$ whose maximum height (with
respect to $\height_e$) is at most $1.534$. \cqfd
\ecoro

%
%
%
\subsection{The inductive construction} 
\label{sec:induction}

Fix arbitrary constants $\epsilon_0\in\RR^*_+\cup\{\infty\}$ and
$\delta_0,\kappa_0\geq 0$, and fix an arbitrary point $\xi_0$ in
$X\cup \partial_\infty X$. Let $(C_n)_{n\in\NN}$ be a family of
$\epsilon_0$-convex subsets of $X$ such that $\xi_0\notin
C_0\cup\partial_\infty C_0$, and let $f_0$ be a $\kappa_0$-penetration
map for $C_0$.

The aim of this section is to construct by induction a sequence
of geodesic rays or lines in $X$, starting from $\xi_0$ and having a
suitable penetration behaviour in the $C_n$'s.

\medskip\noindent{\bf Prescription of constants. }  %
The following constants will appear in the statement, or in the proof,
of the inductive construction:
\begin{enumerate} 
\item[$\bullet$] $c_1 = c'_1(\epsilon_0)>0$ given by Lemma
  \ref{lemma:c1} if $\epsilon_0\neq \infty$ and by Lemma
  \ref{lem:c1pinfty} if $\epsilon_0=\infty$ and
  $(f_0,\delta_0)\neq(\ph_{C_0},0)$; otherwise $c_1=\frac{1}{19}$ ;
\item [$\bullet$] %
  $c_2= c'_2(\epsilon_0)>0$ given by Equation \eqref{eq:c2} if
  $\epsilon_0\neq \infty$ and by \eqref{eq:c2i} otherwise;
\item [$\bullet$] %
  $c_3 =2\sinh c_1+c_2\,e^{2c_1}\sinh c_1$, which is
  positive, and depends  on $\epsilon_0$;
\item [$\bullet$] %
  $c_4=c'_3(\epsilon_0)\sinh(c_1+\delta_0)+
  c_2\;e^{-3c'_3(\epsilon_0)\sinh(c_1+\delta_0)-\log 2}\sinh c_1$,
  where $c'_3(\cdot)$ is given by Equation \eqref{eq:c3} if
  $\epsilon_0\neq \infty$ and by \eqref{eq:c3i} otherwise. Note that
  $c_4$ is positive, and depends on $\epsilon_0,\delta_0$;
\item [$\bullet$] %
  $c_5=c_5(\epsilon_0,\delta_0)=
  2\max\{c_2,c'_3(\epsilon_0)\}\sinh(c_1+\delta_0)$, which is
  positive, and depends on $\epsilon_0,\delta_0$;
\item [$\bullet$] %
  $c_6= 3c_4+\log 2$, which is positive, and depends on
  $\epsilon_0,\delta_0$;
\item [$\bullet$] %
  $h_0=h_0(\epsilon_0,\delta_0,\kappa_0)= \max\{\;\delta_0 +\kappa_0
  \;,\;c_0(\epsilon_0)+\kappa_0\;,\;
  h'(\epsilon_0,\sinh(\delta_0+c_1)) \;, \;\delta_0+2c_1+c_6\;\}$,
  where $c_0(\cdot)$ is given by Equation (\ref{eq:c0}) if
  $\epsilon\neq \infty$ and by (\ref{eq:c0i}) otherwise, and
  $h'(\cdot,\cdot)$ is given by Equation \eqref{eq:h} if
  $\epsilon_0\neq \infty$ and by \eqref{eq:c3i} otherwise;
\item [$\bullet$] %
  For every $h'_0\geq 0$, let 
  $h'_1=h'_1(\epsilon_0,\delta_0,h'_0)=h'_0+2c_5$.
\end{enumerate}
Fix $h'_0\geq h_0$ and $h\geq h'_1$.

\medskip
\noindent 
{\bf Assumptions on the family $(C_n)_{n\in\NN}$.} %
Assume that there exists at least one geodesic ray or line $\ga_0$
starting from $\xi_0$ and meeting $C_0$ with $f_0(\ga_0)=h$, and that
the following conditions are satisfied.
\begin{enumerate}
\item[$(iii)$] {\bf (Almost disjointness property)}
  For every $m,n$ in $\NN$ with $m\neq n$, the diameter of $C_n\cap
 C_m$ is at most $\delta_0$.
\item[$(iv)$] {\bf (Local prescription property)} For every $n$ in
  $\NN-\{0\}$ such that $\xi_0\notin C_n\cup\partial_\infty C_n$, if
  there exists a geodesic ray or line $\alpha$ starting from $\xi_0$
  which meets first $C_0$ and then $C_n$ with $f_0(\alpha)=h$ and
  $\ell_{C_n}(\alpha)\geq h'_0$, then there exists a geodesic ray or
  line $\alpha'$, starting from $\xi_0$ which meets first $C_0$ and
  then $C_n$ with $f_0(\alpha')=h$ and $\ell_{C_n}(\alpha')= h'_0$.
\end{enumerate}
Note that $(iii)$ is satisfied with $\delta_0=0$ if the $C_n$'s have
disjoint interior. In Section \ref{sec:results}, we will use
Proposition \ref{prop:prescribe} to 
check $(iv)$ for various applications,
with $h'_0=\max\{h_0,h_0^{\rm min}\}$ and $h\geq \max\{h'_1,h^{\rm
min}\}$, for the various values of $h_0^{\rm min},h^{\rm min}$ defined
in the various cases of Proposition \ref{prop:prescribe}.

\medskip For every $n$ in $\NN$ such that $\xi_0\notin
C_n\cup\partial_\infty C_n$, define $f_n=\ell_{C_n}:T^1_{\xi_0}X\ra
[0,+\infty]$, and for every geodesic ray or line $\ga$
starting from $\xi_0$ and meeting $C_n$, let $t^-_n(\ga),t^+_n(\ga)\in
\;]-\infty,+\infty]$ be the entrance time and exit time of $\ga$ in and
out of the convex subset $C_n$ respectively. The following remark will
be used later on.

\blemm \label{lem:enterpositiv} %
For every $n>0$, for every geodesic ray or line $\ga$ starting from
$\xi_0$ and entering $C_0$ at time $t=0$, such that $f_0(\ga)=h$ and
$\ga(]\delta_0,+\infty[)$ meets $C_n$, we have $\xi_0\notin
C_n\cup\partial_\infty C_n$ and $t^-_n(\ga)>0$.
\elemm

\dem %
Otherwise, as $\ga(]\delta_0,+\infty[)$ meets $C_{n}$ and by
convexity, there exists $\epsilon>0$ such that the geodesic segment
$\ga([0,\delta_0+\epsilon])$ is contained in $C_{n}$.  By the
Penetration property $(i)$ of $f_0$, the length of $\ga\cap C_0$ is at least
$h-\kappa_0$, which is bigger than $\delta_0$ as $h\geq h'_1> h'_0\geq
h_0\geq\delta_0+\kappa_0$ by the definitions of $h'_1$ and $h_0$. As
$\ga$ enters $C_0$ at time $t=0$, up to taking $\epsilon>0$ smaller,
this implies that the geodesic segment $\ga([0,\delta_0+\epsilon])$ is
also contained in $C_0$. This contradicts the Almost disjointness
property $(iii)$ as $n\neq 0$.
\cqfd

\medskip\noindent{\bf Statement of the inductive construction. } %
We will define by induction an initial segment
$\N$ in $\NN$, and finite or infinite sequences
\begin{enumerate}
\item[$\bullet$] 
  $(\ga_k)_{k\in\N}$ of geodesic rays or lines starting
  from $\xi_0$,
\item[$\bullet$] 
  $(n_k)_{k\in\N}$ of integers such that $\xi_0\notin
  C_{n_k}\cup\partial_\infty C_{n_k}$,
\item[$\bullet$] 
  $(u_k)_{k\in\N}$ of maps $u_k:[0,+\infty[\ra [h'_0,h'_1]$,
\end{enumerate}
such that the following assertions hold, for every $k$ in $\N$, 
where we use $t^\pm_k=t^\pm_{n_k}(\ga_k)$ to simplify notations.
\begin{enumerate}
\item[(1)]
  The geodesic ray or line $\ga_k$ enters $C_0$ at time $t=0$ and
  $f_0(\ga_k)=h$. 
\item [(2)]
  If $k\geq 1$, then $\ga_k$ meets $C_{n_k}$ with $t^-_k\geq 0$
  and $f_{n_k}(\ga_k)=h'_0$. 
\item [(3)]
  If $k\geq 1$, then $d(\ga_k(t),\ga_{k-1}(t))\leq c_3
  \;e^{t-t_k^-}$ for every $t$ in $[0,t_k^-]$.
\item [(4)]
  If $k\geq 1$, then
$$
u_{k}(t) =\sup_{s\,\in\,[0,+\infty[\,\;:\;|s-t|\,\leq 
\,c_4\,e^{t-t^-_k}} u_{k-1}(s) \;+ \;
c_5\;e^{t-t^-_k}
$$ 
for $t\in [0,t^-_k]$ and $u_{k}(t) = h'_0$ if $t>t^-_k$.
\item [(5)]
  If $k\geq 1$, then $t^-_k\geq t^-_{k-1} +c_6$.
\item[(6)] 
  If $k\geq 1$, for every $n$ in $\NN-\{0\}$ such that
  $\ga_k(]\delta_0,+\infty[)$ meets $C_n$ with $t^-_n(\ga_k)\leq
  t^-_k$, we have $f_n(\ga_k)\leq u_k(t^+_n(\ga_k)-\delta_0)$.
\end{enumerate}

Note that by Lemma \ref{lem:enterpositiv} and by (1), if
$\ga_k(]\delta_0, +\infty[)$ meets $C_n$, then $\xi_0\notin
C_{n}\cup\partial_\infty C_{n}$, so that, in particular,
$t^\pm_n(\ga_k)$ are well defined, and (6) does make sense.

\medskip\noindent{\bf Proof of the inductive construction. } %
By the assumptions, let $\ga_0$ be a geodesic ray or line starting from
$\xi_0$ and entering $C_0$ at time $t_0^-=0$, such that
$f_0(\ga_0)=h$. Let $n_0=0$. Let $u_0:[0,+\infty[\ra [h'_0,h'_1]$ be the
constant map with value $h'_0$.  As the conditions (2)--(6) are empty if
$k=0$, the construction is done at step $0$.

Let $k\geq 1$, and assume that $\ga_0,n_0,u_0,\dots, \ga_{k-1},
n_{k-1},u_{k-1}$ are constructed. Note that $u_{k-1}\geq h'_0$ by
induction. If for every $n$ in $\NN-\{0\}$ such that
$\ga_{k-1}(]\delta_0,+\infty[)$ meets $C_n$, we have
$f_n(\ga_{k-1})\leq u_{k-1}(t^+_n(\ga_{k-1})-\delta_0)$, then we stop
and we define $\N=\{0,1\dots,k-1\}$.

Otherwise, let $\tau$ be the greatest lower bound of the
$t^-_{n}(\ga_{k-1})$'s taken over all $n$ in $\NN-\{0\}$ such that
$\ga_{k-1}(]\delta_0,+\infty[)$ meets $C_n$ with $f_n(\ga_{k-1})>
u_{k-1}(t^+_n(\ga_{k-1})-\delta_0)$.

Let us prove that this lower bound is in fact a minimum, attained for
only one such $n$. Let $\epsilon>0$ such that $h'_0> \delta_0 +
\epsilon$, which is possible by the definition of $h_0$, as $h'_0\geq
h_0$. If $t^-_n(\ga_{k-1})$ and $t^-_m(\ga_{k-1})$ belong to
$[\tau,\tau+\epsilon]$ with $f_n(\ga_{k-1})> u_{k-1}(t^+_n(\ga_{k-1})
-\delta_0)$ and $f_m(\ga_{k-1})> u_{k-1}(t^+_m(\ga_{k-1}) -
\delta_0)$, assume for instance that $t^-_n(\ga_{k-1})\leq
t^-_m(\ga_{k-1})$. As $f_n=\ell_{C_n}$, $f_m=\ell_{C_m}$, $u_{k-1}\geq
h'_0$ and $t^-_m(\ga_{k-1}) - t^-_n(\ga_{k-1})\leq \epsilon$, the
subsets $C_n$ and $C_m$ meet along a segment of length at least
$h'_0-\epsilon> \delta_0$. By the Almost disjointness property (iii),
this implies that $n=m$. In particular, we have $\tau=
t^-_n(\ga_{k-1})$ for a unique $n\in\NN-\{0\}$, and we denote this $n$
by $n_k\in\NN-\{0\}$, so that $\ga_{k-1}(]\delta_0,+\infty[)$ meets
$C_{n_k}$ 
with
\begin{equation}\label{eq:defink}
f_{n_k}(\ga_{k-1})> u_{k-1}(t^+_{n_k}(\ga_{k-1})-\delta_0)\geq h'_0\;.
\end{equation}
In particular, $\xi_0\notin C_{n_k}\cup\partial_\infty C_{n_k}$ by
Lemma \ref{lem:enterpositiv} and by Assertion (1) at rank $k-1$.  Note
that $n_k\neq n_{k-1}$, as $f_{n_{k-1}}(\ga_{k-1})=h'_0$ by the
assertion (2) at rank $k-1$, which would contradict Equation
(\ref{eq:defink}) if $n_k=n_{k-1}$.

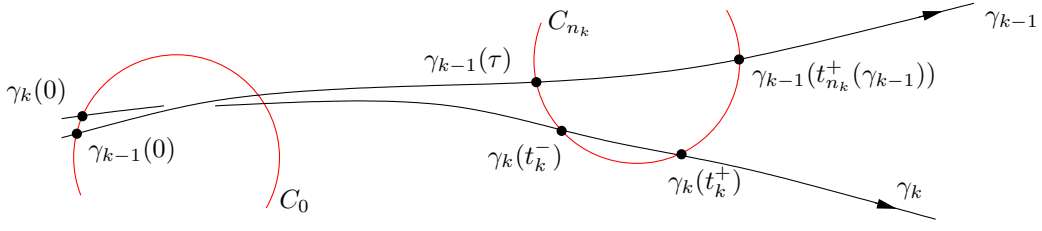
\begin{center}
\input{fig_construcgamk.pstex_t}
\end{center}

\medskip %
By Lemma \ref{lem:enterpositiv}, the geodesic ray or line
$\ga_{k-1}$ first enters $C_0$ and then $C_{n_k}$. Furthermore,
$\ga_{k-1}$ satisfies (1) and $f_{n_k}(\ga_{k-1})\geq h'_0$. Hence, by
the Local prescription property $(iv)$, there exists a geodesic ray or
line $\ga_k$ starting from $\xi_0$ that first enters $C_0$ and then
$C_{n_k}$, with $f_{0}(\ga_{k})=h$ and $f_{n_k}(\ga_{k})= h'_0$.
Choose the parametrization in such a way that $\ga_{k}$ enters $C_0$
at time $0$.  In particular, (1) and (2) hold for $\ga_k$, and
$t_k^-=t^-_{n_k}(\ga_k)>0$.  Define $u_{k}$ by using the induction
formula given in the assertion (4).  Before checking (3)--(6) for
$\ga_k,n_k,u_k$, let us make two preliminary remarks.

\blemm \label{lem:c1prime} %
We have  $d(\ga_{k-1}(\tau),\ga_k(t^-_{k}))\leq c_1$  and
$d(\ga_{k-1}(0),\ga_k(0))\leq c_1$.
\elemm

\dem %
By Lemma \ref{lemma:c1} if $\epsilon_0\neq \infty$ and Lemma
\ref{lem:c1pinfty} otherwise, we have $d(\ga_{k-1}(\tau),
\ga_k(t^-_{k}))\leq c'_1(\epsilon_0)$ and $d(\ga_{k-1}(0),\ga_k(0))
\leq c'_1(\epsilon_0)$.  By the definition of $c_1$, we hence only
have to prove Lemma \ref{lem:c1prime} when $\epsilon_0=\infty$,
$\delta_0=0$ and $f_0=\ph_{C_0}$. In this case, as $c_1=1/19$,
$c_2=5/2$, $c_0(\infty)=4.056$, $\kappa_0 = 2\log(1+\sqrt{2}) =
c'_1(\infty)$, $c'_3(\infty)=5/2$, easy computations show that
$$
h_0=h'(\infty,\sinh c_1)=3\sinh c_1 +c_0(\infty)+c'_1(\infty)\approx
5.9767$$ 
and, for future use,
\begin{equation}\label{eq:valeurhun}
  h'_1(\infty,0,h_0(\infty,0,c'_1(\infty)))\approx 6.5032.
\end{equation}
As $\ph_{C_0}(\ga_k)$ and $\ph_{C_0}(\ga_{k-1})$ are equal to $h\geq
h'_1 \geq h'_0\geq h_0$, and since $h_0/2\geq \log 2$, it follows from
the definition of the map $\ph_{C_0}$ and from Lemma
\ref{lem:petitcalcul} (1) and (2) that $d(\ga_{k-1}(0),\ga_k(0))$ and
similarly $d(\ga_{k-1}(\tau),
\ga_k(t^-_{k}))$  are at most
$\nu(h_0/2)$, where $\nu(.)$ is defined by Equation \eqref{eq:nu}.  An
easy computation shows that $\nu(h_0/2)\leq c_1=1/19$, which proves
the result.  \cqfd

\blemm \label{lem:tau} %
We have  $|\tau -t^-_k|\leq 2c_1$.
\elemm

\dem %
Lemma \ref{lem:enterpositiv}, applied to $n=n_k$ and $\ga=\ga_{k-1}$,
implies that $\tau>0$. We have seen that $t_k^->0$. By the triangular
inequality and the above lemma, we have $|\tau -t^-_k|\leq 2c_1$.
\cqfd

\bigskip
\noindent
{\bf Verification of (5). } 
Note that $\tau=t^-_{n_k}(\ga_{k-1})>t^-_{k-1}$. Otherwise, as 
$$
t^+_{n_k}(\ga_{k-1})=\tau+ f_{n_k}
(\ga_{k-1})\geq \tau+h'_0\geq \tau+h_0>\delta_0
$$ 
by Equation \eqref{eq:defink} and by the definition of $h_0$, we have,
 by the assertion (6) at step $k-1$, the inequality
 $f_{n_k}(\ga_{k-1})\leq u_{k-1}(t^+_{n_k}(\ga_{k-1}) -\delta_0)$,
 which contradicts the definition of $n_k$, see Equation
\eqref{eq:defink}.

Let us first prove that $\tau\geq t^-_{k-1}+h_0-\delta_0$. Assume
first that $\tau\geq t^+_{k-1}$. Then, separating the case $k=1$ where
$f_{n_{k-1}}(\ga_{k-1})=h\geq h'_1\geq h'_0$ from the case $k\geq 2$
where $f_{n_{k-1}}(\ga_{k-1})=h'_0$, we have
\begin{equation}\label{eq:minortau}
\tau-t^-_{k-1}\geq
t^+_{k-1}-t^-_{k-1}= f_{n_{k-1}}(\ga_{k-1})\geq h'_0\geq h_0.
\end{equation}
Hence the result holds. Otherwise, $t^-_{k-1}<\tau<t^+_{k-1}$. By
convexity, $\ga_{k-1}(\tau)$ belongs to $C_{n_{k-1}}$. Note that
$\ga_{k-1}([\tau,\tau + h_0])$ is contained in $C_{n_k}$, since $\tau$
is the entrance time of $\ga_{k-1}$ in $C_{n_k}$, and
$f_{n_k}(\ga_{k-1})\geq h'_0\geq h_0$.  If $\tau + \delta_0<
t^+_{k-1}$, as $\ell_{C_{n_k}}(\ga_{k-1})\geq h_0>\delta_0$ by the
definition of $h_0$, then $C_{n_k}\cap C_{n_{k-1}}$ contains a
geodesic segment of length bigger than $\delta_0$. This contradicts
the Almost disjointness property $(iii)$ since $n_k \neq n_{k-1}$.
Hence
$$
\tau \geq t^+_{k-1}- \delta_0\geq t^-_{k-1}+h_0-\delta_0
$$ 
by Equation (\ref{eq:minortau}), and the result holds.

Now, by Lemma \ref{lem:tau},  
$$
t^-_k-t^-_{k-1}\geq \tau-2c_1- t^-_{k-1}
\geq h_0-\delta_0 -2c_1\geq c_6
$$ 
by the definition of $h_0$. Therefore, the assertion (5) holds at rank
$k$.

\bigskip
\noindent{\bf Verification of (4). } We only have to check that $u_k$
has values in $[h'_0,h'_1]$. We start by proving the following easy but
tedious general lemma.

\blemm \label{lem:un} %
Let $c,c',c'',h_*\geq 0$, let $\M$ be an initial segment in $\NN$, let
$(t_n)_{n\in\M}$ be a sequence of non-negative real numbers, and let
$(u_n:[0,+\infty[\;\ra[0,+\infty[) _{n\in\M}$ be a sequence of maps.
Assume that $u_0$ has constant value $h_*$, and that for every $n$ in
$\M-\{0\}$, we have $t_n-t_{n-1}\geq c''$, $u_n(t)=h_*$ if $t>t_n$ and
if $t\leq t_n$, then
$$u_n(t)= c\;e^{t-t_n}\;+\sup_{s\,\in\,[0,+\infty[\,\;:\;|s-t|\,\leq 
  \,c'\,e^{t-t_n}} u_{n-1}(s)\;.$$ If $c''\geq 3c'+\log 2$, then for
every $t\in[0,+\infty[$, for every $n$ in $\M$, we have
$$h_*\leq u_n(t)\leq h_*+2c\;.$$
\elemm

To prove that $u_k$ has values in $[h'_0,h'_1]$, we apply Lemma
\ref{lem:un} with $c=c_5$, $c'=c_4$, $c''=c_6$, $h_*=h'_0$, $\M=
\{0,1,\dots,k\}$ and $(t_i)_{i\in\M}= (t^-_i)_{1\leq i\leq k}$. Its
hypotheses are satisfied by the definition of the constant $c_6$, by
the assertion (5) at rank less than or equal to $k$, that we just
proved, and by the definition of $u_k$ and the assertion (4) for $u_i$
with $0\leq 1\leq k-1$. Hence the map $u_{k}$ does have values in
$[h'_0,h'_1]$, by the definition of $h'_1$.

\medskip
\noindent{\bf Proof of Lemma \ref{lem:un}. }%
First note that by an easy induction, whatever the value of $c''$ is,
for every $t\in[0,+\infty[$ and $n\in\M$, we have $u_n(t)\geq h_*$.

Let $c''\geq 3c'+\log 2$, $t\in[0,+\infty[$ and $n\in\M$. Let us prove
that $u_n(t)\leq h_*+2c$. We may assume that $t\leq t_n$ and that
$n\geq 1$. Define $t_{-1}=-2c'-1$. Let $m$ be the unique element in
$\NN$ such that $t_{m-1} +2c'<t\leq t_m+2c'$. Let $N=n-m\geq 0$. Note
that for every integer $k$ with $0\leq k\leq N$, we have
$t_{n-k}-t_m\geq (n-m-k)c''$ hence
\begin{equation}
\label{eq:majot}
t-t_{n-k}\leq 2c'-(N-k)c''\;.
\end{equation}

Consider the finite sequence $(x_k)_{0\leq k\leq N}$ defined by
$x_0=0$ and
$$
x_{k+1}=x_k+e^{c'x_k -(N-k)c'' +2c'}
$$ 
for $0\leq k\leq N-1$.  Let us prove by induction on $k$ that $x_k
\leq e^{-(N-k)c''}$, which in particular implies that
\begin{equation}
\label{eq:majxNun}
x_N\leq 1\;.
\end{equation}
Indeed, the result is true for $k=0$. Assume it to be true for some 
$k\leq N-1$. Then
\begin{align*}
x_{k+1}\leq\; & e^{-(N-k)c''} + e^{c'e^{-(N-k)c''} -(N-k)c'' +2c'}\\
\leq \; &
e^{-(N-k-1)c''}\left(e^{-c''}+e^{-c''+3c'}\right)\leq e^{-(N-k-1)c''}
\end{align*}
as $c''\geq 3c'+\log 2\geq \log(1+e^{3c'})$.

Let us now prove by induction on $k$ that, for $0\leq k\leq N$, we
have 
\begin{equation}
\label{form:un}
u_n(t)\leq \sup_{|s-t|\,\leq \,c'\,x_k} u_{n-k}(s)\;+c\;x_k\;.
\end{equation}
This is true if $k=0$, assume it is true for some $k\leq N-1$.  In
particular, $n-k\geq 1$. For every $s\in[0,+\infty[$ such that
$|s-t|\,\leq \,c'\,x_k$, we have
$$ 
u_{n-k}(s) \leq \sup_{|s'-s|\,\leq \,c'\,e^{s-t_{n-k}}}u_{n-k-1}(s') 
\;+c\;e^{s-t_{n-k}}
$$ 
(this is true by definition if $s\leq t_{n-k}$, and also true
otherwise as then $u_{n-k}(s)=h_*$ and $u_{n-k-1}(s')\geq h_*$ for
every $s'$). Hence by the triangular inequality and the equation
\eqref{eq:majot},
\begin{align*}
  u_n(t)\leq & \sup_{|s'-t|\,\leq \,c'\,x_k+c'e^{t+c'x_k -t_{n-k}}}
  u_{n-k-1}(s')\;+
  c\;x_k +ce^{t+c'x_k -t_{n-k}}\\
  \leq & \sup_{|s'-t|\,\leq \,c'\,x_k+c'e^{c'x_k +2c'-(N-k)c''}}
  u_{n-k-1}(s')\;+
  c\;x_k +ce^{c'x_k+2c'-(N-k)c''}\\
  = & \;\sup_{|s'-t|\,\leq \,c'\,x_{k+1}} u_{n-k-1}(s')\;+
  c\;x_{k+1}\;,
\end{align*}
which proves the inductive formula \eqref{form:un}.

Finally, let us prove that $u_n(t)\leq h_*+2c$, which finishes the
proof of the lemma. Take $k=N$ in the inductive formula
\eqref{form:un}, and note that $n-N=m$. For every $\epsilon>0$, let
$s\in[0,+\infty[$ with $|s-t|\,\leq \,c'\,x_N$ such that $\sup_{|s'-t|
  \,\leq \,c'\,x_N} u_{m}(s')\leq u_{m}(s)+\epsilon$. If $s> t_m$ or
$m=0$, then $u_m(s)=h_*$, hence by the inequality \eqref{eq:majxNun},
$$
u_n(t)\leq \sup_{|s'-t|\,\leq \,c'\,x_N} u_{m}(s')\;+c\;x_N\leq 
h_*+\epsilon+c\;,
$$
and the result holds.  Otherwise, $s\leq t_m$ and $m\geq 1$. For every
$s'\in[0,+\infty[$ such that $|s'-s|\,\leq \,c'\,e^{s-t_m}$, we have
$s'\geq s-c'\geq t-2c'> t_{m-1}$. Again,  the definition of $s$ and
the inequality
\eqref{eq:majxNun} gives
\begin{align*}
u_n(t) & \leq  u_m(s)+\epsilon+c\;x_N\\ & = 
\sup_{|s'-s|\,\leq \,c'\,e^{s-t_m}} u_{m-1}(s')
\;+c\;e^{s-t_m}+\epsilon+c\;x_N
\leq \epsilon+h_*+2c\;,
\end{align*}
and the result also holds.
\cqfd

\bigskip
\noindent{\bf Verification of (3). } %
Let $t$ be in $[0,t^-_k]$. Recall that
$d(\ga_{k-1}(\tau),\ga_k(t^-_{k})) \leq c_1$, hence we have
$d(\ga_{k-1}(\tau),\ga_k)\leq c_1$ by Lemma \ref{lem:c1prime}. By
Lemma \ref{lem:expopinch}, we have $d(\ga_{k-1}(0),\ga_k)\leq
e^{-\tau} \sinh c_1 $.  By the Penetration property $(i)$ of $f_0$ and the
definition of $h_0$, we have
$$
\ell_{C_0}(\ga_k)\geq
f_0(\ga_k)-\kappa_0=h-\kappa_0\geq h_0-\kappa_0\geq c_0(\epsilon_0)\;.
$$ 
Thus, by Lemma \ref{lemma:c2} if $\epsilon_0\neq \infty$ and by Lemma
\ref{lem:c2pinfty} if $\epsilon_0=\infty$, and by the definition of
$c_2$, we have 
\begin{equation}
\label{eq:majdistoo}
d(\ga_{k-1}(0),\ga_k(0))\leq c_2\,e^{-\tau}\sinh c_1\;.
\end{equation}
We refer to Section \ref{subsec:notations} for the definition and
properties of the map $\beta_{\xi_0}$. It
follows from the inequality \eqref{eq:majdistoo} that
\begin{equation}
\label{eq:majorebeta}
|\beta_{\xi_0}(\ga_{k-1}(t),\ga_k(t))|=
|\beta_{\xi_0}(\ga_{k-1}(0),\ga_k(0))|\leq
d(\ga_{k-1}(0),\ga_k(0))\leq c_2\,e^{-\tau}\sinh c_1\;.
\end{equation}
For every $s$ in $\RR$, let $\ga_{k-1}(s')$ be the point on the
geodesic line $\ga_{k-1}$ such that the equality $\beta_{\xi_0}
(\ga_{k-1}(s'), \ga_k(s)) =0$ holds.  For every point
$p\in\gamma_{k-1}$, we have
\begin{equation}
\label{eq:astucebeta}
d(p,\gamma_{k-1}(t'))=\big|\beta_{\xi_0}(p,\gamma_{k-1}(t'))\big|=
\big|\beta_{\xi_0}(p,\gamma_{k}(t))\big|\le  
d(p,\gamma_k(t))\;,
\end{equation}
Using the triangle inequality with the point $p$ the closest to
$\ga_k(t)$ on $\ga_{k-1}$, Lemma \ref{lem:expopinch} and Lemma
\ref{lem:c1prime}, we hence have the following inequality
\begin{eqnarray}
d(\ga_k(t),\ga_{k-1}(t')) &\leq &2\, d(\gamma_k(t),\gamma_{k-1})\le 2\,
e^{t-t^-_k}\sinh d(\ga_k(t^-_k),\ga_{k-1}(\tau))\nonumber\\
 &\leq &
2 \,e^{t-t^-_k} \sinh c_1\;.\label{eq:majdittp}
\end{eqnarray}
Note that, using Equation (\ref{eq:astucebeta}) with $p=\ga_{k-1}(t)$
and the inequalities (\ref{eq:majorebeta}),
$$
d(\ga_{k-1}(t),\ga_{k-1}(t'))=|\beta_{\xi_0}(\ga_{k-1}(t),\ga_k(t))|
\leq c_2\,e^{-\tau}\sinh c_1\;..
$$ 
Hence, by the inequality \eqref{eq:majdittp}, we have
\begin{align*}
d(\ga_k(t),\ga_{k-1}(t))\leq & \;\;d(\ga_k(t),\ga_{k-1}(t'))+
d(\ga_{k-1}(t'),\ga_{k-1}(t))\\
\leq  &
\;\;2\,e^{t-t^-_k}\sinh c_1
+c_2\;e^{-\tau} \sinh c_1\;.
\end{align*}
As $\tau \geq t^-_k -2c_1$ by Lemma \ref{lem:tau}, and by the
definition of $c_3$, we get
$$
d(\ga_k(t),\ga_{k-1}(t))\leq c_3 \;e^{t-t^-_k}\;,
$$
which proves the assertion (3) at rank $k$.

\bigskip
\noindent{\bf Verification of (6). } 
By absurd, assume that there exists $n\in\NN-\{0\}$ such that
$\ga_k(]\delta_0,+\infty[)$ meets $C_n$ (so that in particular
$\xi_0\notin C_{n}\cup\partial_\infty C_{n}$ by Lemma
\ref{lem:enterpositiv}), with $t^-_n(\ga_k)\leq t^-_k$ and
\begin{equation}
\label{eq:absurdn}
f_n(\ga_k)> u_k(t^+_n(\ga_k)-\delta_0)\;.
\end{equation} 
To simplify notation, let $s^\pm_{k}=t^\pm_n(\ga_{k})$,
$x=\ga_k(s^-_k), y = \ga_k(s^+_k)$, and, as we will prove later on
that $\ga_{k-1}$ also meets $C_n$, let $s^\pm_{k-1}= t^\pm_n
(\ga_{k-1})$, $x'=\ga_{k-1}(s^-_{k-1}), y' = \ga_{k-1}(s^+_{k-1})$.

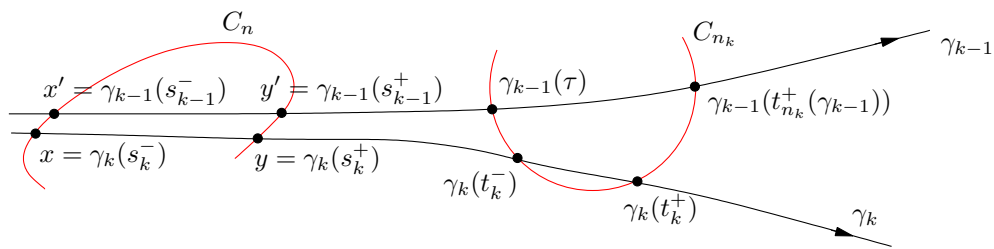
\begin{center}
\input{fig_verifasscinq.pstex_t}
\end{center}

Note that $s^+_k\leq t_k^- +\delta_0$. Otherwise, as $s^-_k\leq t^-_k$
and by convexity, there exists $\epsilon>0$ such that
$\ga_{k}([t^-_k,t^-_k+\delta_0+\epsilon])$ is contained in $C_{n}$.
As $t^+_k-t^-_k = h'_0 \geq h_0> \delta_0$, up to making $\epsilon$
smaller, the geodesic segment $\ga_{k}([t^-_k, t^-_k + \delta_0 +
\epsilon])$ is also contained in $C_{n_k}$. Hence $n$ is equal to
$n_k$ by the Almost disjointness property (iii). But
$f_{n_k}(\ga_k)=h'_0$ and, by Equation \eqref{eq:absurdn}, we have
$f_{n}(\ga_k)> u_k(s^+_k-\delta_0)\geq h'_0$, so that $n$ cannot be
equal to $n_k$.

By Lemma \ref{lem:expopinch} applied to the geodesic triangle with
vertices $\ga_k(t^-_k+\delta_0), \ga_{k-1}(\tau),\xi_0$, and as
$d(\ga_k(t^-_k),\ga_{k-1}(\tau))\leq c_1$ by Lemma \ref{lem:c1prime},
we have
\begin{align}
\label{eq:majdygakplu}
d(y,\ga_{k-1}) & \leq  e^{-d(\ga_k(t^-_k+\delta_0),y)}\,\sinh
d(\ga_k(t^-_k+\delta_0),\ga_{k-1}(\tau))\nonumber\\ 
& \leq  e^{s^+_k-t^-_k-\delta_0}
\,\sinh(\delta_0+c_1)
\end{align}
which is, in particular, at most $\sinh(\delta_0+c_1)$.

Let $q'$ be the closest point to $y$ on $\ga_{k-1}$, and let $p$
(resp.~$q$) be the closest point to $x'$ (resp.~$q'$) on $\ga_{k}$.
Then $d(x',p)\leq d(x,x')\leq c'_1(\epsilon_0)$ by Lemma
\ref{lemma:c1} if $\epsilon_0\neq\infty$ and Lemma \ref{lem:c1pinfty}
otherwise. As closest point maps do not increase distances, we have
$$
d(y,q)\leq d(y,q')=d(y,\ga_{k-1})\leq \sinh(\delta_0+c_1)\;.
$$ 
Note that
\begin{equation}
\label{eq:minohO}
d(x,y)= f_n(\ga_k)> h'_0\geq h_0\geq h'(\epsilon_0,
\sinh(\delta_0+c_1))\geq \sinh(\delta_0+c_1)+ c'_1(\epsilon_0)\;,
\end{equation} 
by the definition of $h_0$ and of $h'(\cdot,\cdot)$ in Equation
  \eqref{eq:h} if $\epsilon_0\neq \infty$ and by \eqref{eq:c3i}
  otherwise. Similarly, we have $d(x,y)\geq h_0\geq
  c_0(\epsilon_0)$. Hence $\xi_0,x',q'$ are in this order on
  $\ga_k$. Therefore, by convexity, 
$$
d(x',\ga_k)\leq d(q',\ga_k)\leq d(q',y)=d(y,\ga_{k-1}) \;.
$$ 
Hence, by  Lemma
\ref{lemma:c2} if $\epsilon_0\neq\infty$ and Lemma \ref{lem:c2pinfty}
otherwise, and by the inequality \eqref{eq:majdygakplu}, we have
\begin{equation}
\label{eq:majdxxprim}
d(x,x')\leq c_2\;d(x',\ga_k)\leq c_2\;e^{s^+_k-t^-_k-\delta_0} 
\,\sinh(\delta_0+c_1)\;.
\end{equation}

Furthermore, as we have seen that $d(x,y)\geq h'(\epsilon_0,
\sinh(\delta_0+c_1))$ and by the inequality \eqref{eq:majdygakplu}, it
follows from Lemma \ref{lem:h0} if $\epsilon_0\neq \infty$ and by
Lemma \ref{lem:c3pinfty} otherwise, that the geodesic line $\ga_{k-1}$
meets $C_n$ and one of the following two assertions hold :
\begin{equation}\label{form:1}
d(y,y')\leq c'_3(\epsilon_0) \;d(x',\ga_k)\leq c'_3(\epsilon_0) 
\;e^{s^+_k-t^-_k-\delta_0}
\;\sinh(\delta_0+c_1)
\end{equation}
or
\begin{equation}
\label{form:2}
d(x',y')\geq d(x,y)\;. 
\end{equation}

Before obtaining a contradiction from both of these assertions, we
prove a technical result.

\blemm\label{lem:smkmun} 
We have $\delta_0<s^-_{k-1}<\tau$, so that $\ga_{k-1}
(]\delta_0,+\infty[)$ meets $C_n$ with $t^-_n(\ga_{k-1})<
\tau$.  
\elemm

\dem Assume first by absurd that $s^-_{k-1}\leq \delta_0$.  If
$s^-_{k-1}\in\;]0,\delta_0]$, we have by the triangular inequality,
Lemma \ref{lem:c1prime} and the inequality \eqref{eq:majdxxprim} ,
\begin{align*}
s^-_{k} & =d(\ga_{k}(0),\ga_{k}(s^-_{k}))\\ 
& \leq 
d(\ga_{k}(0),\ga_{k-1}(0))+d(\ga_{k-1}(0),\ga_{k-1}(s^-_{k-1}))+
d(\ga_{k-1}(s^-_{k-1}),\ga_{k}(s^-_{k}))\\
& \le c_1+ \delta_0+c_2\,\sinh(\delta_0+c_1)\;.
\end{align*}
Let $z_0$ and $z_{s_{k-1}^-}$ be the closest points on $\gamma_k$ to
$\gamma_{k-1}(0)$ and $\gamma_{k-1}(s_{k-1}^-)$, respectively. If
$s^-_{k-1}\leq 0$, then as the closest point projection does not
increase distances, we have
\begin{align*}
s^-_{k}=\; & d(\ga_{k}(0),\ga_{k}(s^-_{k}))\leq 
d\big(\ga_k(0),z_0)\big)+d\big(z_0,\ga_k(s_k^-)\big)\\
\le \; & d\big(\ga_k(0),z_0)\big)+d\big(z_{s_{k-1}^-},\ga_k(s_k^-)\big)\\
\le\; & 
d(\ga_{k}(0),\ga_{k-1}(0))+d(\ga_{k-1}(s^-_{k-1}),\ga_{k}(s^-_{k}))\\
\leq\; &
c_1+c_2\,\sinh(\delta_0+c_1)\;.
\end{align*}
Hence, by the definition of $c_5$ and as $c'_3(\epsilon_0)\geq 1$ (see
the equation \eqref{eq:c3} if $\epsilon_0\neq\infty$ or \eqref{eq:c3i}
otherwise), we have
$$
c_5\geq
c_2\,\sinh(\delta_0+c_1)+c'_3(\epsilon_0)\sinh(\delta_0+c_1)\geq
c_2\,\sinh(\delta_0+c_1)+ \delta_0+c_1\geq s^-_{k}\;.
$$ 
Now $\ell_{C_n}(\ga_k)\geq h_0>\delta_0$ by the definition of $h_0$,
and
$$
\ell_{C_0}(\ga_k)\geq f_0(\ga_k)-\kappa_0=h-\kappa_0\geq h'_1-\kappa_0\geq
h_0+2c_5-\kappa_0>\delta_0+c_5\geq \delta_0+s^-_{k}\;.
$$
As $s^-_{k}\geq 0$ is the entrance time of $\ga_k$ in $C_n$, this
implies that ${\rm diam}(C_0\cap C_n)>\delta_0$. As $n\neq 0$, this
contradicts the Almost disjointness property $(iii)$, hence
$\delta_0<s^-_{k-1}$.

Assume now by absurd that $s^-_{k-1}\geq\tau$. Then as in the case
$s^-_{k-1}\le 0$, we get
$$
t^-_k - s^-_k\leq d(\ga_{k-1}(\tau),\ga_{k}(t^-_k))+
d(\ga_{k-1}(s^-_{k-1}),\ga_{k}(s^-_k))\leq c_1+c'_1(\epsilon_0)\;,
$$
by Lemma \ref{lem:c1prime}, and by Lemma \ref{lemma:c1} if $\epsilon_0
\neq\infty$ and Lemma \ref{lem:c1pinfty} otherwise.  We have seen in
the inequalities \eqref {eq:minohO} that 
$$
h_0\geq \sinh(\delta_0+c_1) +
c'_1(\epsilon_0) > \delta_0 + c_1 + c'_1(\epsilon_0).
$$ 
Hence
$$
t^-_k \geq s^+_k-\delta_0\geq s^-_k+h_0-\delta_0> 
s^-_k+c_1+ c'_1(\epsilon_0)\;,
$$
a contradiction. Hence $s^-_{k-1}<\tau$. 
\cqfd

\medskip
Assume first that the inequality \eqref{form:1} holds. As $s_k^-\geq
0$ by Lemma \ref{lem:enterpositiv} and by the definition of $h_0$, we
have
$$
s^+_k> h'_0+s^-_k \geq h_0\geq \delta_0+2c_1+c_6\;.
$$ 
Hence, as $\tau\geq t^-_k-2c_1$ by Lemma \ref{lem:tau}, we have
$e^{-\tau}\leq e^{-c_6}\;e^{s^+_k-\delta_0-t^-_k}$. By the definition
of $c_6$ and of $c_4$, we have 
$$
c_6=3c_4+\log 2\geq 3\;c'_3(\epsilon_0)
\sinh (c_1+\delta_0)+\log 2\;. 
$$
By the triangular inequality since $s^+_{k-1}\geq 0$ by
\ref{lem:smkmun}, by the equations \eqref{form:1} and
\eqref{eq:majdistoo}, and by the definition of $c_4$, we hence have
\begin{eqnarray}
  |s^+_k-s^+_{k-1}| &\leq & 
              d(y,y')+d(\ga_k(0),\ga_{k-1}(0))\nonumber\\
  &\leq & c'_3(\epsilon_0) \;e^{s^+_k-t^-_k-\delta_0}
  \;\sinh(\delta_0+c_1) + c_2\,e^{-\tau}\sinh c_1 \nonumber\\
& \leq &
  c_4\;e^{s^+_k-\delta_0-t^-_k}\;.\label{eq:majspkmspkm}
\end{eqnarray}
By the Lipschitz property $(ii)$ of $f_n=\ell_{C_n}$ (as $n\neq 0$), by
the inequalities \eqref{eq:majdxxprim} and \eqref{form:1}, and by the
definition of $c_5$, we have
\begin{eqnarray}
|f_n(\ga_{k-1})-f_n(\ga_{k})|&\leq &\!\!
2\max\{d(x,x'),d(y,y')\}\nonumber\\
&\leq &\!\!
2\max\{c_2 \;e^{s^+_k-\delta_0-t^-_k}
\;\sinh(\delta_0+c_1), c'_3(\epsilon_0) \;e^{s^+_k-\delta_0-t^-_k}
\;\sinh(\delta_0+c_1)\}\nonumber\\
&\leq & \!\!c_5\;e^{s^+_k-\delta_0-t^-_k}\;.\label{eq:valabs}
\end{eqnarray}

By the inequalities \eqref{eq:absurdn} and \eqref{eq:valabs}, by the
minimality property of $\tau$, and by Lemma \ref{lem:smkmun}, we have
$$
u_k(s^+_k-\delta_0)< f_n(\ga_{k})\leq f_n(\ga_{k-1})
+c_5\;e^{s^+_k-\delta_0-t^-_k}
\leq u_{k-1}(s^+_{k-1}-\delta_0)+c_5\;e^{s^+_k-\delta_0-t^-_k}\;.
$$

Assume now that the inequality \eqref{form:2} holds instead of the
inequality \eqref{form:1}.  Then $f_n(\ga_k)\leq f_n(\ga_{k-1})$, so
we again have that
$$
u_k(s^+_k-\delta_0)<
u_{k-1}(s^+_{k-1}-\delta_0)+c_5\;e^{s^+_k-\delta_0-t^-_k}.
$$
As $|(s^+_{k-1}-\delta_0)-(s^+_k-\delta_0)|\leq
c_4\;e^{s^+_k-\delta_0-t^-_k}$ by the inequality \eqref{eq:majspkmspkm},
this contradicts the assertion (4) on the map $u_k$. Hence the
assertion (6) at rank $k$ is verified.

\bigskip\noindent{\bf The main corollary of the construction.} %
The above inductive construction will only be used in this paper
through the following summarizing statement.

\bprop\label{prop:main} %
Let $X$ be a proper geodesic ${\rm CAT}(-1)$ metric space. Let
$\epsilon_0$ in $\RR^*_+\cup\{\infty\}$, $\delta_0,\kappa_0\geq 0$ and
$\xi_0\in X\cup\partial_\infty X$. Let $h'_0\geq
h_0(\epsilon_0,\delta_0,\kappa_0)$ and $h\geq h'_1=
h'_1(\epsilon_0,\delta_0,h'_0)$. Let $(C_n)_{n\in\NN}$ be a collection
of $\epsilon_0$-convex subsets of $X$ which satisfies the assertions
$(iii)$ and $(iv)$, and with $\xi_0\notin C_0\cup\partial_\infty C_0$. Let
$f_0:T^1_{\xi_0}X\ra[0,+\infty]$ be a continuous
$\kappa_0$-penetration map in $C_0$.  Assume that there exists a
geodesic ray or line $\ga_0$ starting from $\xi_0$ with
$f_0(\ga_0)=h$. Then there exists a
geodesic ray or line $\ga_\infty$ starting from $\xi_0$, entering
$C_0$ at time $t=0$ with $f_0(\ga_\infty)=h$, such that
$\ell_{C_n}(\ga_\infty)\leq h'_1$ for every $n$ in $\NN-\{0\}$ such
that $\ga_\infty(]\delta_0,+\infty[)$ meets $C_n$.  
\eprop

\dem %
Apply the main construction of the previous subsections with initial
input a geodesic ray or line $\ga_0$ entering $C_0$ at time $t=0$ with
$f_0(\ga_0)=h$, to get finite or infinite sequences
$(\ga_k)_{k\in\N}$, $(n_k)_{k\in\N}$, $(u_k)_{k\in\N}$ satisfying the
assertions (1)--(6).

If $\N$ is finite, with maximum $N$, define $\ga_k=\ga_N$ for $k> N$.
Then the sequence $(\ga_k)_{k\in\NN}$ converges to a geodesic ray or
line $\ga_\infty=\ga_N$ in $T^1_{\xi_0}X$. If $\N$ is infinite, as $X$
is complete, it follows from the assertions (3) and (5), by an easy
geometric series argument, that the sequence $(\ga_k)_{k\in \N}$
converges in $T^1_{\xi_0}X$ to a geodesic ray or line $\ga_\infty$
starting from $\xi_0$ and entering $C_0$ at time $t=0$, as $C_0$ is
closed and convex. By the continuity of $f_0$ and the assertion (1), we have
$f_0(\ga_\infty)=h$.

Suppose by absurd that there exists $n$ in $\NN-\{0\}$ such that
$\ga_\infty(]\delta_0,+\infty[)$ meets $C_n$ and $\ell_{C_n}(\ga_\infty) >
h'_1>0$.  In particular, $\ga_\infty(]\delta_0,+\infty[)$ meets the
interior of $C_n$ and $\xi_0\notin C_n\cup\partial_\infty C_n$ by
Lemma \ref{lem:enterpositiv}. Furthermore, it follows from the
definition of the stopping time, and the fact that $u_k\leq h'_1$ for
every $k$, that $\N$ is infinite. Hence, as the $\ga_k$'s converge to
$\ga_\infty$, and by the continuity of $\ell_{C_n}$, if $k$ is big
enough, then $\ga_k(]\delta_0,+\infty[)$ meets $C_n$ and
$\ell_{C_n}(\ga_k) > h'_1$.

In particular, $t_n^+(\ga_k)> \delta_0$.  Note that $t_n^-(\ga_k)$,
which is at distance at most $c'_1(\epsilon_0)$ from
$t_n^-(\ga_\infty)$ by Lemma \ref{lemma:c1} if $\epsilon_0\neq \infty$
and Lemma \ref{lem:c1pinfty} otherwise, is bounded as $k$ tends to
$\infty$. Hence if $k$ is big enough, then $t_n^-(\ga_k)$ is less than
$t^-_k$, as $t^-_k$ converges to $+\infty$ when $k\ra+\infty$ by the
assertion (5). This contradicts the assertion (6), as $u_k\leq h'_1$.
\cqfd

\medskip 
\rem %
If $X,\epsilon_0,\delta_0,\kappa_0,\xi_0, h'_0,h, (C_n)_{n\in\NN},f_0$
satisfy the hypotheses in the statement of Proposition
\ref{prop:main}, and if for every $n$ such that $\xi_0 \notin
C_n\cup\partial_\infty C_n$, we have a $\kappa$-penetration map $g_n:
T^1_{\xi_0}X\ra[0,+\infty]$ for some constant $\kappa\geq 0$, then
Proposition \ref{prop:main} implies that there exists a geodesic ray
or line $\ga_\infty$ starting from $\xi_0$, entering $C_0$ at time
$t=0$ with $f_0(\ga_\infty)=h$, such that $g_n(\ga_\infty)\leq
h'_1+\kappa$ for every $n$ in $\NN-\{0\}$ such that
$\ga_\infty(]\delta_0,+\infty[)$ meets $C_n$. We will apply this
observation to more general penetration maps
than the $\ell_{C_n}$'s, in Section \ref{sec:results}.

\medskip %
The next corollary yields geodesic lines with the prescribed
penetration in $C_0$, and that essentially avoid the $C_n$'s not only
for positive times, but also for negative ones. The penetration in the
sets $C_n$ for $n\ne 0$ cannot be made quite as small as in
Proposition \ref{prop:main}.

\bcoro\label{coro:main} %
Let $X$ be a proper geodesic ${\rm CAT}(-1)$ metric space. Let
$\epsilon_0$ in $\RR^*_+\cup\{\infty\}$, $\delta_0,\kappa_0\geq 0$.
Let $C_0$ be a proper $\epsilon_0$-convex subset of $X$, and let
$$
f_0 : \bigcup\limits_{\xi\in \partial_\infty X -\partial_\infty{C_0}}
T^1_\xi X \ra[0,+\infty]
$$ 
be a continuous map such that ${f_0}_\mid{}_{T^1_{\xi_0}X}$ is a
$\kappa_0$-penetration map in $C_0$ for every $\xi_0\in
\partial_\infty X -\partial_\infty C_0$. Let $h'_0\geq
h_0=h_0(\epsilon_0,\delta_0, \kappa_0)$, $h\geq
h'_1=h'_1(\epsilon_0,\delta_0,h'_0)$, and
$$
h''_1=h'_1(\epsilon_0,\delta_0,h'_0)+
c'_3(\epsilon_0) (\delta_0+c_1)+c'_1(\epsilon_0)\;
$$
Assume that there exists a geodesic line $\ga_0$ in $X$ with
$f_0(\ga_0)=h$. For every $n$ in $\NN-\{0\}$, let $C_n$ be an
$\epsilon_0$-convex subset of $X$, such that $(C_n)_{n\in\NN}$
satisfies the assertions $(iii)$ and $(iv)$ with respect to every
$\xi_0\in \partial_\infty X -\partial_\infty C_0$.  Then there exists
a geodesic line $\ga'$ in $X$ entering $C_0$ at time $t=0$ with
$f_0(\ga')=h$, such that $\ell_{C_n}(\ga')\leq h''_1$ for every $n$ in
$\NN-\{0\}$.
\ecoro

\dem %
Let $\ga_0$ be a geodesic line in $X$ with $f_0(\ga_0)=h$, and let
$\xi$ be the starting point at infinity of $\ga_0$, which does not
belong to $\partial_\infty C_0$ as $h<\infty$. Applying
Proposition \ref{prop:main} with $\xi_0=\xi$, as $h\geq h'_1$, there
exists a geodesic line $\ga$ starting from $\xi$ and entering $C_0$ at
time $0$, such that $f_0(\ga) =h$ and $\ell_{C_n}(\ga)\leq h'_1$ for
every $n\in\NN-\{0\}$ such that $\ga(]\delta_0,+\infty[)$ meets $C_n$.

Let $\xi'$ be the other endpoint at infinity of $\ga$, which does not
belong to $\partial_\infty C_0$ as $h<\infty$.  Applying Proposition
\ref{prop:main} again with now $\xi_0=\xi'$, we get that there exists
a geodesic line $\gamma'$ starting from $\xi'$ and entering $C_0$ at
time $0$, such that $f_0(\gamma')=h$ and $\ell_{C_n}(\gamma')\leq
h'_1$ for every $n\in\NN-\{0\}$ such that $\ga'(]\delta_0,+\infty[)$
meets $C_n$.

Assume by absurd that there exists $n\in\NN-\{0\}$ such that
$\ell_{C_n}(\gamma') > h''_1>0$. Then $\gamma$ enters $C_n$ at a point
$x'_n$, exiting it at a point $y'_n$ at time at most $\delta_0$, as
$h''_1>h'_1$ by the definition of $h''_1$. In particular, if
$x'=\gamma'(0)$ is the entering point of $\gamma'$ in $C_0$, then
$d(y'_n,x')\leq \delta_0$ if $x',y'_n,x'_n, \xi'$ are not in this
order on $\ga'$.

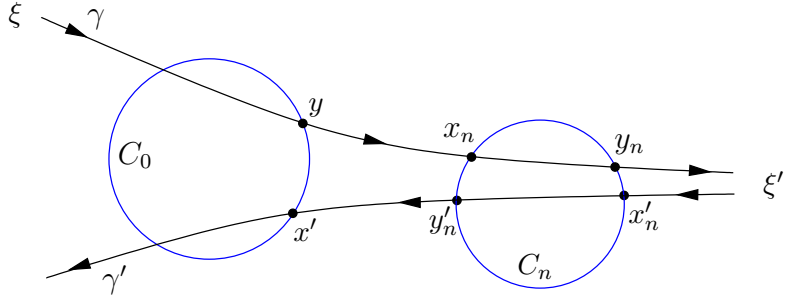
\begin{center} 
\input{fig_horlingene.pstex_t}
\end{center}

Let $y$ be the exiting point of $\gamma$ out of $C_0$. Note that 
\begin{equation}\label{eq:minohunp}
h\geq h'_1 \geq h'_0\geq h_0 \geq 
h'(\epsilon_0,\sinh(\delta_0+ c_1))\geq 
h'(\epsilon_0,\delta_0+ c_1)
\end{equation}
by the definitions of $h'_1, h_0, h'$. By Lemma \ref{lemma:c1} if
$\epsilon_0\neq \infty$ and by Lemma \ref{lem:c1pinfty} if $\epsilon_0
=\infty$ and $(f_0,\delta_0) \neq (\ph_{C_0},0)$, and as in the proof
of Lemma \ref{lem:c1prime} if $(\epsilon_0, f_0, \delta_0) = (\infty,
\ph_{C_0},0)$ since $h\geq h_0$, we have $d(x',y) \leq c_1$. Hence by
convexity, 
$$
d(y'_n,\ga)\leq d(x',\ga) +\delta_0 \leq d(x',y) +
\delta_0 \leq \delta_0+ c_1\;. 
$$
Note that
$$
d(x_n',y'_n)=\ell_{C_n}(\ga')>h''_1\geq h'_1\geq 
h'(\epsilon_0,\delta_0+ c_1)
$$
by the definition of $h''_1$ and by the inequalities \eqref{eq:minohunp}.
Hence, by Lemma \ref{lem:h0} if $\epsilon_0\neq \infty$ and by Lemma
\ref{lem:c3pinfty} otherwise, the geodesic line $\ga$ enters $C_n$ at
a point $x_n$ and exits it at a point $y_n$ such that
\begin{equation}\label{eq:dichobis2}
d(y'_n,x_n)\leq c'_3(\epsilon_0) d(y'_n,\ga)\;\;
{\rm or}\;\;d(x_n,y_n)\geq d(x'_n,y'_n)\;. 
\end{equation}

Let us prove by absurd that $\ga(]\delta_0,+\infty[)$ meets $C_n$.
Otherwise, since $\ga^{-1}(y)\geq 0$, by convexity, and by Lemma
\ref{lemma:c1} if $\epsilon_0\neq \infty$ or Lemma \ref{lem:c1pinfty}
if $\epsilon_0 =\infty$, we have
\begin{equation}\label{eq:adf}
d(x',x'_n)\leq d(x',y)+\delta_0+d(y_n,x'_n)\leq
c_1+\delta_0+c'_1(\epsilon_0)\;.
\end{equation}
By the inequalities \eqref{eq:minohunp} and by the definition of
$h'(\epsilon,\eta)$, we have 
$$
h'_1 \geq h'(\epsilon_0,\sinh(\delta_0+
c_1)) \geq 2\sinh(\delta_0+ c_1) \geq 2\delta_0+c_1\;. 
$$
Hence the inequalities \eqref{eq:adf} contradicts the fact that, by
the definition of $h''_1$,
$$
d(x',x'_n)\geq d(x'_n,y'_n) - \delta_0  
> h''_1-\delta_0\geq h'_1-\delta_0+c'_1(\epsilon_0)
\geq c_1+\delta_0+c'_1(\epsilon_0)\;.
$$
Assume that the second of the inequalities \eqref{eq:dichobis2} holds
true. As $d(x_n',y'_n)> h'_1$, this contradicts the construction of
$\ga$.

Hence we have 
$$
d(y'_n,x_n)\leq c'_3(\epsilon_0) d(y'_n,\ga)\leq
c'_3(\epsilon_0)(\delta_0+ c_1)\;.
$$  
But then, by the triangular
inequality and by the definition of $h''_1$,
$$
d(x_n,y_n)\geq d(x'_n,y'_n)-d(x_n,y'_n)-d(y_n, x'_n)>
h''_1-c'_3(\epsilon_0)(\delta_0+ c_1)-c'_1(\epsilon_0)= h'_1\;,
$$
which contradicts the construction of $\ga$.  \cqfd

\section{Prescribing the penetration of geodesic lines}
\label{sec:results}

%
%
%

In this section, we apply Proposition \ref{prop:main} to
prove a number of results on the geodesic flow of negatively curved
Riemannian manifolds. 

The following constants appear in the theorems, depending on
$\epsilon\in \RR^*_+\cup\{\infty\}, \delta,\kappa\geq 0$.
\begin{itemize}
\item %
  $c''_1=c''_1(\epsilon,\delta,\kappa)=
  \max\big\{2c'_1(\epsilon)+2\delta+\kappa,
  h'_1(\epsilon,\delta,h_0(\epsilon,\delta,c'_1(\infty))\big\}$.
\item %
  $c''_2(\epsilon)=c''_1(\epsilon,0,0)+c'_1(\infty)+2c_1$, where
  $c_1=c'_1(\epsilon)$ if $\epsilon\neq\infty$, and $c_1=1/19$
  otherwise. Note that $ c''_2(\infty)=
  h'_1(\infty,0,h_0(\infty,0,c'_1(\infty))+c'_1(\infty)+2c_1\approx
  8.3712 $ by the definition of $c''_1$ and the approximation
  \eqref{eq:valeurhun}.
\end{itemize}
Recall that the constants $c'_1(\epsilon)$ are given by Lemmas
\ref{lemma:c1} and \ref{lem:c1pinfty}, and that
$h_0(\cdot,\cdot,\cdot)$ and $h'_1(\cdot,\cdot,\cdot)$ are given in
the list of constants in the beginning of the subsection
\ref{sec:induction}.

\subsection{Climbing in balls and horoballs}
\label{sec:ballsandhoroballs}

In this subsection, we construct geodesic rays or lines having
prescribed penetration properties in a ball or a horoball, while
essentially avoiding a family of almost disjoint convex subsets.  Let
us consider the penetration height and inner projection penetration
maps first in horoballs and then in balls. Note that if $C_0$ is
a ball or an horoball, if $f_0=\ph_{C_0}$, then
$\norm{f_0-\ph_{C_0}}=0$ and if $f_0=\ipp_{C_0}$, then
$\norm{f_0-\ph_{C_0}}\leq c'_1(\infty)$ by Section
\ref{subsec:penetrationmaps}.

\btheo\label{theo:line1} %
Let $\epsilon\in \RR^*_+\cup\{\infty\}$, $\delta,\kappa\geq 0$; let
$X$ be a complete simply connected Riemannian manifold with sectional
curvature at most $-1$ and dimension at least $3$; let $\xi_0\in X
\cup \partial_\infty X$; let $C_0$ be a horoball such that $\xi_0
\notin C_0\cup\partial_\infty C_0$; let $f_0=\ph_{C_0}$ or
$f_0=\ipp_{C_0}$; let $(C_n) _{n\in\NN-\{0\}}$ be a family of
$\epsilon$-convex subsets of $X$; for every $n\in\NN-\{0\}$ such that
$\xi_0 \notin C_n\cup\partial_\infty C_n$, let $f_n:T^1_{\xi_0}X\ra
[0,+\infty]$ be a $\kappa$-penetration map in $C_n$.  If
$\diam(C_n\cap C_m)\leq\delta$ for all $n,m$ in $\NN$ with $n\neq m$,
then, for every $h\geq c''_1(\epsilon,\delta,\norm{f_0-\ph_{C_0}})$,
there exists a geodesic ray or line $\gamma$ starting from $\xi_0$ and
entering $C_0$ at time $0$, such that $f_0(\gamma)=h$ and
$f_n(\gamma)\leq c''_1(\epsilon,\delta,\norm{f_0-\ph_{C_0}})+\kappa$
for every $n\geq 1$ such that $\gamma\big(]\delta,+\infty[\big)$ meets
$C_n$.  \etheo

\medskip \dem %
let $h\geq c''_1$. In order to apply Proposition \ref{prop:main},
define $\epsilon_0 =\epsilon,\delta_0=\delta,
\kappa_0=2\log(1+\sqrt{2}) =c'_1(\infty)$ and
$h'_0=h_0(\epsilon_0,\delta_0,\kappa_0)$. Recall that $\ph_{C_0}$ and
$\ipp_{C_0}$ are $\kappa_0$-penetration maps for $C_0$ by Lemma
\ref{lem:proppenehoro}.  For every $n\in\NN-\{0\}$ such that $\xi_0
\notin C_n\cup\partial_\infty C_n$, let us apply Proposition
\ref{prop:prescribe} Case (1) to $C= C_0$, $C'=C_n$, $f=f_0$,
$f'=\ell_{C_n}$, $h'=h'_0$, so that $h^{\rm min}=
2\,c'_1(\epsilon)+2\delta+\norm{f_0-\ph_{C_0}}$ and $h_0^{\rm min} =
2\delta$. Note that $h_0^{\rm min} \leq h'_0$, as
$$
h'_0\geq h'(\epsilon,\sinh(\delta+c_1))\geq
2\sinh(\delta+c_1)\geq 2\delta\;,
$$ 
by the definition of $h_0$ and of $h'(\cdot,\cdot)$. As $h\geq
c''_1\geq h^{\rm min}$ by the definition of $c''_1$, Proposition
\ref{prop:prescribe} (1) hence implies that $(C_n)_{n\in\NN}$ satisfies
the Local prescription property $(iv)$. Thus by Proposition
\ref{prop:main}, there exists a geodesic ray or line $\gamma$ starting
at $\xi_0$ such that $f_0(\gamma)=h$ and $\ell_{C_n}(\gamma)
\leq c''_1$, which implies that $f_n(\gamma)\leq c''_1+\kappa$, for
every $n\geq 1$ such that $\gamma(]\delta,+\infty[)$ meets $C_n$.
\cqfd

\medskip %
The proof of the corresponding result when $C_0$ is a ball of radius
$R\geq\epsilon$ is the same, using Case (2) of Proposition
\ref{prop:main} instead of Case (1). This requires $h\leq h^{\rm
  max}= 2R-2c'_1(\epsilon)-\norm{f_0-\ph_{C_0}}$. To be nonempty, the
following result requires 
$$
R\geq c''_1(\epsilon, \delta,
\norm{f_0-\ph_{C_0}})/2+ \;c'_1(\epsilon) +
\norm{f_0-\ph_{C_0}}/2\;.
$$

\btheo\label{theo:line} %
Let $\epsilon>0$, $\delta,\kappa\geq 0$; let $X$ be a complete simply
connected Riemannian manifold with sectional curvature at most $-1$
and dimension at least $3$; let $C_0$ be a ball of radius $R\geq
\epsilon$; let $\xi_0\in X \cup (\partial_\infty X)-C_0$; let $f_0=
\ph_{C_0}$ or $f_0=\ipp_{C_0}$; let $(C_n) _{n\in\NN-\{0\}}$ be a
family of $\epsilon$-convex subsets of $X$; for every $n\in\NN-\{0\}$
such that $\xi_0 \notin C_n\cup\partial_\infty C_n$, let $f_n:
T^1_{\xi_0} X \ra [0,+\infty]$ be a $\kappa$-penetration map $C_n$.
If $\diam(C_n\cap C_m)\leq\delta$ for all $n,m$ in $\NN$ with $n\neq
m$, then, for every
$$
h\in \Big[c''_1\big(\epsilon,
\delta,\norm{f_0-\ph_{C_0}}\big),\; 2R - 2c'_1(\epsilon)-
\norm{f_0-\ph_{C_0}}\;\Big]\;,
$$ there exists a geodesic ray or line $\gamma$
starting from $\xi_0$ and entering $C_0$ at time $0$, such that
$f_0(\gamma)=h$ and $f_n(\gamma)\leq
c''_1\big(\epsilon,\delta,\norm{f_0-\ph_{C_0}}\big)+\kappa$ for every
$n\geq 1$ such that $\gamma(]\delta,+\infty[)$ meets $C_n$.  \cqfd
\etheo

Varying the family $(C_n) _{n\in\NN-\{0\}}$ of $\epsilon$-convex
subsets appearing in Theorems \ref{theo:line1} and \ref{theo:line},
among balls of radius at least $\epsilon$, horoballs,
$\epsilon$-neighbourhoods of totally geodesic subspaces, etc, we get
several corollaries. We will only state two of them, Corollaries
\ref{coro:horoline} and \ref{coro:ballline}, which have applications
to equivariant families. The proofs of these results are simplified
versions of the proof of Corollary \ref{coro:main}, giving better
(though very probably not optimal) constants.

\bcoro\label{coro:horoline} %
Let $X$ be a complete simply connected Riemannian manifold with
sectional curvature at most $-1$ and dimension at least $3$, and let
$\big(H_n\big)_{n\in\NN}$ be a family of horoballs in $X$ with
disjoint interiors. Then, for every $h\geq c''_1(\infty,0,0)\approx
6.5032$, there exists a geodesic line $\gamma'$ such that
$\ph_{H_0}(\gamma')=h$ and $\ph_{H_n}(\gamma')\leq c''_2(\infty)
\approx 8.3712$ for every $n\geq 1$.  
\ecoro

\dem %
Let $C_0=H_0$ and let $\xi$ be a point in $\partial_\infty
X-\partial_\infty C_0$. We apply Theorem \ref{theo:line1} with
$\epsilon=\infty$, $\delta=0$, $\kappa=0$, $\xi_0=\xi$, $C_n=H_n$ for
every $n$ in $\NN$, $f_0=\ph_{C_0}$, and $f_n=\ell_{C_n}$ for every
$n\neq 0$ such that $\xi_0\notin C_n\cup\partial_\infty C_n$. Note
that for every $n\in\NN$, $f_n$ is a $\kappa$-penetration map in
$C_n$.  As $h\ge c''_1(\epsilon,0,0)$, there exists a geodesic line
$\gamma$ starting from $\xi$ and entering $C_0$ at time $0$, such that
$\ph_{C_0}(\gamma) =h$ and $\ell_{C_n}(\gamma)\leq
c''_1(\epsilon,0,0)$ for every $n\in\NN-\{0\}$ such that $\gamma$
meets $C_n$ at a positive time.

Let $\xi'$ be the other endpoint of $\gamma$. This point is not in
$\partial_\infty C_0$.  Applying Theorem \ref{theo:line1} again, as
above except that now $\xi_0=\xi'$, we get that there exists a
geodesic line $\gamma'$ starting from $\xi'$ and entering $C_0$ at
time $0$, such that $\ph_{C_0}(\gamma')=h$ and
$\ell_{C_n}(\gamma')\leq c''_1(\epsilon,0,0)$ for every $n\in\NN-\{0\}$
such that $\gamma'$ meets $C_n$ at a positive time.

Let $c''_2=c''_2(\epsilon)$. Assume by absurd that there exists
$n\in\NN-\{0\}$ such that $\ph_{C_n}(\gamma')> c''_2>0$. Then $\gamma'$
enters $C_n$ at the point $x'_n$, exiting it at the point $y'_n$ at a
nonpositive time, as $c''_2>c''_1(\epsilon,0,0)$. In particular, if
$x'=\gamma'(0)$ is the entering point of $\gamma'$ in $C_0$, then
$x',y'_n,x'_n,\xi'$ are in this order on $\ga'$ (see the picture in
the proof of Corollary \ref{coro:main}).

Let $y$ be the exiting point of $\gamma$ out of $H_0$.  With
$c_1=1/19$, as in the proof of Lemma \ref{lem:c1prime}, since
$\ph_{C_0}(\ga)$ and $\ph_{C_0}(\ga')$ are equal to
$$
h\geq c''_1(\infty,0,0)\geq h'_1(\infty,0, h_0(\infty,0,c'_1(\infty)))
\geq h_0(\infty,0,c'_1(\infty))=h'(\infty,\sinh c_1)
$$
by the definition of $c''_1, h'_1, h_0$, we have $d(x',y)\leq c_1$.

Let $\xi_n$ be the point at infinity of $H_n$. Let $p'$ be the point
in $[x'_n,y'_n]$ the closest to $\xi_n$, so that 
$$
d(p',y'_n)\geq
\beta_{\xi_n}(y'_n,p') = \ph_{C_n}(\ga')/2 > c''_2/2.
$$ 
Let $p$ be
the point of $\ga$ the closest to $p'$. As, by convexity and the
definition of $c''_2$, we have 
$$
d(p',p)=d(p',\ga)\leq d(x',\ga)\leq
d(x',y)\leq c_1< c''_2/2,
$$ 
it follows that $p$ belongs to the interior
of $C_n$. If $p\in\;]\xi,y]$, then by convexity,
$$
c''_2/2< d(p',y'_n)\leq d(p',x')\leq d(p',p)+d(y,x')\leq 2c_1\;,
$$
a contradiction, as by the definition of $c''_2$, of $c''_1$ and of
$h'(\epsilon,\eta)$ (see the equations (\ref{eq:h}) and
(\ref{eq:c3i})), we have
$$
c''_2\geq c''_1(\epsilon,0,0)+2\,c_1\geq h'(\epsilon,\sinh
c_1)+2\,c_1\geq 2\sinh c_1+ 2c_1 > 4\,c_1\;.
$$ 
Hence $p\in \;]y,\xi'[\;\subset \ga(]0,+\infty[)$, so that $\ga$
meets $C_n$ at a positive time. But, by Proposition
\ref{lem:proppenehoro} and the definition of $c''_2$,
\begin{align*}
\ell_{C_n}(\ga)\geq\;& \ph_{C_n}(\ga)-c'_1(\infty)\geq
2\beta_{\xi_n}(y'_n,p)-c'_1(\infty)\geq
2(\beta_{\xi_n}(y'_n,p')-d(p,p'))-c'_1(\infty)\\ 
> \; &
2(c''_2/2-c_1)-c'_1(\infty)=c''_1(\epsilon,0,0)\;.
\end{align*}
This contradicts the construction of $\ga$.
\cqfd

\bigskip %
Let $M$ be a complete nonelementary geometrically finite Riemannian
manifold with sectional curvature at most $-1$ (see for instance
\cite{Bow} for a general reference). Recall that a {\it cusp} of $M$
is an asymptotic class of minimizing geodesic rays in $M$ along which
the injectivity radius converges to $0$. If $M$ has finite volume,
then the set of cusps of $M$ is in bijection with the (finite) set of
ends of $M$, by the map which associates to a representative of a cusp
the end of $M$ towards which it converges.  Let $\pi:\wt{M}\ra
M$ be a universal Riemannian covering of $M$, with covering group
$\Ga$. If $e$ is a cusp of $M$, and $\rho_e$ a minimizing geodesic ray
in the class $e$, as $M$ is geometrically finite and nonelementary,
there exists (see for instance \cite{BK,Bow,HP2}) a (unique) maximal
horoball $H_e$ in $\wt{M}$ centered at the point at infinity
$\xi_e$ of a fixed lift of $\rho_e$ in $\wt{M}$, such that $\ga
H_e$ and $H_e$ have disjoint interiors if $\ga\in\Ga$ does not fix
$\xi_e$. The image $V_e$ of $H$ in $M$ is called the {\it maximal
  Margulis neighborhood of $e$}. If $\rho_e$ starts from a point in
the image by $\pi$ of the horosphere bounding $H$, then let
$\height_e:M\ra \RR$ be map defined by
$$
\height_e(x)= \lim_{t\ra\infty} \;(t-d(\rho_e(t),x))\;,
$$
called the {\it height function with respect to $e$}. Let
$\maxht_e:T^1M\to\RR$ be defined by
$$
\maxht_e(\gamma)=\sup_{t\in\RR}\height_e(\gamma(t))\;.
$$
The {\it maximum height spectrum} of the pair $(M,e)$ is the subset of
$]-\infty,+\infty]$ defined by
$$
\maxsp(M,e)=\maxht_e(T^1M)\;.
$$

\bcoro\label{coro:maxsp} %
Let $M$ be a complete, nonelementary geometrically finite Riemannian
manifold with sectional curvature at most $-1$ and dimension at least
$3$, and let $e$ be a cusp of $M$. Then $\maxsp(M,e)$ contains
$[c''_2/2,+\infty]$.
\ecoro

Note that $c''_2/2\approx 4.1856$, hence Theorem
\ref{theo:introhallhoro} of the introduction follows. 

\medskip \dem %
With the above notations, let $(H_n)_{n\in\NN}$ be the
$\Gamma$-equivariant family of horoballs in $\widetilde M$ with
pairwise disjoint interiors such that $H_0=H_e$. Apply Corollary
\ref{coro:horoline} to this family to get, for every $h\geq c''_2\geq
c''_1(\infty,0,0)$, a geodesic line $\wt \ga$ in $\widetilde M$ with
$\ph_{H_0}(\wt \ga)=h$ and $\ph_{H_n}(\wt \ga)\leq c''_2 $ for every
$n\geq 1$.  Let $\ga$ be the locally geodesic line in $M$ image by
$\pi$ of $\wt \ga$. Observe that $\height_e\circ\pi=\beta_{H_n}$ in
$H_n$ and that $\ph_{H_n}(\wt \ga)=2\sup_{t\in \RR}\beta_{H_n}(\wt
\ga(t))$ (see Section \ref{subsec:penetrationmaps}).  Hence
$\sup_{t\in\RR} \height_e(\gamma(t))=h/2$ and the result follows.
\cqfd

\medskip Schmidt and Sheingorn \cite{SS} treated the case of
two-dimensional manifolds of constant curvature $-1$ (hyperbolic surfaces)
with a cusp. They showed that in that case $\maxsp(M,e)$ contains the
interval $[\log 100,+\infty]\approx [4.61,+\infty]$. This paper
\cite{SS} was a starting point of our investigations, although the
method we use is quite different from theirs. 

\medskip Let $\P$ be a (necessarily finite) nonempty set of cusps of
$M$. For every $e$ in $\P$, choose a horoball $H_e$, with point at
infinity $\xi_e$ as above Corollary \ref{coro:maxsp}. The horoballs of
the family $(gH_e)_{g\in\Ga/\Ga_{\xi_e}\;,\;e\in \P}$ may have non
disjoint interiors. But as $M$ is geometrically finite and
nonelementary, there exists (see \cite{BK,Bow}) $t\geq 0$ such that
two distinct elements in $(gH_e[t])_{g\in\Ga/\Ga_{\xi_e}\;,\;e\in \P}$
have disjoint interiors. Let $t_\P$ be the lower bound of all such
$t$'s. For every $\ga\in T^1M$, define
$$
\maxht_\P(\gamma)=\max_{e\in \P}\maxht_e(\gamma)\;\;{\rm and}
\;\;\maxsp(M,\P)=\maxht_\P(T^1M)\;.
$$

\medskip \rem %
Let $\C$ be the set of all cusps of $M$. Under the same hypotheses as
in Corollary \ref{coro:maxsp}, the following two assertions hold, by
applying Corollary \ref{coro:horoline} to the family of horoballs
$(gH'_{e'}) _{g\in\Ga/\Ga_{\xi_{e'}} \;,\;e'\in \C}$ with
$H'_{e'}=H_e$ if $e'=e$, and $H'_{e'}=H_{e'}[t]$ for some $t$ big
enough otherwise, for the first assertion, and to the family
$(gH_e[t_\P]) _{g\in\Ga/\Ga_{\xi_e} \;,\;e\in \P}$ for the second one.

(1) For every cusp $e$ of $M$, there exists a constant $t\geq 0$ such
that for every $h\geq t$, there exists a locally geodesic line $\ga$
in $M$ such that $\maxht_e(\ga)=h$ and $\maxht_{e'}(\ga)\leq t$ for
every cusp $e'\neq e$ in $M$.

(2) Let $\P$ be a nonempty set of cusps of $M$.  Then $\maxsp(M,\P)$
contains the halfline $[(c''_2+t_\P)/2,+\infty]$.

\medskip

Now, we prove the analogs of Corollaries \ref{coro:horoline} and
\ref{coro:maxsp} for families of balls with disjoint interiors.
Let 
  $$
  R_0^{\rm min}=7\sinh c'_1(\infty)+\frac{3}{2}c'_1(\infty)\approx
  22.4431.
  $$

\bcoro\label{coro:ballline} %
Let $X$ be a complete simply connected Riemannian manifold with
sectional curvature at most $-1$ and dimension at least $3$, and let
$\big(B_n\big)_{n\in\NN}$ be a family of balls in $X$ with disjoint
interiors such that the radius $R_0$ of $B_0$ is at least $R_0^{\rm
  min}$. For every $h\in \big[c''_1(R_0^{\rm
  min},0,0),2R_0-2\,c'_1(R_0^{\rm min})\big]$, there exists a geodesic
line $\gamma$ in $X$ with $\ph_{B_0}(\gamma) =h$ and
$\ph_{B_n}(\gamma) \le c''_2(R_0^{\rm min})$ for all $n\geq 1$.
\ecoro

\dem %
We start by some computations. Let $\epsilon>0$. With
$c_1=c'_1(\epsilon)$ and $c_5=c_5(\epsilon,0)$ as in Subsection
\ref{sec:induction}, we have $c_5\geq 6\sinh c_1$ since
$c'_3(\epsilon)\geq 3$ by the definition of $c'_3(\epsilon)$
in Equation \eqref{eq:c3}. By the definition of $h_0$ in Subsection
\ref{sec:induction} and of $h'$ in Equation \eqref{eq:h}, we have
$h_0(\epsilon,0,c'_1(\infty))\geq h'(\epsilon,\sinh
c_1)\geq 2\sinh c_1$. Hence, by the definition
of $c''_2, c''_1,h'_1$ and as $\epsilon\mapsto c'_1(\epsilon)$ is
decreasing,
\begin{align*}
c''_2(\epsilon)=\; &
c''_1(\epsilon,0,0)+c'_1(\infty)+2\,c_1 =
\max\{2\,c_1, h_0(\epsilon,0,c'_1(\infty))+2\,c_5\}+
c'_1(\infty)+2\,c_1\\
\geq\;& 2\sinh c_1 + 12\sinh c_1+
c'_1(\infty)+2\,c_1\geq 14\sinh c'_1(\infty)+3\,c'_1(\infty)\;.
\end{align*}

Define now $\epsilon=R_0^{\rm min}$, so that $2\epsilon\leq
c''_2(\epsilon)$ and $R_0\geq \epsilon$. For every $n\neq 0$, let
$R_n$ be the radius of the ball $B_n$. If for some $n\neq 0$ we have
$2R_n\leq c''_2(\epsilon)$, then $\ph_{B_n}(\gamma) \leq
c''_2(\epsilon)$ and the last assertion of Corollary
\ref{coro:ballline} holds for this $n$. Hence up to removing balls, we
may assume that $R_n\geq c''_2(\epsilon)/2\geq\epsilon$ for every
$n\neq 0$, so that the balls in $(B_n)_{n\in\NN}$ are
$\epsilon$-convex.

The end of the proof is now exactly as the proof of Corollary
\ref{coro:horoline}, with the following modifications: 
$\xi$ is any point in $\partial_\infty X$; $\epsilon=R_0^{\rm min}$;
$C_n=B_n$ for every $n$ in $\NN$; we apply Theorem \ref{theo:line}
instead of Theorem \ref{theo:line1}, which is possible by the range
assumption on $h$; we take now $c_1=c'_1(\epsilon)$, so that we still
have $d(x',y)\leq c_1$ by Proposition \ref{lemma:c1}; $\xi_n$ is now
the center of $B_n$, and $\beta_{\xi_n}(u,v) =d(u,\xi_n)- d(v,\xi_n)$
(see Section \ref{subsec:notations}). Besides that, the proof is
unchanged.  \cqfd

\medskip %
A heavy computation shows that 
$$
c'_1(R_0^{\rm min})\approx 1.7627,
c''_1(R_0^{\rm min},0,0) \approx 101.4169\;\;{\rm and}\;\; 
c''_2(R_0^{\rm min})\approx 106.7051\;.
$$ 
Note that the above corollary is nonempty only if $R_0\geq
c''_1(R_0^{\rm min},0,0)/2+c'_1(R_0^{\rm min})\approx 52.4712$.  The
constants in the following corollary are not optimal.  Theorem
\ref{theo:balinjrad} in the introduction follows from it.

\bcoro %
Let $M$ be a complete Riemannian manifold with sectional curvature at
most $-1$ and dimension at least $3$, let $(x_i)_{i\in I}$ be a finite
or countable family of points in $M$ with $r_i=\inj_M x_i$, such that
$d(x_i,x_j)\geq r_i+r_j$ if $i\neq j$ and such that $r_{i_0}\ge 56$
for some $i_0\in I$.  Then, for every $d\in [2,r_{i_0}-54]$, there
exists a locally geodesic line $\ga$ passing at distance exactly $d$
from $x_{i_0}$ at time $0$, remaining at distance greater than $d$
from $x_{i_0}$ at any nonzero time, and at distance at least $r_i-56$
from $x_{i}$ for every $i\neq i_0$. In particular,
$$
\min_{t\in\RR}d(\gamma(t),x_{i_0})=d.
$$ 
\ecoro

\dem %
Let $\pi:\wt M\ra M$ be a universal covering of $M$, with covering
group $\Ga$, and fix a lift $\wt x_i$ of $x_i$ for every $i\in I$. Let
$B_i$ be the ball $B_{\wt M}(\wt x_i,r_i)$.  Apply Corollary
\ref{coro:ballline} to the family of balls $(g\,B_i)_{g\in\Ga\,,\,i\in
  I}$ in $X=\wt M$, which have pairwise disjoint interiors by the
definition of $r_i$ and the assumption on $d(x_i,x_j)$. Note that
$r_{i_0}\geq 56\geq R_0^{\rm min}$ (see the definition of $R_0^{\rm
min}$). Let $h=2(r_{i_0}-d)$, which belongs to $[108,2r_{i_0}-4]$,
which is contained in $[c''_1(R_0^{\rm min},0,0),
2r_{i_0}-2\,c'_1(R_0^{\rm min})]$ by the previous computations. Then
Corollary \ref{coro:ballline} implies that there exists a geodesic
line $\wt \ga$ in $\wt M$ such that $\ph_{B_{i_0}}(\gamma)=h$ and
$\ph_{gB_i}(\gamma) \le c''_2(R_0^{\rm min})< 108$ for all $(g,i)\neq
(1,i_0)$. Parametrize $\wt \ga$ such that its closest point to $\wt
x_{i_0}$ is at time $t=0$. Let $\ga=\pi\circ\wt\ga$, then the result
follows by the definition of $\ph_C$ (see Subsection
\ref{subsec:penetrationmaps}).
\cqfd

%
%
%
\subsection{Spiralling aroung totally geodesic subspaces}
\label{subsec:totgeod}

In this subsection, we apply Proposition \ref{prop:main} and
Corollary \ref{coro:main} when $C_0$ is a tubular neighborhood of a
totally geodesic submanifold. We only give a few of the various
possible applications, others can be obtained by varying the objects
$(C_n)_{n\in\NN-\{0\}}$, as well as the various subcases in
Proposition \ref{prop:prescribe} (3) and (4).

\btheo\label{theo:totgeoddimbig} %
Let $\epsilon>0$, $\delta\geq0$. Let $X$ be a complete simply
connected Riemannian manifold with sectional curvature at most $-1$
and dimension at least $3$. Let $L$ be a complete totally geodesic
submanifold of $X$ with dimension at least $2$, different from $X$,
and $C_0=\N_\epsilon L$. Let $(C_n)_{n\in\NN-\{0\}}$ be a family of
$\epsilon$-convex subsets in $X$ such that $\diam(C_n\cap
C_m)\leq\delta$ for all $n\neq m$ in $\NN$. Let either $f_0=\ftp_L$ or
$f_0=\ell_{\N_\epsilon L}$, with $X$ having constant curvature in this
second case. Let
$$
h'_0=h_0\big(\epsilon_0,\delta_0,\max\{\norm{f_0-\ell_{\N_\epsilon L}},
\norm{f_0-\ftp_L}+2\epsilon-8\,c'_1(\epsilon)\}\big)
$$ 
and $ h\geq h'_1=
h'_1(\epsilon,\delta,h'_0)$. 
\begin{itemize}
\item[$\bullet$] %
  For every $\xi\in(X\cup \partial_\infty X)-(C_0\cup\partial_\infty
  C_0)$, there exists a geodesic ray or line $\gamma$ starting from
  $\xi$ and entering $\N_\epsilon L$ at time $0$ with $f_0(\gamma)
  =h$, and with $\ell_{C_n}(\gamma)\le h'_1$ for every $n\neq 0$ such
  that $\gamma(]\delta,+\infty[)$ meets $C_n$.
\item[$\bullet$] %
  There exists a geodesic line $\gamma$ in $X$ with
  $f_0(\gamma) =h$, and with
  $$
  \ell_{C_n}(\gamma)\le h'_1 + c'_3(\epsilon) \big(\delta+c'_1(\epsilon)\big)+
  c'_1(\epsilon)$$ 
  for all $n\neq 0$.
\end{itemize}
\etheo

Note that if $\ell_{C_n}(\gamma)\le c$, then $f(\ga)\leq c+\kappa$
for any $\kappa$-penetration map $f$ in $C_n$.

\medskip \dem %
We apply Proposition \ref{prop:main} and Corollary \ref{coro:main}
with $\epsilon_0=\epsilon$, $\delta_0=\delta$, 
$$
\kappa_0 =
\max\{\norm{f_0-\ell_{\N_\epsilon L}}, \norm{f_0-\ftp_L}
+2\epsilon-8\, c'_1(\epsilon)\}\;,
$$ 
so that $h'_0=
h_0(\epsilon_0,\delta_0,\kappa_0)$,  and $f_0$ is a continuous
$\kappa_0$-penetration map in $C_0$. As $L$ is a complete totally
geodesic submanifold of dimension and codimension at least $1$, there
does exist a geodesic line $\ga_0$ in $X$ such that $f_0(\gamma_0)
=h$.  Let $h_0^{\rm min}=\delta_0$ and $h^{\rm min}=
4\,c'_1(\epsilon)+2\epsilon+\delta+ \norm{f_0-\ftp_L}$. By the
definitions of $h'_1(\cdot,\cdot,\cdot),
h_0(\cdot,\cdot,\cdot),c_5(\cdot,\cdot)$ in Subsection
\ref{sec:induction}, we have
$$
h'_0=h_0(\epsilon_0,\delta_0,\kappa_0)> \delta_0=h_0^{\rm min}\;,
$$
and
\begin{align*}
h& \geq   h'_1=
h_0(\epsilon_0,\delta_0,\kappa_0)+2\,c_5(\epsilon_0,\delta_0)
\geq \kappa_0+12\sinh(c'_1(\epsilon_0)+\delta_0)\\ &\geq 
\kappa_0+12\, c'_1(\epsilon)+\delta\geq h^{\rm min}\;.
\end{align*}
The family $(C_n)_{n\in \NN}$ hence satisfies the Local prescription
property $(iv)$ by Proposition \ref{prop:prescribe} (3). Therefore,
the result follows from Proposition \ref{prop:main} and Corollary
\ref{coro:main} \cqfd

\brema {\rm  
%
If the $C_n$'s are disjoint from $\N_\epsilon L$ (and $\delta=0$),
then the same result as Theorem \ref{theo:totgeoddimbig} also holds
when $L$ has dimension $1$, by replacing Proposition
\ref{prop:prescribe} (3) by Proposition \ref{prop:prescribe} (4) in
the above proof and $h_0^{\rm min}=\delta_0$ by $h_0^{\rm min}=0$.

%
}
\erema

\btheo\label{theo:totgeoddimun} %
Let $\epsilon>0$, $\delta\geq0$. Let $X$ be a complete simply
connected Riemannian manifold with sectional curvature at most $-1$
and dimension at least $3$. Let $(L_n)_{n\in\NN}$ be a family of
geodesic lines in $X$, such that $\diam(\N_\epsilon L_n\cap\N_\epsilon
L_m)\leq\delta$ for all $n\neq m$ in $\NN$. Let either $f_0=\ftp_{L_0}$, or
$f_0=\ell_{\N_\epsilon L_0}$ if $X$ has constant curvature, or
$f_0=\cp_{L_0}$ if the metric spheres for the Hamenst\"adt distances (on
$\partial_\infty X-\{\xi\}$ for any $\xi\in\partial_\infty X$) are
topological spheres. Let
$$
h'_0=\max\{5\,c'_1(\epsilon)+5\epsilon+\delta\,,\, h_0(\epsilon,
\delta, \max\{\norm{f_0-\ell_{\N_\epsilon L_0}},
\norm{f_0-\ftp_{L_0}}+2\epsilon-8\,c'_1(\epsilon)\})\}
$$ 
and $h\geq h'_1=
h'_1(\epsilon, \delta, h'_0)$. 
\begin{itemize}
\item[$\bullet$]
For every $\xi\in(X\cup \partial_\infty X)-(\N_\epsilon L_0 \cup
\partial_\infty L_0)$ (and $\xi\in\partial_\infty X-\partial_\infty L_0$
if $f_0=\cp_{L_0}$), there exists a geodesic ray or line $\gamma$ starting
from $\xi$ and entering $\N_\epsilon L_0$ at time $0$ with
$f_0(\gamma) = h$, such that $\ell_{\N_\epsilon L_n} (\gamma)\le h'_1$
for every $n\neq 0$ such that $\gamma(]\delta,+\infty[)$ meets
$\N_\epsilon L_n$.
\item[$\bullet$]
There exists a geodesic line $\gamma$ in $X$ such that $f_0 (\gamma)
=h$, and, if $n\neq 0$, then $\ell_{\N_\epsilon L_n}(\gamma)\le h'_1 +
c'_3(\epsilon) \big(\delta+c'_1(\epsilon)\big)+c'_1(\epsilon)$.
\end{itemize}
\etheo

Note that if $\ell_{\N_\epsilon L_n}(\gamma)\le c$, then $f(\ga)\leq c+\kappa$
for any $\kappa$-penetration map $f$ in $\N_\epsilon L_n$.

\medskip \dem %
As in the previous proof, we apply Proposition \ref{prop:main} and
Corollary \ref{coro:main} with $C_n=\N_\epsilon L_n$, $\epsilon_0=
\epsilon$, $\delta_0= \delta$, 
$$
\kappa_0 =\max\{\norm{f_0-\ell_{\N_\epsilon L_0}}, 
\norm{f_0-\ftp_{L_0}} +2\epsilon-8\, c'_1(\epsilon)\}\;.
$$
For every $n\neq 0$, let $h_0^{\rm
  min}=3\,c'_1(\epsilon) +3\epsilon+\delta+ \norm{\ell_{\N_\epsilon
  L_n} - \ftp_{L_n}}$ and $h^{\rm min}=
  4\,c'_1(\epsilon)+2\epsilon+\delta+ \norm{f_0-\ftp_L}$. In
  particular, $h'_0=\max\{5\,c'_1(\epsilon)+5\epsilon+\delta\,,\,
  h_0(\epsilon, \delta,\kappa_0)\}\geq h_0^{\rm min}$ by Lemma
  \ref{fellowtravelpenemap}. As in the end of the previous proof, the
  family $(C_n)_{n\in\NN}$ hence satisfies the property $(iv)$ by
  Proposition \ref{prop:prescribe} (4), and the result follows.  \cqfd

\medskip %
Let $M$ be a complete Riemannian manifold with sectional curvature at
most $-1$ and dimension $n\geq 3$. Fix a universal cover $\wt M\ra M$
of $M$. For $\epsilon>0$, $\delta\geq 0$, a (possibly not connected,
but any two components having equal dimension) immersed complete
totally geodesic submanifold $L$ (of dimension at least $1$ and at
most $n-1$) will be called {\it $(\epsilon,\delta)$-separated} if the
diameter of the intersection of the $\epsilon$-neighbourhoods of two
lifts to $\wt M$ of two components of $L$ is at most $\delta$.

\medskip
\noindent {\bf Examples. } (1) If $L$ is compact and embedded, then
there exists $\epsilon>0$ such that $L$ is $(\epsilon,0)$-separated.
For instance, a finite family of disjoint simple closed geodesics is
$(\epsilon,0)$-separated for $\epsilon$ small enough.

(2) If $L$ is compact, and if $L$ is {\it self-transverse} (i.e.~if
the tangent spaces at every double point of $L$ are transverse), then
for every $\epsilon>0$ small enough, $L$ is $(\epsilon,1)$-separated.
In particular, a finite family of closed geodesics (possibly
nonsimple) is $(\epsilon,1)$-separated for $\epsilon$ small enough.

(3) The lift of a  locally geodesic line $\ga:\RR\to M$ to the unit
tangent bundle $T^1M$ is the map $\tilde\ga:\RR\to T^1M$ (or by abuse
its image) given by $\tilde\ga(t)=(\ga(t),\ga'(t))$ for every
$t\in\RR$.  For every $\rho>0$, if the $\rho$-neighbourhood (for the
standard Riemannian metric of $T^1M$) of the lift of $\ga$ to $T^1M$
is a tubular neighbourhood, then there exists $\delta(\rho)\geq 0$
such that $\ga$ is $(\rho,\delta(\rho))$-separated. Indeed, if the
intersection of the $\rho$-neighbourhoods of two different lifts to a
universal cover of $\ga$ has diameter big enough (depending only on
$\rho$), then by arguments similar to the ones in the proof of Lemma
\ref{lem:midpointclose}, two subsegments of the two lifts will follow
themselves closely for some time, hence the tangent vectors at two
points on these two lifts will be closer than $\rho$.

\medskip %
Let $L$ be an $(\epsilon,\delta)$-separated immersed complete totally
geodesic submanifold. Let $({\wt L}_\alpha)_{\alpha\in\A}$ be the
family of (connected) complete totally geodesic submanifolds of ${\wt
M}$, that are the lifts to ${\wt M}$ of the components of $L$. Note
that in particular, the family $(\N_\epsilon ({\wt L}_\alpha))
_{\alpha\in\A}$ is locally finite.

Let $f$ be one of the symbols $\ell, \bp,\ftp,\cp$ and assume that $L$
has dimension $1$ if $f= \cp$. Let $\kappa_f$ be respectively $0$,
$2\,c'_1(\epsilon)$, $2\,c'_1(\epsilon)+2\epsilon$,
$2\,c'_1(\epsilon)+2\,c'_1(\infty)+2\epsilon$. For every locally
geodesic line $\ga$ in $M$, consider a lift $\wt\ga$ of $\ga$ to $\wt
M$.  For every $\alpha\in\A$, let $f_\alpha=\ell_{\N_\epsilon
L_\alpha}, \bp_{L_\alpha},\ftp_{L_\alpha},\cp_{L_\alpha}$
respectively, which is a $\kappa_f$-penetration map in $\N_\epsilon
L_\alpha$ by Subsection \ref{subsec:penetrationmaps}.

The family $(f_\alpha(\wt\ga))_{\alpha\in\A}$ will be called the
family of {\it spiraling times} of $\ga$ along $L$ with respect to $f$
(and {\it length spiraling times}, {\it fellow-traveling times} or
{\it crossratio spiraling times} if $f=\ell,\ftp,\cp$ respectively).
Up to permutation of $\A$, it does not depend on the choice of the
lift $\wt\ga$ of $\ga$.  The entering times of $\tilde\ga$ in the sets
$\N_\epsilon L_\alpha$ with $f_\alpha(\tilde\ga)>\delta+\kappa_f$,
where $\alpha$ varies in $\A$, form a discrete subset (with
multiplicity one) of $\RR$, as $\N_\epsilon L_\alpha\cap \N_\epsilon
L_\beta$ has diameter at most $\delta$ if $\alpha\neq \beta$. We will
only be interested in the corresponding spiraling times.  It is also
then possible to order these spiraling times using the order given by
the parametrisation on $\wt\ga$, but we will not need this here. When
$L$ is embedded, and $\epsilon$ is small enough so that the
$\epsilon$-neighborhood of $L$ is a tubular neighborhood, then the
(big enough) fellow-traveling times are the ones defined in the
introduction, see the picture below.


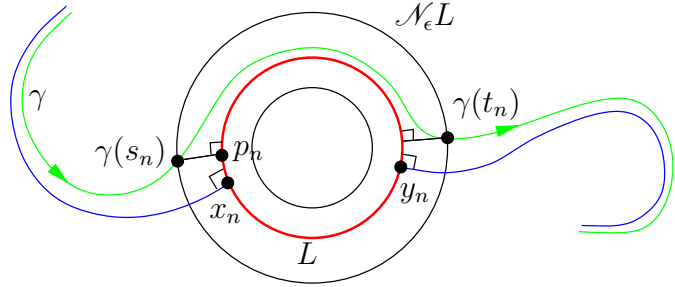
\begin{center}
\input{fig_fellotravtim.pstex_t}
\end{center}


\bcoro\label{coro:separated} %
Let $M$ be a complete Riemannian manifold with sectional curvature at
most $-1$ and dimension $n\geq 3$. Let $\epsilon>0$, $\delta\geq 0$.
Let $L$ be an $(\epsilon,\delta)$-separated immersed complete totally
geodesic submanifold (of dimension at least $1$ and at most $n-1$).
Let $f$ be one of the symbols $\ell,\ftp,\cp$, and
$\kappa'_f=\max\{0,4\epsilon-6c'_1(\epsilon)\}$, $\;2c'_1(\epsilon) +
2\epsilon$, $\;2c'_1(\epsilon) + 2\epsilon+2c'_1(\infty)$ respectively. If
$f=\ell$, assume that $M$ has constant curvature. If $f=\cp$, assume
that $L$ has dimension $1$ and that the metric spheres for the
Hamenst\"adt distances (on the punctured boundary of a universal cover
of $M$) are topological spheres.

For every 
$$
h\ge h'_1=h'_1\big(\epsilon, \delta,
\max\{5\,c'_1(\epsilon)+5\epsilon+\delta\,,\, h_0(\epsilon, \delta,
\kappa'_f)\}\big),
$$ 
there exists a locally geodesic line $\ga$ in $M$ having one
spiraling time with respect to $f$ exactly $h$, and all others
being at most $h'_1+ c'_3(\epsilon)
\big(\delta+c'_1(\epsilon)\big)+c'_1(\epsilon)$.

If furthermore $M$ is nonelementary and geometrically finite, then for
every cusp $e$ of $M$, we may also assume that the locally
geodesic line $\ga$ does not enter too much into the maximal Margulis
neighborhood of $e$, i.e.~ $\gamma$ satisfies
$$
\maxht_e(\ga)\leq \sup_{x\in L} \height_e(x)+\epsilon +
\frac{1}{2}\big(h'_1+ c'_3(\epsilon)(\delta+c'_1(\epsilon)) +
c'_1(\epsilon)\big)\;.
$$
\ecoro

\dem %
Let $\pi:\wt M\ra M$ be a universal cover of $M$, with covering group
$\Ga$.
With $\kappa_0$ the constant in the proofs of the theorems
\ref{theo:totgeoddimbig} and \ref{theo:totgeoddimun}, it is easy to
check, using Section \ref{subsec:penetrationmaps}, that $\kappa'_f\geq
\kappa_0$ for every case of $f$.  The first assertion follows from
Theorem \ref{theo:totgeoddimun} applied to the family $(L_n)_n$ of the
lifts of the components of $L$ to $\wt M$, if the dimension of $L$ is
$1$, and from Theorem \ref{theo:totgeoddimbig} otherwise.

To prove the last assertion, with the notations of Section
\ref{sec:ballsandhoroballs}, let $t_e=\sup_{x\in L} \height_e(x)
+\epsilon$. We add to the family of convex subsets in Theorem
\ref{theo:totgeoddimbig} if $\dim L\geq 2$, and in the proof of
Theorem \ref{theo:totgeoddimun} otherwise, the family of horoballs
$\ga H_e[t_e]$ for $\ga$ in $\Ga$ (modulo the stabilizer
$\Ga_{\xi_e}$). Note that these horoballs have pairwise disjoint
interiors, and that their interiors are disjoint from the
$\epsilon$-neighborhood of every lift of a component of $L$.  \cqfd

\medskip
Theorem \ref{theo:introftp} in the introduction follows
from this one, by the above example (1).

\medskip\noindent{\bf Remark.} %
(1) If we wanted to have the same locally geodesic line $\ga$ for every
cusp $e$ of $M$ in the second assertion of Corollary \ref{coro:separated},
we should add the bigger family of horoballs $(\ga 
H_e[t_e])_{\ga\in\Ga/\Ga_{\xi_e}\,,\:e\in\C}$, and replace there $t_e$
by $\max_{e\in \C}\{t_e,t_\C\}$, where $\C$ is the set of cusps of
$M$, and $t_\C$ is the lower bound of $t\geq 0$ such that two distinct
elements in $(\ga H_e[t])_{\ga\in\Ga/\Ga_{\xi_e}\,,\:e\in\C}$ have
disjoint interiors (see the definition above Corollary
\ref{coro:ballline}), in order for the new horoballs to have disjoint
interiors.

\medskip
(2) With $M$ and $L$ as above, let $f$ be one of the symbols
$\ell,\bp,\ftp,\cp$. Define, for every locally geodesic line $\ga$ in $M$,
$$
\operatorname{maxspt}_{L,f}(\ga)=\sup_{\alpha\in\A}f_\alpha(\tilde\ga)\;,
$$
the least upper bound of spiraling times of $\ga$ around $L$ with respect
to $f$.
Let 
$$
\operatorname{MaxSp}_{L,f}(M)=\{\operatorname{maxspt}_{L,f}(\ga): \ga\in T^1M\}
$$
be the {\it maximum spiraling
  spectrum} 
$\operatorname{MaxSp}_{L,f}(M)$ around $L$ with respect to $f$.
Theorem
\ref{coro:separated} gives, in particular, sufficient conditions for the
maximum spiraling spectrum to contain a ray $[c,+\infty]$.

%
%
%

\subsection{Recurrent geodesics and related results}
\label{subsec:recur}

In this section, when $M$ is geometrically finite, we construct
locally geodesic lines that have a prescribed height in a cusp
neighbourhood of $M$, and furthermore satisfy some recurrence
properties. We will use the notation introduced in Section
\ref{sec:ballsandhoroballs} concerning the cusps $e$, and the objects
$\height_e, V_e, H_e,\xi_e$.

\bcoro 
\label{coro:totgeodbound} %
Let $M$ be a complete, nonelementary, geometrically finite Riemannian
manifold with compact totally geodesic boundary, with sectional
curvature at most $-1$ and dimension at least $3$. Let $e$ be a cusp
of $M$. Then there exists a constant $c''_{3}=c''_{3}(e,M)$ such that
for every $h'\ge c''_{3}$, there exists a locally geodesic line
$\ga$ in $M$ with $\maxht_e(\ga)=h'$, such that the spiraling
times of $\ga$ along the boundary $\partial M$ are at most $c''_{3}$.
\ecoro

Up to changing the constant $c''_{3}$, we may also assume that $\ga$
stays away from some fixed (small enough) cusp neighbourhood of every
cusp different from $e$. Note that, up to changing the constant
$c''_{3}$, the last assertion of the corollary does not depend on the
choice of $f=\ell,\bp,\ftp,\cp$, with respect to which the spiraling
times are computed, and we will use $f=\ell$.

\medskip \dem As $\partial M$ is compact, there exists
$\epsilon'\in\;]0,1[$ such that the $\epsilon'$-neighbourhood of the
geodesic boundary $\partial M$ is a tubular neighbourhood of $\partial
M$. By definition of manifolds with totally geodesic boundary, there
exists a complete simply connected Riemannian manifold $\wt M$, a
nonelementary, torsion-free, geometrically finite discrete subgroup
$\Gamma$ of isometries of $\wt M$, a $\Gamma$-equivariant collection
$(L^+_k)_{k\in\NN}$ of pairwise disjoint open halfspaces with totally
geodesic boundary $(L_k)_{k\in\NN}$, such that $M$ is isometric with
$\Gamma\backslash (\wt M-\bigcup_{k\in\NN}L^+_k)$. We will identify
$M$ and $\Gamma\backslash (\wt M-\bigcup_{k\in\NN}L^+_k)$ by such an
isometry from now on. Note that $(\N_{\epsilon'}L^+_k)_{k\in\NN}$ is a
family of pairwise disjoint $\epsilon'$-convex subsets in $\wt M$.

Let $t_{e,\partial M}=\max_{x\in\partial M}\height_e(x)$, which exists
since $\partial M$ is compact. Note that the family $(g
H_e[t_{e,\partial M}+1]) _{g\in \Ga/\Ga_{\xi_e}}$ is a
$\Gamma$-equivariant family of pairwise disjoint horoballs in $\wt M$,
which are disjoint from $\N_{\epsilon'} L^+_n$ for all $n\in\NN$. Let
us relabel this family of horoballs as $(H_k)_{k\in\NN}$ such that
$H_0= H_e[t_{e,\partial M}+1]$. Note that the horoballs $H_k$,
$k\in\mathbb N$, are $\epsilon'$-convex.

Define 
$$  
c''_3= 
\max\big\{h'_1\big(\epsilon',0,h_0(\epsilon',0,c'_1(\infty))\big), 
c'_1(\epsilon')\big\}+ t_{e,\partial M}+ 1+ 
c'_1(\epsilon')(c'_3(\epsilon')+1) \;.
$$ 
and let $h'\geq c''_3$. We apply Corollary \ref{coro:main} with $X=\wt
M$; $\epsilon_0 = \epsilon'$; $\delta_0=0$; $\kappa_0=c'_1(\infty)$;
$C_0=H_0$; $f_0= \ph_{C_0}$; $h'_0=h_0(\epsilon_0,\delta_0,\kappa_0)$;
$C_{2k+1}=\N_{\epsilon'}L^+_k$; $C_{2k}=H_k$; $h=2h'- 2(t_{e,\partial
  M}+ 1)$.  Note that $f_0$ is a $\kappa_0$-penetration map in $C_0$
by Lemma \ref{lem:proppenehoro} and $h\ge
h'_1(\epsilon_0,\delta_0,h'_0)$, as $h'\geq c''_3$.  As
$\widetilde{M}$ is a manifold of dimension at least $2$, there does
exist a geodesic line $\ga_0$ in $X$ with $f_0(\ga_0)=h$. The family
$(C_n)_{n\in\NN}$, whose elements have pairwise disjoint interiors,
satisfies the assertion $(iii)$. It also satisfies $(iv)$, by the same
proof as the one of Theorem \ref{theo:line1} as $h\geq
2c'_1(\epsilon')$. Hence, by Corollary \ref{coro:main}, there exists a
geodesic line $\wt\gamma$ in $X$ with $\ph_{H_0}(\wt\gamma)=h$ and
$$
\ell_{C_n}(\wt\gamma)\le
h''_1=h'_1(\epsilon_0,\delta_0,h'_0)+
c'_1(\epsilon_0)(c'_3(\epsilon_0)+1)
$$ 
for all $n\neq 0$.

As $\ell_{C_{2n+1}}(\wt\gamma)$ is finite, the geodesic $\wt \ga$
doesn't cross the boundary of $L_n^+$, hence stays in $\wt M-
\bigcup_{k\in\NN}L^+_k$.  Let $\pi:\wt M- \bigcup_{k\in\NN}L^+_k\ra M$
be the canonical projection, and $\ga=\pi\circ\wt\ga$. 
Hence, the length traveling times of $\ga$ are at most $c''_3$.

Note that 
$$
\ph_{H_e}(\wt\ga)=\ph_{C_0}(\wt\ga)+2(t_{e,\partial M}+1)=
h+2(t_{e,\partial M}+1)= 2h',
$$ 
by the paragraph above Lemma \ref{lem:proppenehoro}.  Furthermore,
if $g\in(\Ga-\Ga_{\xi_e})/\Ga_{\xi_e}$, then there exists $k$ in
$\NN-\{0\}$ such that
\begin{align*}
\ph_{gH_e}(\wt\ga)&=\ph_{C_{2k}}(\wt\ga)+2(t_{e,\partial M}+1)\leq
\ell_{C_{2k}}(\wt\ga)+c'_1(\infty)+2(t_{e,\partial M}+1)\\ &
\leq h''_1+c'_1(\infty)+2(t_{e,\partial M}+1)\leq 2c''_3\leq 2h'\;.
\end{align*}
Therefore
$\maxht_e(\ga)=h'$ by the same proof as in the end of the proof of
Corollary \ref{coro:maxsp}.  
\cqfd
\bigskip

Let $M$ be a compact, connected, orientable, irreducible, acylindrical,
atoroidal, boundary incompressible $3$-manifold with nonempty boundary
(see for instance \cite{MT} for references on $3$-manifolds and
Kleinian groups).  A {\it hyperbolic structure} on a manifold is a
complete Riemannian metric with constant sectional curvature $-1$. A
cusp $e$ of a hyperbolic structure is {\it maximal} if the maximal
Margulis neighborhood of $e$ is a neighborhood of an end of the
manifold. Let $P$ be the union of the torus components of $\partial
M$, and $\gf(M)=\gf(M,P)$ be the (nonempty) space of complete
geometrically finite hyperbolic structures in the interior of $M$
whose cusps are maximal, up to isometries isotopic to the identity.
Recall that $\gf(M)$ is homeomorphic to the Teichm\"uller space of
$\partial_0M=\partial M-P$.

For every $\sigma$ in $\gf(M)$, the cusps of $\sigma$ are in
one-to-one correspondence with the torus components of $\partial M$,
as any minimizing geodesic ray representing a cusp converges to an end
of the interior of $M$ corresponding to a torus component of $\partial
M$. If $e$ is a torus component of $\partial M$, let
$\maxht_{\sigma,e} (\gamma)$ denotes the maximum height of a locally
geodesic line $\gamma$ in $\sigma$ with respect to the cusp
corresponding to $e$. The {\it convex core} of a structure $\sigma$ in
$\gf(M)$ is the smallest closed convex subset of the interior of $M$,
whose injection in the interior of $M$ induces an isomorphism on the
fundamental groups.

The following result generalizes Theorem \ref{theo:intromaxrec} in the
introduction to the case of several cusps.

\bcoro \label{coro:max+recurrent} %
Let $M$ be a compact, connected, orientable, irreducible,
acylindrical, ator\-oi\-dal, boun\-dary incompressible $3$-manifold with
boundary, and let $e$ be a torus component of $\partial M$.  For every
compact subset $K$ in $\gf(M)$, there exists a constant
$c''_{4}=c''_{4}(K)$ such that for every $h\ge c''_{4}$ and every
$\sigma\in K$, there exists a locally geodesic line $\ga$ contained in
the convex core of $\sigma$ such that $\maxht_{\sigma,e}(\ga)=h$, and
$\maxht_{\sigma,e'}(\ga)\leq c''_{4}$ for every torus component
$e'\neq e$ of $\partial M$.  \ecoro

\dem For a subset $A$ of $\partial_\infty\HH^3_\RR$, we denote by
$\conv A$ the hyperbolic convex hull of $A$ in $\HH^3_\RR$.  A
subgroup $\Gamma$ of $\pi_1M$ is called a {\it boundary subgroup} if
there are an element $\gamma\in\pi_1M$, a component $C$ of
$\partial_0M$, and a point $x\in C$ such that $\Gamma=\gamma\,
\Im\big(\pi_1(C,x)\to \pi_1(M,x)\big)\,\gamma^{-1}$.  Let
$(\Gamma_n)_{n\in\NN}$ be the collection of boundary subgroups of
$\pi_1M$. Let $(\Gamma'_n)_{n\in\NN}$ be the collection of maximal
(rank $2$) abelian subgroups of $\pi_1M$, with $\Gamma'_0$ conjugated
to $\pi_1 e$.

Let $\rho_\sigma:\pi_1M\to\Isom(\HH^3_\RR)$ be a holonomy
representation corresponding to $\sigma\in K$.  By assumption,
$\Ga=\rho_\sigma (\pi_1M)$ is a (particular) web group (see for
instance \cite{AM}). More precisely, for all $n\in\NN$,
$\rho_\sigma(\Gamma_n)$ is a quasifuchsian subgroup of $\Ga$
stabilizing a connected, simply connected component
$\Omega_{n,\sigma}$ of the domain of discontinuity of
$\rho_\sigma\pi_1M$, such that $\Omega_{n,\sigma}$ and
$\Omega_{m,\sigma}$ have disjoint closures if $n\ne m$, and that
$\partial\Omega_{n,\sigma}$ contains no parabolic fixed points of
$\Ga$. Let $(H_{k,\sigma}) _{k\in\NN}$ be a maximal family of
horoballs with pairwise disjoint interiors such that $H_{k,\sigma}$ is
$\rho_\sigma(\Gamma'_k)$-invariant (such a family is unique if $M$ has
only one torus component). To make it canonical over $\gf(M)$, we may
fix an ordering $e_1=e,e_2,\dots, e_m$ of the torus components of
$\partial M $, and take by induction $H_{k,\sigma}$, for the $k$'s in
$\NN$ such that $\Ga'_k$ is conjugated to $\pi_1(e_i)$, to be
equivariant and maximal with respect to having pairwise disjoint
interiors as well as having their interior disjoint with the interior of
$H_{k_*,\sigma}$, for the $k_*$'s in $\NN$ such that $\Ga'_{k_*}$ is
conjugated to $\pi_1(e_j)$ with $j<i$.  Note that the $H_{k,\sigma}$'s,
besides the ones such that $\Ga'_k$ is conjugated to $\pi_1e$, are not
the maximal horospheres that allow to define the height functions, but
this changes their values only by a constant (uniform on $K$).

Hence, as $K$ is compact, there exists $\delta>0$ such that for every
$\sigma\in K$, the $1$-convex subsets $\N_1(\conv
\;\Omega_{n,\sigma})$ and $H_{k,\sigma}$ for $n,k\in\NN$ meet pairwise
with diameter at most $\delta$.

The claim follows as in Corollary \ref{coro:totgeodbound} by applying
Corollary \ref{coro:main} to $X=\HH^3_\RR$, $\epsilon_0=1$, $\delta_0
=\delta$, $\kappa_0=c'_1(\infty)$, $C_0=H_{0,\sigma}$, $f_0=
\ph_{C_0}$, $h'_0=h_0(\epsilon_0,\delta_0,\kappa_0)$, $C_{2n+1}=
\N_1(\conv\;\Omega_{n,\sigma})$, $C_{2n}=H_{n,\sigma}$ to get a
geodesic line $\wt\ga$ in $X$ with prescribed penetration in $C_0$,
and penetration bounded by a constant in $C_n$ for $n\neq 0$. The
finiteness of the intersection lengths $\ell_{C_{2n+1}}(\wt \ga)$ for
$n\in\NN$ implies that $\wt \ga$ stays in the convex hull of the limit
set of $\Ga$.  \cqfd \bigskip

\rem The fact that a locally geodesic line stays in the convex core of
the manifold and does not converge (either way) to a cusp is
equivalent with the locally geodesic line being two-sided recurrent.

%
%
%

\subsection{Prescribing the asymptotic penetration behavior}
\label{subsec:limsupspec}

Let $X$ be a proper geodesic ${\rm CAT}(-1)$ metric space and let
$\xi\in X\cup \partial _\infty X$. Let
$\epsilon\in\RR^*_+\cup\{+\infty\}$, $\delta,\kappa\geq 0$.  Let
$(C_\alpha)_{\alpha\in\A}$ be a family of $\epsilon$-convex subsets of
$X$ which satisfies the Almost disjointness condition $(iii)$ with
parameter $\delta$. For each $\alpha\in\A$, let $f_\alpha$ be a
$\kappa$-penetration map. Let $\gamma$ be a geodesic ray or line, with
$0$ in the domain of definition of $\ga$ (as we are only interested in
the asympotic behavior, the choice of time $0$ is unimportant). These
assumptions guarantee that the set $\E_\ga$ of times $t\geq 0$ such
that $\ga$ enters in some $C_\alpha$ at time $t$ with
$f_\alpha(\ga)>\delta+\kappa$ is discrete in $[0,+\infty[$, and that
$\alpha=\alpha_t$ is then unique.  The set $\E_\ga$ is finite if
$f_\beta(\ga)=+\infty$ for some $\beta$.  Hence
$\E_\ga=(t_i)_{i\in\N}$ for some initial segment $\N$ in $\NN$, with
$t_i<t_{i+1}$ for $i,i+1$ in $\N$.  With
$a_i(\ga)=f_{\alpha_{t_i}}(\ga)$, the (finite or infinite) sequence
$\big(a_i(\ga)\big)_{n\in\N}$ will be called the (nonnegative) {\em
penetration sequence} of $\gamma$ with respect to
$(C_\alpha,f_\alpha)_ {\alpha\in\A}$. In this section, we study the
asymptotic behavior of these penetration sequences. We will only state
some results when the $C_\alpha$'s are balls or horoballs, but similar
ones are valid, for instance for $\epsilon$-neighborhoods of geodesic
lines in $X$ (see for instance \cite{HPP}.  We may also prescribe the
asymptotic penetration in one cusp, while keeping the heights in the
other cusps (uniformly) bounded.

In the following results, we show how to prescribe the asymptotic
behaviour of the penetration sequence of a geodesic ray or line with
respect to horoballs and their penetration height functions.
First, we prove a general result, and we give the more explicit result
for Riemannian manifolds as  Corollary \ref{coro:horolimsup}.

\btheo \label{theo:limsuppenseq} %
Let $X$ be a proper geodesic ${\rm CAT}(-1)$ metric space, with
$\partial_\infty X$ infinite.  Let $(H_\alpha)_{\alpha\in\A}$ be a
family of horoballs with pairwise disjoint interiors.  Assume that
there exists $K\in[0,+\infty[$ and a dense subset $Y$ in
$\partial_\infty X$ such that, for every geodesic ray $\ga$ in $X$
with $\ga(+\infty)\in Y$, we have $\liminf_{t\ra+\infty}d(\ga(t),
\bigcup_{\alpha\in\A} H_\alpha) \le K$.  Let $\xi\in
X\cup\partial_\infty X$ and $c,c'\geq 0$. Assume that for every $h\geq
c$ and $\alpha\in\A$ such that $\xi\notin H_\alpha\cup
H_\alpha[\infty]$, there exists a geodesic ray or line $\ga$ starting
from $\xi$ and entering $H_\alpha$ at time $t=0$ with
$\ph_{H_\alpha}(\ga) = h$, and with $\ph_{H_\beta}(\ga)\leq c'$ for
every $\beta$ in $\A-\{\alpha\}$ such that $\ga(]0,+\infty[)$ meets
$H_\beta$. Let $\big(a_i(\ga')\big)_{n\in\N}$ be the penetration
sequence of a geodesic ray or line $\ga'$ with respect to
$(H_\alpha,\ph_{H_\alpha})_ {\alpha\in\A}$.

Then, for every 
$$
h\geq h_*=\max\big\{c\,,\;c'+3\,c'_1(\infty)+
10^{-5}\big\},
$$
there exists a geodesic ray or line $\ga$ starting from $\xi$ such
that
$$
\limsup_{i\to+\infty} \;a_i(\ga)=h\;.
$$
\etheo

\medskip \dem %
To simplify notation, let $f_\alpha=\ph_{H_\alpha}$, $ c_*=
c'+3\,c'_1(\infty)+ 10^{-5} $, so that $h_*=
\max\{c_*,c\}$. If a geodesic ray or line $\ga$ starting from $\xi$
meets $H_\alpha$, let $t^-_\alpha(\ga)$ and $t^+_\alpha(\ga)$ be the
entrance and exit times.

Let $h\geq h_*$, and let $\alpha_0\in\A$ such that $\xi\notin
H_{\alpha_0}\cup H_{\alpha_0}[\infty]$, which exists by the
assumptions.  As $h\geq h_*\geq c$, there exists a geodesic ray or
line $\gamma_0$ starting from $\xi$, entering $H_{\alpha_0}$ at time
$0$, such that $f_{\alpha_0}(\ga_0)=h$, and $f_\alpha(\ga_0)\le c'$
for every $\alpha\neq \alpha_0$ such that $\ga_0(]0,+\infty[)$ meets
$H_\alpha$.

We construct, by induction, sequences $(\ga_k)_{k\in\NN}$ of geodesic
rays or lines starting from $\xi$, $(\alpha_k)_{k\in\NN}$ of elements
of $\A$, and $(t_k)_{k\in\NN-\{0\}}$ of elements in $[0,+\infty[$
converging to $+\infty$, such that for every $k\in\NN$,
\begin{enumerate}
\item %
  $\ga_k$ enters the interior of $H_{\alpha_0}$ at time $0$, with
  $d(\ga_k(0),\ga_{k-1}(0))\leq \frac1{2^k}$ if $k\geq 1$;
\item %
  $\ga_k$ enters $H_{\alpha_k}$ and $f_{\alpha_k}(\ga_k)=h$;
\item %
  if $0\leq
  j\leq k-1$,  
  then $\ga_k(]0,+\infty[)$ enters the interior of
  $H_{\alpha_j}$ before entering $H_{\alpha_k}$ with
  $t^-_{\alpha_j}(\ga_k)< t_k= t^+_{\alpha_{k-1}}(\ga_k)$;
\item %
   if $k\geq 1$, then for every $\alpha$ such that $\ga_k(]0,+\infty[)$
  meets $H_\alpha$, we have $\big| f_{\alpha}(\ga_k)- f_{\alpha}
  (\ga_{k-1}) \big| < \frac1{2^k}$ if $t^-_\alpha(\ga_k)< t_k$, and
  $f_{\alpha} (\ga_k)\leq c_*$ if $ t_k\leq t^-_\alpha (\ga_k)<
  t^-_{\alpha_k}(\ga_k)$, and $f_{\alpha} (\ga_k) \leq c'$ if $
  t^-_\alpha(\ga_k)\geq t^+_{\alpha_k}(\ga_k)$.
\end{enumerate}

Let us first prove that the existence of such sequences implies
Theorem \ref{theo:limsuppenseq}.  By the assertion (1), as $\ga_k(0)$
stays at bounded distance from $\ga_0(0)$, up to extracting a
subsequence, the sequence $(\ga_k)_{k\in\NN}$ converges to a geodesic
ray or line $\ga_\infty$ starting from $\xi$, entering in
$H_{\alpha_0}$ at time $t=0$, by continuity of the entering point in
an $\epsilon$-convex subset. Let us prove that $\limsup_{i\to+\infty}
\;a_i(\ga_\infty)=h$.

The lower bound $\limsup_{i\to+\infty} \;a_i(\ga_\infty)\geq h$ is
immediate by a semicontinuity argument. Indeed, for every $k>i$ in
$\NN$, we have by the assertions (2), (3) and (4),
$$
\big| f_{\alpha_i}(\ga_k)-h\big|=
\big |f_{\alpha_i}(\ga_k)-f_{\alpha_i}(\ga_i)\big|
\leq \sum_{j=i}^{k-1}
\big|f_{\alpha_i}(\ga_{j+1})-f_{\alpha_i}(\ga_j)\big| 
\leq \sum_{j=i}^{k-1}\frac1{2^{j+1}}\leq \frac1{2^i}\;.
$$
Hence by continuity of $f_{\alpha_i}$, we have the inequality $f_{\alpha_i}
(\ga_\infty) \geq h-\frac1{2^i}$, whose right side converges to $h$ as $i$ tends
to $+\infty$, which proves the lower bound, as $h>c'_1(\infty)$ and
$f_{\alpha_i}$ is a $c'_1(\infty)$-penetration map in $H_{\alpha_i}$ (see
Section \ref{subsec:penetrationmaps}).

To prove the upper bound, assume by absurd that there exists $\epsilon
>0$ such that for every $\lambda>0$, there exists
$\alpha=\alpha(\lambda)\in\A$ such that $\ga_\infty$ enters $H_\alpha$
with $f_{\alpha}(\ga_\infty) \geq h+ \epsilon$ and
$t^-_\alpha(\ga_\infty) > \lambda+2\,c'_1(\infty)$. Take 
$\lambda_0=\max\big\{t_{i+1}:\frac 1{2^i}\ge\frac\epsilon 2\big\}$, and 
$\alpha= \alpha(\lambda_0)$

By continuity of $f_{\alpha}$, if $k$ is big enough, we have
$f_{\alpha}(\ga_k)\geq h+ \frac{\epsilon}{2}\geq h_*\geq c_*\geq
c'$. Thus,
$\ga_k$ meets $H_\alpha$ as $h_*>0$. The entry time is positive, as 
$d(\ga_k(0),\ga_\infty(0))\le c'_1(\infty)$ and the entrance
points of $\ga_k$ and $\ga_\infty$ in $H_\alpha$ are at distance at
most $c'_1(\infty)$, both by Lemma \ref{lemma:c1}, and as the entrance time of
$\ga_\infty$ in $H_\alpha$ is strictly bigger than $2c'_1(\infty)$.  Hence, by
the assertion (4), we have $t^-_\alpha(\ga_k)< t_k$. Let $i\leq k-1$
be the minimum element of $\NN$ such that for $j=i,\dots, k-1$, the
geodesic $\ga_{j+1}$ meets $H_\alpha$ at a positive time with
$t^-_\alpha(\ga_{j+1})< t_{j+1}$. By the triangular inequality, we
have
$$
\big |t^-_\alpha(\ga_{i+1})-t^-_\alpha(\ga_\infty) \big|\leq
d\big(\ga_{i+1}(t^-_\alpha(\ga_{i+1})),
\ga_\infty(t^-_\alpha(\ga_\infty))\big)+
d\big(\ga_{i+1}(0),\ga_\infty(0)\big)\leq 2\,c'_1(\infty)\;.
$$
Hence 
$$
t_{i+1}> t^-_\alpha(\ga_{i+1})\geq t^-_\alpha(\ga_\infty)
-2\,c'_1(\infty)> \lambda_0+2\,c'_1(\infty)-2\,c'_1(\infty)=\lambda_0.
$$  
By the definition of $\lambda_0$, we hence have $\frac{1}{2^{i}}<
\frac{\epsilon}{2}$.  By the definition of $i$ and by the assertion
(4), we have
\begin{align*}
f_{\alpha}(\ga_i)=\; & f_{\alpha}(\ga_k)+\sum_{j=i}^{k-1}
\big(f_{\alpha}(\ga_{j})-f_{\alpha}(\ga_{j+1})\big)
\geq h+ \frac{\epsilon}{2} -
\sum_{j=i}^{k-1} \frac{1}{2^{j+1}}\\ \geq \;&
h+ \frac{\epsilon}{2} - \frac{1}{2^{i}}\geq h\geq h_*\;,
\end{align*}
and in particular by the same argument as for $\ga_k$ above, $\ga_i$
enters $H_\alpha$ at a positive time and $t^-_\alpha(\ga_i)< t_i$.
This contradicts the minimality of $i$. This completes the proof,
assuming the existence of a sequence with properties (1)--(4).

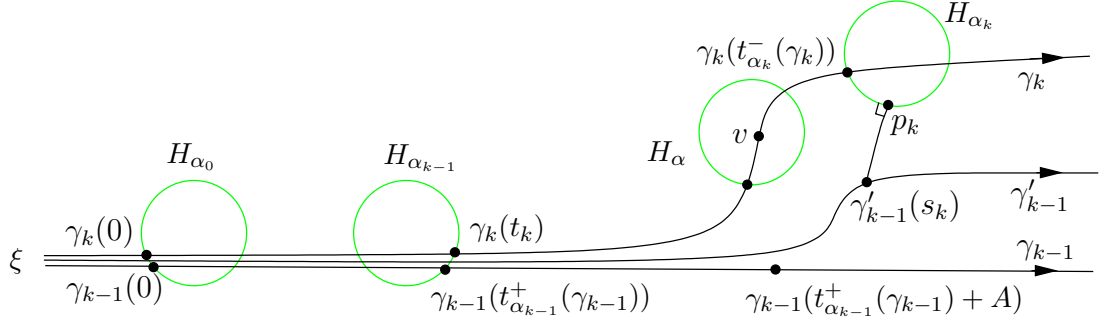
\begin{center} 
\input{fig_limsup_horo.pstex_t}
\end{center}

Let us now construct the sequences $(\ga_k)_{k\in\NN}$,
$(\alpha_k)_{k\in\NN}$, $(t_k)_{k\in\NN-\{0\}}$. We have defined
$\ga_0$, $\alpha_0$, and they satisfy the properties (1)--(4). Let
$k\geq 1$, and assume that $\gamma_{k-1}$, $\alpha_{k-1}$, as well as
$t_{k-1}$ if $k\geq 2$, have been constructed.

As $Y$ is dense in $\partial_\infty X$, for every $A>0$, there exists
a geodesic ray or line $\ga'_{k-1}$ starting from $\xi$ with
$\ga'_{k-1}(+\infty)\in Y$, entering in $H_{\alpha_0}$ at time $t=0$,
which is very close to $\ga_{k-1}$ on $[0,t^+_{\alpha_{k-1}}(\gamma_{k-1})+A]$.
By the definition of $K$, let $s_k$ be the first time $t\geq
t^+_{\alpha_{k-1}}(\ga_{k-1})+A$ such that there exists $\alpha$ in $\A$ with
$d(\ga'_{k-1}(t), H_\alpha)\leq K+1$, and let $\alpha_k$ be such an
$\alpha$ with $d(\ga'_{k-1}(s_k), H_\alpha)$ minimum.  Let $p_k$ be
the closest point of $H_{\alpha_k}$ to $\ga'_{k-1}(s_k)$.  
Note that $\xi\notin H_{\alpha_k}\cup H_{\alpha_k}[\infty]$, if $A$ is
big enough (in particular compared to $K$), as $H_{\alpha_0}$ and
$H_{\alpha_k}$ have disjoint interiors.

By the hypothesis, let $\ga_k$ be a geodesic ray or line starting from
$\xi$ with $f_{\alpha_k}(\ga_k)=h$ (which proves the assertion (2) at
rank $k$ as $h>0$) and $f_{\alpha}(\ga_k)\le c'$ for every $\alpha$
such that $\ga_k ([t_{\alpha_k}^+(\ga_k),+\infty[)$ enters $H_\alpha$.  As a
${\rm CAT}(-1)$ metric space is $\log(1+\sqrt{2})$-hyperbolic, the
geodesic $]\xi,p_k]$ is contained in the $\log(1+\sqrt{2})$-
neighbourhood of the union $]\xi,\ga'_{k-1}(s_k)]\cup
[\ga'_{k-1}(s_k),p_k]$.  By Lemma \ref{lemma:c1}, we have
$d(\ga_k(t^-_{\alpha_k}(\ga_k)),]\xi,p_k])\leq c'_1(\infty)$, and therefore
$]\xi,\ga_k(t^-_{\alpha_k}(\ga_k))]$ is contained in the
$(c'_1(\infty)+\log(1+\sqrt{2}))$-neighbourhood of $]\xi,\ga'_{k-1}(s_k)]\cup
[\ga'_{k-1}(s_k),p_k]$.  Up to choosing $A$ big enough, we may hence
assume that $\ga_k$ is very close to $\ga_{k-1}$ between the times $0$
and $t^+_{\alpha_{k-1}}(\ga_{k-1})+1$.  Using this and properties (1)
and (3) at rank $k-1$, we have
\begin{enumerate}
\item[$\bullet$] $\ga_k$ does enter the interior of $H_{\alpha_0}$, at
  a time that we may assume to be $0$, with
  $d(\ga_k(0),\ga_{k-1}(0))\leq \frac1{2^k}$ (this proves the
  assertion (1) at rank $k$);
\item[$\bullet$]  for $0\leq j\leq k-1$, as $\ga_{k-1}$ passes in the interior of
  $H_{\alpha_{j}}$ at a time strictly between $0$ and
  $t^+_{\alpha_{k-1}}(\ga_{k-1})$, by the
  inductive assertions (3) if $k\neq 1$ and $j\leq k-2$, or (1) if
  $k=1$ or (2) if $j=k-1$, so does the geodesic ray or line $\ga_k$; this
  allows, in particular, to define $t_k= t^+_{\alpha_{k-1}}(\ga_k)$, and
  proves the assertion (3) at rank $k$;
\item[$\bullet$] for every $\alpha$ such that $\ga_k(]0,+\infty[)$
  meets $H_\alpha$ and $t^-_\alpha(\ga_k)< t_k$, we may assume by
  continuity that $\big| f_{\alpha}(\ga_k)-f_{\alpha}(\ga_{k-1})\big|
  < \frac1{2^k}$.
\end{enumerate}
Hence, to prove the assertion (4) at rank $k$, we consider
$\alpha\in\A$ such that $\ga_k$ meets $H_\alpha$ with $t_k \leq
t^-_\alpha(\ga_k) < t^-_{\alpha_k}(\ga_k)$, and we prove that
$f_{\alpha} (\ga_k)\leq c_*$. We may assume that $f_{\alpha}
(\ga_k)>0$. Let $v$ be the highest point of $\ga_k$ in $H_\alpha$,
which, by disjointness, belongs to $]\ga_k(t^-_\alpha(\ga_k)),
\ga_k(t^-_{\alpha_k}(\ga_k))[$. Let $u$ be a point in
$]\xi,\ga'_{k-1}(s_k)]\cup [\ga'_{k-1}(s_k),p_k]$ at distance at most
$c'_1(\infty)+\log(1+\sqrt{2})$ from $v$. Assume first that
$u\in[\ga'_{k-1}(s_k),p_k]$. Note that by the minimality assumption on
$\alpha_k$, the point $u$ then does not belong to $H_\alpha$.  As
$c_*\geq 2\,c'_1(\infty)+2\log(1+\sqrt{2})=3\,c'_1(\infty)$, this
implies that $f_{\alpha} (\ga_k)\leq 2\,d(u,v)\leq c_*$.  Assume now
that $u=\ga'_{k-1}(t)$ with $t\in[t^+_{\alpha_{k-1}}(\gamma_{k-1})+A,s_k[$.
Then by the minimality of $s_k$, the point $u$ again does not belong
to $H_\alpha$ (it is in fact at distance at least $K+1$ from
$H_\alpha$).  Hence similarly $f_{\alpha} (\ga_k)\leq c_*$.  Finally,
assume that $u=\ga'_{k-1}(t)$ with
$t\in[0,t^+_{\alpha_{k-1}}(\gamma_{k-1})+A[$.  Let $u'$ be a point on
$\ga_{k-1}([t^+_{\alpha_{k-1}}(\gamma_{k-1}),t^+_{\alpha_{k-1}}(\gamma_{k-1})+A])$ with
$d(u,u')\leq 10^{-5}/2$ (as $\gamma_{k-1}$ and $\gamma'_{k-1}$ were
assumed to be very close on that range). The point $u'$ is at distance
at most $c'/2$ from a point in the complement of $H_\alpha$, as if it
belongs to the interior of $H_\alpha$, then
$f_\alpha(\gamma_{k-1})\leq c'$ by the inductive hypothesis (4) on
$\gamma_{k-1}$. Hence
$$
f_\alpha(\ga_k)\leq 2\,d(v,u')+2(c'/2)\leq 
2c'_1(\infty)+10^{-5}=c_*\;.
$$
This proves the result.  
\cqfd 

\medskip\noindent{\bf Remark. } %
The proof when $(H_\alpha)_{\alpha\in\A}$ is a family of balls of
radius $R>0$, replacing $c'_1(\infty)$ by $c'_1(R)$, and assuming both
in the hypothesis and in the conclusion that $h\leq c''$ for some $c''$, is
the same.

\bcoro\label{coro:horolimsup} %
Let $X$ be a complete simply connected Riemannian manifold with
sectional curvature at most $-1$ and dimension at least $3$, and let
$\big(H_\alpha\big)_{\alpha\in\A}$ be a family of horoballs in $X$
with disjoint interiors. Assume that there exists $K\in[0,+\infty[$
and a dense subset $Y$ in $\partial_\infty X$ such that, for every
geodesic ray $\ga$ in $X$ with $\ga(+\infty)\in Y$, we have
$\liminf_{t\ra+\infty}d(\ga(t), \bigcup_{\alpha\in\A} H_\alpha) \le
K$. Then, for every $\xi\in X\cup\partial_\infty X$ and
$$
h\geq c''_1(\infty,0,0)+
4\,c'_1(\infty) + 10^{-5} \approx 13.5542,
$$ 
there exists a geodesic ray or line $\gamma$ starting from $\xi$ such
that, with $(a_i(\ga))_{n\in\N}$ the penetration sequence of $\ga$
with respect to $(H_\alpha,\ph_{H_\alpha})_ {\alpha\in\A}$, we have
$$
\limsup_{i\to+\infty} \;a_i(\ga)=h\;.
$$
\ecoro

\dem %
Let  $c= c''_1(\infty,0,0)$, $c'= c''_1(\infty,0,0)
+ c'_1(\infty)$. We apply Theorem \ref{theo:line1} with $\epsilon=
\infty$, $\delta=0$, $\kappa=c'_1(\infty)$, $\xi_0=\xi$,
$C_0=H_\alpha$ where $\alpha\in\A$ satisfies $\xi\notin
H_\alpha\cup\H_\alpha[\infty]$, $f_0=\ph_{C_0}$, $(C_n)_{n\geq 1}$ is
$(H_\beta)_{\beta\in\A-\{\alpha\}}$ (up to indexation). Then the
assumptions of Theorem
\ref{theo:limsuppenseq} are satisfied. An easy computation of $h_*$
in Theorem \ref{theo:limsuppenseq} then yields the result.
\cqfd

\medskip\noindent{\bf Remark. } %
Using Theorem \ref{theo:line} instead of Theorem \ref{theo:line1}, the
same statement when $(H_\alpha)_{\alpha\in\A}$ is a family of balls of
radius $R>0$, for $h\in[c''_1(R,0,0)+4\,c'_1(R) + 10^{-5},2R-c'_1(R)]$
holds true.

\medskip %
As in Section \ref{sec:ballsandhoroballs}, we consider a complete,
nonelementary, geometrically finite Riemannian manifold $M$, and $e$
an 
end of $M$.  The {\it asymptotic height spectrum} of the pair $(M,e)$
is
$$
\limsupsp(M,e)=\big\{\limsup_{t\to\infty}
\height_e(\gamma(t)):\gamma\in T^1M\big\}\;.
$$

In classical Diophantine approximation, the {\it Lagrange spectrum} is
the subset of $[0,+\infty[$ consisting of the {\it approximation constants}
$c(x)$ of an irrational real number $x$ by rational numbers $p/q$,
defined by
$$
c(x)=-2\log\mu=\liminf_{q\ra\infty}|q|^2\big|x-\frac{p}{q}\big|.
$$
Using the well known connection between the Diophantine approximation
of real numbers by rational numbers and the action of the modular
group $\PSLZ$ on the upper halfplane model of the real hyperbolic
plane, the asymptotic height spectrum of the modular orbifold
$\PSLZ\backslash\hdr$ is the image of the Lagrange spectrum by the map
$t\mapsto -2\log t$ (see for instance \cite[Theo.~3.4]{HP3}).  Hall
\cite{Hall1, Hall2} showed that the Lagrange spectrum contains an
interval $[0,c]$ for some $c>0$. The maximal such interval $[0,\mu]$
(which is closed as the Lagrange spectrum is closed, by Cusick's
result, see for instance \cite{CF}), called {\it Hall's ray}, was
determined by Freiman \cite{Freiman} (see also
\cite{Slo} where the map $t\mapsto 1/t$ has to be applied). The
geometric interpretation of Freiman's result in our context is that
$\limsupsp(\PSLZ\backslash\hdr)$ contains the maximal interval
$[c,+\infty]$ with
$$
c=-2\log\left(\frac{491993569}{2221564096+283748\sqrt{462}}\right)\approx 
3.0205\;.
$$

The following result is the asymptotic analog of Corollary
\ref{coro:maxsp}, and has a completely similar proof.  Theorem
\ref{theo:introhallquoi} in the introduction follows, since
$(c''_1(\infty,0,0)+ 4\,c'_1(\infty) + 10^{-5})/2\approx 6.7771$. The
result proves the existence of Hall's ray in our geometric context,
which is much more general (there is no assumption of arithmetic
nature, nor of constant curvature nature), and with a universal
constant (though we do not know the optimal one) $6.7771$ which is not
too far from the geometric Freiman constant $3.0205$ from the above
particular case.

\bcoro\label{coro:limsupsp} %
Let $M$ be a complete, nonelementary, geometrically finite Riemannian
manifold with sectional curvature at most $-1$ and dimension at
least $3$, and let $e$ be a cusp of $M$. Then $\limsupsp(M,e)$
contains the interval 
$$
[(c''_1(\infty,0,0)+ 4\,c'_1(\infty) + 10^{-5})/2,+\infty]\;.
$$  
\cqfd 
\ecoro

In the next section \ref{sec:diophapp}, we consider a number of
arithmetically defined examples, illustrating this last result. But we
need first to recall some properties and do some computations in the
real and complex hyperbolic spaces.

\section{Applications to Diophantine approximation in
  negatively curved manifolds}
\label{sec:diophapp}

%
%
%

\subsection{On complex hyperbolic geometry and the Heisenberg group}
\label{subsec:comphypgeom}

To facilitate computations, we identify elements in
$\CC^{n-1}$ with their coordinate column matrices. We will denote by
$A^*=\mbox{}^t\overline{A}$ the adjoint matrix of a complex matrix
$A$. 
In particular, the standard hermitian scalar product of 
$w,w'\in\CC^{n-1}$ is $
\overline{w^*\;w'}= \sum_{i=1}^{n-1}w_i\overline{w'_i}$. We also use
the notation $|w|^2=w^*\;w$.

\medskip 
Let
$\HH^n_\CC$ be the Siegel domain model of complex hyperbolic
$n$-space whose
underlying set is
$$
\HH^n_\CC=\{(w_0,w)\in
\CC\times\CC^{n-1}\;:\; 2\,{\rm Re}\;w_0 -|w|^2>0\}\;,
$$
and whose Riemannian metric is 
$$
ds^2_\CC=\frac{4}{(2\,{\rm Re}\;w_0 -|w|^2)^2}
\big((dw_0-dw^*\;w)((\overline{dw_0}-w^*\;dw)+
(2\,{\rm Re}\;w_0 -|w|^2)\;dw^*\;dw\big)\;,
$$
(see  for instance  \cite[Sect.~4.1]{Gol}). 
Complex hyperbolic space has constant holomorphic sectional curvature
$-1$, hence
its real sectional curvatures are bounded between $-1$ and $-\frac14$.
Its boundary at infinity is 
$$
\partial_\infty\HH^n_\CC=\{(w_0,w)\in
\CC\times\CC^{n-1}\;:\; 2\,{\rm Re}\;w_0 -|w|^2=0\}
\cup\{\infty\}\;.
$$ 
The horoballs centered at $\infty$ in $\HH^n_\CC$ are the subsets
$$
\H_s=\{(w_0,w)\in
\CC\times\CC^{n-1}\;:\; 2\,{\rm Re}\;w_0 -|w|^2\geq s\}\;,
$$ 
for $s>0$. Note that the subset $\HH^1_\CC= \{(w_0,w)\in\HH^n_\CC\;:\;
w=0\}$ is the right halfplane model of the real hyperbolic plane with
constant curvature $-1$, and it is totally geodesic in $\HH^n_\CC$.
In particular, the (unit speed) geodesic line starting from $\infty$,
ending at $(0,0)\in \partial_\infty\HH^n_\CC$ and meeting the
horosphere $\partial \H_2$ at time $t=0$ is the map $c_0:\RR\ra
\HH^n_\CC$ defined by $c_0:t\mapsto (e^{-t},0)$.

Let $q$ be the nondegenerate Hermitian form $-z_0\overline{z_n}
-z_n\overline{z_0} + |z|^2$ of signature $(1,n)$ on
$\CC\times\CC^{n-1}\times\CC$ with coordinates $(z_0,z,z_n)$. This is
not the form considered in \cite[page 67]{Gol}, hence we need to do
some computations with it, but it is better suited for our purposes.
The Siegel domain $\HH^n_\CC$ embeds in the complex projective
$n$-space $\PP_n(\CC)$ by the map (using homogeneous coordinates)
$$
(w_0,w)\mapsto [w_0:w:1]\;.
$$
Its image is the negative cone of $q$, that is
$\{[z_0:z:z_n]\in\PP_n(\CC)\;:\;q(z_0,z,z_n)<0\}$. This embedding
extends continuously to the boundary at infinity, by mapping
$(w_0,w)\in
\partial_\infty\HH^n_\CC-\{\infty\}$ to $[w_0:w:1]$ and $\infty$ to
$[1:0:0]$, so that the image of $\partial_\infty\HH^n_\CC$ is the null
cone of $q$, that is $\{[z_0:z:z_n]\in\PP_n(\CC)\;:\;
q(z_0,z,z_n)=0\}$.  We use matrices by blocks in the decomposition
$\CC\times\CC^{n-1}\times\CC$.

Let 
$$
Q=\left(\begin{array}{ccc} 0 & 0 & -1\\
    0 & I & 0 \\ -1 & 0 & 0\end{array}\right)
$$ 
be the matrix of $q$. If
$$
X=\left(\begin{array}{ccc} a & \gamma^* & b\\
    \alpha & A & \beta \\ c & \delta^* & d\end{array}\right),
$$ 
then
$$
Q^{-1}X^*Q=\left(\begin{array}{ccc} \overline{d} & -\beta^* &
    \overline{b}\\ -\delta & A^* & -\gamma \\ \overline{c} & -\alpha^*
    & \overline{a}\end{array}\right).
$$  
If ${\rm U}_{Q}$ is the group of invertible matrices with
complex coefficients preserving the hermitian form $q$, then $X$
belongs to ${\rm U}_{Q}$ if and only if $X$ is invertible with inverse
$Q^{-1}X^*Q$.  In particular, if $X$ belongs to ${\rm U}_{Q}$, then
\begin{equation}\label{eq:equationsUq}
  \left\{\begin{array}{c}c\overline{d}-\delta^*\delta+d\overline{c}   
 = 0 \\ 
      a\overline{b}-\ga^*\ga+b\overline{a}   = 0 \\ 
      -\alpha\beta^*+A A^* - \beta\alpha^*  =I \\
      c\overline{b}-\delta^*\ga+d\overline{a}  =1\\
      \overline{d}\alpha-A\delta+\overline{c}\beta =0\\
      \overline{b}\alpha-A\ga+\overline{a}\beta  =0 \;.
\end{array}\right.
\end{equation}

The group ${\rm U}_{Q}$ acts projectively on $\PP_n(\CC)$, preserving
the negative cone of $q$, hence it acts on $\HH^n_\CC$. We will denote in
the same way the action of $U_Q$ on $\HH^n_\CC$ and the action of
$U_Q$ on the image of $\HH^n_\CC$ in $\PP_n(\CC)$. It is well known
(see for instance \cite{Gol}) that ${\rm U}_{Q}$ preserves the
Riemannian metric of $\HH^n_\CC$.

The Heisenberg group ${\rm Heis}_{2n-1}$ is 
the real Lie group with underlying space
$\CC^{n-1}\times\RR$ and law
$$
(\zeta,v)(\zeta',v')=(\zeta+\zeta',v+v'-2\;{\rm Im}\;\zeta^*\zeta')\;.
$$
It has a Lie group embedding in ${\rm U}_{Q}$, defined by 
$$
(\zeta,v)\mapsto u_{\zeta,v}=\left(\begin{array}{ccc} 1 & \zeta^* & 
\frac{|\zeta|^2}{2}-i\frac{v}{2}\\   0 & I & \zeta \\ 0 & 0 & 1
\end{array}\right),
$$
whose image preserves the point $\infty$ as well as each horoball
centered at $\infty$, as an easy computation shows.

The {\it Cygan distance} (see \cite[page 160]{Gol}) on ${\rm
  Heis}_{2n-1}$ is the unique left-invariant distance $d_{\rm Cyg}$
such that 
$$
d_{\rm Cyg}((0,0),(\zeta,v))=(|\zeta|^4+v^2)^{1/4}\;.
$$
We introduce the {\it modified Cygan distance} $d'_{\rm Cyg}$ as the
unique left-invariant distance $d'_{\rm Cyg}$ such that
$$
d'_{\rm Cyg}((0,0),(\zeta,v))=
((|\zeta|^4+v^2)^{1/2} +|\zeta|^2)^{1/2}\;.
$$
It is straightforward to check that $d'_{\rm Cyg}$ is indeed a
distance, in the same way as the Cygan distance, see for instance \cite{KR},
and that it is equivalent
to the Cygan distance, 
$$
d_{\rm Cyg}\leq d'_{\rm Cyg}\leq \sqrt{2}\;d_{\rm Cyg}\;.
$$
Hence, its induced length distance is
equivalent to the Carnot-Carathéo\-do\-ry distance on the Heisenberg
group ${\rm Heis}_{2n-1}$ (see \cite[page 161]{Gol}).

As the action of ${\rm Heis}_{2n-1}$ on $\partial_\infty\HH^n_\CC-
\{\infty\}$ is simply transitive, $d_{\rm Cyg}$ and $d'_{\rm Cyg}$
define distances on $\partial_\infty\HH^n_\CC-\{\infty\}$, which are
invariant under the action of ${\rm Heis}_{2n-1}$. We also call these
distances the {\it Cygan distance} and the {\it modified Cygan
distance}, and again denote them by $d_{\rm Cyg}$ and $d'_{\rm
Cyg}$. Explicitly, these distances are given by
$$
d_{\rm Cyg}(u_{\zeta,v}(0,0),u_{\zeta',v'}(0,0))=
d_{\rm Cyg}((\zeta,v),(\zeta',v'))\;,
$$
and the similar expression for the modified Cygan distance.

\blemm %
The distance $d_{\rm Cyg}$ (resp.~$d'_{\rm Cyg}$) is the unique
distance on 
$\partial_\infty\HH^n_\CC-\{\infty\}$ invariant under ${\rm
  Heis}_{2n-1}$ such that $d'_{\rm Cyg}((w_0,w),(0,0))=
\sqrt{2|w_0|}$ (resp.~
$d'_{\rm Cyg}((w_0,w),(0,0))=
\sqrt{2|w_0|+|w|^2}$).
\elemm

\dem For every $(w_0,w)$ in $\partial_\infty\HH^n_\CC-\{\infty\}$,
note that $(w_0,w)=u_{\zeta,v}(0,0)$ if and only if $v=-2\;{\rm Im}\;
w_0$ and $\zeta=w$, and that $2\;{\rm Re}\; w_0=|w|^2$. Hence
$$
d'_{\rm Cyg}(u_{\zeta,v}(0,0),(0,0))=
((4\;{\rm Re}^2\; w_0+4\;{\rm Im}^2\;
w_0)^{1/2} +|w|^2)^{1/2}=\sqrt{2|w_0|+|w|^2}\;.
$$
A similar proof gives the result for the Cygan distance.\cqfd
\medskip

In particular, if $n=2$, then $d'_{\rm Cyg}$ is indeed defined as in
the statement of Theorem \ref{theo:applicygan} in the introduction.

\medskip
Let $d_{\HH^n_\CC}$ be the Riemannian distance on $\HH^n_\CC$, and
$d'_{\HH^n_\CC}=\frac{1}{2} \,d_{\HH^n_\CC}$ be the Riemannian
distance of the Riemannian metric of $\HH^n_\CC$ renormalized to have
maximal real sectional curvatures $-1$.

\bprop \label{prop:cuspdistcalccomphyp} %
For every $\xi,\xi'$ in $\partial_\infty\HH^n_\CC-\{\infty\}$, for
every $s_0>0$, the distance $\ell'$ for the renormalized Riemannian
distance $d'_{\HH^n_\CC}$ between the horoball $\H_{s_0}$ and the
horoball centered at $\xi$ and tangent to the geodesic line between
$\infty$ and $\xi'$ is, if these horoballs are disjoint,
$$
\ell'=-\log d'_{\rm Cyg}(\xi,\xi') +\frac12\log(\frac{s_0}{2})\;.
$$
\eprop

\dem By invariance of the modified Cygan distance, of each horoball
centered at $\infty$, and of the normalized Riemannian distance, by
the action of the Heisenberg group, we may assume that $\xi=(0,0)$.
Let $(\zeta,v)\in{\rm Heis}_{2n-1}$ such that $\xi'=u_{\zeta,v}(\xi)$.
As $u_{\zeta,v}$ sends geodesic lines to geodesic lines, and fixes
$\infty$, the geodesic line (for $d_{\HH^n_\CC}$) starting from
$\infty$ and ending at $\xi'$ is $u_{\zeta,v}\circ c_0$, which by an
easy computation is
$$
u_{\zeta,v}\circ c_0:t\mapsto (e^{-t}+(|\zeta|^2-iv)/2,\zeta)\;.
$$
The matrix $X_0=\left(\begin{array}{ccc} 0 & 0 & 1\\
    0 & I & 0 \\ 1 & 0 & 0\end{array}\right)$ belongs to $U_Q$, as
$X_0^{-1}=X_0$ and $Q^{-1}X_0^*Q=X_0$, and the corresponding isometry of
$\HH^n_\CC$ sends $\infty\in\partial_\infty\HH^n_\CC$ to $(0,0)\in
\partial_\infty\HH^n_\CC$. Hence $X_0$ sends the horoballs centered at
$\infty$ to the horoballs centered at $(0,0)$. Let $s>0$, an easy
computation shows that
$$
X_0\H_s=\{(w_0,w)\in \CC\times\CC^{n-1}\;:\; 
2\,{\rm Re}\;w_0 -|w|^2\geq s |w_0|^2\}\;.
$$
For every $t$ in $\RR$, the point $u_{\zeta,v}\circ c_0(t)$ belongs to
the horosphere $X_0\partial\H_s$ if and only if
$$
2\,{\rm Re}\;(e^{-t}+(|\zeta|^2-iv)/2) -|\zeta|^2 = 
s |e^{-t}+(|\zeta|^2-iv)/2|^2\;,
$$
that is, if and only if 
$$
s\;e^{-2t}+(s|\zeta|^2-2)e^{-t}+\frac{s}{4}\,(|\zeta|^4 +v^2)=0\;.
$$
The horoball $X_0\H_s$ is hence tangent to the geodesic line
$u_{\zeta,v}\circ c_0$ if and only if the above quadratic equation
with unknown $e^{-t}$ has a double solution, that is, if and only if
its discriminant $\Delta$ is $0$. An easy computation gives $
-\Delta=s^2v^2+4s|\zeta|^2-4 $. Thus,  the horoball $X_0\H_s$ is
tangent to $u_{\zeta,v}\circ c_0$ if and only if
\begin{equation}\label{eq:value_s}
s=\frac{2}{\sqrt{|\zeta|^4 +v^2}+|\zeta|^2}\;.
\end{equation}
As the geodesic line $c_0$ passes through the point at infinity of
both horoballs $\H_{s_0}$ and $X_0\H_s$ (which have disjoint interiors
if $s_0$ is big enough), the Riemannian distance between them is the
length of the subsegment of $c_0$ joining them. Note that $c_0$ meets
$X_0\partial\H_s$ at $(\frac{2}{s},0)$. Hence, by an easy computation in
$\HH^1_\CC$,
$$
\ell'=d'_{\HH^n_\CC}(\H_{s_0},X_0\H_s)=
\frac12\;d_{\HH^n_\CC}(\H_{s_0},X_0\H_s)=
\frac12\;d_{\HH^n_\CC}((\frac{s_0}{2},0),(\frac{2}{s},0))=
\frac12(\log\frac{s_0}{2}-\log\frac{2}{s})\;.
$$
By Equation (\ref{eq:value_s}), the result follows.
\cqfd

\medskip %
For every $X$ in ${\rm U}_Q$, we will denote by $c=c(X)$ its
$(3,1)$-coefficient in its matrix by blocks. 
Note that $X$ fixes
$\infty$ if and only if $c=0$, by the  equations
(\ref{eq:equationsUq}).
 Equivalently, by the same set of equations, a
matrix fixes $\infty$ if and only if it is upper triangular (this is the
main reason why we chose the hermitian form $q$ rather than the one in
\cite{Gol}). The following lemma is completely analogous to
Proposition 5.14 of \cite{HP4}, but as we are using a different
quadratic form, we need to give a proof.

\blemm \label{lem:depthcalccomphyp} %
For every $X$ in ${\rm U}_Q$ and every $s>0$ such that the horoballs
$\H_s$ and $X\H_s$ have disjoint interiors, we have
$$
d'(\H_s,X\H_s)=\log |c|+\log\frac{s}{2}\;.
$$
\elemm

\dem As $\H_s$ and $X\H_s$ have disjoint interiors, $X$ does not fix
$\infty$, hence $c\neq 0$.  Left and right multiplication of $X$ by an
element $u_{\zeta,v}$ for some $(\zeta,v)$ in ${\rm Heis}_{2n-1}$ does
not change the coefficient $c$ of $X$, nor does it change
$d'(\H_s,X\H_s)= \frac12\;d(\H_s,X\H_s)$, as $u_{\zeta,v}$ preserves
the distance $d$ and each horosphere centered at $\infty$. 
Hence, as
${\rm Heis}_{2n-1}$ acts transitively on $\partial_\infty\HH^n_\CC-
\{\infty\}$, we may assume that $X\infty=(0,0)$ and that
$X^{-1}\infty=(0,0)$.  
As $X\infty=(0,0)$,
the coefficients $a,\alpha$ of $X$ are $0$, 
and hence by the second equation
of (\ref{eq:equationsUq}), the coefficient $\gamma$ is $0$. 
As $X^{-1}\infty=(0,0)$,
the coefficients $d,\delta$ of $X$ are $0$, and
hence by the fifth equation of (\ref{eq:equationsUq}), the coefficient
$\beta$ is 
$0$. Therefore, by the third and fourth equation of
(\ref{eq:equationsUq}), the matrix $X$ has the form
$\left(\begin{array}{ccc} 0 & 0 & \frac{1}{\overline{c}}\\ 0 & A & 0
    \\ c & 0 & 0 \end{array}\right)$, with $A$ unitary. An easy
computation, similar to the one we already did with $X_0$, shows that
$$
X\H_s=\{(w_0,w)\in \CC\times\CC^{n-1}\;:\; 
2\,{\rm Re}\;w_0 -|w|^2\geq s  |c|^2|w_0|^2\}\;.
$$
Hence, as above,
$$
d'(\H_s,X\H_s)=\frac12\;d(\H_s,X\H_s)=
\frac12\;d((\frac{s}{2},0),(\frac{2}{s|c|^2},0))=
\log |c|+\log\frac{s}{2} \;.\;\;\mbox{\cqfd}
$$

\medskip 
Let $m$ be a squarefree positive integer, let
$K_{-m}=\QQ(i\sqrt{m})$ be the corresponding 
imaginary quadratic number field, and let $\O_{-m}$ be the ring of integers
of $K_{-m}$.
An {\it order} $\O$ in $K_{-m}$ is a unitary subring of 
$\O_{-m}$ which is a free $\ZZ$-module of rank $2$.
We 
use for instance \cite[chap.~7]{Cox} for a general reference on these
objects. 
An example of an order in $K_{-m}$ is $\ZZ[i\sqrt m]$, 
and $\O_{-m}$ is the maximal order of $K_{-m}$. 
In particular, $\O$ contains a $\QQ$-basis of $K_{-m}$, and the
field of fractions of $\O$ is $K_{-m}$.
Let $\omega$ be the element of
$\O$ with ${\rm Im}\,\omega>0$ such that $\O=\ZZ[\omega]=\ZZ+\omega\ZZ$.

As $\O$ is stable by complex conjugation, the subset
$$
{\rm SU}_Q(\O)=
{\rm SU}_Q\cap\M_{n+1}(\O)
$$ 
is a discrete subgroup of the
semi-simple connected real Lie group ${\rm SU}_Q={\rm U}_Q \cap{\rm
  SL}_{n+1}(\CC)$. 

Let $\I$ be a non-zero
ideal of $\O$.
We denote by $\Ga_{\CC,\I}$ the preimage, by the
group morphism 
${\rm SU}_Q(\O)\ra {\rm SL}_{n+1}(\O/\I)$ of reduction
modulo $\I$, of the parabolic subgroup of matrices whose first column
has all its coefficients $0$ except the first one.
As $\O/\I$ is finite  
($\I$ is nonzero),
$\Ga_{\CC,\I}$ is a finite
index subgroup of 
${\rm SU}_Q(\O)$ 
as

Recall that a horoball $H$ centered at a point $\xi$ in a ${\rm
  CAT}(-1)$ metric space $X$ is {\it precisely invariant} under a
group of isometries $\Ga$ if for every $g\in \Ga$ that does not fix
$\xi$, the intersection $g\stackrel{\circ}{H}\cap \stackrel{\circ}{H}$
is empty.

\blemm \label{lem:horoprecinvcomphyp} %
For $v=2\;{\rm Im}\;\omega$ if ${\rm Re}\;\omega\in\ZZ$, and
$v=4\;{\rm Im}\;\omega$ otherwise, the horoball $\H_v$ is precisely
invariant under $\Ga_{\CC,\I}$. Furthermore, if $\I=\O=\O_{-1}$,
then $\H_2$ is the maximal horoball centered at $\infty$ which is
precisely invariant under $\Ga_{\CC,\I}$.  \elemm

\dem %
With $v$ as in the statement, the element $-iv/2$ belongs to $\O$, as
$i\;{\rm Im}\;\omega=\omega-\;{\rm Re}\;\omega$ belongs to $\O$ if
${\rm Re}\;\omega\in\ZZ$, and $2i\;{\rm Im}\;\omega =\omega-
\overline{\omega}$ belongs to $\O$ (which is stable by conjugation).
Hence $u_{0,v}$ belongs to $\Ga_{\CC,\I}$. It follows for instance
from \cite[Prop.~5.7]{HP00} (which is an easy consequence of the
complex hyperbolic Shimizu inequality of Kamiya \cite{Kam} and Parker
\cite{Par}) that the 
horoball $\H_v$ is precisely invariant (the hermitian form $q$ in
\cite{HP00} is not the same one as the above, but it is equivalent by a
permutation of coordinates, hence we may indeed apply
\cite[Prop.~5.7]{HP00}).

If $\I=\O=\O_{-1}$, then $X_0$ defined above belongs to $\Ga_{\CC,\I}$
and $\omega=i$. By Lemma \ref{lem:depthcalccomphyp}, we have
$d(\H_2,X_0\H_2)=0$, hence the last assertion follows.  \cqfd

\medskip For every $(a,\alpha,c)\in\O\times\O^{n-1}\times \O$, let
$\langle a,\alpha,c\rangle$ be the ideal of $\O$ generated by $a$, $c$
and the components of $\alpha$.

\bprop \label{prop:identiforbit} %
If $n=2$ and $\O=\O_{-m}$, then
\begin{enumerate}
\item for every $\I$, the set of parabolic fixed points of $\Ga_{\CC,\I}$
  is exactly 
  the set of points in $\partial_\infty\HH^n_\CC$ having  homogeneous
  coordinates in $\PP_n(\CC)$ that are elements in $K_{-m}$;
\item
the orbit $\Ga_{\CC,\I}\cdot\infty$ 
is exactly the set of points in
$\partial_\infty\HH^n_\CC$ having homogeneous coordinates in $\PP_n(\CC)$ of
the form $[a:\alpha:c]$ with
$(a,\alpha,c)\in\O\times\I^{n-1}\times\I$, $2\;{\rm Re}\;
a\overline{c}=|\alpha|^2$ and $\langle a,\alpha,c\rangle=\O$;
\item
if $m=1,2,3,7,11,19,43,67,163$ and $\I=\O$, then 
$\Ga_{\CC,\I}$ has only one orbit of parabolic fixed points.
\end{enumerate}
\eprop

\dem (1) If $\I=\O$, then the first result is due to Holzapfel
\cite{Hol1}, see \cite[page 280]{Hol2}. As
$\Ga_{\CC,\I}$ has finite index in 
${\rm SU}_Q(\O)$, and as a discrete
group and a finite index subgroup have the same set of parabolic fixed
points, the first claim follows.

(2) 
A result of Feustel \cite{Feu} (see \cite[page
280]{Hol2}, \cite{Zin}) says that the map which associates 
to a parabolic fixed point of 
${\rm SU}_Q(\O)$
the fractional ideal
generated by its homogeneous coordinates in $\O_{-m}$ induces a
bijection from the set of orbits under 
${\rm SU}_Q(\O)$
of
parabolic fixed points of 
${\rm SU}_Q(\O)$  to the set of ideal
classes of $K_{-m}$. As $\infty$ corresponds to $[1:0:0]$ whose
coordinates generate the trivial fractional ideal, the second claim
follows if $\I=\O$, as well as claim (3).

If $M\left(\begin{array}{c}1\\0\\0\end{array}\right)=
\left(\begin{array}{c}a\\\alpha\\c\end{array}\right)$, then
$\left(\begin{array}{c}a\\\alpha\\c\end{array}\right)$ is the first
column of the matrix $M$, so that the second claim if $\I\neq\O$
follows by the definition of $\Ga_{\CC,\I}$.  \cqfd

\subsection{Quaternions and $5$-dimensional real hyperbolic geometry}
\label{subsec:quatrealhypgeom}

Let $\HH$ be Hamilton's quaternion algebra over $\RR$, generated
as a real vector space by the standard basis $1,i,j,k$, with products
$k=ij=-ji$, $i^2=-1$, $j^2=-1$ and unit $1$. Recall that the {\it
  conjugate} of the quaternion $z=x_1+x_2i+x_3j+x_4k$ is
$\overline{z}=x_1-x_2i-x_3j-x_4k$, which satisfies
$\overline{z\,w}=\overline{w}\,\overline{z}$, and that its {\it
  absolute value} (or square root of its reduced norm) is
$$
|z|=\sqrt{N(z)}=\sqrt{z\overline{z}}=
\sqrt{\overline{z}z}=\sqrt{x_1^2+x_2^2+x_3^2+x_4^2}\;.
$$
The {\it Dieudonné determinant} (see \cite{Die} and \cite{Asl}) $\Delta$ is
the group 
morphism from the group ${\rm GL}_2\HH$ of invertible $2\times 2$ matrices
with coefficients in $\HH$ to $\RR^*_+$, given by
$$
\Delta(\left(\begin{array}{cc} a & b \\ c & d\end{array}\right))=
\left\{\begin{array}{cl} |\,ad-aca^{-1}b\,| & {\rm if}\; a\neq 0\\
|\,cb-cac^{-1}d\,| & {\rm if}\; c\neq 0\end{array}\right.\;.
$$
We will denote by ${\rm SL}_2(\HH)$ the group of $2\times 2$ quaternionic
matrices with Dieudonné determinant $1$ (this notation is different,
hence should not be confused with the notation ${\rm SL}(2,C_n)$ for
$n=3$ and $C_3=\HH$ of Vahlen and Ahlfors \cite{Ahl}, see also
\cite{MWW}, giving a description of the isometry group of the real
hyperbolic $(n+1)$-space using the $2^n$-dimensional real Clifford
algebra $C_n$). We refer for instance to \cite{Kel2} for more
information on ${\rm SL}_2\HH$.

The group ${\rm SL}_2(\HH)$ acts on the Alexandrov compactification
$\HH\cup\{\infty\}$ of $\HH$ by
$$
\left(\begin{array}{cc} a & b \\ c & d\end{array}\right)\cdot z =
\left\{\begin{array}{ll} (az+b)(cz+d)^{-1} & 
{\rm if}\; z\neq \infty,-c^{-1}d \\ 
ac^{-1} & {\rm if}\; z=\infty, c\neq 0\\ 
\infty & {\rm otherwise.}\end{array}\right.
$$
It is well known (see for instance \cite{Kel2}) that ${\rm
  PSL}_2(\HH)={\rm SL}_2(\HH)/\{\pm{\rm Id}\}$ is the orientation
preserving conformal group of the $4$-sphere $\HH\cup\{\infty\}$ with
its standard conformal structure defined by the $4$-dimensional
Euclidean space $(\HH,|\cdot|)$. In the upper halfspace model
$\HH^5_\RR$ of the $5$-dimensional real hyperbolic space with constant
curvature $-1$, consider the coordinates $(z,t)$ with $z\in\HH$ and
$t>0$ (called the {\it vertical coordinate}), so that
$\partial_\infty\HH^5_\RR$ identifies with the union of $\HH$ (for
$t=0$) and of $\{\infty\}$.  By the Poincaré extension procedure (see
for instance \cite[Sect.~3.3]{Bea}), the group ${\rm PSL}_2(\HH)$ hence
identifies to the group of orientation preserving isometries of
$\HH^5_\RR$.  We will denote 
the Riemannian
distance on $\HH^5_\RR$ by $d_{\HH^5_\RR}$.

\blemm \cite[Theo.~1.-2)]{Hel} \label{lem:calvertcomp}%
For every $g=\left(\begin{array}{cc} a & b \\ c & d\end{array}\right)$
in ${\rm SL}_2(\HH)$, and $(z,t)$ in $\HH^5_\RR$, the vertical
coordinate of $g(z,t)$ is
$$
\frac{t}{|cz+d|^2 + |c|^2t^2}\;.
$$
\elemm

\dem As \cite{Hel} is an announcement, we give a proof for the sake of
completeness. The proof is  an adaptation of the proof for ${\rm
SL}_2(\CC)$ in \cite[page
58]{Bea}, the main problem consists of being careful with the
noncommutativity of $\HH$. We may assume that $c\neq 0$, as the map
$z\mapsto \alpha z\beta+\ga$ for $\alpha,\beta,\ga$ in
$\HH^*\times\HH^*\times\HH$ is an Euclidean similitude of ratio
$|\alpha\beta|$. Define the {\it isometric sphere} of $g$ to be the
sphere $S_g$ of center $-c^{-1}d$ and radius $\frac1{|c|}$ in the
Euclidean space $(\HH,|\cdot|)$. By the definition of an Euclidean
reflection with respect to a sphere in this Euclidean space, the map
$$
\sigma: z\mapsto -c^{-1}d + 
\frac{1}{|c|^2}\frac{z+c^{-1}d}{|z+c^{-1}d|^2} 
$$
is the Euclidean reflection with respect to the sphere $S_g$.
An easy computation shows that the map $\varphi=g\circ\sigma$ is
$$
z\mapsto (b-ac^{-1}d)(\overline{z}\,\overline{c}  + 
\overline{d}) + ac^{-1}\;,
$$
which is an Euclidean isometry, as $z\mapsto \overline{z}$ is, and
$|cb-cac^{-1}d|=1$. The Poincaré extension of $\varphi$ preserves the
vertical coordinates, and the Poincaré extension of $\sigma$ is the
Euclidean reflection with respect to the sphere in $\HH^5_\RR$ whose
equator is $S_g$. As $g=\varphi\circ\sigma$, the result follows by an
easy computation. 
\cqfd

\medskip %
The horoballs centered at $\infty$ in $\HH^5_\RR$ are the
subsets $\H_s$ for $s>0$, where
$$
\H_s=\{(z,t)\in\HH\times]0,+\infty[\;:\; t\geq s\}\;.
$$

\blemm \label{lem:depthcalc5hyp} %
For every $g=\left(\begin{array}{cc} a & b \\ c & d\end{array}\right)$
in ${\rm SL}_2\HH$, and every $s>0$ such that the horoballs $\H_s$
and $g\H_s$ have disjoint interiors, we have
$$
d(\H_s,X\H_s)=2\log |c|+2\log s\;.
$$
\elemm

\dem %
As $\H_s$ and $g\H_s$ have disjoint interiors, we have $c\neq 0$.  The
map $g$ sends the geodesic line between $-c^{-1}d$ and $\infty$ to the
geodesic line between $\infty$ and $ac^{-1}$, hence the point
$(-c^{-1}d,s)$ of intersection of the first line with $\H_s$ to the
point $g(-c^{-1}d,s)$ of intersection of the second line with $g\H_s$.
The vertical coordinate of $g(-c^{-1}d,s)$ is $\frac{1}{|c|^2s}$ by
the previous lemma \ref{lem:calvertcomp}. Hence the result follows by
an easy computation of hyperbolic distances. \cqfd

\medskip %
We will use \cite{Vig} and \cite[Section 2]{MR} as 
general references on quaternion algebras.
Let $A(\QQ)$ be a quaternion algebra over $\QQ$, which is {\it
  ramified} over $\RR$, that is, the real algebra $A(\QQ)\otimes_\QQ
\RR$ is isomorphic to Hamilton's algebra $\HH$. We identify $A(\QQ)
\otimes_\QQ\RR$ and $\HH$ by any such isomorphism. Let $\O'$ be an
{\it order} of $A(\QQ)$, that is an unitary subring which is a
finitely generated $\ZZ$-module generating the $\QQ$-vector space
$A(\QQ)$. For instance, if 
$$
A(\QQ)=\{x_1+x_2i+x_3j+x_4k\in\HH \;:\;
x_1,x_2,x_3,x_4\in\QQ\},
$$ 
we can take 
$$
\O'=\{x_1+x_2i+x_3j+x_4k \in\HH
\;:\;x_1,x_2,x_3,x_4\in\ZZ\},
$$ 
or the Hurwitz ring
$$
\O'=\big\{x_1\frac{1+i+j+k}{2}+x_2i+x_3j+x_4k\in\HH\;:\;
x_1,x_2,x_3,x_4\in\ZZ\big\}\;,
$$
which is a maximal order. Let $\I'$ be a non-zero two-sided ideal in
the ring $\O'$.

We denote by $\Ga_{\I'}$ the preimage in the group morphism ${\rm
  SL}_2(\O')\ra {\rm GL}_{2}(\O'/\I')$ of reduction modulo $\I'$ of the
subgroup of upper triangular matrices.
As $\O'/\I'$ is finite ($\I'$ is nonzero), $\Ga_{\I'}$ is a finite index
subgroup of ${\rm SL}_2(\O')$.

\blemm \label{lem:horoprecinv5hyp} %
The horoball $\H_1$ is precisely invariant under $\Ga_{\I'}$.
Furthermore, if $\I'=\O'$, then $\H_1$ is the maximal horoball
centered at $\infty$ which is precisely invariant under $\Ga_{\I'}$.
\elemm

\dem %
The element $\left(\begin{array}{cc} 1 & 1 \\ 0 & 1\end{array}\right)$
belongs to $\Ga_{\I'}$. It follows 
from \cite[page 1091]{Kel2} 
that the horoball $\H_1$ is
precisely invariant.

If $\I'=\O'$, then $g=\left(\begin{array}{cc} 0 & 1 \\ -1 &
    0\end{array}\right)$ belongs to $\Ga_{\I'}$, and by Lemma
\ref{lem:depthcalc5hyp}, $d(\H_1,g\H_1)=0$, hence the last assertion
follows.  
\cqfd

\subsection{On arithmetic lattices}
\label{subsec:arithsubgrou}

The following result follows from the work of Borel and Harish-Chandra
\cite{BHC} and of Borel 
\cite[Theo.~1.10]{Bor2} 
(see \cite{Bor1} for 
an elementary presentation of semisimple algebraic groups). Two
subgroups $A$ and $B$ of a group $C$ are said to be {\it
commensurable} in this theorem if $A\cap B$ has finite index in both
$A$ and $B$.

\btheo \cite{BHC,Bor2}\label{theo:aritgroulatt} %
Let $\underline{G}$ be a connected semisimple algebraic group defined
over $\QQ$ of $\RR$-rank one and $\underline{P}$ be a minimal
parabolic subgroup of $\underline{G}$ defined over $\QQ$, let
$G=\underline{G}(\RR)_0$ and $P=G\cap\underline{P}(\RR)$, let $\Ga$ be
a subgroup of $G$ commensurable to $\underline{G}(\ZZ)\cap G$, then
$\Ga$ is a lattice in $G$, and the set of parabolic fixed points of
$\Ga$ on $G/P$ is $\underline{G}(\QQ)P$.
\cqfd
\etheo

Such a subgroup $\Ga$ will be called an {\it arithmetic lattice} in
$G$. Note that the $\RR$-rank assumption is equivalent to the fact
that for every (or equivalently any) maximal compact subgroup $K$ of
the Lie group $G$, the associated symmetric space $X=G/K$ may be
endowed with a $G$-invariant Riemannian metric with sectional
curvature at most $-1$. Such a metric is then unique up to
multiplication by a positive constant, and $P$ is the stabilizer of a
point in the boundary at infinity $\partial_\infty X$. The orbital map
at this point hence induces a $G$-equivariant homeomorphism
between $G/P$ and $\partial_\infty X$. Note that there is a
terminology problem: by a parabolic element, we mean an isometry of $X$
having a unique fixed point (called a {\it parabolic fixed point}) on
$X\cup\partial_\infty X$, that belongs to $\partial_\infty X$, but the
set of real points of a parabolic subgroup of $\underline{G}$ also
contains non parabolic elements~!

\medskip
\noindent{\bf Examples. } (1) Let $m$ be a squarefree positive
integer, and let $\I$ be a non-zero ideal in an order $\O$ in the ring of
integers $\O_{-m}$ of the imaginary quadratic number field
$K_{-m}=\QQ(i\sqrt{m})$. Let $(1,\omega)$ be a basis of $\O$ as a
$\ZZ$-module. It is also a basis of $K_{-m}$ as a $\QQ$-vector space,
and of $\CC$ as an $\RR$-vector space. If $x_1,\dots,x_n, y_1,\dots,
y_n$ are real numbers, as $\omega$ is a quadratic integer, note that
$$
\prod_{i=1}^n (x_i+\omega\,y_i)=P(x_1,\dots,x_n,y_1,\dots,y_n)+
\omega\,Q(x_1,\dots,x_n,y_1,\dots,y_n) 
$$ 
where $P$ and $Q$ are polynomials with integer coefficients in
$x_1,\dots,x_n,y_1,\dots,y_n$.

By writing each coefficient of a $n\times n$ complex matrix $X$ in the
basis $(1,\omega)$ over $\RR$, the equation ${\rm det}\;X=1$ gives a
system of two polynomial equations with integer coefficient, with
unknown the coordinates of the coefficients of $X$ in $(1,\omega)$.

Hence, there exists an algebraic group $\underline{G}$ defined over
$\QQ$ such that $\underline{G}(\ZZ)={\rm SL}_2(\O)$,
$\underline{G}(\QQ)={\rm SL}_2(K_{-m})$ and $\underline{G}(\RR)= {\rm
  SL}_2(\CC)$. As the Lie group $\underline{G}(\RR)$ is connected and
semisimple, with associated symmetric space the real hyperbolic
$3$-space, the algebraic group $\underline{G}$ is connected,
semisimple with $\RR$-rank one. Let $\underline{P}$ be the algebraic
subgroup of $\underline{G}$ corresponding to the upper triangular
subgroup of $2\times 2$ matrices, so that $P$ is the stabilizer of the
point at infinity $\infty$ in the upper halfspace model of
$\HH_\RR^3$.

Let $\Ga_{\RR,\I}$ be the finite index subgroup of the group ${\rm
  SL}_2(\O)$, which is the preimage in  the group morphism ${\rm
  SL}_2(\O)\ra {\rm SL}_{2}(\O/\I)$ of reduction modulo $\I$  of the
subgroup of upper triangular matrices.
By Theorem \ref{theo:aritgroulatt}, the subgroup $\Ga_{\RR,\I}$
is a lattice in ${\rm SL}_2(\CC)$, and its set of parabolic fixed
points is 
$$
\P_{\Ga_{\RR,\I}}={\rm
  SL}_2K_{-m}\cdot\infty=K_{-m}\cup\{\infty\},
$$ 
as
$\left(\begin{array}{cc} 0 & -1 \\ 1 & 0\end{array}\right)$ and
$\left(\begin{array}{cc} x & 0 \\ 1 & x^{-1}\end{array}\right)$, for
every $x$ in $K_{-m}-\{0\}$, are elements of ${\rm SL}_2(K_{-m})$ sending
$\infty$ to $0$ and $x$ respectively.

Note that if $\I=\O=\O_{-m}$, then $\Ga_{\RR,\I}={\rm PSL}_2(\O_{-m})$
is a Bianchi group, which is well-known to be a lattice in ${\rm
  PSL}_2(\CC)$ (see for instance \cite{MR}). The fact that
$\P_{\Ga_{\RR,\O_{-m}}}= K_{-m}\cup\{\infty\}$ is also proven in
\cite[Prop.~2.2, page 314]{EGM}.

\medskip (2) Recall that $N(\omega)=\omega\ov\omega$ and ${\rm
  Tr}(\omega)=\omega+\ov\omega=2\;{\rm Re}\;\omega$ are integers, as $\omega$
is an algebraic integer. 
If $x,y,x',y'$ are real numbers, 
note that
$$
(x+\overline{\omega}y)(x'+\omega y)=
(xx'+N(\omega)\;yy'+{\rm Tr}(\omega)\;yx')+ \omega(xy'-yx')
\;.
$$
Recall that the matrix $Q$ (introduced in Section
\ref{subsec:comphypgeom}) has integer coefficients. Hence by writing
each coefficient of a $(n+1)\times (n+1)$ complex matrix $X$ in the
basis $(1,\omega)$ over $\RR$, the system of equations given by ${\rm
det}\;X=1$ and $X^*\;QX=Q$ becomes a system of $2((n+1)^2+1)$
polynomial equations with integer coefficients, with unknown the
coordinates of the coefficients of $X$ in $(1,\omega)$.

Therefore there exists an algebraic group $\underline{G}$ defined over
$\QQ$ such that $\underline{G}(\ZZ)={\rm SU}_Q(\O)$,
$\underline{G}(\QQ)={\rm SU}_Q(K_{-m})$ and $\underline{G}(\RR)= {\rm
  SU}_Q$. As the Lie group $\underline{G}(\RR)$ is connected and
semisimple, with associated symmetric space the complex hyperbolic
$n$-space, the algebraic group $\underline{G}$ is connected,
semisimple with $\RR$-rank one. Let $\underline{P}$ be the algebraic
subgroup of $\underline{G}$ corresponding to the upper triangular
subgroup of $(n+1)\times (n+1)$ matrices, so that $P$ is the
stabilizer of the point at infinity $\infty$ in the Siegel domain
model of $\HH_\CC^n$, or of the point $[1:0:0]$ in the projective
model.

By Theorem \ref{theo:aritgroulatt}, the group $\Ga_{\CC,\I}$ defined
in Section \ref{subsec:comphypgeom} is a lattice in ${\rm SU}_Q$, and
its set of parabolic fixed points is $\P_{\Ga_{\CC,\I}}={\rm
  SU}_Q(K_{-m})\cdot\infty$. By Witt's theorem, 
${\rm SU}_Q(K_{-m})$ acts transitively on the isotropic lines in
${K_{-m}}^{n-1}$ for the hermitian form $q$. Hence $\P_{\Ga_{\CC,\I}}$
is exactly the set of rational points of the quadric over $K_{-m}$
with equation $q=0$ in $\PP_n(\CC)$. If $n=2$ and $\O=\O_{-m}$, we
recover Proposition \ref{prop:identiforbit} (1).

\medskip (3) Let $\I'$ be a non-zero two-sided ideal in an order $\O'$
of a quaternion algebra $A(\QQ)$ over $\QQ$ such that
$A(\QQ)\otimes_\QQ \RR=\HH$. For
every field $K$ containing $\QQ$, define $A(K)=A(\QQ)\otimes_\QQ K$.
Let $(e_1,e_2,e_3,e_4)$ be a basis of $\O'$ as a $\ZZ$-module. It is
also a basis of $A(K)$ as a $K$-vector space for every field $K$
containing $\QQ$.

If $x=x_1+x_2i+x_3j+x_4k$ is an element in $\HH$, written in the
standard basis $(1,i,j,k)$, let ${\rm Tr}\;x=2x_1$ be its reduced
trace, and $N(x)=x_1^2+x_2^2+x_3^2+x_4^2$ be its reduced norm. 
A $2\times 2$ matrix $X$ with coefficients in $\HH$ has Dieudonné
determinant $1$ if and only if
\begin{eqnarray}\label{eq:detdieudun}
N(ad)+N(bc)+{\rm Tr}(a\overline{c}d\overline{b})=1 
\end{eqnarray}
(see for instance \cite[page 1084]{Kel2}). 
The maps  $\RR^4\ra\RR$ defined by $(x_1,x_2,x_3,x_4)\mapsto
N(x_1e_1+x_2e_2+x_3e_3+x_4e_4)$ and $(x_1,x_2,x_3,x_4)\mapsto
{\rm Tr}(x_1e_1+x_2e_2+x_3e_3+x_4e_4)$ are polynomial maps in
$x_1,x_2,x_3,x_4$ with rational coefficients. 

By writing each coefficient $a,b,c,d$ of a $2\times 2$ matrix $X$ with
coefficients in $A(K)$ in the basis $(e_1,e_2,e_3,e_4)$ for any field
$K$, the equation (\ref{eq:detdieudun}) becomes a polynomial equation
with coefficients in $\QQ$, with unknown the coordinates of the
coefficients of $X$ in $(e_1,e_2,e_3,e_4)$.

Hence there exists an algebraic group $\underline{G}$ defined over
$\QQ$ such that $\underline{G}(\ZZ)={\rm SL}_2(\O')$ and
$\underline{G}(K)={\rm SL}_2(A(K))$ for every field $K$ containing
$\QQ$.  As the Lie group $\underline{G}(\RR)={\rm SL}_2(\HH)$ is
connected and semisimple, with associated symmetric space the real
hyperbolic $5$-space, the algebraic group $\underline{G}$ is
connected, semisimple with $\RR$-rank one. Let $\underline{P}$ be the
algebraic subgroup of $\underline{G}$ corresponding to the upper
triangular 
matrices, so that $P$ is the
stabilizer of the point at infinity $\infty$ in the upper halfspace
model of $\HH_\RR^5$.

Let $\Ga_{\I'}$ be the group introduced in Section
\ref{subsec:quatrealhypgeom}, which has finite index in ${\rm
  SL}_2(\O')$. By Theorem \ref{theo:aritgroulatt}, the subgroup
$\Ga_{\I'}$ is a lattice in ${\rm SL}_2(\HH)$, and its set of
parabolic fixed points is $\P_{\Ga_{\I'}}={\rm
  SL}_2(A(\QQ))\cdot\infty$. As $A(\QQ)$ is a division algebra, the
same argument as for  example (1) shows that
$\P_{\Ga_{\I'}}=A(\QQ)\cup\{\infty\}$.

\subsection{The ubiquity of Hall rays}
\label{subsec:ubiqhallways}

In this subsection, we give applications of our geometric results from
Section \ref{sec:results} to the framework of Diophantine approximation
in negatively curved manifolds, introduced in \cite{HP3,HP4}, to which
we refer for notation and background.  In particular, we will consider
arithmetically defined examples. See also the previous works of
\cite{For,Ser}, among many others.

Let $M$ be a complete, nonelementary, geometrically finite Riemannian
manifold with sectional curvature at most $-1$ and
dimension at least $3$. Let $\pi:\wt M\ra M$ be a universal Riemannian
covering, with covering group $\Ga$. Let $e$ be a cusp of $M$, and, as
in Section \ref{sec:ballsandhoroballs}, let $V_e$ be a fixed Margulis
neighborhood of $e$, $H_e$ a horoball in $\wt M$ with $\pi(H_e)=V_e$
and $\xi_e$ the point at infinity of $H_e$. Note that requiring $V_e$
to be a Margulis neighborhood is equivalent to requiring $H_e$ to be
precisely invariant under $\Ga$. In the previous works
\cite{HP3,HP4}, it was required that $V_e$ is the maximal Margulis
neighborhood, as this makes the constructions independent of the
choice of $V_e$. But as it is not always easy to determine the maximal
Margulis neighborhood of a cusp, and as it is not necessary for the
statements, we will fix some choice of $V_e$ (or equivalently $H_e$)
which is not necessarily maximal.

\medskip\noindent{\bf Three (classes of) examples. } 
These examples could in fact be orbifolds rather than manifolds, but
the extension to this context is obvious.  
We use the same notation as in the examples of Subsection
\ref{subsec:arithsubgrou}. 

(1) Let $\Ga_{\RR,\I}$ be the finite index subgroup of the group ${\rm
  SL}_2(\O)$, which is the preimage, by the group morphism ${\rm
  SL}_2(\O)\ra {\rm SL}_{2}(\O/\I)$ of reduction modulo $\I$, of the
subgroup of upper triangular matrices.
The quotient
$M=\Ga_{\RR,\I}\backslash\HH^3_\RR$ is a finite volume real hyperbolic
orbifold.  Let $\pi:\HH^3_\RR\ra M$ be the canonical projection, $e$
the cusp of $M$ corresponding to $\xi_e=\infty$, and $H_e$ be the
horoball of points of Euclidean height at least $1$. As
$\left(\begin{array}{cc} 1 & 1 \\ 0 & 1\end{array}\right)$ belongs to
$\Ga_{\RR,\I}$, it is well known that $H_e$ is precisely invariant
under $\Ga_{\RR,\I}$.  Furthermore, if $\I=\O$, then
$\left(\begin{array}{cc} 0 & 1 \\ -1 & 0\end{array}\right)$ belongs to
$\Ga_{\RR,\I}$, hence $H_e$ is maximal. For more details, see
\cite{HP3}, end of Section 5.

(2) Let $\Ga_{\CC,\I}$ be the group of isometries of the Siegel domain
model $\HH^n_\CC$ of the complex hyperbolic $n$-space with (constant)
holomorphic sectional curvature $-1$, that was introduced in Section
\ref{subsec:comphypgeom}.
Let $M$ be the finite volume complex
hyperbolic orbifold $\Ga_{\CC,\I} \backslash \HH^n_\CC$, which is
endowed with the quotient of the renormalized Riemannian distance
$d'_{\HH^n_\CC}$ in order for its sectional curvatures to be at most
$-1$. Let $\pi: \HH^n_\CC\ra M$ be the canonical projection, $e$ be
the cusp of $M$ corresponding to $\xi_e=\infty$, and $H_e$ the
horoball $\H_{2\;{\rm Im}\,\omega}$ if ${\rm Re}\,\omega\in\ZZ$, and
$\H_{4\;{\rm Im}\,\omega}$ otherwise, which is precisely invariant
under $\Ga_{\CC,\I}$ by Lemma \ref{lem:horoprecinvcomphyp} (and
maximal if $\I=\O=\O_{-1}$).

(3) Let $\I'$ be a non-zero two-sided ideal in an order $\O'$ of a
quaternion algebra $A(\QQ)$ over $\QQ$ such that $A(\QQ)\otimes_\QQ
\RR=\HH$, and $\Ga_{\I'}$ the group of isometries of the upper
halfspace model $\HH^5_\RR$ of the real hyperbolic $5$-space with
(constant) sectional curvature $-1$, that was introduced in Section
\ref{subsec:quatrealhypgeom}. 
Let $M$ be the
finite volume real hyperbolic orbifold $\Ga_{\I'} \backslash
\HH^5_\RR$, $\pi: \HH^5_\RR\ra M$ be the canonical projection, $e$ be
the cusp of $M$ corresponding to $\xi_e=\infty$, and $H_e$ the
horoball $\H_1$, which is precisely invariant under $\Ga_{\I'}$ by
Lemma \ref{lem:horoprecinv5hyp} (and maximal if $\I'=\O'$).

\medskip %
Let the {\em link} of $e$ in $M$, ${\rm Lk}_e={\rm Lk}_e(M)$,  be the
space of locally geodesic lines (up to 
translation at the source) starting from $e$ in $M$ that are
nonwandering (i.e.~such that each of them accumulates in some compact
subset of $M$). Let ${\rm Rat}_e$ be the space of locally geodesic lines
starting from $e$ and converging to $e$. Normalize the locally geodesic lines
in ${\rm Lk}_e\cup {\rm Rat}_e$ so that their first intersection with
$\partial V_e$ is at time $0$. Endow ${\rm Lk}_e\cup {\rm Rat}_e$ with
the compact-open topology. Let $\Lambda\Gamma\subset
\partial_\infty\widetilde M$ be the limit set of $\Ga$,
$\P_\Ga\subset\Lambda\Gamma$ be the set of parabolic fixed points of
$\Ga$, and $\Ga_\infty$ be the stabilizer of $\xi_e$ in $\Ga$. Then
the maps $\Ga\xi_e-\{\xi_e\}\ra {\rm Rat}_e$ and
$\Lambda\Gamma-\P_\Ga\ra {\rm Lk}_e$, which associate to $x$ the
projection in $M$ by $\pi$ of the geodesic line starting from
$\xi_e$ and ending at $x$, induce a bijection $\Ga_\infty\backslash
(\Ga\xi_e-\{\xi_e\})\ra {\rm Rat}_e$ and a homeomorphism
$$
\Ga_\infty\backslash (\Lambda\Gamma-\P_\Ga)\ra {\rm Lk}_e\;.
$$
We identify these spaces by these maps. Note that ${\rm Lk}_e\cup
{\rm Rat}_e$ is compact if and only if $M$ has only one cusp, and that
${\rm Rat}_e$ is dense in ${\rm Lk}_e\cup {\rm Rat}_e$ (as $\Ga\xi_e$
is dense in $\Lambda\Gamma$).  
Diophantine approximation in $M$ (see
\cite{HP3,HP4,HP2}) studies the rate of convergence of sequences of
points in ${\rm Rat}_e$ to given points in ${\rm Lk}_e$.

For every $r$ in ${\rm Rat}_e$, let $D(r)$, called the {\it depth} of
$r$, be the length of the subsegment of $r$ between the first and the
last meeting points with $\partial V_e$.

\medskip
\noindent{\bf Examples. }%
(1) Consider $M=\Ga_{\RR,\I}\backslash\HH^3_\RR$. Then 
$\P_{\Ga_{\RR,\I}}\subset \CC\cup\{\infty\}$ is exactly
$K_{-m}\cup\{\infty\}$, by the example (1) of Section
\ref{subsec:arithsubgrou}. Thus, ${\rm
  Lk}_e(M)=(\Ga_{\RR,\I})_\infty\backslash(\CC-K_{-m})$.  In a
commutative unitary ring $R$, 
we denote by $\langle p_1,\dots,p_k\rangle$ the ideal
generated by 
$p_1,\dots,p_k\in R$. It is easy to
prove (see for instance \cite[Lem.~2.1, page 314]{EGM}) that ${\rm
  Rat}_e$ is the set of elements $r=p/q$ (modulo
$(\Ga_{\RR,\I})_\infty$) with $(p,q) \in\O\times\I$ such that $\langle
p,q\rangle=\O$. Furthermore (see \cite[Lem.~2.10]{HP3})
$$
D(r)=2\log|q|\;.
$$

(2) Consider $M=\Ga_{\CC,\I}\backslash \HH^n_\CC$. Let $\Q(\RR)$ be
the real quadric $\partial_\infty\HH^n_\CC-\{\infty\}$. By
considering a basis of $K_{-m}$ over $\QQ$, it is easy to see that
$\Q(\RR)$ is the set of $\RR$-points of a quadric $\Q$ defined over
$\QQ$ (which depends on $m$), whose set $\Q(\QQ)$ of $\QQ$-points is
$\Q(\RR)\cap(K_{-m}\times K_{-m}^{n-1})$. 
We have ${\rm Lk}_e=(\Ga_{\CC,\I})_\infty
\backslash (\Q(\RR)- 
\P_{\Ga_{\CC, \I}})$. By the example (2) of Section
\ref{subsec:arithsubgrou}, we have $\P_{\Ga_{\CC,\I}}= \Q(\QQ)\cup
\{\infty\}$.

Then ${\rm Rat}_e$ is the quotient modulo $(\Ga_{\CC,\I})_\infty$ of
the subset of $\Q(\QQ)$ of points of the form $(a/c,\alpha/c)$ with
$(a,\alpha,c)\in\O\times \I^{n-1}\times \I$ such that there exist
$b,c,\beta,\gamma,\delta,A$ matrices of the appropriate size such that
$\left(\begin{array}{ccc}a & \ga^* & b\\ \alpha & A & \beta \\
c&\delta^* & d \end{array}\right)$ belongs to $\Ga_{\CC,\I}$. By
Proposition \ref{prop:identiforbit} (2), this existence requirement is
equivalent to the requirement 
that $q(a,\alpha,c)=0$ and $\langle a,\alpha,c\rangle=\O$, if $n=2$
and $\O=\O_{-m}$.  By Proposition \ref{prop:identiforbit} (3), ${\rm
Rat}_e= (\Ga_{\CC,\I})_\infty\backslash\Q(\QQ)$ if $n=2$,
$\I=\O=\O_{-m}$ and $m=1,2,3,7,11,19,43,67,163$.  

If $r\in{\rm Rat}_e$ is of the form $(a/c,\alpha/c)$ (modulo
$(\Ga_{\CC,\I})_\infty$) as above, then by Lemma
\ref{lem:depthcalccomphyp}, we have
$$
D(r)=\log|c|+\left\{\begin{array}{ll} \log\;{\rm Im}\;\omega & 
    {\rm if}\;\; {\rm Re}\;\omega\in\ZZ, \\ 
\log(2\;{\rm Im}\;\omega) & {\rm otherwise.}\end{array}\right.
$$

(3) Consider $M=\Ga_{\I'} \backslash \HH^5_\RR$. We have
$\P_{\Ga_{\I'}} = A(\QQ)\cup \{\infty\}$, by the example (3) of
Section \ref{subsec:arithsubgrou}. Hence
${\rm Lk}_e= (\Ga_{\I'})_\infty \backslash (\HH-A(\QQ))$.

It is easy to see
that ${\rm Rat}_e$ is the set of elements $r=pq^{-1}$ (modulo
$(\Ga_{\I'})_\infty$) with $(p,q)\,\in\,\O'\times(\I'-\{0\})$ such that there
exists $r,s\in\O'$ with $|qr-qpq^{-1}s|=1$. Furthermore, by Lemma
\ref{lem:depthcalc5hyp}, we have
$$
D(r)=2\log|q|\;.
$$

\bigskip %
The {\it cuspidal distance}
 $d'_e(\ga,\ga')$ 
of $\ga,\ga'$ in ${\rm Lk}_e\cup {\rm Rat}_e$ 
is the minimum of the $\wt{d'_e}(\wt{\ga},\wt{\ga'})$
for $\wt{\ga},\wt{\ga'}$ two lifts of $\ga,\ga'$ to $\wt M$ starting
from $\xi_e$, where $\wt{d'_e}(\wt{\ga},\wt{\ga'})$ is the greatest lower bound
of $r>0$ such that the horosphere centered at $\wt{\ga}(+\infty)$, at
signed distance $-\log 2r$ from $\partial H_e$ on the geodesic line
$]\xi_e,\wt{\ga}(+\infty)[$, meets $\wt{\ga'}$  (see
\cite[Sect.~2.1]{HP3}). Though not necessarily 
an actual distance, $\wt{d'_e}$ is equivalent to the Hamenst\"adt
distance (see Subsection \ref{subsec:penetrationmaps} and
\cite[Rem.~2.6]{HP3}).

\medskip
\noindent{\bf Examples. } %
(1) If $M$ has constant curvature $-1$, if one identifies ${\rm
  Lk}_e\cup {\rm Rat}_e$ with a subset of $\partial H_e$ by the first
intersection point, then $d'_e$ is the induced Riemannian distance on
$\partial H_e$, which is Euclidean (see \cite[Sect.~2.1]{HP3}); in
particular, if $M=\Ga_{\RR,\I}\backslash \HH^3_\RR$ or
$M=\Ga_{\I'}\backslash \HH^5_\RR$, then $d'_e$ is the quotient of the
standard Euclidean distance on ${\rm Lk}_e\cup {\rm Rat}_e$ identified
with a subset of $(\Ga_{\RR,\I})_\infty\backslash\CC $ or
$(\Ga_{\I'})_\infty\backslash\HH$.

(2) If $M$ is Hermitian with constant holomorphic sectional curvature
$-1$, then $d'_e$ is no longer Riemannian, but by Proposition
\ref{prop:cuspdistcalccomphyp} is a multiple of the modified Cygan
distance $d'_{\rm Cyg}$. In particular, if $M=\Ga_{\CC,\I}\backslash
\HH^n_\CC$, then $d'_e$ is the quotient by $(\Ga_{\CC,\I})_\infty$
acting on $\partial_\infty\HH^n_\CC$ of the distances
$$
\left\{\begin{array}{ll} \frac{1}{2\sqrt{{\rm Im}\;\omega}}\;
d'_{\rm Cyg} &
    {\rm if}\;\; {\rm Re}\;\omega\in\ZZ \\
    \frac{1}{2\sqrt{2\;{\rm Im}\;\omega}}\;d'_{\rm Cyg} & {\rm
      otherwise~.}
\end{array}\right.
$$

\medskip \rem %
The claim in the first paragraph of Section 3.11 in \cite{HP4} (where
we only considered $m=1$ and $\I=\O=\O_{-1}$) that the cuspidal
distance coincides with the Hamenst\"adt distance is incorrect; as the
Cygan distance and the modified Cygan distance are equivalent, this
does not change the statement of the main results, Theorems 3.1, 3.2,
3.4, and 3.5 of \cite{HP4}. As $d'_{\rm Cyg}\leq \sqrt{2}\;d_{\rm Cyg}$ so
that $d_{\rm 
  Cyg}\geq \sqrt{2}d'_e$, the constant $\frac{1}{\sqrt[4]{5}}$
appearing in Theorem 3.6 of \cite{HP4} has to be replaced by
$\frac{1}{2\sqrt[4]{5}}$.

\medskip %
With $M$ as in the beginning, for every $x$ in ${\rm Lk}_e$, define
the {\it approximation constant} $c(x)$ of $x$ as
$$
c(x) =\liminf_{r\in {\rm Rat}_e\,,\; 
D(r)\ra\infty} d'_e(x,r)\;e^{D(r)}\;.
$$
The {\it Lagrange spectrum} of $M$ with respect to $e$ is the subset
${\rm Sp}_{\rm Lag}(M,e) $ of $\RR$ consisting of the constants
$c(x)$ for $x$ in ${\rm Lk}_e$. It is shown in \cite{HP3} that 
\begin{itemize}
\item $c(x)$ is well defined (as ${\rm Rat}_e$ is dense in ${\rm
  Lk}_e\cup{\rm Rat}_e$ and $\{D(r)\;:\;r\in {\rm Rat}_e\}$ is a
discrete subset of $\RR$ with finite multiplicities),
\item  $c(x)$ is finite (as $x$ is nonwandering), (Note that if $y$
is a locally geodesic line starting from $e$ in $M$ that converges
into a cusp of $M$, then the same formula would yield $c(y)=+\infty$.)
\item  the least upper bound of ${\rm Sp}_{\rm Lag}(M,e)$, denoted by
$K_{M,e}$ and called the {\it Hurwitz constant}, is finite. 
\end{itemize}
In particular, ${\rm Sp}_{\rm Lag}(M,e)\subset [0,K_{M,e}]$. The
following result tells us that the Lagrange spectrum contains a
nontrivial initial interval $[0,c]$, with a universal lower bound on
$c$ (whose optimal value we do not know).

\btheo \label{theo:applidioappox}%
Let $M$ be a complete, nonelementary, geometrically finite Riemannian
manifold with sectional curvature at most $-1$ and
dimension at least $3$, and let $e$ be a cusp of $M$.
The Lagrange spectrum ${\rm Sp}_{\rm Lag}(M,e)$ contains the
interval $[0,0.0337]$. In particular, $K_{M,e}\geq 0.0337$.
\etheo

\dem %
By \cite{HP3}, the map $x\mapsto -2\log x$ sends bijectively the
Lagrange spectrum onto the asymptotic height spectrum. We  apply
Corollary \ref{coro:limsupsp} and the computation above it.  
\cqfd

\medskip %
A precise version of this theorem is stated as corollaire 5 in
\cite{PP2} when $M$ is a real or complex hyperbolic manifold. Theorem
\ref{theo:applibianchi} in the introduction follows immediately, by
the first example discussed in this section.  By varying the (non
uniform) arithmetic lattices in the isometry group of a negatively
curved symmetric space (see for instance
\cite{MR,MWW}), other arithmetic applications are possible. We only
state two of them in what follows.

\medskip %
Let $\I'$ be a non-zero two-sided ideal in an order $\O'$ of a
quaternion algebra $A(\QQ)$ over $\QQ$ ramifying over $\RR$, and $N$
the reduced norm on $A(\RR)=A(\QQ)\otimes_\QQ\RR$ (see for instance
\cite{Vig}, and Section \ref{subsec:quatrealhypgeom}). For every $x\in
A(\RR)-A(\QQ)$, define the {\it approximation constant} of $x$ by
$$
c(x)=\liminf_{(p,q)\,\in\,\O'\times\I'\;:\; 
\exists\; r,s\in\O'\; N(qr-qpq^{-1}s)=1\,,
\;N(q)\ra\infty} \;\;  N(q)N(x-pq^{-1})^{\frac{1}{2}}\;,
$$ 
and the {\it Hamilton-Lagrange spectrum} for the approximation of
elements of $\HH$ by elements of $\O'\I'^{-1}$ as the subset of $\RR$
consisting of the $c(x)$ for $x\in A(\RR)-A(\QQ)$. Note that
$c(x)$ is finite if $x\notin A(\QQ)$, as then $x$ is not a parabolic
fixed point of $\Ga_{\I'}$. Apply Theorem \ref{theo:applidioappox} to
$M=\Ga_{\I'}\backslash \HH^5_\RR$ with the above discussions of the
third example to get the following result.

\btheo \label{theo:quaternionentiers} %
The Hamilton-Lagrange spectra  contain the interval $[0,0.0337]$.
\cqfd \etheo

In the case when  $\I'=\O'$ and $\O'$ is the Hurwitz maximal order in
Hamilton's 
quaternion algebra $A(\QQ)\subset\HH$ (see Subsection
\ref{subsec:quatrealhypgeom}), A.~Schmidt \cite{Sch2} proved that the
Hamilton-Lagrange spectrum contains $\sqrt 2\;{\rm Sp}_\QQ$ where ${\rm
Sp}_\QQ$ is the classical Lagrange spectrum for the approximation of
real numbers by rational numbers. As ${\rm Sp}_\QQ$ contains $[0,\mu]$
where $\mu$ is Freiman's constant (see the end of Subsection
\ref{subsec:limsupspec}), this proves that the Hamilton-Lagrange
spectrun in this case contains the interval $[0,0.312]$, which is
reasonably close to our general estimate. Note that the fact that our
approximation constant coincides with the inverse of A.~Schmidt's
approximation constant follows from \cite[Thm.5]{Sch1}
\medskip

Let $m$ be a squarefree positive integer, $\I$ be a non-zero ideal in
an order $\O$ in the ring of integers $\O_{-m}$ of the imaginary
quadratic number field $\QQ(i\sqrt{m})$, and $\omega$ be the element
of $\O_{-m}$ with ${\rm Im}\,\omega>0$ such that $\O=\ZZ+\omega\ZZ$.
Let $\E_{\O,\I}$ be the set of $(a,\alpha,c)$ in
$\O\times\I^{n-1}\times\I$ such that there
exists a matrix of the form $\left(\begin{array}{ccc}a & \ga^* & b\\
    \alpha & A & \beta \\ c&\delta^* & d \end{array}\right)$ that
belongs to $\Ga_{\CC,\I}$. If $n=2$ and $\O=\O_{-m}$, then, as seen
previously,  
$$
\E_{\O,\I}=\{(a,\alpha,c)\in\O\times\I^{n-1}\times\I\;:\;
q(a,\alpha,c)=0,\;\langle a,\alpha,c\rangle=\O\}\;.
$$ 
We do not know if this is  the case for every $n,\I,\O$.
For every $x\in\Q(\RR)-\Q(\QQ)$, define the {\it
  approximation constant} of $x$ by
$$
c(x)=\liminf_{(a,\alpha,c)\in\,\E_{\O,\I}\,,
\;|c|\ra\infty} \;\;  |q|\;d'_{\rm Cyg}(x,(a/c,\alpha/c))\;,
$$ 
and the {\it Heisenberg-Lagrange spectrum}, for the approximation of
elements of $\Q(\RR)$ by elements of $\{(a/c,\alpha/c)\;:\;
(a,\alpha,c)\in\E_{\O,\I}\}\subset \Q(\QQ)$, as the subset of $\RR$
consisting of the $c(x)$ for $x\in\Q(\RR)-\Q(\QQ)$. Note that $c(x)$
is finite if $x\notin\Q(\QQ)$, as then $x$ is not a parabolic fixed
point of $\Ga_{\CC,\I}$. Our last result follows from Theorem
\ref{theo:applidioappox} and the previous discussions of the second
example.

\btheo \label{theo::applicygangene} %
The Heisenberg-Lagrange spectra contain the interval
$[0\,,\;0.0674\,({\rm Im}\;\omega)^{-\frac12}]$ if ${\rm Re}
\;\omega\in\ZZ$ and $[0\,,\;0.0476\,({\rm Im} \;\omega)^{-\frac12}]$
otherwise.  \cqfd 
\etheo

Theorem \ref{theo:applicygan} in the introduction follows from this
one by taking $m=1$, $n=2$, $\I=\O=\O_{-1}$, as then $\omega=i$.

\bigskip
{\small\noindent \begin{tabular}{l} 
Department of Mathematics and Statistics, P.O. Box 35\\ 
40014 University of Jyv\"askyl\"a, FINLAND.\\
{\it e-mail: parkkone@maths.jyu.fi}
\end{tabular}
\medskip

\noindent \begin{tabular}{l}
D\'epartement de Math\'ematique et Applications, UMR 8553 CNRS\\
Ecole Normale Sup\'erieure, 45 rue d'Ulm\\
75230 PARIS Cedex 05, FRANCE\\
{\it e-mail: Frederic.Paulin@ens.fr}
\end{tabular}
}

\end{document}